\documentclass[aos]{imsart}

\RequirePackage{amsthm,amsmath,amsfonts,amssymb}
\RequirePackage[numbers,sort&compress]{natbib}
\RequirePackage[colorlinks,citecolor=blue,urlcolor=blue]{hyperref}
\RequirePackage{graphicx}

\startlocaldefs
\usepackage{bm}
\usepackage{multirow}
\usepackage{bbm}
\usepackage{ulem}
\usepackage{dsfont}
\usepackage{makecell}
\usepackage{booktabs}
\usepackage{mathrsfs}
\usepackage{float}
\usepackage{graphicx}
\usepackage{amssymb}
\usepackage{amsmath}
\usepackage{amsthm}
\usepackage{empheq} 
\usepackage{subfigure}
\usepackage{epstopdf}
\usepackage{mdframed}
\usepackage{mathtools}
\usepackage{booktabs}
\usepackage{lipsum}
\usepackage{graphicx}
\usepackage{mathtools}  
\mathtoolsset{showonlyrefs}
\usepackage{anyfontsize}
\usepackage{tikz}
\usetikzlibrary{backgrounds,automata}
\usetikzlibrary{shapes.geometric}
\usetikzlibrary{shapes.arrows}
\usetikzlibrary{calc}
\usetikzlibrary{intersections}
\usepackage{color}
\usepackage{psfrag}
\usepackage{pgfplots}
\usepackage{hyperref}
\usepackage[ruled,linesnumbered]{algorithm2e}
\graphicspath{ {./images/}}
\hypersetup{
  colorlinks,
  linkcolor={red!50!black},
  citecolor={blue},
  urlcolor={blue!80!black}
}
\usepackage[toc,page]{appendix}
\usepackage{enumitem}
\usepackage{comment}
\definecolor{ocre}{RGB}{243,102,25}
\definecolor{mygray}{RGB}{243,243,244}
\definecolor{deepGreen}{RGB}{26,111,0}
\definecolor{shallowGreen}{RGB}{235,255,255}
\definecolor{deepBlue}{RGB}{61,124,222}
\definecolor{shallowBlue}{RGB}{235,249,255}
\definecolor{americanrose}{rgb}{1.0, 0.01, 0.24}
\definecolor{airforceblue}{rgb}{0.36, 0.54, 0.66}
\definecolor{aliceblue}{rgb}{0.94, 0.97, 1.0}
\definecolor{lightcyan}{rgb}{0.88, 1.0, 1.0}
\definecolor{amber(sae/ece)}{rgb}{1.0, 0.49, 0.0}
\definecolor{amber}{rgb}{1.0, 0.75, 0.0}
\definecolor{amaranth}{rgb}{0.9, 0.17, 0.31}

\newcommand{\overbar}[1]{\mkern 1.5mu\overline{\mkern-1.5mu#1\mkern-1.5mu}\mkern 1.5mu}

\DeclareMathOperator*{\argmin}{arg\,min}

\newcommand{\bb}[1]{\mathbb{#1}}
\SetKwComment{Comment}{/* }{ */}
\newcommand{\norm}[1]{\left\lVert#1\right\rVert}

\newcommand{\eq}[1]{\begin{equation} #1\end{equation}}
\newcommand{\ip}[1]{\left\langle#1 \right\rangle}
\newcommand{\mc}[1]{\mathcal{#1}}
\newcommand{\longeq}[1]{\begin{equation}
    \begin{split}
        #1
    \end{split}
\end{equation}}

\def\hat{\widehat}
\def\tilde{\widetilde}
\def\t{^{\top}}
\def\bar{\overbar}

\def\bo{\boldsymbol }

\def\ret{\mathsf{Retr}}

\newtheorem{theorem}{Theorem}[section]
\newtheorem{lemma}[theorem]{Lemma}
\newtheorem{claim}[theorem]{Claim}
\newtheorem{definition}[theorem]{Definition}
\newtheorem{assumption}[theorem]{Assumption}
\newtheorem{corollary}[theorem]{Corollary}
\newtheorem{remark}[theorem]{Remark}
\newtheorem{proposition}[theorem]{Proposition}

\newtheorem{example}[theorem]{Example}
\newcommand{\mb}{\boldsymbol}
\def\hat{\widehat}

\def\Log{\mathsf{Log}}
\def\Exp{\mathsf{Exp}}
\def\dist{\mathsf{dist}}

\def\tilde{\widetilde}
\def\t{^{\top}}
\usepgfplotslibrary{colorbrewer}
\def\spacingset#1{\renewcommand{\baselinestretch}
{#1}\small\normalsize}

\spacingset{1}
\pgfplotsset{compat=1.18}
\endlocaldefs

\begin{document}
\begin{frontmatter}
\title{High-order Accurate Inference on Manifolds}
\runtitle{High-order Accurate Inference on Manifolds}

\begin{aug}
\author[A]{\fnms{Chengzhu}~\snm{Huang}\ead[label=e1]{ch3786@columbia.edu}} and
\author[B]{\fnms{Anru R.}~\snm{Zhang}\ead[label=e2]{anru.zhang@duke.edu}}
\address[A]{Department of Statistics, Columbia University \printead[presep={ ,\ }]{e1}}

\address[B]{Department of Biostatistics \& Bioinformatics and Department of Computer Science, Duke University\printead[presep={,\ }]{e2}}
\end{aug}

\begin{abstract}
We present a new framework for statistical inference on Riemannian manifolds that achieves high-order accuracy, addressing the challenges posed by non-Euclidean parameter spaces frequently encountered in modern data science. Our approach leverages a novel and computationally efficient procedure to reach higher-order asymptotic precision. In particular, we develop a bootstrap algorithm on Riemannian manifolds that is both computationally efficient and accurate for hypothesis testing and confidence region construction. Although locational hypothesis testing can be reformulated as a standard Euclidean problem, constructing high-order accurate confidence regions necessitates careful treatment of manifold geometry. To this end, we establish high-order asymptotics under an appropriate coordinate representation induced by a second-order retraction, thereby enabling precise expansions that incorporate curvature effects. We demonstrate the versatility of this framework across various manifold settings, including spheres, the Stiefel manifold, fixed-rank matrix manifolds, and rank-one tensor manifolds; for Euclidean submanifolds, we also introduce a class of projection-like coordinate charts with strong consistency properties. Finally, numerical studies confirm the practical merits of the proposed procedure.
\end{abstract}



\end{frontmatter}
\begin{sloppypar}
\section{Introduction}\label{sec:intro}
This paper focuses on the statistical inference for parameters defined on a Riemannian manifold. Let $\mc X_n \coloneqq \{X_1, \ldots, X_n\} \subset \mc S$ be observations in a sample space $\mc S$ independently sampled from a distribution $\bb P_{\theta_0}$, where $\theta_0$ resides on a $p$-dimensional Riemannian manifold $\mc M$. Our aim is to make inference on $\theta$ based on $\mc X_n$. For any positive integer $N$, we write $[N] \coloneqq \{1, \ldots, N\}$. Given an objective function $L(\theta, x): \mc M \times \mc S \rightarrow \bb R$ and its sample average analog $L_n(\theta,\mc X_n) \coloneqq \sum_{i\in[n]} L(\theta, X_i)/n$, an M-estimator, an important instance of extremum estimator, can be defined as
\eq{\label{eq: argmin}
\tilde \theta_n \coloneqq    \argmin_{\theta \in \mc M} L_n(\theta,\mc X_n).  
}
For brevity, we often write $L_n(\theta, \mc X_n)$ in place of $L_n(\theta)$.

A natural approach to assess the location of $\tilde \theta_n$ is to exploit a suitable coordinate representation. Under the assumption that the likelihood function is geodesically convex but possibly non-smooth, \cite{brunel2023geodesically} explored the first-order asymptotics of $\tilde \theta_n$ and proved that the coordinate of $\tilde{\theta}_n$ under the logarithmic mapping centered at $\theta_0$ and a basis $\{e_i\}$ of $\mathrm T_{\theta_0}\mc M$ is asymptotically distributed as
\eq{
\sqrt{n} \cdot \pi_{\{e_i\}} \circ \Log_{\theta_0}(\tilde \theta_n) \Rightarrow \mc N(\mb 0, \mb \Sigma),
}
where $\pi_{\{e_i\}}$ denotes a canonical mapping from $\mathrm T_{\theta_0}\mc M$ to $\bb R^p$ induced by $\{e_i\}$,  $\Log_{\theta_0}$ represents the logarithmic mapping, and 
\longeq{
& \mb \Sigma \coloneqq \Big(\big( \bb E\left[\nabla^2 L(\theta_0, X_1)[e_i,e_j]\right]\big)_{i,j\in[p]}\Big)^{-1} \big(\bb E\big[e_i L(\theta_0, X_1)e_j L(\theta_0, X_1)\t\big]\big)_{i,j\in[p]} \\ 
& \qquad \cdot \Big(\big( \bb E\left[\nabla^2 L(\theta_0, X_1)[e_i,e_j]\right]\big)_{i,j\in[p]}\Big)^{-1}. \label{eq: definition of Sigma}
}

The primary goal of this paper is to go beyond this asymptotic normality of the exact extremum estimators and establish computationally feasible inference methods that achieve high-order accuracy in manifold settings. Specifically, we aim to design a hypothesis test that ensures a high-order accurate control of type-I error. 
Additionally, for a confidence level $1 - \alpha$, our goal is to construct a confidence region $I\subseteq \mc M$ that satisfies
\eq{
\big\vert\bb P_{\mc X_n}[\theta_0 \in I] - (1 - \alpha) \big\vert = o(n^{-1/2}).
}

To achieve this goal, we employ a resampling procedure for inference on $\theta$, so as to mimic the distribution of the targeted quantity as well as to construct the region $I$. Toward the theoretical guarantees for the inference tasks, we focus on demonstrating the distributional closeness between the approximate solution and the resampled approximate solution after appropriate studentification, in line with the Edgeworth expansions of (approximate) M-estimators \citep{hall1996bootstrap,andrews2002higher}.

\paragraph*{Challenges in Manifold Settings}\label{sec:challenges}

We describe a smooth manifold, without delving into concepts of differential geometry, as a space that locally behaves like Euclidean space. A smooth manifold endowed with a metric structure is a Riemannian manifold. The primary challenge in dealing with such an abstract structure lies in the locality of the coordinates; specifically, we must account for the different representations of a point under various charts. The core idea of this paper is to consider exponential mapping, where the relationship between different charts is determined by the Riemannian connection defined on the manifold. Leveraging results from differential geometry \citep{gavrilov2007double}, we are able to upper bound the curvature's effect on the coordinates and explicitly interpret the distribution of the proposed estimator. Furthermore, exponential mappings are often computationally expensive in many applications \cite{edelman1998geometry,absil2012projection}, prompting us to explore alternative methods, such as retraction mappings, to achieve second-order accuracy.

\subsection{Extremum Estimators on Manifolds}\label{sec:extreme}

An extremum estimator is obtained by optimizing an objective function that depends on observed data \citep{newey1994large,fisher1922mathematical,wald1949note}. Maximum likelihood estimation (MLE) is a prominent example and has been widely used for inference in parametric models. Notably, the framework of extremum estimators naturally lends itself to a manifold-based formulation in several important statistical analysis scenarios.

\paragraph*{Constrained Parameter Space}

The asymptotic properties of statistical models with a fixed-dimensional parameter in Euclidean space have been extensively studied \citep{van2000asymptotic}. However, there are many statistical problems that do not exhibit the desirable properties associated with Euclidean spaces. For example, the behavior of the MLE and the associated information metric in curved exponential families exhibit notable complexity \citep{efron1975defining,efron1978geometry}. The nonlinearity constraints in such problems offer intriguing interpretability, while also posing challenges in terms of coordinate representation, naturally falling into the differential geometric viewpoint. Given a parametric space $ \Theta = \{\theta: f_i(\theta) = 0,~i\in [r]\} \subset \bb R^p$ ($r<p$) where $\{f_i\}_{i\in[r]}$ are smooth functions whose derivatives are linearly independent, the constant-rank level set theorem \citep{lee2013smooth} shows $\Theta$ is an embedded submanifold with dimension $(p-r)$. However, the approach to enable high-order accurate inference in these settings is as of yet undeveloped. 

\paragraph*{Nonconvex Optimization} Nonconvex optimization problems often arise in modern machine learning and signal processing,  often accompanied by complex geometric structures that cannot be fully captured using standard Euclidean tools. In particular, the matrix/tensor-valued parameter estimation problems with low-rank constraints naturally lead to formulations on specialized manifolds such as the fixed-rank matrix manifold, Stiefel manifold, and the Grassmannian manifold. Recent developments have shown the significance of such manifold-based approaches, either theoretically or algorithmically (e.g., \cite{boumal2019global,li2021weakly,chen2020proximal,huang2016riemannian}). 

\paragraph*{Examples} We present several key Riemannian manifolds below:
\begin{itemize}
\item The Stiefel manifold $\mathrm{St}(p,r) \coloneqq \{\mb X \in \bb R^{p\times r}: \mb X\t \mb X = \mb I_r\}$ consists of all $p\times r$ orthonormal matrices. This manifold generalizes the spheres to higher dimensions and is pivotal in applications requiring orthonormal bases or dimensionality reduction, such as in factor analysis \citep{zheng2025sofari}. 
\item The fixed-rank matrix manifold, denoted as $\mc R_{r, p_1, p_2}$, represents the set of all $p_1 \times p_2$ matrices with rank $r$ and possesses a smooth manifold structure. This manifold is commonly explored in the study of low-rank matrix completion and approximation problems \citep{vandereycken2013low,mishra2014fixed}. For example, when the observed data consists of noisy measurements of an underlying element in $\mc R_{r, p_1, p_2}$, a critical challenge lies in recovering the true underlying parameter and performing uncertainty quantification.

\end{itemize}
In short, viewing estimation and inference problems through the lens of extremum estimators on manifolds can reveal new insights, particularly in scenarios where traditional Euclidean frameworks fall short.

\subsection{Bootstrap and Studentization on Manifolds}
As the sample size increases in a fixed model, it is widely recognized that the distribution of the extremum estimator exhibits certain asymptotic properties. However, a significant challenge lies in the lack of direct access to the actual asymptotic distribution, unless the sample distribution is known \cite{brunel2023geodesically} or the underlying manifold possesses particular geometric structures \cite{lee2024huber}. To address this, a notable idea is to {\it studentize} the target statistic while resampling the data, namely, to divide the (resampled) statistic by a sample-based estimate of its standard deviation (e.g., $t$-statistic) to achieve a more accurate approximation of the statistic’s distribution. Over the past decades, computationally efficient bootstrap procedures have been explored extensively in works such as \cite{davidson1999bootstrap,andrews2002higher}. However, these multi-step estimators are not readily applicable to manifold settings, where unique geometrical structures introduce additional complexities. This limitation motivates our investigation into bootstrap procedures specifically tailored for manifold-based problems.

In this paper, we propose a fine-grained inference procedure to address the locality challenges in manifold settings, enabling the derivation of high-order accurate confidence regions through appropriate studentization techniques. To ensure high-order accuracy in the distributional approximation obtained via bootstrap, we carefully control each resampled term, incorporating the effects of the Riemannian curvature of the underlying manifold.

\subsection{Our Contributions}

The contributions of this paper include the following three points.

First, for a class of well-formalized M-estimators on Riemannian manifolds, we develop a bootstrapping algorithm that stands out for its computational efficiency and high-order accuracy in hypothesis testing and confidence region construction. 

Secondly, we discover that, while performing a locational hypothesis testing on a Riemannian manifold can be naturally transformed into a well-studied testing problem on a Euclidean space, constructing a high-order accurate confidence region poses additional challenges due to a manifold’s inherent geometrical structure. To ensure both the computational efficiency and the coverage accuracy of the proposed confidence interval construction methods, we establish high-order asymptotics of the specific statistics in our algorithms. To address the ambiguity of the coordinate expression, we adopt a fixed chart centered at the true parameter in our analysis. This approach allows us to expand the functions of interest through a Euclidean lens while accounting for the impact of manifold curvature when transitioning across different charts.

Finally, our methodology and theoretical results extend to a wide range of statistical settings. We specifically study the inference on spheres, Stiefel manifolds, fixed-rank matrix manifolds, and rank-one tensors manifolds using our developed framework. Notably, for Euclidean submanifolds, we introduce a class of projection-like coordinate charts equipped with the necessary consistency properties, providing a robust foundation for further applications.

\subsection{Preliminaries on Riemannian Manifolds}\label{sec:notation-preliminaries}

To ensure this paper is self-contained, we provide a concise introduction to Riemannian manifolds, particularly in the context of numerical optimization \citep{absil2009optimization}.

\paragraph*{Smooth Structure} A smooth manifold $\mathcal{M}$ of dimension $p$ is a connected paracompact Hausdorff space for which every point has a neighborhood $U$ that is homeomorphic to an open subset $\Omega$ of $\bb R^p$. Such a homeomorphism $\phi: U \rightarrow \Omega$ is called a (coordinate) chart. A tangent vector $X$ to the manifold $\mc M$ at $x$ is a mapping from the set of real-valued functions defined on $\mc M$ (i.e., $C^{\infty}(\mc M)$) to $\bb R$ such that there exists a smooth curve $\gamma(t)$ on $\mc M$ with $\gamma(0)= x$ and $X: f \mapsto \frac{\mathrm d}{\mathrm d t}f\circ \gamma(t) \big\vert_{t=  0} \eqqcolon X f $ for all $f\in C^{\infty}(\mc M)$. The tangent space $\mathrm T_x \mc M$ of a smooth manifold $\mc M$ is the linear space of all tangent vectors at $x$. The tangent bundle $\mathrm T \mc M$ is the disjoint union of tangent spaces over all points $x \in \mc M$, given by $\mathrm T\mc M \coloneqq \bigsqcup_{x\in \mc M}\mathrm T_x \mc M$.  A vector field on $\mc M$ is a smooth function that assigns a tangent vector $v \in \mathrm T_{\theta}\mc M$ to each point $\theta\in \mc M$. $\Gamma(\mc M)$ denotes the collection of possible vector fields on $\mc M$. Naturally, given a vector field $X\in \Gamma(\mc M)$ and a smooth function $f\in C^\infty(\mc M)$ on $\mc M$, $Xf$ is again a smooth function on $\mc M$. For further use, we define the mapping $XY:C^\infty(\mc M) \rightarrow \mathbb{R}$ for $X, Y\in \Gamma(\mc M)$ such that $XY f\mapsto X(Yf)$ for $f \in C^\infty(\mc M)$. Given two smooth manifolds $\mc M, \mc N$ and a smooth function $f: \mc N \rightarrow \mc M$, we define a mapping called the differential of $f$ at $x\in \mc N$ as a linear mapping $\mathrm T_x \mc N \rightarrow \mathrm T_{f(x)}\mc M$ such that $\mathrm d f(X) g= X(g\circ f)$ for every $X \in \mathrm T_x\mc N$ and every $g\in C^\infty(\mc M)$. Given a smooth curve $\gamma \subset \mc M$, we define the velocity of $\gamma$ at $s$, denoted by $\frac{\mathrm d}{\mathrm dt} \gamma \big\vert_{s}$ or $\dot \gamma(s)$, to be the tangent vector $\mathrm d \gamma (\frac{\mathrm d}{\mathrm dt}\big\vert_s)$, where $\frac{\mathrm d}{\mathrm dt}\big\vert_s$ is the standard coordinate basis in $\mathrm T_{s} \bb R$. 

\paragraph*{Riemannian Manifold} Further, a Riemannian manifold $\mc M$ is a smooth manifold equipped with a smooth positive-definite inner product $\langle \cdot, \cdot \rangle$ on the tangent space $\mathrm T_p\mc M$ at each point $p$.
With an inner product at each point, the length of a piece-wise smooth curve $ \gamma(t),t\in[0,1]$ from $x$ to $x'$ can be defined as 
$L( \gamma) = \int_{0}^1 \ip{\dot{ \gamma}(t), \dot{ \gamma}(t)}^{\frac{1}{2}}dt
$.
We also denote the distance of two points $x, x'$ on manifold $\mathcal{M}$ by $\dist(x,x')$, which is defined by the infimum of $L(\gamma)$ over all possible piece-wise smooth curves from $x$ to $x'$. Further, given two vector fields on $\mc M$, we note that $\ip{X, Y}$ is a smooth function on $\mc M$ by the smoothness of $\ip{\cdot, \cdot}$. 

\paragraph*{Levi-Civita Connection} To enable the comparison of tangent vectors across disparate tangent spaces, we introduce the notion of connection $\nabla: \Gamma(\mc M)\times \Gamma(\mc M) \rightarrow \Gamma(\mc M)$, as a generalization of the directional derivative on a vector field to Riemannian manifolds, satisfying that (i) $\nabla_{f_1 X_1 + f_2 X_2} Y = f_1\nabla_{X_1}Y + f_2\nabla_{X_2}Y $, (ii) $\nabla_{X}(\alpha_1 Y_1 + \alpha_2 Y_2) = \alpha_1 \nabla_XY_1 + \alpha_2 \nabla_X Y_2$, (iii) $\nabla_X(fY) = (Xf)Y + f\nabla_XY$ with $X, Y, X_1, X_2, Y_1, Y_2\in \Gamma(\mc M)$, $f_1, f_2, f\in C^\infty(\mc M)$, and $\alpha_1, \alpha_2\in \bb R$. Among all possible choices of affine connections, the Levi-Civita connection is the unique one that preserves the metric and is torsion-free; to be specific, an affine connection is called the Levi-Civita connection if (i) $X\ip{Y, Z} = \ip{\nabla_XY, Y} + \ip{Y, \nabla_X Z}$ holds for any vector fields $X, Y, Z$\footnote{
    Since $\ip{Y, Z}$ is a smooth function on $\mc M$, $X\ip{Y, Z}\coloneqq X(\ip{Y,Z})$ is also a smooth function given an additional tensor field $X$ on $\mc M$.}; (ii) $\nabla_X Y - \nabla_Y X - (XY - YX)= 0$ holds for any vector fields $X, Y, Z$ (see the details of the Levi-Civita connection's interpretation in \cite{jost2008riemannian} and \cite{lee2018introduction}). We \emph{only} consider the Levi-Civita connections throughout this paper. Given a smooth curve $\gamma\subset \mc M$ with velocity $\dot \gamma$, we denote $\nabla_{\dot \gamma(s)}\dot \gamma$ by $\frac{\mathrm d^2}{\mathrm dt^2} \gamma\big\vert_{s}$. 

\paragraph*{Parallel Transport}
Parallel transport is a local realization of a connection that intrinsically defines an isometry between tangent spaces. 
Letting $\mu \in \Gamma(\mc M)$ be a vector field on $\mc M$ and $\gamma(t)$ be a smooth curve, $\mu(\gamma(t))$ is called parallel transport of $\mu(\gamma(0))$ along $\gamma$ if $\nabla_{\dot{\mb \gamma}(t)}\mu = 0$ holds for arbitrary $t$. 
If the shortest curve connecting $x$ with $x'$ is unique, we denote the tangent vector at $x'$ from the parallel transport of $v \in \mathrm T_x\mc M$ in the tangent space of $x$ along the shortest path (also a geodesic by \cite[Theorem 6.4]{lee2018introduction}) by $T_{x \rightarrow x'}(v)$; in fact, Theorem 6.17 in \cite{lee2018introduction} ensures the uniqueness of the geodesic connecting two arbitrary points in a small neighborhood. This local viewpoint will frequently be used in the subsequent discussion.

\paragraph*{Retraction} The notion of retraction mapping is a crucial tool in developing efficient and practical optimization algorithms on manifolds. In much of the previous literature, a retraction mapping $R$ refers to a smooth mapping from $\mathrm T \mc M$ to $\mc M$. However, since in this paper we are mainly interested in quantities related to specific coordinate representations, we adopt a slightly refined notion. We define a retraction $R \coloneqq \{R_\theta\}_{\theta\in \mc M}$ as a collection of retraction mappings, where each retraction mapping $R_\theta$ centered at $\theta$ is a diffeomorphism from an open neighborhood $U \subset \mathbb{R}^p$ of $\mb 0$ to $\mathcal{M}$. The diffeomorphism satisfies the conditions: (i) $R_\theta(\mathbf{0}) = \theta$, (ii) $\{\mathrm dR_\theta \delta_i\}$ is an orthogonal basis of $\mathrm T_\theta\mc M$ for each $\theta\in \mc M$, where $\delta_i$ denotes the $i$-th canonical basis of $\bb R^p$, (iii) the mapping $R_\theta \circ \mathrm dR_\theta^{-1}( y)$, with respect to $\theta$ and $y$, is \emph{a smooth mapping from an open neighborhood of the zero section of $\mathrm T\mc M$ to $\mc M$}. The Euclidean realization of a retraction presented here also leads to the subsequent discussion on the necessity of smoothness in coordinate bases across different charts. To impose regularity conditions on retractions, we further introduce the notion of a second-order retraction, defined as follows.
\begin{definition}\label{definition:second order retraction}
A retraction $R$ is second-order if it holds for every $\theta\in \mc M$ and $\mb y$ in an open neighborhood of $\mb 0 \in \bb R^p$ that $\frac{\mathrm d^2}{\mathrm d t^2} R_\theta(t\mb y)\big\vert_{0} = 0$. 
\end{definition}

\paragraph*{Derivatives}
    Given a Riemannian manifold $\mathcal{M}$ equipped with the Levi-Civita connection $\nabla$ and a smooth function $f$ on $\mathcal{M}$, the $k$-th derivative $\nabla^k \subset \Gamma(\otimes_{i=1}^k \mathrm {T}\mathcal{M}^*)$\footnote{
    $\Gamma(\otimes_{i=1}^k \mathrm {T}\mathcal{M}^*)$ denotes the collection of smooth functions that assign an element of $\underbrace{\mathrm T_x\mc M^* \otimes \cdots \otimes \mathrm T_x\mc M^*}_{\text{$k$ times}}$ to each $x$ in $\mc M$ where $\mathrm T_x\mc M^*$ denotes the dual space of $\mathrm T_x\mc M$.  
    } of $f$ ($k \in \bb N^+$) is defined recursively as follows
    \begin{align}
        &\nabla f[X_1] = \nabla^1 f[X_1] = X_1 f,\\ 
        &\nabla^2 f[X_1,X_2] = X_1(\nabla^1 f[X_2]) - \nabla^1 f[\nabla_{X_1}X_2],\\
        &\nabla^k f[X_1, \ldots, X_k] = X_1\nabla^{k-1}f[X_2,\ldots, X_k]- \nabla^{k-1}f[\nabla_{X_1} X_2,\ldots, X_k] 
        \\&\quad -\nabla^{k-1}f[X_2,\nabla_{X_1} X_3,\ldots, X_k]- \cdots  - \nabla^{k-1}f[X_2,X_3,\ldots,\nabla_{X_1}X_k].
    \end{align}
Note that a retraction $R = \{R_x\}_{x\in \mc M}$ is a second-order retraction if and only if $\nabla^2 R_x = \mb  0$ for every $x\in \mc M$.

\begin{remark}
	In a way, a retraction typically suggests the topological correspondence between the derivative of mapping and the tangent space of the manifold but does not involve any metric structure defined on the manifold. However, a second-order retraction goes beyond the topological structure; namely, it takes the metric structure of a manifold into account. As a consequence, different metric structures on a manifold usually lead to different second-order retractions. 
\end{remark}

\subsection{Notation}\label{sec:notation}
 We denote the $i$-th canonical basis of $\bb R^p$ by $\delta_i$. Given a matrix $\mb A\in \bb R^{N\times N}$, we denote its Moore–Penrose inverse by $\mb A^{\dagger}$. The $\ell_2$-norm of a vector $\mb v$ is denoted by $\norm{\mb v}_2$, and the Frobenius norm of a matrix $\mb A \in \bb R^{n_1\times n_2}$ is denoted by $\norm{\mb A}_F$. Given a vector space $\mathcal T$ with norm $\norm{\cdot}$, we denote a sphere $\{y\in \mathcal T: \norm{y - x} \leq \rho\}$ centered at $x \in \mathcal T$ with radius $\rho$ by $B(x, \rho)$. 

Consider a Riemannian manifold $\mc M$ of dimension $p$. The norm of a tangent vector $v \in \mathrm T \mc M$ is written as $\norm{\mb v} \coloneqq \ip{v, v}^{\frac12}$.  
We denote by $\Exp_x(\cdot)$ and $\Log_x(\cdot)$ the exponential mapping and the logarithmic mapping at $x\in \mc M$. 
Given a basis $\{e_i\}$ of $\mathrm T \mc M$, we define the canonical mapping $\pi_{\{e_i\}}:\mathrm T_x\mc M \rightarrow \bb R^p$ at $x$ by $\pi_{\{e_i\}}(y) = \mb v$ with $y = \sum_{i=1}^p v_i e_i$ and $\mb v = (v_1, \ldots, v_p)\in \bb R^p$. 
For a function $f: \mc M \rightarrow \mathbb R$ defined on $\mc M$, we denote its composition with $R_{\theta}$ by $ f^{\ret(\theta)}$. 

Throughout this paper, coordinate expressions and abstract tensor field expressions will be used interchangeably. To avoid confusion, we will use the notation $\nabla^k f$ for $f: \mc M \rightarrow \bb R$ to denote the $k$-th derivative of $f$ in the abstract tensor field expression, while we use $\bar\nabla^{\mb \nu} f(\mb x)\coloneqq \frac{\partial^{\nu_1}}{\partial x_1^{\nu_1}}\cdots\frac{\partial^{\nu_p}}{\partial x_p^{\nu_p}}f(\mb x) $ for $f: \bb R^p \rightarrow \bb R$ to denote the partial derivative of $f$ related to $\mb \nu \in \bb N^p$ in the coordinate expression. Moreover, given a sufficiently smooth function $f: \bb R^p \rightarrow \bb R$, we denote the gradient column vector by $\bar \nabla f(\mb x) = (\frac{\partial }{\partial x_j}f(\mb x))_{j\in[p]}$ and the Hessian matrix by $\bar \nabla^2 f(\mb x) = (\frac{\partial^2 }{\partial x_i \partial x_j} f(\mb x ))_{i,j\in[p]} \in \bb R^{p\times p}$. Further, given a tangent basis $\{e_i\}$ of $\mathrm T_\theta\mc M$ at $\theta\in \mc M$ and a multi-index $\mb \nu \in \bb N^p$, we denote $$\nabla^{|\mb \nu|}f(\theta)[\underbrace{e_1,\ldots,e_1}_{\nu^{(1)}},\underbrace{e_2,\ldots,e_2}_{\nu^{(2)}}, \ldots, \underbrace{e_p,\ldots, e_p}_{\nu^{(p)}}]$$ by $\nabla^{\mb \nu}f(\theta)$, referred to as the $\mb \nu$-derivative of $f$ under the basis $\{e_i\}$; in the notion of differential geometry, this is also equivalent to $\frac{\partial^{v^{(1)}}}{\partial x_1^{v^{(1)}}}\cdots \frac{\partial^{v^{(p)}}}{\partial x_p^{v^{(p)}}}f\circ \Exp \circ \pi^{-1}_{\{e_i\}}$.

For an objective function $L(\theta, x): \mc M \times \mc S \rightarrow \bb R$ or a composition $L(R_\theta(\mb \eta),x)$ with a mapping $R_\theta(\mb \eta)$ from $\bb R^p$ to $\mc M$, we slightly abuse the notation by performing all differentiation operations with respect to the first argument, without explicitly indicating it.

\section{High-order Accurate Inference Using Extremum Estimators}\label{sec:extremum-estimators}

Inference with high-order accuracy typically utilizes the Edgeworth expansions for both the statistics of interest and their resampled counterparts. We first define the resampled objective function and its corresponding (heuristic) resampled M-estimator as
\begin{align}
	&L_n^*(\theta) \coloneqq   \sum_{i\in[n]}L(\theta, X_i^*) /n,\quad \tilde{\theta}_n^* = \argmin_{\theta\in \mc M} L_n^*(\theta),\label{eq: empirical argmin}
\end{align}
where $\mc X_n^* = \{X_i^*\}_{i\in[n]}$ is the resampled collection with replacement from $\mc X_n$. However, the minimizers in \eqref{eq: argmin} and \eqref{eq: empirical argmin} may suffer from non-uniqueness, a subtle but important geometric complication. To circumvent this, we instead consider arbitrary solutions in a neighborhood $U$ of $\theta_0$ to first-order conditions, which serve as relaxed formulations of the global optimizers: 
\begin{align}
& \tilde \theta_n \in \{\theta \in U:~\nabla L_n(\theta) = 0\} \label{eq: first-order condition of M estimation},\\ 
& \tilde \theta_n^* \in \{\theta \in U:~\nabla L_n^*(\theta) = 0 \}\label{eq: resampled first-order condition of M estimation}.
\end{align}

Despite the validity of the high-order asymptotics related to the exact solution to the stationary conditions, the explicit solution is often infeasible. To address this, we employ the Riemannian Newton method to proceed with estimation and inference. 

\subsection{Riemannian Newton Iteration}\label{section: newton's method}

The classic Newton's method \citep{dennis1996numerical} is an iterative procedure that incorporates second-order information to seek a critical point of the cost function. It is well-established that in deterministic settings, Newton's method achieves a quadratic rate of convergence if an initial point is close enough to a minimizer with a nonsingular second derivative in Euclidean settings \citep{dennis1996numerical} as well as in general Riemannian manifold settings \citep{absil2009optimization}. In the former case, letting $\mb \theta_{k-1}$ be the iterative point at the $(k-1)$-th step and $f$ be an objective function, Newton's update is written as 
\eq{
\mb \theta_{k} \coloneqq    \mb \theta_{k-1} - \big[\bar\nabla^2 f({\mb\theta}_{k-1}) \big]^{-1} \bar \nabla f({\mb\theta}_{k-1}) \in \bb R^p.
}
Given a Riemannian manifold $\mathcal{M}$, Newton's update\footnote{The Moore–Penrose inverse is introduced to address the potential singularity issue of the Hessian matrix.} takes form
    \eq{\label{manifoldestimator}
    {\theta}_k\coloneqq    R_{{\theta}_{k-1}}\left( - \big[\bar \nabla^2 f\circ R_{{\theta}_{k-1}}\left(\mb 0 \right)\big]^{\dagger }\bar \nabla f\circ R_{{\theta}_{k-1}}\left(\mb 0 \right)\right) \in \mc M, 
    }
    where $R$ is a retraction defined in Definition 1 and $f: \mc M\rightarrow \bb R$ is an objective function. 
    
    In what follows, we consider an objective function $L$ composed with a random data set $\mc Y$ and introduce a Riemannian Newton iteration procedure as described in Algorithm~\ref{algorithm: update}\footnote{One may notice that the domain of $R_\theta$ is an open neighborhood of $\mb 0$, but not necessarily $\bb R^p$. To ensure the well-definedness, we simply let $R_{\theta}(\mb y)$ be $\theta$ for any $\mb y$ outside the domain of $R_{\theta}$. }, primarily for constructing confidence regions in subsequent steps. A key difference between optimization on a Riemannian manifold and in Euclidean space is that there is no fixed coordinate representation (or coordinate chart) for the manifold. Instead, at each iteration of Algorithm~\ref{algorithm: update}, we operate on the local chart $R^{-1}_{\theta^{(s-1)}}(\cdot)$, centered at the current estimate $\theta^{(s-1)}$. 
        \begin{algorithm}[htp]
        \caption{Riemannian Newton Iteration on a Manifold 
        }\label{algorithm: update}
        \KwIn{Initial point ${ \theta}^{\text{initial}}$ and data collection $\mathcal Y$, iteration times $t$}
            \KwOut{$\theta^{(t)}$}
            Let $\theta^{(0)} ={\theta}^{\text{initial}}$
        \For{$s = 1, \ldots, t$}{
        Update: $$
        \theta^{(s)} = R_{\theta^{(s-1)}}\Big( - \big[\bar \nabla^2 \big(\sum_{ X\in \mc Y} L( R_{ \theta^{(s-1)}}(\cdot), X) \big)\big|_{\mb 0}/ |\mc Y|\big]^{\dagger}\big[\bar \nabla \big(\sum_{ X\in \mc Y} L( R_{ \theta^{(s-1)}}(\cdot), X) \big)\big|_{\mb 0}/ |\mc Y|\big]\Big). 
        $$
        }
    \end{algorithm}
   
    \begin{remark}
    As discussed in \cite{adler2002newton,absil2009optimization}, local quadratic convergence is maintained even if Algorithm~\ref{algorithm: update} employs only a first-order retraction. However, as we will demonstrate in Section~\ref{subsection: curvature effect}, a second-order retraction is essential to explicitly incorporate the high-order asymptotics of the constructed statistics, thereby enabling high-order accurate inference.
    \end{remark}

	Since the empirical solution to \eqref{eq: first-order condition of M estimation} is likely to stay close to the resampled solution to \eqref{eq: resampled first-order condition of M estimation} with high probability, it naturally serves as an ideal initialization. We therefore propose to integrate the estimation with the resampling procedures, allowing the stochastic optimization process to produce a sequence of replicates based on a similar underlying mechanism as the original estimate. Our procedure generates the estimate $\hat\theta_n$ and the resampled estimates $\{\theta_{n,i}^*\}_{i\in [b]}$, as summarized in Algorithm~\ref{algorithm: resampled newton iteration}.
    
    \begin{algorithm}[htp]
        \caption{Resampled Riemannian Newton Iteration on a Manifold}
        \label{algorithm: resampled newton iteration}
        \KwIn{Initial estimate $\hat{ {\theta}}_n^{(0)}$; number of refinement steps $\mathsf{burn}$; number of resamplings $b$. }
        \KwOut{Approximate extremum estimate $\hat \theta_n$ and the resampled approximate extremum estimate series $\{\hat \theta_{n}^{*[i]}\}_{i=1}^{b}$}
        Implement Algorithm~\ref{algorithm: update} with initial point $\theta^{\text{initial}} = \hat{\theta}_{n}^{(0)} $, $t=\mathsf{burn}$, and the dataset $\mc X_n$ to obtain $\hat{ \theta}_n$. \\
        \For{$i = 1, \ldots, b$}{
        Resample the dataset $\mc X_n$ to obtain $\mc X_n^{[i]}$ (i.e., draw $n$ samples with replacement from $\mc X_n$ to form $\mc X_n^{[i]}$), and implement Algorithm~\ref{algorithm: update} with \emph{the initial point $\hat{\theta}_n$, $t=2$, and the dataset $\mathcal X_n^{[i]}$} to obtain $\hat{\theta}_{n}^{*[i]}$.
        }
      \end{algorithm}

At the beginning, we apply a $\mathsf{burn}$-step Newton to compute an approximation $\hat \theta_n$ to $\tilde\theta_n$. 
In the bootstrapping stage of Algorithm~\ref{algorithm: resampled newton iteration}, we perform a two-step Riemannian Newton initialized from the point estimate $\hat \theta_n$, provided the resampled dataset $\mc X_n^{[i]}$. Heuristically speaking, both $\hat \theta_n$ and $\tilde\theta_n$ are likely to lie within a radius-$O(n^{-\frac{1}{2}}(\log n)^{\frac{1}{2}})$ spherical neighborhood around $\theta_0$ with high probability, and $\hat \theta_n$ usually serves as a nice approximation of $\tilde\theta_n$. Provided the quadratic convergence rate of Newton's method, the two-step update is expected to provide a good approximation of the solution to the first-order condition given the dataset $\mc X^{[i]}_n$, which is further confirmed by our high-order asymptotic analysis in Section~\ref{sec: theoretical guarantees}.
 
\subsection{Hypothesis Testing on Manifold}\label{subsection: hypothesis testing}

We consider a hypothesis testing problem on a Riemannian manifold, specifically focusing on the following locational testing problem:
\eq{
H_0: \theta_0 = \theta_1 \in \mc M \quad \text{ versus }\quad H_1: \theta_0 \neq \theta_1, \quad \text{with significance level $\alpha$.}
}
We focus on controlling the test rejection probability when the null hypothesis holds. To achieve this, we note that by restricting attention to a neighborhood $U\subseteq \mc M$ of $\theta_1$, where $(\phi, U)$ is a smooth chart with $\phi(\theta_1) = \mb 0\in \bb R^p$, the testing problem reduces to 
\eq{
H_0: \phi(\theta_0)= \mb 0 \text{\quad versus\quad }H_1: \phi(\theta_0) \neq \mb 0, \quad \text{with significance level $\alpha$.}
}
In light of this, a natural approach is to employ an alternative loss function $\sum_{i\in[n]}L(\phi^{-1}(\cdot), X_i) / n: \phi( U) \rightarrow \bb R$. 
Moreover, we define the $t$-statistics $T_{\phi,j}$, $j \in [p]$, the Wald statistic $W_\phi$, and their bootstrap alternatives $T_{\phi,j}^*$, $W_\phi^*$ as follows\footnote{In the case of singular issues, such as $(\mathbf{\Sigma}_\phi)_{j,j} = 0$, $\mathbf{\Sigma}_\phi$ being non-invertible, or $\hat{\theta}_n$ lying outside the domain of $\phi$, the corresponding statistic is set to zero. }:
	\begin{align}
		& T_{\phi,j} = \frac{  \phi(\hat \theta_n )_j}{(\mb \Sigma_\phi)_{j,j}^{1/2}}, \quad T_{\phi,j}^* = \frac{  \phi(\hat \theta_n^*)_j -   \phi(\hat \theta_n)_j}{(\mb   \Sigma_\phi^*)_{j,j}^{1/2}}, \quad \text{for $j\in [p]$}, \label{eq: t stat for testing}
        \\ 
		& W_\phi =   \phi (\hat\theta_n)\t    \mb \Sigma_{\phi}^{-1}   \phi (\hat\theta_n), \quad W_\phi^* = \big(  \phi (\hat\theta_n^*) -   \phi (\hat\theta_n)\big)\t  {{}  \mb \Sigma_{\phi}^*}^{\dagger} \big(  \phi (\hat\theta_n^*) -   \phi (\hat\theta_n)\big). \label{eq: wald stat for testing}
	\end{align}
	Here $\hat \theta_n^*$ denotes a generic $\hat{\theta}_{n,1}^*$ in Algorithm~\ref{algorithm: resampled newton iteration} and $\mc X^*_n = \{X_i^*\}_{i\in[n]}$ denotes a generic resampled dataset $\mc X^{[1]}_n$. The matrices $\mb \Sigma_\phi$ and ${\mb \Sigma}_{\phi}^*$ are defined as 
	\begin{align}
		&  \mb \Sigma_\phi \coloneqq \left(\bar \nabla^2 L_n (\phi^{-1}(\cdot ) )\big|_{ \phi(\hat \theta_n)}\right)^{-1}\\
  & \qquad \Big( \frac{1}{n}\sum_{i\in[n]}\bar \nabla L(\phi^{-1}(\cdot), X_i)\big|_{ \phi(\hat \theta_n)}\bar \nabla L(\phi^{-1}(\cdot), X_i)\big|_{ \phi(\hat \theta_n)}\t \Big) \left(\bar \nabla^2 L_n (\phi^{-1}(\cdot ) )\big|_{ \phi(\hat \theta_n)}\right)^{-1},\\
		&   \mb \Sigma_\phi^* \coloneqq \left(\bar \nabla^2 L_n^* (  \phi^{-1}(\cdot) )\big|_{ \phi(\hat \theta_n^*)}\right)^{-1}\\ 
  & \qquad \Big( \frac{1}{n}\sum_{i\in[n]}\bar \nabla L ( \phi^{-1}(\cdot ) , X_i^*)\big|_{\phi(\hat \theta_n^*)}\bar \nabla L ( \phi^{-1}(\cdot ) , X_i^*)\big|_{\phi(\hat \theta_n^*)}\t \Big) \left(\bar \nabla^2 L_n^* (  \phi^{-1}(\cdot) )\big|_{ \phi(\hat \theta_n^*)}\right)^{-1}.
	\end{align}
Note that the coordinate derivative operations, $\bar{\nabla}$ and $\bar{\nabla}^2$, refer specifically to the first argument of $L$, i.e., the coordinates under the retraction mapping $\phi^{-1}$.

Intuitively, $\phi(\hat\theta_n)$ and $  \phi(\hat\theta_n^*)$ are close to the exact solutions $\phi(\tilde \theta_n)$ and $ \phi(\tilde \theta_n^*)$, where $\tilde \theta_n$ and $\tilde \theta_n^*$ are defined in \eqref{eq: first-order condition of M estimation} and \eqref{eq: resampled first-order condition of M estimation}. In view of the delta method (Lemma~\ref{lemma:deltamethod} in the Supplementary Material), it thus suffices to focus on the high-order asymptotics of the latter ones.
	Moreover, a useful observation is that $ \phi(\tilde \theta_n)$ and $  \phi(\tilde \theta_n^*)$ are also the solutions to 
	\begin{align}
		& \bar \nabla   L_n(\phi^{-1}(\mb x)) = 0 ~~ \text{and}~~
		\bar \nabla   L_n^*(\phi^{-1}(\mb x)) = 0,
	\end{align}
respectively, where $L_n(\phi^{-1}(\cdot))$ could be viewed as an empirical loss function defined on a neighborhood of $\bb R^p$. Since well-established high-order asymptotic results exist for M-estimation \citep{hall1996bootstrap,andrews2002higher} for samples in a Euclidean space, we restrict our focus to $\mc S = \bb R^m$ so that these can be applied in our context. We have the following formal result, where the assumptions mainly stem from those of \cite[Theorem~2]{andrews2002higher}.
	\begin{proposition}\label{proposition: hypothesis testing}
		Suppose the following assumptions hold:
		\begin{enumerate}
        \item  The zero vector is the unique solution in a sufficiently small spherical neighborhood $U$ of $\theta_0$ to $\bb E\big[\bar \nabla L(\phi^{-1}(\cdot), X_1)\big] = \mb 0$. There exists a function $C_g(\mb x)$ such that $\norm{\bar \nabla L(\phi^{-1}(\mb \eta_1),\mb  x) - \bar \nabla L(\phi^{-1}(\mb \eta_2),\mb  x)}_2 \leq C_g(\mb x)\norm{\mb \eta_1 - \mb \eta_2}_2$ for all $\mb x$ in the support of $ X_1$ and all $\mb \eta_1, \mb \eta_2 \in \phi(U)$. $\bb E[C_g^{q_1}(X_1)] < \infty$ and $\bb E\norm{\bar \nabla L(\phi^{-1}(\mb \eta), X_1)}_2^{q_1}< \infty $ for all $0< q_1 < \infty$. 
        \item The $p$-by-$p$ Hessian matrix $\bb E[\bar\nabla^2 L(\phi^{-1}(\mb 0), X_1)]$ is of full rank $p$. $\bar \nabla  L(\phi^{-1}(\mb \eta), \mb x)$ is nine-times differentiable with respect to $\mb \eta$ on $\phi(U)$ for all $\mb x$ in the support of $X_1$. Let $f(\mb \eta,\mb x)$ be a vector containing $\bar \nabla L(\phi^{-1}(\mb \eta),\mb x)$, $\bar \nabla L(\phi^{-1}(\mb \eta),\mb x)\bar \nabla L(\phi^{-1}(\mb \eta),\mb x)\t$, and their derivatives through order $5$. There is a function $C_{\partial f}(\mb x)$ such that $\norm{ \bar\nabla^{\mb \nu} f(\mb 0, \mb x) -  \bar\nabla^{\mb \nu} f(\mb \eta, \mb x)}_2 \leq C_{\partial f}(\mb x) \norm{\mb \eta}_2$ for all $\mb \eta \in \phi(U)$, all $\mb \nu\in \bb N^p$ with $\sum_{i\in[p]}v_i \leq 4$, and all $\mb x$ in the support of $X_1$. $\bb E\big[ C_{\partial f}(X_1)^{q_2}\big] < \infty$ and $\bb E\big[ \big(\bar \nabla^{\mb \nu} f(\mb 0,  X_1)\big)^{q_2}\big] < \infty$ for all $q_2 >0 $ and all $\mb \nu$ with $\sum_{i\in[p]} \nu_i  \leq 4$. 
        The function $f(\mb 0, \mb x)$  is once differentiable with respect to $\mb x$ with uniformly continuous first derivative. 
        \item The random vector $f(\mb 0, X_1)$ satisfies Cram\'{e}r's condition, that is, \eq{\limsup_{\norm{\mb t} \rightarrow \infty} \big|\bb E\big[e^{i\mb t\t \mb f(\mb 0, X_1)}\big]\big|<1.} 
		\end{enumerate}
		If $\theta_0 = \theta_1$ and the initial estimate satisfies $\dist(\theta^{\text{initial}},\theta_0) \leq c_0$ for some sufficiently small $c_0$, then it holds for every $\xi \in (0,\frac{1}{2}]$ that  
		\begin{align}
			&  \bb P_{\mc X_n}\big[T_{\phi,j} >  z^*_{T_{\phi,j},2, \alpha}\big]  = \alpha + o(n^{-1 + \xi}),\\ 
            & \bb P_{\mc X_n}\big[T_{\phi,j} <  z^*_{T_{\phi,j},2, 1-\alpha}\big] = \alpha + o(n^{-1 + \xi}), \\ 
			& \bb P_{\mc X_n} \big[W_\phi > z^*_{W_{\phi},2,\alpha}\big] = \alpha + o(n^{-\frac{3}{2} + \xi}),
		\end{align}
    where $z^*_{T_{\phi,j},k, \alpha}$ is defined as a value that minimizes $\big|\bb P\big[T_{\phi,j}^* \leq z \mid \mc X_n\big] - (1- \alpha) \big|$ over $z\in \bb R$ and $z^*_{W_\phi,2,\alpha}$ is defined analogously. 
	\end{proposition}

Proposition~\ref{proposition: hypothesis testing} indicates that by repeatedly resampling the data, we can obtain a consistent estimate of \( z^*_{T_{\phi,j},2,\alpha} \) or \( z^*_{W_\phi,2,\alpha} \) to serve as a critical value. According to the Glivenko–Cantelli theorem \cite[Theorem 2.4.9]{durrett2019probability}, this approach yields a more accurate approximation of the distribution of \( T_{\phi,j} \) or \( W_\phi \), ensuring that the actual rejection probability closely aligns with the intended significance level \(\alpha\).

\subsection{Confidence Region Constructions on Manifolds}
In this section, we delve into the construction of confidence regions on Riemannian manifolds. A natural approach, drawing on classic ideas from statistics, would be to convert the hypothesis test introduced in Section \ref{subsection: hypothesis testing} into a confidence region. However, this approach may be computationally infeasible since we no longer have access to the presumed location information of $\theta_0$, which previously allowed for a convenient reduction from the manifold setting to a Euclidean one. Additionally, the choice of coordinate representations introduces further randomness, making the construction of confidence regions notably distinct from that of hypothesis testing. 
\subsubsection{Algorithms for Confidence Region Constructions}
\label{subsubsec: algorithms for confidence region constructions}
Selecting appropriate coordinate representations is essential for performing inference on Riemannian manifolds. For hypothesis testing, it is natural to adopt a fixed coordinate system at $\theta_0$. However, this approach no longer applies to confidence interval construction, necessitating transitions between different charts. These coordinate transformations are inherently nonlinear and influenced by the curvature of the local region. Consequently, when resampling the data and constructing studentized statistics, curvature effects may cause the distribution function of the studentized statistics for $\hat{\theta}_{n}^{*[i]}$ to deviate from being second-order consistent with that related to $\hat{\theta}_n$. 

\paragraph*{Algorithm for Wald-Type Statistic}\label{sec:Wald}
The Wald-type statistic is often advantageous for constructing an elliptically shaped confidence region within a fixed chart or directly on the parameter manifold $\mathcal{M}$. By leveraging the information gathered during the resampling process, we can approximate the distribution of the following Wald-type statistic:	
\eq{\label{eq: Wald Statistic}
    W \coloneqq  R_{\hat \theta_n}^{-1}( \theta_0)\t \hat{\mb \Sigma}^\dagger R^{-1}_{\hat \theta_n}( \theta_0),
	}
	where
	\eq{\label{eq: hatSigma}
	\hat{\mb \Sigma}\coloneqq \big(\bar \nabla^2 L_n (R_{\hat\theta_n}(\cdot ))\big|_{\mb 0}\big)^{-1}\Big(\sum_{j\in[n]} \bar \nabla L(R_{\hat \theta_n}(\cdot),X_j)\big|_{\mb 0}\bar \nabla L(R_{\hat \theta_n}(\cdot),X_j)\big|_{\mb 0}\t /n \Big)\big(\bar \nabla^2 L_n (R_{\hat\theta_n}(\cdot ))\big|_{\mb 0}\big)^{-1}.
	}
    As later elaborated in Section~\ref{subsection: curvature effect}, one can show that $- R^{-1}_{\theta_0}(\hat\theta_n)$ is approximately equal to $ R^{-1}_{\hat\theta_n}(\theta_0)$ up to a rotation as noted in \eqref{eq: Wald Statistic}, and a natural approximation of $R^{-1}_{\hat \theta_n}(\theta_0)$ is $R_{\hat \theta_n^{*[i]}}^{-1}(\hat \theta_n)$, which is closely approximated by $R_{\theta_0}^{-1}(\hat \theta_n) - R_{\theta_0}^{-1}(\hat \theta_n^{*[i]})$. Therefore, when working with the $i$-th resampled dataset $\mc X^{[i]}_n= \{X_j^{[i]}\}_{j\in[n]}$ in Algorithm~\ref{algorithm: resampled newton iteration}, replacing $R_{\hat \theta_n}(\theta_0)$ and $\hat{\mb \Sigma}$ with $R_{\hat \theta_n^{*[i]}}(\hat \theta_n)$ and $\check{\mb \Sigma}^{[i]}$, respectively, yields the counterpart resampled statistic $W^{*[i]}$ to $W$ as follows: 
    \eq{\label{eq: empirical Wald statistic}
     W^{*[i]}\coloneqq 
    R_{\hat \theta_{n}^{*[i]}}^{-1}\big(\hat \theta_n\big)\t 
    {{}\check {\mb \Sigma}^{[i]}}^{\dagger}
    R_{\hat \theta_{n}^{*[i]}}^{-1}\big(\hat \theta_n\big), 
    }
    where $\check{\mb \Sigma}^{[i]}$ is defined as 
    \eq{\label{eq: check Sigma}
    \check{\mb \Sigma}^{[i]} \coloneqq   
    \big(\bar \nabla^2 L_n^{[i]}( R_{\hat \theta_{n}^{*[i]} }(\cdot))\big|_{\mb 0}\big)^{-1}
    \big( \sum_{j\in[n]} \bar \nabla L( R_{\hat \theta_{n}^{*[i]}} (\cdot),X^{[i]}_j)\big|_{\mb 0}\bar \nabla L( R_{\hat \theta_{n}^{*[i]}} (\cdot),X^{[i]}_j)\big|_{\mb 0}\t/ n \big)
        \big(\bar \nabla^2 L_n^{[i]}( R_{\hat \theta_{n}^{*[i]} }(\cdot))\big|_{\mb 0}\big)^{-1}
    }
    with the shorthand notation $L_n^{[i]}$ for $\sum_{j\in[n]} L(\cdot, X_j^{[i]})/n$ provided the $i$-th resampled dataset $\mc X_n^{[i]} = \{X_j^{[i]}\}_{j\in[n]}$.
    
    The algorithm is summarized in Algorithm~\ref{algorithm:wald based statistic}\footnote{Similar to the handling of the locality of $R_\theta$ in Algorithm~\ref{algorithm: resampled newton iteration}, if $\hat \theta_n^{*[i]}$ is outside the domain of $R^{-1}_{\hat \theta_n}$, we set the corresponding quantity $W^{*[i]}$ to be zero. Similar rules apply to the quantity $T^{*[i]}$ in Algorithm~\ref{algorithm:t based statistic}. }. Heuristically, since the statistic of interest $W$ and its bootstrapped counterpart are largely unaffected by the choice of chart, as discussed in Section~\ref{subsection: curvature effect}, the conditional distribution of $W^{*[i]}$ given $\mc X_n$ is expected to closely approximate the distribution of $W$ with high probability. This will be further justified in Section~\ref{sec: theoretical guarantees}. 
    
    \begin{algorithm}[htp]
        \caption{Approximate Resampled M-Estimates on a Manifold with Wald-based Statistic}\label{algorithm:wald based statistic}
        \KwIn{Approximate extremum estimate $\hat \theta_n$ and the resampled approximated extremum estimate series $\{\hat \theta_{n}^{*[i]}\}_{i=1}^{b}$ derived from Algorithm~\ref{algorithm: resampled newton iteration}}
        \KwOut{confidence region for $\theta_0$ on the chart $R_{\hat{\theta}_n}$}
   Compute the $\alpha$-level quantile 
   $w_{\alpha} \coloneqq \argmin_{x\in \bb R}\big| \sum_{i\in[b]}1\{ W^{*[i]} \leq x\} / b - (1-\alpha)\big|$ with $W^{*[i]}=
    R_{\hat \theta_{n}^{*[i]}}^{-1}\big(\hat \theta_n\big)\t 
    {{}\check {\mb \Sigma}^{[i]}}^{\dagger}
    R_{\hat \theta_{n}^{*[i]}}^{-1}\big(\hat \theta_n\big)
   $\;
        Construct the $(1-\alpha)$-level confidence region:
        \eq{
                    I_{1-\alpha}^{\mathsf{Wald}} = \left\{\theta\in \mathrm{range}(R_{\hat \theta_n}): R^{-1}_{\hat \theta_n}(\theta)\t \hat{\mb \Sigma}^{\dagger} R^{-1}_{\hat \theta_n}(\theta)\leq w_{\alpha}\right\} 
                 }
        with $\hat{\mb \Sigma}$ defined in \eqref{eq: hatSigma}. 
    \end{algorithm}

\paragraph*{Algorithm for Intrinsic $t$-statistics}

We now move on to investigate how to quantify the uncertainty given a single coordinate of $\theta_0$ under the chart $R_{\hat \theta_n}$. Before presenting our algorithm for constructing confidence intervals on a single coordinate, here we would like to discuss further the insight into why inference on a single coordinate of the true parameter is an even more difficult task beyond establishing certain asymptotics in Euclidean cases.

\subparagraph*{Intractability of Local Charts on Riemannian Geometry}

A natural approach to inference is to establish specific asymptotics of the targeted quantities under a specific coordinate chart $\tilde R$, namely, the limiting distribution of $\tilde R(\theta_0)$ where $\tilde{R}$ represents a coordinate chart centered at some location $\theta'$ that could depend on the data. However, this local Euclidean viewpoint introduces some extra intractability, as the basis $\{\tilde e_i\}$ induced by the chart $\tilde R$ we use is usually related to the samples themselves, which are intertwined with some complicated dependency, and consequently deter the expected limit distribution. Notably, the cooperated dependency varies across different manifolds, making it challenging to derive a general inference procedure for a single coordinate. 

\begin{figure}
		\centering
		\includegraphics[width = 0.5\linewidth]{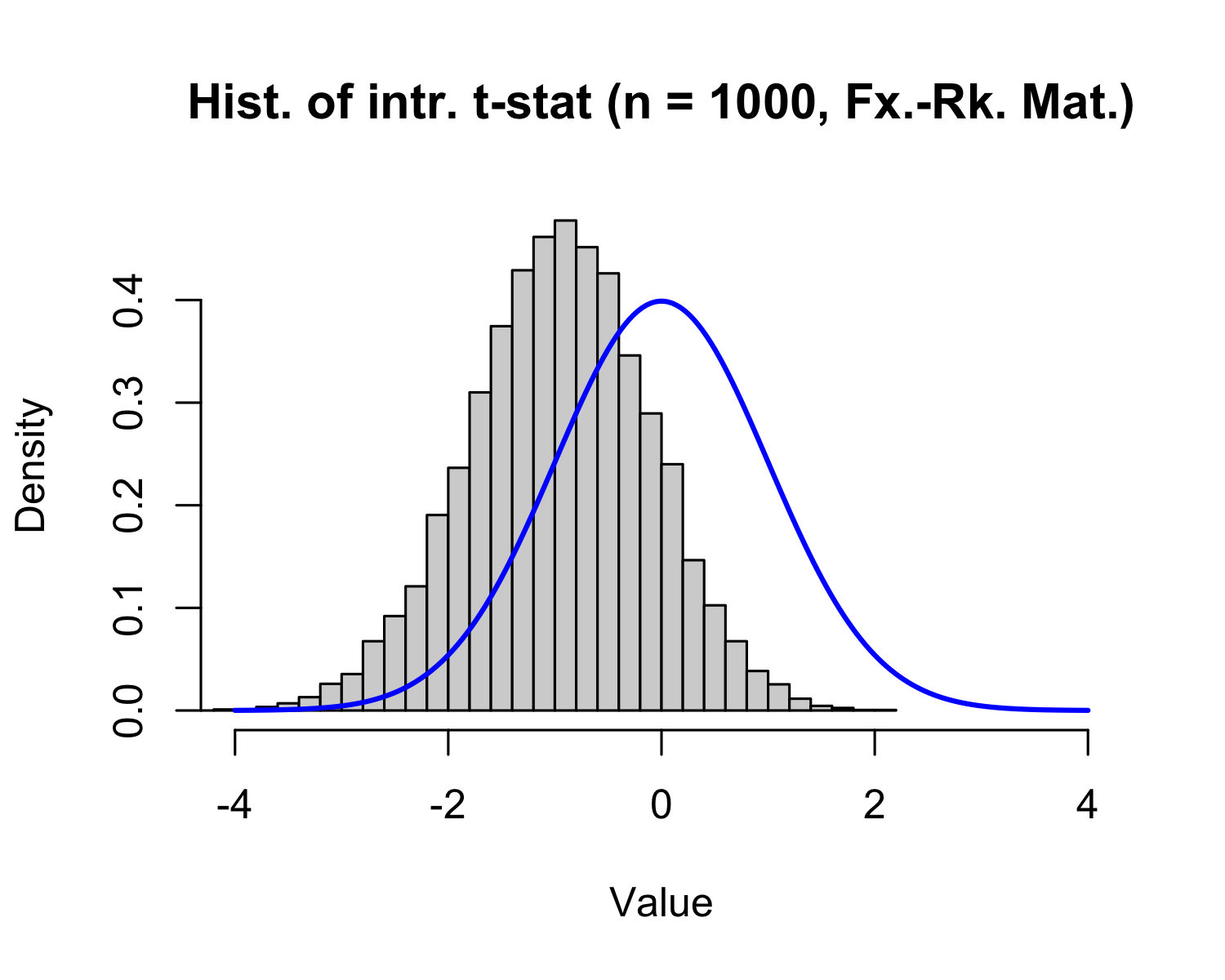}
		\caption{A Counterexample: Irregular Empirical Distribution Using an Irregular Chart. We consider the fixed-rank matrix manifold with the sample size $n = 1000$ and the ground truth parameter $\boldsymbol \theta_0 = \textrm{diag}(1,1,0,0)$. The histogram represents the studentized version of $\big(R_{\hat\theta_n}^{-1}(\theta_0)\big)_1$ while the blue curve depicts the probability density function of the standard normal distribution. }
		\label{fig: counterexample}
	\end{figure}

Consider, for instance, the fixed-rank matrix manifold and the SVD-based retraction in Example \ref{example: fixed-rank matrices manifold}, which is elaborated on later. Suppose that all the nonzero singular values of the ground truth $\theta_0$ are identical. In this scenario, the vector field induced by the first coordinate\footnote{The first coordinate refers to the top left entry of $\mb A$ for a tangent vector $\mb Z = \mb U\left[\begin{matrix}
            \mb A & \mb C\\ 
            \mb B& \mb 0
        \end{matrix}\right]\mb V\t$.} of the charts is not even differentiable. This non-differentiability results in the limiting distribution of the studentized version of $\big(R_{\hat\theta_n}^{-1}(\theta_0)\big)_1$ deviating significantly from the standard Gaussian distribution, even with a large sample size, as shown in Figure~\ref{fig: counterexample}.

\subparagraph*{A Differentiable Chart Helps} However, if specific smoothness conditions, for instance, the additional assumption in Theorem~\ref{theorem:bootstrapmanifold}.\ref{item: thm1.2}, on charts locally around $\theta_0$ are imposed, we are able to construct confidence intervals with respect to a linear combination of coordinates on the chart centered at $\hat \theta_n$. In particular, given a fixed vector $\mb a\in \bb R^p$, we aim to approximate the distribution of the studentized $t$-statistic 
\eq{
T = \frac{\mb a\t R^{-1}_{\hat \theta_n}(\theta_0)}{\big( \mb a\t \hat{\mb \Sigma}\mb a\big)^{\frac12}}. 
}
Similar to the construction of the resampled Wald statistic, we consider the counterpart resampled statistic $T^{*[i]}$: 
\eq{
T^{*[i]} = \frac{\mb a\t R^{-1}_{\hat \theta_n^{*[i]}}(\hat \theta_n)}{\big( \mb a\t \check{\mb \Sigma}^{[i]}\mb a\big)^{\frac12}}. 
}
The algorithm is outlined in Algorithm~\ref{algorithm:t based statistic}.

      \begin{algorithm}[htbp]
        \caption{Approximate Resampled M-Estimates on a Riemannian Manifold with Intrinsic $t$-statistic}
        \label{algorithm:t based statistic}
        \KwIn{Direction $\mb a$, the approximate extremum estimate $\hat \theta_n$, and the resampled approximate extremum estimate series $\{\hat \theta_{n}^{*[i]}\}_{i=1}^{b}$ derived from Algorithm~\ref{algorithm: resampled newton iteration}}
        \KwOut{Confidence region for $\theta_0$ on the chart $R_{\hat{\theta}_n}$ along the direction $\mb a$}
    Compute the $(1 - \alpha)$-quantile $w_{1 - \alpha} \coloneqq \argmin_{x\in \bb R}\big| \sum_{i\in[b]}1\{ T^{*[i]} \leq x\} / b - \alpha\big|$ with $T^{*[i]}=
    \frac{\mb a\t R_{\hat \theta_{n}^{*[i]}}^{-1}(\hat \theta_n) }{
    \big(\mb a\t {{}\check {\mb \Sigma}^{[i]}}^{\dagger} \mb a\big)^{1/2}}
   $\; 
        Construct the $(1-\alpha)$-level confidence region:
        \eq{I_{1-\alpha}^{\mathsf{Intr.t}} =\left\{ \theta \in \mathrm{range}(R_{\hat \theta_n}): \mb a\t R^{-1}_{\hat \theta_n}(\theta) \in  \big(-\infty, - w_{1 - \alpha} \cdot  \big(\mb a\t \hat {\mb \Sigma}\mb a\big)^{\frac{1}{2}}\big]  \right\}
        }
        with $\hat{\mb \Sigma}$ defined in \eqref{eq: hatSigma}.
    \end{algorithm}

\paragraph*{Algorithm for Statistics Related to Extrinsic Coordinates}

For problems involving a submanifold embedded in the Euclidean space, our framework can also accommodate a broad range of extrinsic representations of the true parameter. Specifically, let $\mc M$ be a submanifold embedded in $\bb R^{d}$. Suppose that the true parameter is $ \theta_0 \in \mc M \subset \bb R^{d}$ and $f:\bb R^{d} \rightarrow \bb R$ is a smooth function with $\nabla f\circ R_{\theta_0}\neq  0$. Under these conditions, Algorithm~\ref{algorithm:extrinsic statistic} yields a second-order accurate confidence interval for $f(\theta_0)$. 

\begin{algorithm}[htbp]
    \caption{Approximate Resampled M-Estimates on a Riemannian Manifold with Extrinsic $t$-Statistic}
    \label{algorithm:extrinsic statistic}    
\KwIn{Smooth function $f$, the approximate extremum estimate $\hat \theta_n$, and the resampled approximate extremum estimate series $\{\hat \theta_{n}^{*[i]}\}_{i=1}^{b}$ obtained from Algorithm~\ref{algorithm: resampled newton iteration}}
    \KwOut{Confidence interval for $f(\theta_0)$}
Compute the $(1 - \alpha)$-quantile $w_{1 - \alpha} \coloneqq \argmin_{x\in \bb R}\big| \sum_{i\in[b]}1\{ T^{*[i]} \leq x\} / b - (1-\alpha)\big|$ with 
$
T^{*[i]}=
\big(f(\hat{\theta}_{n}^{*[i]})- f(\hat{\theta}_n)\big) /  \sqrt{\big(\bar \nabla f(R_{\hat \theta_n^{*[i]}}(\cdot))\mid_{\mb 0}\t  
    \check {\mb \Sigma}^{[i]}\bar \nabla f(R_{\hat \theta_n^{*[i]}}(\cdot))\mid_{\mb 0} \big)^{\frac{1}{2}}}
$\; 

    Construct the $(1-\alpha)$-level confidence interval for $f(\theta_0)$: 
    \eq{
        I_{1-\alpha}^{\mathsf{Extr.t}} = \Big(-\infty, f(\hat{\theta}_n) - w_{\alpha}\cdot  \big(\bar \nabla f(R_{\hat{\theta}_{n}}(\cdot))\mid_{\mb 0}\t 
        \hat {\mb \Sigma}\bar \nabla f(R_{\hat{\theta}_{n}}(\cdot))\mid_{\mb 0} \big)^{\frac{1}{2}}\Big]
    }
    with $\hat{\mb \Sigma}$ defined in \eqref{eq: hatSigma}.
\end{algorithm}

\subsubsection{Confidence Regions Construction: An Insight into Coordinate Representations}\label{subsection: curvature effect}

To clarify the importance of employing a second-order retraction and selecting appropriate coordinate representations across different charts in Algorithm~\ref{algorithm:wald based statistic} and Algorithm~\ref{algorithm:t based statistic}, we begin by examining the geometry of a two-dimensional submanifold $\mc M$ embedded in $\bb R^3$. Recall that $\theta_0$ represents the true parameter, $\hat{\theta}_n$ denotes the outputs of Algorithm~\ref{algorithm: resampled newton iteration}, and $\hat \theta^*_{n}$ refers to a generic version of $\hat \theta^{*[i]}_{n}$.

One appealing but impractical choice to represent $\theta_0, \hat \theta_n, \hat \theta_n^*$ is to utilize the chart centered at $\theta_0$ so that the problem is reduced to the well-studied Euclidean M-estimation problem by considering the objective function $L( R_{\theta_0}(\cdot ), \cdot ):\bb R^p \times \mc S \rightarrow \bb R$. Although in practice the mapping $R_{\theta_0}$ is inaccessible without the knowledge of $\theta_0$, it enlightens us to bridge the actual coordinate representations with the oracle coordinates under $R_{\theta_0}$; indeed, the roadmap starts with relating $R_{\hat\theta_n}^{-1}(\theta_0)$ and $-R_{\hat \theta_n^*}^{-1}(\hat \theta_n)$ with $-R_{\theta_0}^{-1}(\hat\theta_n)$ and $R_{\theta_0}^{-1}(\hat\theta_n^*) - R_{\theta_0}^{-1}(\hat\theta_n)$, respectively. Consider the case where $\hat\theta_n$ and $\hat \theta_n^*$ are contained in a sufficiently small neighborhood of $\theta_0$, which can be justified under some regularity conditions as $n$ goes to infinity. 
   	
\begin{figure}[htbp]
	\centering
	\subfigure[Tangent Spaces $\mathrm T_{\theta_0}\mc M$ and $\mathrm T_{\hat\theta_n}\mc M$]{
	\begin{tikzpicture}[scale=.5]
	\filldraw[color=americanrose!10] (1,5.5) -- (-1.5,3.5) -- (8.5,2.5) -- (10.5,4.8) -- (1,5.5);
	\filldraw[color = amber!10] (5.42,2.82) -- (2.3,6.2) -- (6.5,8.0) -- (9.36,4.88);
	\filldraw[color = amber!10] (7.8,0.2) -- (12,2) -- (10,4.24) -- (8.5,2.5) -- (5.42,2.82);
	\draw [densely dotted] (1,5.5) -- (-1.5,3.5) -- (8.5,2.5) -- (10.5,4.8) -- (1,5.5);
	\draw [densely dotted] (2.3,6.2) -- (6.5,8.0) -- (12,2) -- (7.8,0.2) -- (2.3,6.2);
	\coordinate [label=-170:{\color{americanrose}$\hat{\theta}_n$}] (thetahat) at (4,4);
	\coordinate [label=right:{\color{amber}$\theta_0$}] (theta0) at (7.5,3.2);
	\coordinate [label=left:{\color{airforceblue}$\hat \theta_n^*$}] (thetahatstar) at (6,2);
	\coordinate [label=left:{\color{americanrose}$\textrm T_{\hat\theta_n}\mc M$}] (T_thetahat) at (0,5);
	\coordinate [label=left:{\color{amber}$\textrm T_{\theta_0}\mc M$}] (T_theta0) at (10,7);
	\coordinate (C) at ($ (thetahat) ! .5 ! (theta0) $);
	\coordinate [label = left:$\mc M$] (M) at (5,0.5);
	\coordinate (A) at ($ (thetahat) ! .5 ! (thetahatstar) $);
	\coordinate [label={[shift={(0.15,-0.05)}]\tiny{\color{americanrose}$\mathrm dR_{\hat \theta_n}\circ R^{-1}_{\hat \theta_n}(\theta_0)$}}](D) at ($(thetahat)+(3.5,0.3)$);
	\coordinate [label={[shift={(0.1,-0.05)}]\tiny{\color{amber}$dR_{\theta_0}\circ R^{-1}_{\theta_0}(\hat\theta_n)$}}](E) at ($(theta0)+(-3, 2)$);
	\node [fill=americanrose,inner sep=.8pt,circle] at (thetahat) {};
	\node [fill=amber,inner sep=.8pt,circle] at (theta0) {};
	\node [fill=airforceblue,inner sep=.8pt,circle] at (thetahatstar) {};
	\draw (0,2)	.. controls (0.5,5) and (3,5) .. (6,-1);
	\draw (0,2)	.. controls (1,2.5) and (2,3) .. (3.34,3);
	\draw (6,-1)	.. controls (7.5,0) and (9.5,1) .. (12,0);
	\draw (0.5,3.4)	.. controls (1.7,5.5) and (11,8) .. (12,0);	
	\draw [dashed] (thetahat)	.. controls ($(C)+(-0.5,0.5)$) and ($(C)+(0.5,0.5)$) .. (theta0);	
	\draw [->,color = americanrose](thetahat) -- (D);
	\draw [densely dotted,color=amber(sae/ece)] (D) .. controls ($ (D) ! .5 ! (theta0) + (0.05,-0.1)$) and ($ (D) ! .5 ! (theta0) + (0.05,-0.1)$) .. (theta0);
	\draw [->, color = amber](theta0) -- (E);
	\draw [densely dotted,color=amber(sae/ece)] (E) .. controls ($ (E) ! .5 ! (thetahat) + (-0.05,-0.1)$) and ($ (E) ! .5 ! (thetahat) + (-0.05,-0.1)$) .. (thetahat);
	\end{tikzpicture}
}	\hfill
	\subfigure[Tangent Spaces $\mathrm T_{\theta_0}\mc M$ and $\mathrm T_{\hat\theta_n^*}\mc M$]{
	\begin{tikzpicture}[scale=.5]
	\filldraw[name path = line1, color = amber!10] (2.3,6.2) -- (6.5,8.0) -- (12,2) -- (7.8,0.2) -- (2.3,6.2);
	\path[name path = line2] (4.5,-0.5) -- (9,0) -- (7.5,8.5) -- (3,8) -- (4.5,-0.5);
	\draw [densely dotted] (2.3,6.2) -- (6.5,8.0) -- (12,2) -- (7.8,0.2) -- (2.3,6.2);
	\draw[name intersections={of=line1 and line2}]
		(intersection-1) node {}
    	(intersection-2) node {}
    	(intersection-3) node {}
    	(intersection-4) node {}
    	;
    \path[name path = line4] (7.8,0.2) -- (2.3,6.2);
    \path [name path = line3] (3,2.2) -- (10, 3.6);
   	\path [name path = line5] (7.5,8.5) -- (9,0);
    \path[name intersections={of = line4 and line3, by = {E}}];
    \path[name intersections={of = line5 and line3, by = {F}}];
    \filldraw[color = airforceblue!10] (7.5,8.5) -- (3,8) -- (4.5,-0.5) -- (9,0) -- (intersection-3) -- (7.8,0.2) -- (E) -- (F);
    \draw [densely dotted] (7.5,8.5) -- (3,8) -- (4.5,-0.5) -- (9,0) -- (intersection-3) -- (7.8,0.2) -- (E) -- (F) -- (7.5,8.5);
   	\draw [densely dotted] (2.3,6.2) -- (6.5,8.0) -- (12,2) -- (7.8,0.2) -- (2.3,6.2);
	\coordinate [label=-170:{\color{americanrose}$\hat{\theta}_n$}] (thetahat) at (4,4);
	\coordinate [label=right:{\color{amber}$\theta_0$}] (theta0) at (7.5,3.2);
	\coordinate [label=left:{\color{airforceblue}$\hat \theta_n^*$}] (thetahatstar) at (6,2);
	\coordinate [label=left:{\color{airforceblue}$\textrm T_{\hat\theta_n^*}\mc M$}] (T_thetahat) at (7,9);
	\coordinate [label=left:{\color{amber}$\textrm T_{\theta_0}\mc M$}] (T_theta0) at (10,7);
	\coordinate (C) at ($ (thetahat) ! .5 ! (theta0) $);
	\coordinate [label = left:$\mc M$] (M) at (5,0.5);
	\coordinate (A) at ($ (thetahat) ! .5 ! (thetahatstar) $);
	\coordinate [label={[shift={(0.1,-0.05)}]\tiny{\color{amber}$dR_{\theta_0}\circ R^{-1}_{\theta_0}(\hat\theta_n)$}}](E) at ($(theta0)+(-3, 2)$);
	\coordinate [label = {[shift = {(-1,-0.5)}]\color{airforceblue}\tiny $\mathrm d R_{\hat\theta_n^*} \circ R^{-1}_{\hat \theta_n^*}(\hat \theta_n)$}](F) at ($(thetahatstar) + (-0.9,2.7)$); 
	\coordinate [label = right:{\color{amber}\tiny $\mathrm d R_{\hat\theta_n^*} \circ R^{-1}_{\hat \theta_n^*}(\hat \theta_n)$}](G) at ($(theta0) + (-1.2,-1.2)$); 
	\node [fill=americanrose,inner sep=.8pt,circle] at (thetahat) {};
	\node [fill=amber,inner sep=.8pt,circle] at (theta0) {};
	\node [fill=airforceblue,inner sep=.8pt,circle] at (thetahatstar) {};
	\draw (0,2)	.. controls (0.5,5) and (3,5) .. (6,-1);
	\draw (0,2)	.. controls (1,2.5) and (2,3) .. (3.34,3);
	\draw (6,-1)	.. controls (7.5,0) and (9.5,1) .. (12,0);
	\draw (0.5,3.4)	.. controls (1.7,5.5) and (11,8) .. (12,0);
	\draw [->, color = airforceblue] (thetahatstar) -- (F);
	\draw [->, color = amber] (theta0) -- (G);
	\draw [->, color = amber](theta0) -- (E);
	\draw [dashed] (thetahat)	.. controls ($(A)+(-0.5,1)$) and ($(A)+(0.5,0.7)$) .. (thetahatstar);
	\draw [dashed] (thetahat)	.. controls ($(C)+(-0.5,0.5)$) and ($(C)+(0.5,0.5)$) .. (theta0);
	\draw [dashed] (thetahatstar)	.. controls ($ (theta0) ! .5 ! (thetahatstar) + (-0.05,-0.1)$) and ($ (thetahatstar) ! .5 ! (theta0) + (0.05,-0.1)$) .. (theta0);
	\draw[densely dotted, color=amber(sae/ece)] (thetahat) .. controls ($ (F) ! .5 ! (thetahat) + (-0.05,0.1)$) and ($ (thetahat) ! .5 ! (F) + (0.05,0.1)$) .. (F);
	\draw[densely dotted, color=amber(sae/ece)] (thetahatstar) .. controls ($ (G) ! .5 ! (thetahatstar) + (-0.05,-0.1)$) and ($ (thetahatstar) ! .5 ! (G) + (0.05,-0.1)$) .. (G);
	\draw [densely dotted,color=amber(sae/ece)] (E) .. controls ($ (E) ! .5 ! (thetahat) + (-0.05,-0.1)$) and ($ (E) ! .5 ! (thetahat) + (-0.05,-0.1)$) .. (thetahat);
	\end{tikzpicture}
	}
	\caption{Tangent spaces and coordinates of a dimension-$2$ submanifold $\mc M$ in $\bb R^3$. }
	\label{fig:tangent space}
\end{figure}
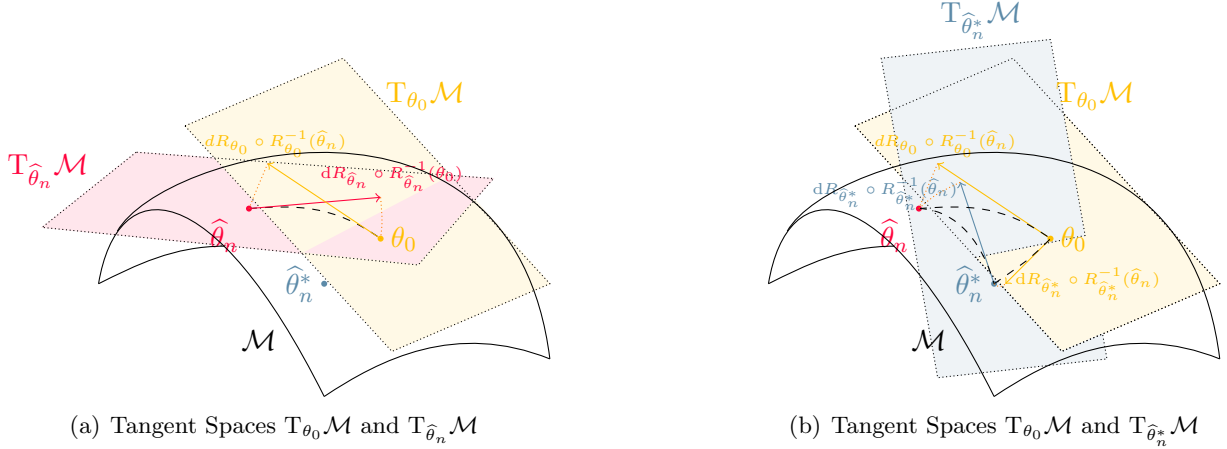

We note that the Edgeworth expansions of the studentized versions of $R^{-1}_{\theta_0}(\hat{\theta}_n)$ and $R^{-1}_{\theta_0}(\hat{\theta}_n^*)$ are readily accessible, owing to the previous results (e.g., \cite{bhattacharya1978validity, andrews2002higher}) on Edgeworth expansions for (approximate) M-estimations. Consequently, it is natural to pursue equating $R^{-1}_{\theta_0}(\hat \theta_n)$ with $-R^{-1}_{\hat\theta_n}(\theta_0)$ and also $R_{\theta_0}^{-1}(\hat \theta_n^*) -R_{\theta_0}^{-1}(\hat \theta_n)$ with $-R_{\hat \theta_n^*}(\hat \theta_n)$, as illustrated in Figure~\ref{fig:tangent space}. This analysis is conducted in two distinct parts: 
   	
	\begin{itemize}
		\item \textbf{Identify $R_{\hat\theta_n}^{-1}(\theta_0)$ with $-R_{\theta_0}^{-1}(\hat\theta_n)$.} We relate the coordinate of $\hat\theta_n$ under $R_{\theta_0}$ to the exponential mapping by letting
		\begin{equation}\label{eq:triangle1}
					R^{-1}_{\theta_0}(\hat \theta_n) = \pi_{\{\mathrm dR_{\theta_0} \delta_i \}}\circ \Log_{\theta_0}(\hat \theta_n) + \triangle_1. 
		\end{equation}
		
		Similarly, for $R_{\hat\theta_n}(\theta_0)$ one has 
		\begin{equation}\label{eq:triangle1'}
			R^{-1}_{\hat\theta_n}(\theta_0) = \pi_{\{\mathrm dR_{\hat\theta_n} \delta_i\}}\circ \Log_{\hat\theta_n}(\theta_0) + \triangle_1'.
		\end{equation}
		
		Moreover, the invariance property of parallel transport yields that 
		\begin{equation}
			\pi_{\{\mathrm dR_{\theta_0} \delta_i\}}\circ \Log_{\theta_0}(\hat \theta_n)  = - \pi_{\{T_{\theta_0\rightarrow \hat \theta_n}\mathrm dR_{\theta_0} \delta_i\}}\circ \Log_{\hat\theta_n}(\theta_0).
		\end{equation}
		Combining the above pieces together implies that
		\begin{equation}\label{eq:implication1}
			\norm{R^{-1}_{\theta_0}(\hat \theta_n)+  \pi_{\{T_{\theta_0\rightarrow \hat \theta_n}\mathrm dR_{\theta_0} \delta_i\}}\circ\pi_{\{\mathrm dR_{\hat\theta_n} \delta_i\}}^{-1} \circ R^{-1}_{\hat\theta_n}(\theta_0)} \leq \norm{\triangle_1} + \norm{\triangle_1'}. 
		\end{equation}

		\item \textbf{Identify $-R_{\hat\theta_n^*}^{-1}(\hat \theta_n)$ with $R_{\theta_0}^{-1}(\hat \theta_n^*) -R_{\theta_0}^{-1}(\hat \theta_n)$. } We denote that 
		\begin{align}
			& R^{-1}_{\theta_0}(\hat\theta_n^*) =  \pi_{\{\mathrm dR_{\theta_0} \delta_i\}}\circ \Log_{\theta_0}(\hat \theta_n^*) + \triangle_2,\\
			& R^{-1}_{\hat\theta_n^*}(\hat\theta_n) = \pi_{\{\mathrm dR_{\hat\theta_n^*} \delta_i\}}\circ \Log_{\hat\theta_n^*}(\hat\theta_n) + \triangle_3. 
		\end{align}
		
		Again, leveraging parallel transports, we derive that 
		\begin{align}
			& \pi_{\{\mathrm dR_{\theta_0} \delta_i\}}\circ \Log_{\theta_0}(\hat \theta_n^*)  = - \pi_{\{T_{\theta_0\rightarrow \hat \theta_n^*}\mathrm dR_{\theta_0} \delta_i\}}\circ \Log_{\hat\theta_n^*}(\theta_0).
		\end{align}
		
		Then the difference between $- \pi_{\{T_{\theta_0\rightarrow \hat \theta_n^*}\mathrm dR_{\theta_0} \delta_i\}} \circ \pi^{-1}_{\{\mathrm dR_{\theta_0} \delta_i\}}\circ R_{\hat\theta_n^*}^{-1}(\hat \theta_n)$ and $R_{\theta_0}^{-1}(\hat \theta_n^*) -R_{\theta_0}^{-1}(\hat \theta_n)$ could be written as 
		\begin{equation}\label{eq:implication2}
			\norm{R_{\theta_0}^{-1}(\hat \theta_n^*) -R_{\theta_0}^{-1}(\hat \theta_n) + \pi_{\{T_{\theta_0\rightarrow \hat \theta_n^*}\mathrm dR_{\theta_0} \delta_i\}}\circ \pi^{-1}_{\{\mathrm dR_{\theta_0} \delta_i\}}\circ R_{\hat\theta_n^*}^{-1}(\hat \theta_n)} \leq \triangle_1 + \triangle_2 + \triangle_3 + \tilde \triangle,
		\end{equation}
		where $\tilde \triangle$ is defined as 
		\begin{equation}
			\tilde \triangle \coloneqq    \norm{\Log_{\theta_0}(\hat \theta_n^*) - \Log_{\theta_0}(\hat \theta_n) +T_{\hat\theta_n^* \rightarrow\theta_0} \Log_{\hat\theta_n^*}(\hat\theta_n)} .
		\end{equation}
	\end{itemize}
   	
    Informally speaking, since the magnitudes of $\dist(\hat\theta_n, \theta_0)$ and $\dist(\hat\theta_n^*, \theta_0)$ are both $O(n^{-\frac{1}{2}}(\log n)^{\frac12})$, to prove that the distributions of $-\pi_{\{T_{\theta_0\rightarrow \hat \theta_n}\mathrm dR_{\theta_0} e_i\}}\circ\pi_{\{\mathrm dR_{\hat\theta_n} \delta_i\}}^{-1} \circ R^{-1}_{\hat\theta_n}(\theta_0)$ and $-\pi_{\{T_{\theta_0\rightarrow \hat \theta_n^*}\mathrm dR_{\theta_0} \delta_i\}}\circ \pi^{-1}_{\{\mathrm dR_{\theta_0} \delta_i\}}\circ R_{\hat\theta_n^*}^{-1}(\hat \theta_n)$ align with those of $R_{\theta_0}(\hat \theta_n)$ and $R_{\theta_0}^{-1}(\hat \theta_n^*) -R_{\theta_0}^{-1}(\hat \theta_n)$ in the sense of high order, we need to show that the magnitudes of $\triangle_1,\triangle_1', \triangle_2, \triangle_3, \tilde\triangle$ are small enough, to be specific, of order $O(n^{-\frac32}(\log n)^{\frac32})$.  
	
	By the second-order retraction's property (see Lemma~\ref{lemma: uniform control for second-order retraction} in the supplement) along with the condition $\dist(\hat\theta_n, \theta_0) \vee \dist(\hat\theta_n^*, \theta_0) = O(n^{-\frac12} (\log n)^{\frac12})$, one can show that $\triangle_1, \triangle_1', \triangle_2$, and $\triangle_3$ satisfy the above requirement. What remains to be shown is the upper bound for $\tilde \triangle$ which relies on the so-called \emph{double exponential mapping} \citep{brewin2009riemann,gavrilov2006algebraic,gavrilov2007double}.

\paragraph*{Double Exponential Mapping} Given $x\in \mc M$, the double exponential mapping is defined as $\Log_x\big(\Exp_{\Exp_x(v)}(T_{x\rightarrow \Exp_x v}(w))\big)$ with $v,w\in \mathrm T_x\mc M$. Intuitively, the difference $\Log_x\big(\Exp_{\Exp_x(v)}(T_{x\rightarrow \Exp_x v}(w))\big) - v - w \in \mathrm T_{x}\mc M$ represents a measure of how curvature affects the triangle relation among $v$, $w$, and $v + w$.  This difference captures the deviation from the expected vector space addition $v+w$ when $w$ is transported from $x$ to $\Exp_{x}(v)$. Informally, it could be shown that 
	\eq{
	\norm{\Log_x\big(\Exp_{\Exp_x(v)}(T_{x\rightarrow \Exp_x v}(w))\big) - v - w} \leq C_1 \norm{v}\norm{w} + C_2 \norm{v}\norm{w}^2,
	}
    where the constants $C_1$ and $C_2$ depend on the local curvature around $x$. A formal result is stated in Lemma~\ref{lemma: normal coordinate transomation} in the Supplementary Material.
	
\begin{figure}[htbp]
	\centering
	\begin{tikzpicture}[scale = 0.3]		
		\draw (0,0) .. controls (4,-0.6) and (8,-0.2) .. (12,1); 
		\draw (12,1) .. controls (8,6) .. (6,12); 
		\draw (0,0) .. controls ($(0,0) ! 0.5 ! (6, 12) + (0,1)$) .. (6,12); 
		\coordinate [label=left:$\hat \theta_n^*$] (thetahatstar) at (0,0); 
		\coordinate [label=left:$\theta_0$] (theta0) at (12,1); 
		\coordinate [label=left:$\hat{\theta}_n$] (thetahat) at (6,12); 
		\coordinate [label=left:{\color{airforceblue}$\Log_{\hat \theta_n^*}(\hat \theta_n)$}] (logthetahatstar) at ($(thetahatstar) + (1.5, 4.5)$);
		\coordinate [label=above:{\color{americanrose}$T_{\hat \theta_n^* \rightarrow \theta_0}\Log_{\hat \theta_n^*}(\hat \theta_n)$}] (logtheta0) at ($(theta0) + (1.5, 4.5)$);
		\coordinate [label=below:{\color{americanrose}$\Log_{\theta_0}(\hat \theta_n^*)$}] (logtheta02) at ($(theta0) + (-4.5, -1.3)$);
		\coordinate [label=above:{\color{americanrose}$\Log_{\theta_0}(\hat \theta_n)$}] (logtheta03) at ($(theta0) + (-4, 4)$);
		\node [fill=black,inner sep=.8pt,circle] at (thetahat) {};
		\node [fill=black,inner sep=.8pt,circle] at (theta0) {}; 
		\node [fill=black,inner sep=.8pt,circle] at (thetahatstar) {};
		\draw [color=deepBlue,thick,->] (thetahatstar) -- (-60:1);
		\draw [color=deepBlue,thick,->] (thetahatstar)-- (30:1);
		\draw [color=deepBlue,thick,->] (12,1) -- (12.8,0.4);
		\draw [color=deepBlue,thick,->] (12,1)-- (12.6,1.8);
		\draw [color=airforceblue,thick, ->] (thetahatstar) -- (logthetahatstar);
		\draw [color=americanrose,thick, ->] (theta0) -- (logtheta0);
		\draw [color=americanrose,thick, ->] (theta0) -- (logtheta02);
		\draw [color=americanrose,thick, ->] (theta0) -- (logtheta03);
		\draw [color = amber(sae/ece),dashed] (30:1) .. controls ($(30:1) ! 0.5 ! (12.6, 1.8) + (0,-0.5)$) .. (12.6,1.8); 
		\draw [color=amber(sae/ece), dashed] (logthetahatstar) .. controls ($(logthetahatstar) ! 0.5 ! (logtheta0) + (0,-0.5)$) .. (logtheta0);
		\draw [color=amber(sae/ece),dashed] (-60:1) .. controls (6,-1) .. (12.8,0.4); 
		\coordinate [label= above:$ e_1'$] (e1') at (30:1); 
		\coordinate [label= below:$ e_2'$] (e2') at (-60:1); 
		\coordinate [label= right:$ e_1$] (e1) at (12.6,1.8); 
		\coordinate [label= below:$ e_2$] (e2) at (12.8,0.4); 
		\coordinate [label= below:{\color{amber(sae/ece)}$ \text{Parallel Transport}$}] (parallel transport) at (6,2.7); 
	\end{tikzpicture}
	\caption{Double exponential mapping.}
	\label{fig:coordinates}
\end{figure}
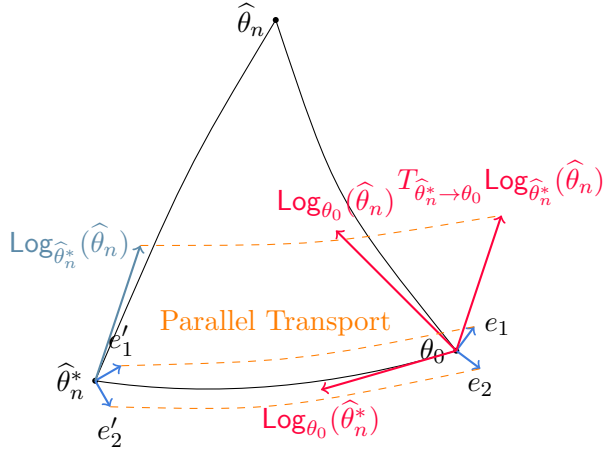

	Back to the problem above, we let $x$ be $\theta_0$, $v$ be $\Log_{\theta_0}(\hat \theta_n^*)$ and $w$ be $T_{\hat\theta_n^* \rightarrow\theta_0} \Log_{\hat\theta_n^*}(\hat\theta_n)$ so that $\Log_x\big(\Exp_{\Exp_x(v)}(T_{x\rightarrow \Exp_x v}(w))\big)$ becomes $\Log_{\theta_0}(\hat \theta_n)$, as depicted in Figure~\ref{fig:coordinates}. Then we arrive at
	\eq{
	\tilde \triangle = O(n^{-3/2}(\log n)^{\frac32}),
	}
	which satisfies the requirement. 
	
\paragraph*{Consistency after Studentization} 
Studentizing the differences is crucial for achieving high-order consistency between the target distribution and its resampling counterpart. However, the studentization term based on the chart $\pi_{\{e_i\}} \circ \Log_{\theta_0}(\cdot )$ in our theoretical analysis does not directly appear in the algorithmic implementation; in practice, we naturally use the studentization terms centered at $\hat{\theta}_n$ or $\hat{\theta}_n^*$. These terms, originating from different perspectives, can be reconciled via Lemma~\ref{lemma: Hessian consistency} in the Supplementary Material, which establishes that the approximation error of a Hessian operator depends on both the distance between the centers of the two charts and the magnitude of the gradient.

\paragraph*{Rotation Issue}
Finally, the approximations in \eqref{eq:implication1} and \eqref{eq:implication2} depend on certain orthogonal transformations, specifically $ \pi_{\{T_{\theta_0\rightarrow \hat \theta_n}\mathrm dR_{\theta_0} \delta_i\}}\circ\pi_{\{\mathrm dR_{\hat\theta_n} \delta_i\}}^{-1}$ and $\pi_{\{T_{\theta_0\rightarrow \hat \theta_n^*}\mathrm dR_{\theta_0} \delta_i\}}\circ \pi^{-1}_{\{\mathrm dR_{\hat \theta_n^*} \delta_i\}}$, which are typically intractable. To address this issue, we adopted the following two approaches.
\begin{itemize}
    \item \textbf{Rotation-Insensitivity of Wald Statistics.} When computing the Wald-type statistic related to $R^{-1}_{\hat{\theta}_n}(\theta_0)^\top \hat{\mb{\Sigma}}^{-1} R^{-1}_{\hat{\theta}_n}(\theta_0)$ and $R^{-1}_{\hat{\theta}_n^*}(\hat{\theta}_n)^\top \hat{\mb{\Sigma}}^{*\dagger} R_{\hat{\theta}_n^*}^{-1}(\hat{\theta}_n)$, the orthogonal mappings are canceled out regardless of the choice of retraction. Consequently, these statistics are insensitive to potential rotations.
		\item \textbf{Imposing Differentiability on Retractions.} For $t$-statistics, the example in Figure~\ref{fig: counterexample} has informed us that the non-differentiability of a chart may lead to irregularities. To mitigate this, we impose certain differentiability conditions and focus on proving that the proposed bootstrap procedure can effectively capture the deviation tendency of a chart.
	\end{itemize}

In summary, the key distinction between Euclidean spaces and general Riemannian manifolds lies in the locality of coordinate expressions, particularly when transitioning across different charts. To address this, we bridge gaps between idealized and realizable coordinate systems. The imposition of differentiability conditions and the adoption of rotation-insensitive techniques further address irregularities arising from rotation issues. These efforts culminate in a rigorous framework that achieves high-order consistency between the target distribution and its resampling counterpart, providing a foundation for high-order accurate inference on Riemannian manifolds.

\section{Theoretical Distributional Guarantees}\label{sec: theoretical guarantees}

In this section, we show that the resampling procedures implemented by the aforementioned algorithms yield high-order accurate approximations of the distributions for the constructed statistic. 
\subsection{Assumptions}
Our approach requires the following regularity conditions to perform Taylor's expansion of the objective function and to ensure the validity of Edgeworth expansion, part of which is inherited from \cite{bhattacharya1978validity,hall1996bootstrap} from the Euclidean viewpoint, while further modifications are introduced to adapt the method to manifold settings. We make the following assumption:
\begin{assumption}
\label{assumption:manifold}
We let $\{e_i\} $ be an orthogonal basis of $\mathrm T_{\theta_0}\mc M$. 
    \begin{enumerate}[label=\Alph*]
        \item \label{item:A}
       	The loss function $L(\theta, x)$ is $4$-times differentiable with respect to $\theta$ for $(\theta,  x) \in \mc  M \times \mc S$. 
        \item \label{item:B} There exists a compact subset $K$ of $\mc M$ such that $\theta_0$ is an interior point in $K$. 
        For each $\mb \nu$ with $1\leq |\mb \nu|\leq 4$, $\bb E_{\theta_0}\big[| \nabla^{\mb \nu}L(\theta_0, X_1)|^{32}\big]<+\infty$ holds. 
         For every $\theta_1\in \mc M$, there exists a function $C(\mb x)$ such that $\big|\nabla^{\mb \nu} L(\theta_0, \mb x) - \bar\nabla^{\mb \nu}L\circ \Exp_{\theta_0}\circ \pi^{-1}_{\{e_i\}}\big(\pi_{\{e_i\}} \circ \Log(\theta_1), \mb x\big)\big| \leq C(\mb x)\cdot  \dist(\theta_0, \theta_1)$ for $x \in \mc S$, $\mb \nu\in \bb N^p$ with $1\leq |\mb \nu| \leq 4 
    $ and $\mathbb E_{\theta_0}\big[C(X_1)^{16}\big] \leq \infty$. 

    \item \label{item:C} $\bb E_{\theta_0}\big[\nabla_{e_i} L(\theta_0, X_1)\big] = 0$ holds for $i\in[p]$, and $
        \big(\bb E_{{\theta}_0}\left[\nabla_{e_i} L({\theta}_0, X_1)\nabla_{e_j}  L({\theta}_0, X_1)\right]\big)_{i,j \in[p]}
    $
    is nonsingular. 
    \item \label{item:E} We let $\mb Z_i$ be a vector that collects $\nabla^{\mb \nu}L(\theta_0,X_i)$, $\nabla^{\mb \nu_1} L(\theta_0, X_i)\nabla^{\mb \nu_2} L(\theta_0, X_i)$ for $1\leq |\mb \nu|\leq 3$,$ 1 \leq |\mb \nu_1|+ |\mb \nu_2| \leq 2$. The random vector $\mb Z_1$ satisfies the Cr\'amer condition, that is, 
    \begin{equation}
    	\limsup_{\norm{\mb t}_2 \rightarrow \infty} \big|\bb E\big[e^{i\mb t\t \mb Z_1}\big]\big|<1. 
    \end{equation}
    \end{enumerate}
\end{assumption}
We briefly remark on the assumptions that (i) the moment conditions on the derivatives of $L$ up to order four, as well as on $C(X_1)$, are primarily required for applying the moderate deviation inequality (see Lemma~\ref{lemma: thm1 in von1967central} in the Supplement), and (ii) the Cr\'amer condition is a standard assumption in establishing the validity of the Edgeworth expansion \cite{bhattacharya1978validity,hall2013bootstrap}. 

\subsection{Point Estimation Convergence Guarantees}
We first showcase our convergence guarantees for the output sequence of Algorithm~\ref{algorithm: resampled newton iteration} based on Assumption~\ref{assumption:manifold}. 
\begin{theorem}[Convergence guarantees: Algorithm~\ref{algorithm: update} with full data]\label{thm:convergence rate of newton}
	If Assumption~\ref{assumption:manifold} holds and the initial estimate satisfies $\dist(\theta^{\text{initial}},\theta_0) \leq c_0$ for some sufficiently small $c_0$, then with $\mathsf{burn} = c\log \log n $, where $c$ is a constant determined by $c_0$, the objective function $L$, and $\mc Y = \mc X_n$, the output of Algorithm~\ref{algorithm: update} satisfies
	\longeq{
	& \dist(\hat \theta_n, \tilde \theta_n ) = O(n^{-2})
	}
	with probability at least $1 - O(n^{-1}(\log n)^{-2})$. 
\end{theorem}

\begin{theorem}[Convergence guarantees: Algorithm~\ref{algorithm: update} with bootstrapped data]\label{thm: convergence rate of resampled newton}
    Suppose that Assumption~\ref{assumption:manifold} holds. Then given an initial estimate $\theta^{\text{initial}}$ satisfying 
        $\dist({\theta}^{\text{initial}}, \tilde{ \theta}_{n}) = O(n^{-2})$ and $i\in[b]$, the output of Algorithm~\ref{algorithm: update} with $t = 2$ and $\mc Y = \mc X_n^{[i]}$, $\hat \theta_{n}^{*[i]}$, satisfies
    \eq{
        \dist\big(\hat{\theta}_{n}^{*[i]}, \tilde\theta_{n}^{[i]}\big) = O\big(n^{-2}(\log n)^2\big)
    } 
    with probability $1 - O(n^{-1}(\log n)^{-2})$ conditional on $\mc X_n$, where $\tilde \theta_{n}^{[i]}$ represents a solution in \eqref{eq: resampled first-order condition of M estimation} given the resampled dataset $\mc X^{[i]}_n$. 
\end{theorem}

\begin{remark}
    Given a deterministic objective function, the quadratic convergence of Newton's algorithms on Riemannian manifolds has been studied in a sequence of previous work \citep{smith1993geometric,smith2014optimization,absil2009optimization}. Theorem~\ref{thm:convergence rate of newton} and Theorem~\ref{thm: convergence rate of resampled newton} extend the convergence guarantees to the stochastic settings with some additional analysis of the fluctuation of derivatives.  
\end{remark}
\begin{remark}
    As demonstrated in previous works \citep{absil2009optimization}, quadratic convergence does not depend on the second-order property of the retraction being used. Nonetheless, the second-order retraction is primarily for the sake of high-order asymptotics upon the coordinate calculation. 
\end{remark}

\subsection{High-order Accurate Distributional Guarantees on Riemannian Manifolds}

This subsection introduces our main theorem on the distributional guarantees for the approximate M-estimators on Riemannian manifolds.

From the statement of Algorithm~\ref{algorithm: resampled newton iteration}, it is clear that, conditioning on the original dataset $\mc X_n$, $\{\hat\theta_{n}^{*[i]}\}_{i\in[b]}$ are identically distributed and conditionally independent of each other. Therefore, we omit the superscript $[i]$ to denote by $\hat{\theta}_n^*$ the generic resampled output. For the intrinsic $t$-statistic and the intrinsic Wald-statistic, the main result is as follows.

\begin{theorem}\label{theorem:bootstrapmanifold}\footnote{Consistent with the convention in Algorithm~\ref{algorithm:wald based statistic} and Algorithm~\ref{algorithm:t based statistic}, we let $R^{-1}_x(y)$ be $\mb 0$ for any $y$ outside the range of $R$. }
    Suppose Assumption~\ref{assumption:manifold} holds and the initial estimate satisfies $\dist(\theta^{\text{initial}},\theta_0) \leq c_0$ for some sufficiently small $c_0$.
    \begin{enumerate}
        \item For every $0< \xi < \frac{1}{2}$, it holds that  
        \longeq{
        & \limsup_{n\rightarrow \infty} n\bb P_{\mc X_n}\bigg[\sup_{-\infty < x< \infty} \Big| \bb P\left[R^{-1}_{\hat{\theta}_n^*}(\hat\theta_n)\t {{}\check{\mb \Sigma}}^{\dagger}R^{-1}_{\hat{\theta}_n^*}(\hat \theta_n) \leq x \mid \mc X_n\right]\\ 
        &\qquad  - \bb  P_{\mc X_n} \left[R^{-1}_{\hat{\theta}_n}(\theta_0)\t {{}\hat{\mb \Sigma}}^{\dagger}R^{-1}_{\hat{\theta}_n}(\theta_0) \leq x \right]\Big|  > n^{-1  +\xi} \bigg] < +\infty ,\label{eq: thm1.1}
        }  
        where $\hat{\mb \Sigma}$ and $\check{\mb \Sigma}$ are defined in \eqref{eq: hatSigma} and \eqref{eq: check Sigma}, respectively.\label{item: thm1.1}
    \item If the retraction mappings $R_{\theta}$ for $\theta$'s in a neighborhood of $\theta_0$ satisfy that the function $\mathrm{d}R_{\theta}^{-1}(T_{\theta_0 \rightarrow \theta}(\mathrm{d} R_{\theta_0}(\delta_i)))$ is twice differentiable with respect to $\theta$ for every $i \in[p]$, then for an arbitrary $\mb a\in \bb R^{p}$ and every $0< \xi< \frac{1}{2}$ it holds that 
	\begin{equation}\label{eq: thm1.2}
     \limsup_{n\rightarrow \infty} n\bb P_{\mc X_n}\bigg[\sup_{-\infty< x < \infty} \Big|\bb P\Big[\frac{\mb a\t R^{-1}_{\hat \theta_n^*}(\hat \theta_n)}{\big(\mb a\t\check{\mb \Sigma}\mb a\big)^{\frac{1}{2}}} \leq  x \mid \mc X_n \Big] - \bb P_{\mc X_n}\Big[\frac{\mb a\t R^{-1}_{\hat \theta_n}(\theta_0)}{\big(\mb a\t\hat{\mb \Sigma}\mb a\big)^{\frac{1}{2}}} \leq x  \Big]  \Big| >n^{-1 +\xi} \bigg] <+\infty .
	\end{equation}\label{item: thm1.2}
\end{enumerate}
\end{theorem}

Moreover, for the extrinsic $t$-statistic, we establish the following distributional guarantee, which can be viewed as a byproduct of the Euclidean Edgeworth expansion theory:  
\begin{proposition}\label{proposition:edgeworth for extrinsic t}
If $\mc M$ is a submanifold embedded in $\bb R^{d'}$, Assumption~\ref{assumption:manifold} holds, and the initial estimate satisfies $\dist(\theta^{\text{initial}},\theta_0) \leq c_0$ for some sufficiently small $c_0$, then for any smooth and Lipschitz-continuous function $f:\bb R^{d'} \rightarrow \bb R$ with $\nabla f\circ \Exp\mid_{\theta_0} \neq \mb 0$ it holds that 
    \eq{\label{eq: thm1.3 }
        \begin{split}
        & \limsup\limits_{n\rightarrow \infty} n\bb P_{\mc X_n}\bigg[\sup_{-\infty< x<\infty}\Big| \bb P\Big[\frac{f(\hat{\theta}_n^*) - f(\hat{\theta}_n)}{\big(\bar \nabla f(R_{\hat \theta_n^*}(\cdot))\mid_{\mb 0} \t 
    \check {\mb \Sigma}\bar \nabla f(R_{\hat \theta_n^*}(\cdot))\mid_{\mb 0}  \big)^{\frac{1}{2}}}\leq x \mid \mc X_n\Big] \\ 
    & \qquad - \bb P_{\mc X_n}\big[\frac{f(\hat{\theta}_n) - f(\theta_0)}{\big(\bar \nabla f(R_{\hat \theta_n}(\cdot))\mid_{\mb 0}\t  
    \hat {\mb \Sigma}\bar \nabla f(R_{\hat \theta_n}(\cdot))\mid_{\mb 0} \big)^{\frac{1}{2}}} \leq x\Big]\Big| \geq n^{-1+\xi} \bigg] \leq \infty
    \end{split}}
\end{proposition}

In practice, it is impossible to resample the data an infinite number of times. However, by increasing the number of resamples, we can approximate the ideal case with arbitrary precision. Formally, as $b \rightarrow\infty$, we can control the coverage probability of the confidence regions constructed in Algorithms~\ref{algorithm:wald based statistic},~\ref{algorithm:t based statistic}, and~\ref{algorithm:extrinsic statistic}, as stated in the following corollary.
\begin{corollary}\label{corollary: coverage}
    Suppose Assumption~\ref{assumption:manifold} holds and the initial estimate satisfies $\dist(\theta^{\text{initial}},\theta_0) \leq c_0$ for some sufficiently small $c_0$. 
    \begin{enumerate}
        \item For every $0 < \xi < \frac12 $, it holds for Algorithm~\ref{algorithm:wald based statistic} that 
        \eq{
        \limsup_{b\rightarrow \infty}\Big|\bb P_{\mc X_n}\big[\theta_0\in I_{1-\alpha}^{\mathsf{Wald}}\big] - (1 - \alpha) \Big| \leq O(n^{-1+ \xi}).
        }
        \item  If the retraction mappings $R_{\theta}$ for $\theta$s in a neighborhood of $\theta_0$ satisfy that the function $\mathrm{d}R_{\theta}^{-1}(T_{\theta_0 \rightarrow \theta}(\mathrm{d} R_{\theta_0}(\delta_i)))$ is twice differentiable with respect to $\theta$ for every $i \in[p]$, then for the confidence region in Algorithm~\ref{algorithm:t based statistic}, an arbitrary $\mb a\in \bb R^{p}$, and every $0< \xi< \frac{1}{2}$, it holds that
        \eq{
         \limsup_{b\rightarrow \infty}\Big|\bb P_{\mc X_n}\big[\theta_0\in I_{1-\alpha}^{\mathsf{Intr.t}} \big] - (1 - \alpha) \Big|   \leq O(n^{-1+ \xi}).
        }
        \item It holds for Algorithm~\ref{algorithm:extrinsic statistic} and every $0<\xi<\frac12$ that 
        \eq{
         \limsup_{b\rightarrow \infty}\Big|\bb P_{\mc X_n}\big[f(\theta_0)\in I_{1-\alpha}^{\mathsf{Extr.t}} \big] - (1 - \alpha) \Big|  \leq O(n^{-1+\xi}).
        }
    \end{enumerate}
\end{corollary}

\section{Specific Applications}
\label{sec:applications}
    
    In this section, we showcase our methodology on several commonly encountered Riemannian manifolds that are crucial in various statistical applications. We will demonstrate the implementation of the proposed algorithms across different Riemannian manifolds under a range of statistical scenarios. Section~\ref{subsec: riemannian manifold examples} provides detailed expressions for Riemannian metrics, gradients, Hessian operators, retractions, and inverse retractions. Section~\ref{subsec: stats applications} specifies the loss functions used in the Gaussian location estimation and Barycenter problems.

    Before proceeding, we first introduce a simple method for constructing a \emph{smooth} mapping $f^{\mathsf{ortho},\mb A_0}$ over a neighborhood of a rank-$p_2$ matrix $\mb A_0\in \bb R^{p_1 \times p_2}$. This mapping maps a rank-$p_2$ matrix $\mb A \in \bb R^{p_1\times p_2}$ with $p_1 > p_2>0$ to an orthonormal matrix $f^{\mathsf{ortho},\mb A_0}(\mb A) \in \mathrm{St}(p_1, p_1-p_2)$ such that $\mb A\t f^{\mathsf{ortho},\mb A_0}(\mb A) = \mb 0_{p_2\times(p_1 - p_2)}$. The procedure is outlined as follows:
    \begin{enumerate}
        \item Compute the top-$p_2$ left singular vectors $\mb U_0 = (\mb u_{\mb A_0,1}, \cdots, \mb u_{\mb A_0, p_2})$ of $\mb A_0$. Fix an orthonormal matrix $\mb U_{0,\perp} \in \bb R^{p_1\times(p_1 - p_2)}$ such that $\mb A_0\t \mb U_{0,\perp} = \mb 0_{p_2\times(p_1 - p_2)}$. 
        \item For any $\mb A$ in a sufficiently small neighborhood of $\mb A_0$, compute the top-$p_2$ left singular vectors $\mb U_{\mb A} = (\mb u_{\mb A,1}, \cdots, \mb u_{\mb A, p_2}) \in \bb R^{p_1 \times p_2}$ of $\mb A$. 
        \item Define $f^{\mathsf{ortho},\mb A_0}(\mb A)$ as a matrix that collects the top-$(p_1 - p_2)$ left singular vectors of $(\mb I_{p_1} - \mb U_{\mb A}\mb U_{\mb A}\t )(\mb U_0, \mb U_{0,\perp})\mathrm{diag}(p_1, p_1-1, \cdots, 1)$.
    \end{enumerate}
    From the construction, the column space of $f^{\mathsf{ortho},\mb A_0}(\mb A)$, which is identical to that of $(\mb I_{p_1} - \mb U_{\mb A}\mb U_{\mb A}\t )(\mb U_0, \mb U_{0,\perp})\mathrm{diag}(p_1, p_1-1, \cdots, 1)$, is orthogonal to the column space of $\mb A$ since $\mb A\t (\mb I_{p_1} - \mb U_{\mb A} \mb U_{\mb A}\t) = \mb 0$.
    This mapping will be used in the following examples to satisfy the differentiability requirements necessary for constructing intrinsic $t$-statistics as discussed in Theorem~\ref{theorem:bootstrapmanifold}. 

    Furthermore, as one may deduce from the proofs of the results in Section~\ref{sec: theoretical guarantees}, the coordinate calculation does not necessarily require the inverse of a mapping in the adopted retraction; in fact, those results hold if we replace $\{R^{-1}_x\}$ with an atlas\footnote{An atlas refers to a collection of charts which covers $\mc M$. } $\{\phi_x\}$ of $\mc M$, such that for every $x$,
    \eq{\label{eq: definition of second-order inverse retraction}
    \norm{\phi_{x}(y) - \pi_{\{\mathrm d \phi^{-1} \delta_i\}}\circ \Log_x(y)}_2 = O(\dist(x, y)^3)
    } for every $y$ in the domain of $\phi_x$. In the subsequent discussion, the term "inverse retraction" may refer to either the inverse of a given retraction or a chart in an atlas.
\subsection{Riemannian Manifold Examples}\label{subsec: riemannian manifold examples}
    We begin with a straightforward yet fundamental example. 
    \begin{example}[Sphere Manifold]\label{example: sphere}
        Consider the unit sphere $\bb S^{p-1}$ endowed with the Euclidean metric in $\bb R^p$. This induced metric is identical to the quotient metric. 
        For every $x\in \bb S^{p-1}$, there exists a canonical isometry between the tangent space $\mathrm T_{\mb x} \bb S^{p-1}$ at $\mb x$ and the subspace $\{\mb z\in \bb R^p: \mb z\t \mb x = 0\}$ of $\bb R^p$, therefore we usually represent the tangent vectors in terms of coordinates in $\bb R^p$ with a slight abuse of notations hereinafter. For an arbitrary smooth function $f:\bb S^{p-1} \rightarrow \bb R$, the gradient $\nabla f$ at $\mb x \in \bb S^{p-1}$ is written as 
    \begin{equation}
    	\nabla f(\mb x)= \left( \mb I_p - \mb x \mb x \t \right) \nabla \bar f(\mb x)\in \{\mb z\in \bb R^p: \mb z\t \mb x = 0\} \cong \mathrm T_{\mb x}\bb S^{p-1},
    \end{equation}
    where $\bar f$ is a smooth function on a neighborhood $U$  of $\mb x$ in $\bb R^p$  whose restriction on $\bb S^{p-1}\cap U$ is $f$\footnote{The existence of such an extension is always guaranteed by Lemma 5.34 in \cite{lee2013smooth}, for example.}. The result in \cite{absil2013extrinsic} yields the Hessian operator as
    \begin{equation}
    	\nabla^2 f(\mb x)[\mb \xi] =\big( \mb I_p - \mb x \mb x\t \big)  \nabla^2 \bar{f}(
    	\mb x)[\mb \xi] \in \{\mb z\in \bb R^p: \mb z\t \mb x = 0\} \cong \mathrm T_{\mb x}\bb S^{p-1},
    \end{equation}
    where we again employ the canonical identification between the tangent space and a subspace in $\bb R^p$. 
    
    Finally, we employ the projection retraction $R_{\mb x}: \bb R^{p-1} \rightarrow \bb S^{p-1}$, defined as 
    \eq{\label{eq: retraction on sphere}
    R_{\mb x}(\mb v) = \frac{\mb x + f^{\mathsf{ortho}, \hat \theta_n}(\mb x)\mb v}{\norm{\mb x + f^{\mathsf{ortho}, \hat \theta_n}(\mb x)\mb v}_2},} which has been proved to be a second-order retraction in \cite{gawlik2018high}. The inverse of $R_{\mb x}$ is given by
    \eq{\label{eq: inverse retraction on sphere}
    R^{-1}_{\mb x}(\mb y) = f^{\mathsf{ortho}, \hat \theta_n}(\mb x)\t (\frac{\mb y}{\mb y\t \mb x} - \mb x).}
    The role of $\hat \theta_n$ here is to ensure the smoothness of the retraction in the neighborhood of the true parameter, leveraging the previously established convergence rate of $\hat \theta_n$. 
    By substituting this retraction and its inverse into our framework, we derive the specific algorithms.
    \end{example}

    A generalization of a sphere manifold is the Stiefel manifold. We consider the canonical metric induced by a quotient formulation \citep{tagare2011notes,gawlik2018high}. 
    \begin{example}[Stiefel Manifold] The Stiefel manifold is defined by $\mathrm{St}(p,r) = \{\mb U\in \bb R^{p\times r}: \mb U\t \mb U = \mb I_r\}$ with $p\geq r$. 
    For an arbitrary $\mb U \in \mathrm{St}(p,r)$, we can fix a perpendicular matrix $\mb U_\perp\in \bb R^{p\times (p-r)}$ such that $(\mb U, \mb U_\perp)\t(\mb U, \mb U_\perp) = \mb I_p$ and identify the tangent space at $\mb U$ to the subspace
    \eq{
    \mathrm T_{\mb U}\mathrm{St}(p,r) = \{\mb U \mb A + \mb U_\perp \mb B: \mb A \in \bb R^{r\times r}, \mb A + \mb A\t = \mb 0_{r\times r}, \mb B \in \bb R^{(p-r)\times r)} \}.
    }
    We define a metric at $\mb U$ as $
        \ip{\mb \xi_1, \mb \xi_2} = \frac{1}{2}\mathrm{tr}(\mb A_1\t \mb A_2) + \mathrm{tr}(\mb B_1\t\mb B_2)
    $ with $\mb  \xi_1 = \mb U\mb A_1 + \mb U_{\perp}\mb B_1,\mb \xi_2= \mb U\mb A_2 + \mb U_{\perp}\mb B_2\in \mathrm T_{\mb U}\mathrm St(p,r)
     $. Given a smooth function $f: \mathrm St(p,r)\rightarrow \bb R$, the gradient at $\mb U$ is written as 
     \eq{
      \nabla f_{\mb U} = \nabla \bar f_{\mb U} - \mb U\nabla \bar f_{\mb U}\t \mb U, 
     }
     where $\bar f$ is a smooth extension of $f$ to a local neighborhood of $\mb U$ in $\bb R^{p\times r}$ and $\nabla \bar f$ is written in a corresponding $p$-by-$r$ matrix form. According to \cite{edelman1998geometry}, the Hessian operator is expressed as 
    \eq{
    \nabla^2 f_{\mb U}[\mb \xi]= \nabla^2 f_{\mb U}[\mb \xi] + \frac{1}{2}\big(\mb U \mb \xi\t \nabla \bar f_{\mb U} + \nabla \bar f_{\mb U}\mb \xi\t \mb U) - \frac{1}{2}\big(\mb I_{p} - \mb U \mb U\t\big)\mb \xi\big(\nabla 	\bar f_{\mb U}\t \mb U + \mb U\t \nabla \bar f_{\mb U} \big) \in \mathrm T_{\mb U}\mathrm{St}(p,r)
    }
    for $\mb \xi \in \mathrm T_{\mb U}\mathrm{St}(p,r)$. 
    
    To define retraction mappings, we introduce a second-order retraction on $\mathrm{St}(p, r)$ based on Theorem 5 from \cite{gawlik2018high}. We first introduce a mapping $\mathsf{skew}(\mb v)$ that maps a vector $\mb v\in \bb R^{r(r-1)/2}$ to a skew-symmetric matrix with its upper-triangle part being $\mb v$. Additionally, we let $\mathsf{mat}_{p_1\times p_2}(\mb v) \in \bb R^{p_1\times p_2}$ be the matricization of a vector $\mb v \in \bb R^{p_1p_2}$ 
    Given $\mb U \in \mathrm{St}(p, r)$, the second-order retraction $R_{\mb U}$ is expressed as:
      \eq{
        R_{\mb U}(\mb v) = \mc P\left(\mb U - \frac{1}{3}\mb H\t \mb H -\frac{1}{2}\mb U\t \mb H \mb U \mb H +\mb H\left(\mb I + \frac{1}{2}\mb U\t \mb H\right)  \right)
    }
    for $\mb v= (\mb v_1\t, \mb v_2\t)\t  \in \bb R^{r(r-1)/2 + (p-r)r}$ with $\mb v_1 \in \bb R^{r(r-1)/2}$ and $\mb v_2 \in \bb R^{(p-r)r}$. Here, the matrix $\mb H$ is defined as  
    \eq{\label{eq: tangent vector of stiefel manifold}
    \mb H = \mb U\mathsf{skew}(\mb v_1) + f^{\mathsf{ortho}, \hat \theta_n}(\mb U)\mathsf{mat}_{r\times(p-r)}(\mb v_2)
    }
    and the unifying operator $\mc P:\bb R^{p\times r}\rightarrow \mathrm{St}(p,r)$ is defined as $\mc P: \mb X \mapsto \mb L\mb R\t$ where the singular value decomposition of $\mb X$ is written as $\mb X = \mb L\mb D \mb R\t$ with $\mb L \in \mathrm{St}(p, r)$ and $\mb R \in \mathrm{St}(r,r)$.  For the inverse retraction, we recommend using the logarithmic mapping approximation proposed in \cite{mataigne2025efficient}, in conjunction with the expression in \eqref{eq: tangent vector of stiefel manifold}.

    \end{example}

    It is worth noting that, while computationally tractable retractions exist for the examples discussed above, such explicit forms may not be always available for general Riemannian manifolds. Nonetheless, in what follows we will focus on Euclidean submanifolds, as a convenient framework that offers the (inverse) retraction mappings we require.   
    \paragraph*{Euclidean Submanifold}
    A Euclidean submanifold refers to a submanifold of a Euclidean space, whose importance could be revealed from the following two distinct perspectives: 
    	\begin{itemize}
    		\item\textbf{Level Sets of Constant Rank Function.}
    		    Suppose that the parameter space $\mc M$ is determined by $\mc M = \{\mb x\in \bb R^{d}:g_i(\mb x) = 0,i=1,\ldots, d- p\}$ with $d >p$. If the mappings $g_i$ are smooth and $\mb g = ( g_1, \ldots, g_{d - p})$ has a constant rank $d - p$, by the constant rank level set theorem (Theorem 5.12 in \cite{lee2013smooth}) we deduce that $\mc M$ is an embedded submanifold of $\bb R^{d}$. 
    
    		\item\textbf{Whitney's and Nash's Embedding Theorems.} Even though we are working with a smooth manifold with an abstract and complicated differential structure, the Whitney embedding theorem shows that an arbitrary smooth manifold of dimension $n$ always admits a proper embedding into $\bb R^{2n+1}$. Furthermore, the Nash embedding theorem establishes that a Riemannian manifold can be embedded into a Euclidean space in a manner that preserves its metric structure. 
            In view of uncertainty quantification concerning parameter location, it thus suffices to establish inference procedures under an embedding mapping and its Euclidean metric. 
    	\end{itemize} 
    	
        For ease of presentation, we use the canonical representation to identify a tangent vector in the tangent space of a submanifold with a vector in $\bb R^d$. This identification similarly extends to normal vectors.
    		We consider a Riemannian manifold of dimension $p$ embedded in $\bb R^{d}$. By Theorem 5.8 in \cite{lee2013smooth}, for every $x \in M$ there exists a chart $(U, \phi)$ for $\bb R^{d}$ centered at $x$ such that $\phi(U\cap \mc M)=\{(x_1, \dots, x_p, 0, \ldots,0) \in U\}$. 
    	We denote by $\tilde \phi$ the last $d - p$ coordinates of $\phi$ and then define the projection operator on $\bb R^d$ as 
        \eq{\mc P_{\mb x} \coloneqq    \mb I_{d} - \nabla\tilde{\phi} \left(\nabla \tilde \phi\t \nabla \tilde \phi\right)^{-1
            }\nabla \tilde \phi\t,
            \label{eq: projection operator}
            } where $\nabla \tilde \phi \coloneqq (\frac{\partial }{\partial x_i}\tilde \phi_j)_{ij}\in \bb R^{d\times(d - p)}$ with a slight abuse of notation. 	    
    		     Given a smooth function $\bar f $ on $\bb R^{d}$, we denote the function whose domain is constrained on $\mc M$ by $f$. Then, following equations (3) and (7) of \cite{absil2013extrinsic}, the gradient and the Hessian of $f$ on $\mc M$ are expressed as
                 \begin{align}
                 & \nabla f = \mc P_{\mb x}\nabla \bar f,
                 \label{eq: gradient expression} \\ 
                 & \nabla^2 f(\mb x)[v] = \mc P_{\mb x}(\nabla^2\bar{f}(\mb x)[v])+ \mc P_{\mb x}\big((\nabla_{v} \mc P_{\mb x}) [\nabla \bar f]\big), 
            \label{eq: hessian expression}
                 \end{align}
                 where $ \nabla  \bar f$ denotes the gradient of $\bar f$ over $\bb R^{d}$.

    	    \paragraph*{Second-Order inverse retraction of Euclidean Submanifold}
    	Now we are in a position to acquire a viable method to approximate the logarithmic mapping so as to obtain a coordinate representation for Euclidean submanifolds. We refer to the inverse retractions in an atlas that satisfy \eqref{eq: definition of second-order inverse retraction} as the second-order inverse retractions.
        Inspired by the projection-like retraction proposed in \cite{absil2012projection}, we introduce projection-like inverse retractions as a means of approximating logarithmic mappings for general Euclidean submanifolds.
        
        We first introduce the concept of inverse retractor. Firstly, two subspaces of dimension $p$ and $(d-p)$ in $\bb R^d$ are said to be transverse, if their direct sum is $\bb R^d$. 
        \begin{definition}[Inverse Retractor]
        	Let $\mc M$ be a $p$-dimensional smooth submanifold of a $d$-dimensional Euclidean space. An inverse retractor on $\mc M$ is a smooth mapping from a neighborhood $U \times U \in \mc M \times \mc M$ to the Grassmannian manifold $\mathrm{Gr}(d - p, \bb R^d)$, which consist of all $(d-p)$-dimensional linear subspaces of $\bb R^d$, such that the subspace $D( x, x)$ is transverse to $\mathrm{T}_x\mc M$ for all $x\in U$. 
        \end{definition}
        
        An inverse retraction $L_x$ at $x$ is called second-order if it satisfies any of the following equivalent conditions:
        \begin{itemize} 
            \item  $\frac{\mathrm d^2}{\mathrm dt^2} L_x(\Exp_x(tv)) = 0$ for every $v \in \mathrm T_x\mc M$;
        \item $\Log_x(y) = L_x(y) + O(\dist(x,y)^3)$ as $\dist(x,y) \rightarrow 0$. 
    \end{itemize}
        Given an inverse retractor $D$ defined on $U \times U$ and $x \in U$, we define an inverse retraction centered at $x$ as follows:
        \begin{equation}   
        \label{eq:inverser-retraction}
        	L_{x}: U \rightarrow \mathrm T_{x}\mc M; \quad L_{x}(y) = v\text{ such that }x + v - y\in \mathrm  D(x, y).
        \end{equation}      
        Building on the work \cite{absil2012projection}, which established the second-order consistency of a retraction induced by a retractor \cite[Definition 14]{absil2012projection}, we justify the second-order consistency of $L_{x}$ induced by an easy-to-compute inverse retractor in the following proposition.  This can be viewed as an inverse counterpart to the projection-like retraction in \cite{absil2012projection}. 
        
        \begin{proposition}\label{proposition:inverse retraction}
        	We let $L_{x}$ be a mapping defined above with an inverse retractor $D$ satisfying that $D(x, y) = \mathrm N_{x} \mc M$ for all $x,y\in U$, where $\mathrm N_{x} \mc M$ denotes the normal space of $\mc M$ at $x$. Then $L_{x}$ is a second-order inverse retraction. 
        \end{proposition}

    By combining the elements established above, we can proceed to calculate the necessary derivative and coordinate expressions across a broad spectrum of Euclidean submanifolds.
    	\begin{example}[Example: Manifold of Fixed-Rank Matrices]\label{example: fixed-rank matrices manifold}
    	The fixed-rank matrix manifold is a well-studied matrix manifold \citep{mishra2012riemannian,luo2024geometric}, which could be viewed as a submanifold of the Euclidean space. We consider the matrix manifold $\mc R_{r,p_1,p_2}$ composed of all rank-$r$ matrices in $\bb R^{p_1\times p_2}$, whose dimension is $\left(p_1r+p_2r-r^2\right)$. For a matrix $\mb X \in \mc R_{r,p_1,p_2}$ with the singular value decomposition 
        $
        \mb X = \mb U\left[\begin{matrix}
            \mb \Sigma & \mb 0\\ 
            \mb 0 & \mb 0
        \end{matrix}\right]\mb V\t
        $, every element $\mb Z$ of $\mathrm T_{\mb X}\mc R_{r,p_1,p_2}$ could be written as 
        \eq{
        \mb Z = \mb U\left[\begin{matrix}
            \mb A & \mb C\\ 
            \mb B& \mb 0
        \end{matrix}\right]\mb V\t, \quad \mb U = \left[ \mb U_1,\mb U_2\right], \mb \quad \mb V = \left[ \mb V_1, \mb V_2\right],
        }
        where $\mb U\t \mb U = \mb I_{p_1}$ and $\mb V\t \mb V = \mb I_{p_2}$. Note that, if identical nonzero singular values exist in the above decomposition, the choice of $\mb U_1$, $\mb V_1$ becomes non-unique. Moreover, the decomposition of the matrices in the neighborhood of such a matrix is not smooth, potentially leading to the irregular asymptotic distribution shown in Figure~\ref{fig: counterexample}. Therefore, we assume the nonzero singular values of $\mb X$ are distinct hereinafter.

        Given a smooth function $f\in \mc R_{r, p_1, p_2}$ with an extension $\bar f$ over $\bb R^{p_1\times p_2}$, the gradient of $f$ at $\mb X$ is given by 
        \eq{
        \nabla f_{\mb X} = \nabla \bar f(\mb X) - \mb U_2\mb U_2\t \nabla \bar f(\mb X) \mb V_2 \mb V_2\t. 
        }
        By Section 4.5 in \cite{absil2013extrinsic}, the Hessian operator of $f$ at $\mb X$ is given by 
        \eq{
        \nabla^2 f_{\mb X}[\mb \xi] = \nabla^2 \bar f_{\mb X}[\mb \xi] - \mb U_2 {
        \mb U_2}\t  \nabla^2 \bar f_{\mb X}[\mb \xi] \mb V_2{\mb V_2}\t +\nabla \bar f_{\mb X}\mb \xi\t(\mb X^\dagger)\t + (\mb X^\dagger)\t \mb\xi\t \nabla \bar f. 
        }
        By Proposition 4.11 in \cite{absil2012projection}, the mapping $R_{\mb X}$ centerd at $\mb X$ given by 
        \eq{
        R_{\mb X}(\mb v) = \mb U\left[\begin{matrix}
            \mb A & \mb C\\ 
            \mb B& \mb B(\mb \Sigma + \mb A)^{-1} \mb C
        \end{matrix}\right]\mb V\t
        }is a second-order retraction at ${\mb X}$. Here, $\mb v = (\mb v_1\t, \mb v_2\t, \mb v_3\t)\t \in \bb R^{p_1p_2 - (p_1 - r)(p_2-r)}$ is transformed into $\mb A = \mathsf{mat}_{r\times r}(\mb v_1)$, $\mb B = \mathsf{mat}_{(p_1 - r)\times r}(\mb v_2)$, and $\mb C = \mathsf{mat}_{r\times (p_2- r)}(\mb v_3)$. And the matrices $\mb U_2$ and $\mb V_2$ in the definitions of $\mb U$, $\mb V$ are specified to be $\mb U_2 = f^{\mathsf{ortho}, \mb U_{\hat \theta_n}}(\mb U_1)$ and  $\mb V_2 = f^{\mathsf{ortho}, \mb V_{\hat \theta_n}}(\mb V_1)$, where $\mb U_{\hat \theta_n}$ and $\mb V_{\hat \theta_n}$ refer to the left and right top-$r$ singular vectors of $\hat \theta_n\in \mc R_{r,p_1,p_2}$, respectively. 
        
        On the other hand, for any $\mb Y \in \mc R_{r,p_1,p_2}$ close to $\mb X$ enough, the inverse of $\mb R_X$ is 
        \eq{
        R_{\mb X}^{-1}(\mb Y) = (\mathsf{vec}(\mb A)\t, \mathsf{vec}(\mb B)\t, \mathsf{vec}(\mb C)\t)\t, 
        }given $\mb Y = \mb U\left[\begin{matrix}
            \mb A & \mb C\\ 
            \mb B& \mb D
        \end{matrix}\right]\mb V\t$, where $\mathsf{vec}$ denotes the vectorization operator. 
        By Proposition~\ref{proposition:inverse retraction}, $R_{\mb X}^{-1}$ is a second-order inverse retraction.

        \end{example}
    \begin{example}[Example: Manifold of Rank-one Tensors]
    The set $\mc M_{p_1,p_2, \ldots	, p_k}^{(1)}$ of Tucker tensors \cite{kasai2016low} of fixed rank $\mb r = \left(1, \ldots, 1\right)$ in $\bb R^{p_1\times p_2\times \cdots\times  p_k}$ also form a Riemannian submanifold in the Euclidean space.  The tangent space of $\mc M_{p_1,p_2, \ldots, p_k}^{(1)}$ at $\mb X = x_0\times\mb u_1\times\mb u_2 \times\cdots \times \mb u_k $ can be expressed as 
    \begin{equation}
    \label{eq:rank-one tensor tangent space}
    	\mathrm{T}_{\mb X}\mc M_{p_1,p_2, \ldots	, p_k}^{(1)} = \Big\{\mb v\big| \mb v = a_0 \underset{i\in[k]}{\times}\mb u_i +\sum_{i\in[k]} 1 \underset{j\neq i}{\times}\mb u_j \times_i\mb u_i',\quad \text{with }{\mb u_i'}\t \mb u_i = 0,a_0\in \bb R\Big\}.
    \end{equation}
    Provided with a smooth function $f$ and $\xi \in \mathrm T_{\mb X}\mc M^{(1)}_{p_1, \ldots, p_k}$, the gradient and the Hessian operator at $\mb X$ can be expressed as follows \citep{heidel2018riemannian}
    \begin{align}
    	& \nabla f_{\mb X} = \sum_{i=1}^k \nabla \bar f_{\mb X} \underset{i\in[k]}{\times}\mb u_i \mb u_i\t +\sum_{i\in[k]} \nabla \bar f_{\mb X}\underset{j\neq i}{\times}\mb u_j\mb u_j\t \times_i\mb U_{i,\perp}\mb U_{i, \perp}\t  \\ 
    	& \nabla^2 f_{\mb X}[\mb \xi] =  \sum_{i=1}^k \nabla^2 \bar f_{\mb X}[\mb \xi] \underset{i\in[k]}{\times}\mb u_i \mb u_i\t +\sum_{i\in[k]} \nabla^2 \bar f_{\mb X}[\mb \xi]\underset{j\neq i}{\times}\mb u_j\mb u_j\t \times_i\mb U_{i,\perp}\mb U_{i, \perp}\t \\ 
    	 & + \sum_{i\in[k]}a_0 \times_i \Big(\mb U_{i, \perp}\mb U_{i, \perp}\t\sum_{k \neq i}\big(\nabla \bar f_{\mb X}\times_k \mb u_k' \underset{k\neq j\neq i}{\times}\mb  u_j \big)  \Big)\underset{j\neq i}{\times} \mb u_j \in \mathrm T_{\mb X}\mc M^{(1)}_{p_1, \ldots, p_k}.
    \end{align}
    We then consider a projective-like retraction, inspired by the estimator proposed in \cite{zhang2019cross}, defined as follows:
    \begin{equation}
    	R_{\mb X}(\mb v) = (x_0 + a_0)\underset{i\in[k]}{\times}\left(\mb u_i + \frac{x_0}{x_0+a_0}\mb u_i'\right),
    \end{equation}
    where $\mb v = (a_0, \mb v_1\t, \cdots, \mb v_k\t)\t$, $\mb u_i = \mb U_{i,\perp}\mb v_i$, and $\mb U_{i,\perp}$, $i\in[k]$ are specified as $\mb U_{i,\perp} = f^{\mathsf{ortho}, \mb u_{i,\hat\theta_n}}(\mb u_i)$, $i \in [k]$, given $\hat \theta_n = a_{0, \hat \theta_n} \times_{i\in[k]}\mb u_{i,\hat \theta_n}$. In view of Theorem 22 in \cite{absil2012projection}, this retraction is a second-order retraction. 
    We then have the following second-order inverse retraction by Proposition~\ref{proposition:inverse retraction} 
    \begin{equation}
    	R^{-1}_{\mb X}(\mb Y) = (a_0, \mb v_1\t, \cdots, \mb v_k\t)\t,
    \end{equation}
    where
    $$
     a_0 = \mb Y\underset{i\in[k]}{\times} \mb u_i\t, \quad \mb v_i =  \frac{1}{x_0}\mb U_{i,\perp}\t\mb Y\underset{j\neq i}{\times} \mb u_i\t.
    $$

    \end{example}
    \subsection{Statistical Applications}
    \label{subsec: stats applications}
    We now turn to two classic statistical problems that naturally align with manifold formulations. Note that all these manifolds can be written as submanifolds in $\bb R^d$s, with the $d$s being self-evident from the preceding discussion.
    \subsubsection{Maximum Likelihood Estimation under Curved Exponential Family}\label{subsubsec: gaussian location estimation}
    The notion of curved exponential family provides a classic example of how geometric structrues naturally arise in statistical modelling \cite{efron1975defining,efron1978geometry,ay2017information,amari1983differential}. A family of distributions $\tilde{\mc M}$ is said to form a curved exponential family when it constitutes a submanifold of an exponential family of the form $\big\{p(\bo \theta, \mb x) = \exp\big(\bar \gamma(\mb x) + \bo \theta\t \mb x - \bar \psi(\bo \theta)\big);~\bo \theta\in \bb R^d \big\}$. Since each parameter $\boldsymbol{\theta}$ corresponds uniquely to a distribution in the family, we represent the parameter submanifold (equipped with the Euclidean metric) by $\mc M$. The objective function is then given by $L(\bo \theta, \mb X) \coloneqq -\bar \gamma(\mb X) - \bo \theta\t \mb X + \bar \psi(\bo \theta)$ with $\bo \theta\in \mc M$. Using the expressions \eqref{eq: gradient expression} and \eqref{eq: hessian expression}, we obtain
    \begin{align}
        & \nabla L(\bo \theta, \mb X) =\mc P_{\bo \theta} (\nabla \bar \psi(\bo \theta) - \mb X) , \\ 
        & \nabla L^2(\bo \theta, \mb X)[v] =   \mc P_{\bo \theta}(\nabla^2\bar{\psi}(\bo \theta)[\mb v])+ \mc P_{\bo \theta}\big((\nabla_{v} \mc P_{\bo \theta}) [\nabla \bar \psi - \mb X ]\big) 
    \end{align}
    for $v \in \mathrm T_{\bo \theta}\mc M$, where $\mc P_{\bo \theta}$ denotes the projection operator in \eqref{eq: projection operator}. 

    As a specific instance of the curved exponential family, consider the estimation of a location parameter in additive Gaussian noise models. Suppose we observe $\boldsymbol X_i = \boldsymbol\theta_0 + \boldsymbol\varepsilon_i$ where the true parameter $\boldsymbol{\theta}_0 \in \mc M \subset \bb R^d$ is unknown and $\boldsymbol{\varepsilon}_i \sim N(0, \mb I_d)$ are independent Gaussian noise terms. 
    The resulting objective function, derived from the maximum likelihood estimation, is given by
    \eq{
        L(\boldsymbol \theta, \mb X) \coloneqq \frac{1}{2}\norm{\mb X- \boldsymbol \theta}_F^2 \text{\quad with \quad }  \boldsymbol \theta \in \mathcal M. 
    }
    The sample average version $L_n$ is written as 
    $$
    L_n(\boldsymbol \theta) = \frac{1}{2n} \sum_{i\in[n]} \norm{\mb X_i - \boldsymbol \theta}_F^2 \text{\quad with \quad }  \boldsymbol \theta \in \mathcal M. 
    $$
    Extending $L_n$ to the whole Euclidean coordinate space, we naturally obtain that $\bar L_n(\boldsymbol \theta) = \frac{1}{2n} \sum_{i\in[n]} \norm{\mb X_i - \boldsymbol \theta}_F^2$ with $\boldsymbol \theta \in \bb R^d$. 

    With the formulations of gradients and Hessian operators in place of Section~\ref{subsec: riemannian manifold examples}, we substitute $\bar L_n$ into those expressions and collectively list the expressions in Table~\ref{tab:gradient and Hessian}, where the notations follow the conventions in Section~\ref{subsec: riemannian manifold examples}. 

    \begin{table}[htbp]
        \centering
        \resizebox{\columnwidth}{!}{
        \begin{tabular}{c|c|c}
        \toprule
        Parameter Space& Gradient $\nabla L_n(\boldsymbol \theta)$ & Hessian $\nabla^2 L_n(\boldsymbol\theta)[\mb \xi]$ \\ 
        \midrule
         Sphere Manifold $\bb S^{p-1} $ & $(\mb I_p - \mb x \mb x\t)\big(\boldsymbol \theta - \bar{\mb X}_n \big)$ & $(\mb I_p - \mb x \mb x\t) \mb \xi$ \\ 
         \midrule
        Stiefel Manifold $\mathrm{St}(p,r) $ & $\big(\boldsymbol \theta - \bar{\mb X}_n \big) - \mb U\big(\boldsymbol \theta - \bar{\mb X}_n \big)\t\mb U$ & \makecell{$\mb \xi + \frac{1}{2}\big(\mb U \mb \xi\t \big(\boldsymbol \theta - \bar{\mb X}_n \big) + \big(\boldsymbol \theta - \bar{\mb X}_n \big)\mb \xi\t \mb U)$ \\$ - \frac{1}{2}\big(\mb I_{p} - \mb U \mb U\t\big)\mb \xi\big(\big(\boldsymbol \theta - \bar{\mb X}_n \big)\t \mb U + \mb U\t \big(\boldsymbol \theta - \bar{\mb X}_n \big) \big)$ }\\ 
        \midrule
        \makecell{fixed-rank matrix manifold\\ $\mathcal R_{r,p_1,p_2}$} & $\big(\boldsymbol \theta - \bar{\mb X}_n \big) - \mb U_2\mb U_2\t \big(\boldsymbol \theta - \bar{\mb X}_n \big) \mb V_2 \mb V_2\t$ & \makecell{$\mb \xi - \mb U_2 {
        \mb U_2}\t  \mb \xi \mb V_2{\mb V_2}\t + \big(\boldsymbol \theta - \bar{\mb X}_n \big)\mb \xi\t(\bo \theta^\dagger)\t$ \\ 
        $+ (\bo \theta^\dagger)\t \mb\xi\t\big(\boldsymbol \theta - \bar{\mb X}_n \big)$} \\ 
        \midrule
        \makecell{Rank-$1$ Tensors Manifold\\ $\mathcal M^{(1)}_{p_1,p_2, \cdots, p_k} $} &  \makecell{$\sum_{i=1}^k \big(\boldsymbol \theta - \bar{\mb X}_n \big) \underset{i\in[k]}{\times}\mb u_i \mb u_i\t$ \\ 
        $+\sum_{i\in[k]} \big(\boldsymbol \theta - \bar{\mb X}_n \big)\underset{j\neq i}{\times}\mb u_j\mb u_j\t \underset{i}{\times}\mb U_{i,\perp}\mb U_{i, \perp}\t $} 
        & \makecell{$\sum_{i=1}^k \mb \xi \underset{i\in[k]}{\times}\mb u_i \mb u_i\t +\sum_{i\in[k]} \mb \xi\underset{j\neq i}{\times}\mb u_j\mb u_j\t \underset{i}{\times} \mb U_{i,\perp}\mb U_{i, \perp}\t$ \\ 
    	 $ + \sum_{i\in[k]}a_0 \underset{i}{\times}  \Big(\mb U_{i, \perp}\mb U_{i, \perp}\t\sum_{k \neq i}\big(\big(\boldsymbol \theta - \bar{\mb X}_n \big)\underset{k}{\times} \mb u_k' \underset{k\neq j\neq i}{\times} \mb u_j \big)  \Big)\underset{j\neq i}{\times} \mb u_j $} \\ 
        \bottomrule
        \end{tabular}
        }
        \caption{Gradient and Hessian Expressions for Different Manifolds. }
        \label{tab:gradient and Hessian}
    \end{table}

    \subsubsection{Barycenter Problem} \label{subsubsection: barycenter}
    Barycenter problems, serving as a specific case of M-estimations on manifolds, attracted much attention in the past decades \citep{frechet1948elements,karcher1977riemannian,bhattacharya2003large,bhattacharya2005large,pennec2019curvature,lee2025general}. To be specific, given $n$ manifold-valued samples $\{x_i\}$ independently drawn from a distribution on a Riemannian manifold $\mc M$, the barycenter is defined as the minimizer $x\in \mc M$ of 
    \begin{equation}
    	L_n(x) = \frac{1}{2n}\sum_{i\in[n]}\dist(x_i,x)^2. 
    \end{equation}
    While previous studies have primarily examined the geometrical and statistical properties of the barycenter, we revisit the barycenter problem from an optimization perspective. Algorithm~\ref{algorithm: resampled newton iteration}  together with the algorithms in Section~\ref{subsubsec: algorithms for confidence region constructions} yields a computationally feasible approach to construct confidence regions for the Barycenter with high accuracy.

        	Next, we focus on the sphere case $\mc M = \bb S^{p-1}$. According to \cite{pennec2018barycentric}, the expressions of gradient and Hessian operator at  $\mb x$ are given by
       	\longeq{
        & \nabla {L_n}(\mb x) = \frac{-\sum_{i\in[n]} \Log_{\mb x}\mb x_i}{n} = -\frac{1}{n}\sum_{i\in[n]} \frac{\arccos(\mb x\t \mb x_i)}{\sqrt{1 - (\mb x\t \mb x_i)^2}}(\mb I_p - \mb x \mb x\t ) \mb x_i,\\ 
        & \nabla^2 {L_n}(\mb x)[\mb v] =  \frac{1}{n} \sum_{i\in[n]}\frac{(\mb I_p - \mb x \mb x\t ) \mb x_i\mb x_i\t(\mb I_p - \mb x \mb x\t )\mb v }{(1 - (\mb x\t \mb x_i)^2)} \\+& \frac{1}{n}\sum_{i\in [n]} \frac{\arccos(\mb x_i\t \mb x)}{\sqrt{1 - (\mb x\t \mb x_i)^2}}(\mb x_i\t \mb x) \Big(\mb I_p - \mb x \mb x\t - \frac{(\mb I_p - \mb x \mb x\t ) \mb x_i\mb x_i\t(\mb I_p - \mb x \mb x\t ) }{(1 - (\mb x\t \mb x_i)^2)}\Big) \mb v
        }
        for $\mb v \in \mathrm T_{\mb x}\bb S^{p-1}$.   As outlined in Example \ref{example: sphere}, we employ the second-order retraction \eqref{eq: retraction on sphere} in the optimization procedure and use its inverse \eqref{eq: inverse retraction on sphere} to calculate the coordinate statistics.

\section{Real Data Analysis}
\begin{figure}[hbtp]
            \centering 
            \subfigure{\includegraphics[width = 0.45\textwidth]{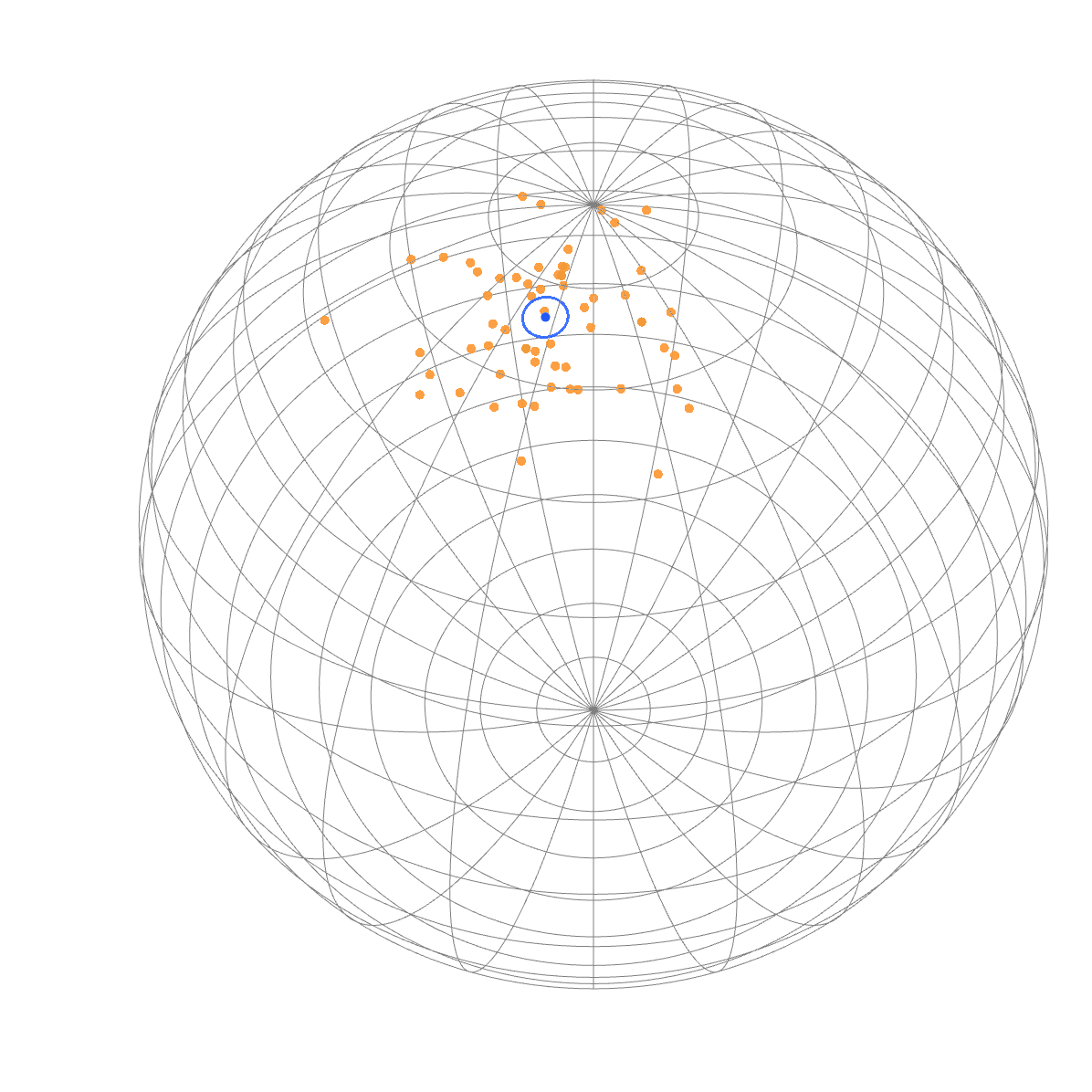}}
            \subfigure{\includegraphics[width = 0.45\textwidth]{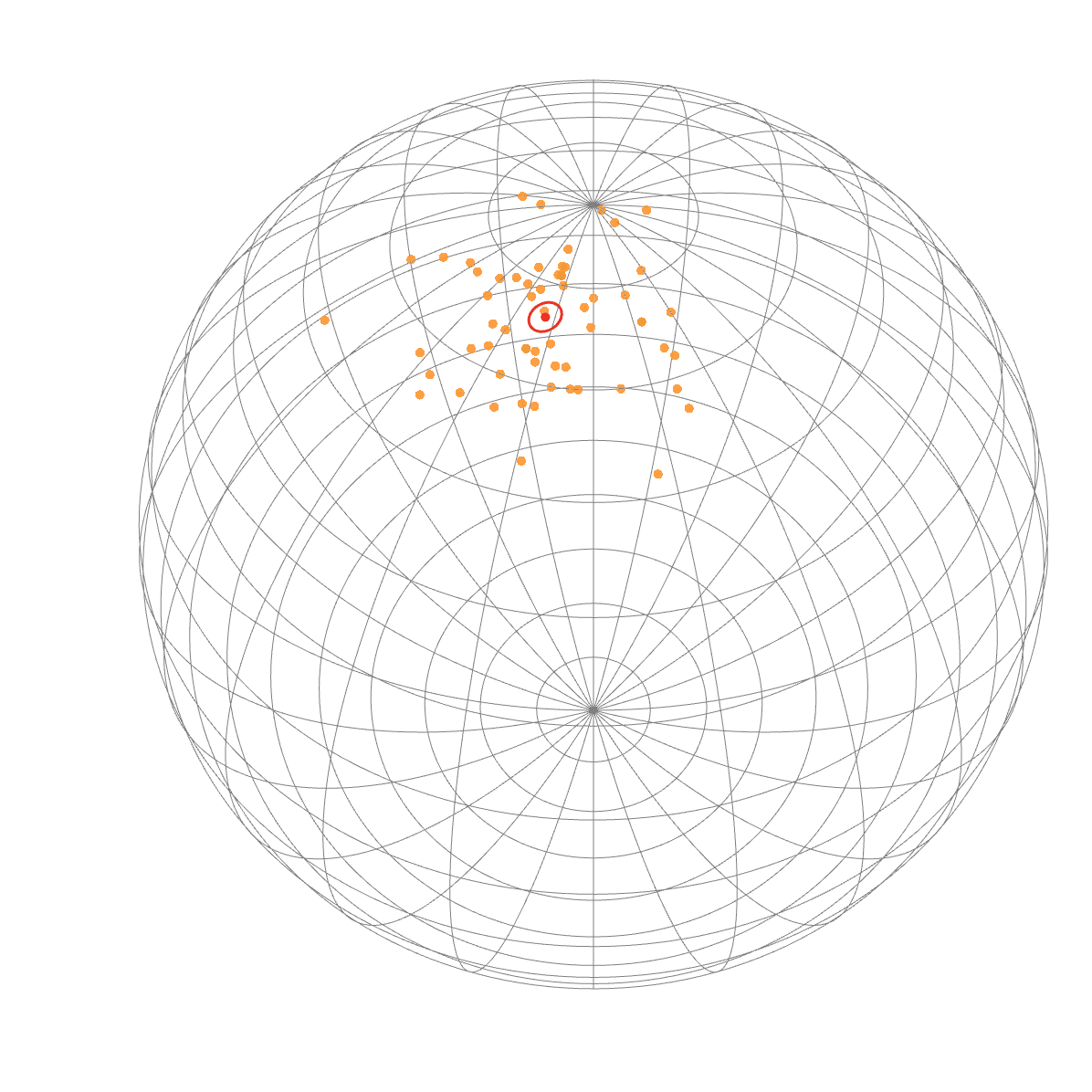}}
    \caption{Estimates and Confidence Regions for the Paleomagnetic Pole Position. The orange points represent the magnetic remanence directions of the samples.  In the left-hand panel, the blue point denotes the estimated pole position, and the surrounding blue circle represents the confidence region constructed in \cite{paindaveine2020inference}. The right-hand panel displays our corresponding estimate and bootstrap-based confidence region in red. }
            \label{fig: real data}
    \end{figure}

We analyze the proposed confidence region for a spherical dataset representing the magnetic remanence of 62 rock specimens collected from Prospect, New South Wales. Each remanence vector records the direction and intensity of Earth’s magnetic field at the time the rocks acquired their magnetization. By examining these remanence directions, we can infer the orientation of the ancient geomagnetic field and estimate the corresponding paleomagnetic pole position. This dataset was provided in \cite[Appendix~B8]{Fisher_Lewis_Embleton_1987}. Its inference procedure with respect to the center parameter $\bo \theta$ was studied in \cite{paindaveine2020inference} under the celebrated Fisher-von Mises-Langevin (FvML) distribution. 
Here we study the barycenter of the underlying distribution instead, following the procedures in Section~\ref{subsubsection: barycenter}, which does not rely on the prior assumption about the distribution form. The resulting confidence region, shown in Figure~\ref{fig: real data}, is compared with that of \cite{paindaveine2020inference}. Our regions appear narrower and better adapted to the rotational asymmetry exhibited by the data.

\section{Conclusion and Discussions}
In this study, we presented a novel framework for high-order accurate inference on Riemannian manifolds, addressing the challenges faced in statistical inference when the parameter spaces are not Euclidean. This framework encompasses the development of computationally friendly bootstrapping algorithms, high-order asymptotics for specific statistics under curvature effect, and the exploration of a comprehensive treatment across various modern statistical settings, including spheres, the Stiefel manifold, Euclidean submanifolds, and specific applications like fixed-rank matrices and rank-one tensors manifolds. Compared with existing works on high-order asymptotics of specific statistics, this work elucidates some essential ingredients for treating nonconvex problems in the scope of manifold language, especially shedding light on the coordinate transformation issue. 

Looking ahead, multiple avenues for future research emerge. Firstly, regarding high-order asymptotics, the symmetric case achieves higher accuracy than the asymmetric case. However, this advantage diminishes in the manifold setting due to the triangle error introduced by the double exponential mapping. Therefore, a correction accounting for local curvature might be needed to achieve the same accuracy as in Euclidean spaces. 

Our current focus is on the fixed manifold case, where the sample size $n$ tends to infinity but dimension $p$ remains fixed. When $p$ increases along with $n$, the necessary conditions on $p$ and the loss function $L$ become significantly more complex and stringent, and the bootstrap may fail as indicated in \cite{el2018can}. Therefore, this study concentrates on the fixed-dimensional case with $n\to \infty$. Given the recent advances in high-dimensional statistical theory, it would be of broad interest to investigate the curvature effect when the dimensions of a sequence of manifolds also grow with the sample size. 

Moreover, our study highlights an additional application of the double exponential mapping, previously employed to investigate curvature effects in the barycenter problem \citep{pennec2019curvature}. Our technique may offer insights to various phenomena in Riemannian optimization with higher precision.

\bibliographystyle{imsart-number} 
\bibliography{ref}       

\begin{appendix}

\section{Related Literature}\label{sec:related-literature}

This work intersects with several key areas of research, including Riemannian optimization, high-order asymptotics, and manifold inference.

Firstly, the extension of classic optimization methods \citep{smith1993geometric,smith2014optimization,huper2004newton,helmke2000jacobi,dennis1996numerical} to Riemannian manifold settings \citep{absil2009optimization} has attracted increasing attention, driven by the emergence of large-scale matrix data \citep{luo2024recursive,luo2022nonconvex} and tensor data \citep{luo2023low}.

Next, the study of high-order asymptotics dates back to the efforts in refining the central limit theorem for fixed-dimensional cases. Foundational results such as \cite{bhattacharya1978validity,babu1984one,bhattacharya1987some,hall1996bootstrap,hall2013bootstrap} established the Edgeworth expansion and its bootstrap adaptations for empirical averages and functions of weighted sums. Additionally, \cite{andrews2002higher} introduced computationally efficient algorithms combined with studentized statistics to achieve high-order asymptotics for M-estimation problems. These developments serve as a foundation for extending such methods to manifold settings. There is another series of research focusing on the Edgeworth expansion of U-statistics \citep{callaert1980edgeworth,bickel1986edgeworth,lai1993edgeworth}, especially on network moment \citep{zhang2022edgeworth,shao2022higherorder}. 

Moreover, research in manifold inference covers diverse geometrical approaches, with a focus on Barycenter problems \citep{bhattacharya2003large,bhattacharya2005large,eltzner2019smeary}. Recent advancements have addressed uncertainty quantification for low-rank models, offering perspectives tailored to high-dimensional settings \citep{chen2019inference,xia2021statistical,xia2022inference,xie2024higherorder}.

    \section{Simulation Studies}
    \label{sec: numerical simulations}
    In this section, we validate our methods and theoretical results through a range of numerical simulations. 
    \paragraph*{Simulation Settings Overview}
    We set out by providing an overview of all simulation settings. The first three simulations proceeded in an additive Gaussian noise framework as discussed in Section~4.2.1 and the last simulation investigates the application of our methods to the Barycenter problem on a sphere manifold. The precise parameter setups are as follows: 
    \begin{itemize}
    \item \textbf{Setting 1: $\bb S^2$. } We consider the sphere manifold $\bb S^2$ embedded in $\bb R^3$ and set the ground-truth parameter as $\boldsymbol \theta_0 = (0, 1, 0)^{\top}$. 
        \item 	\textbf{Setting 2: $\mathrm{St}(4,2)$.}
        We set the ground-truth parameter to be $\boldsymbol \theta_0 = \left(\begin{matrix}
            1/\sqrt{2}& 1/\sqrt{2} \\ 
            1/\sqrt{2} & -1/\sqrt{2} \\ 
            0 & 0 \\ 
            0 & 0
        \end{matrix}\right) \in \mathrm{St}(4,2)$. 
    \item
    \textbf{Setting 3: $\mathcal R_{2,4,4}$. }
	We consider the fixed-rank matrix manifold to be the collection of $4$-by-$4$ matrices with rank $2$ and let the ground truth $\mb \Theta_0$ to be the matrix with all zero entries except for the first two diagonal entries, which are set to $5$ and $1$, respectively.
    \item \textbf{Setting 4: $\mathcal M^{(1)}_{3,3,3}$. }
	We consider the rank-one tensors manifold embedded in $\bb R^{3\times 3 \times 3}$. The ground truth parameter is set to be $\mb \Theta = \mb u_1\otimes \mb u_1\otimes \mb u_1$ with $\mb u_1 = (1, 0, 0)^{\top}$. 
    \end{itemize}

\begin{figure}[hbtp]
	\centering
	\subfigure{
		\includegraphics[width = 0.31\textwidth]{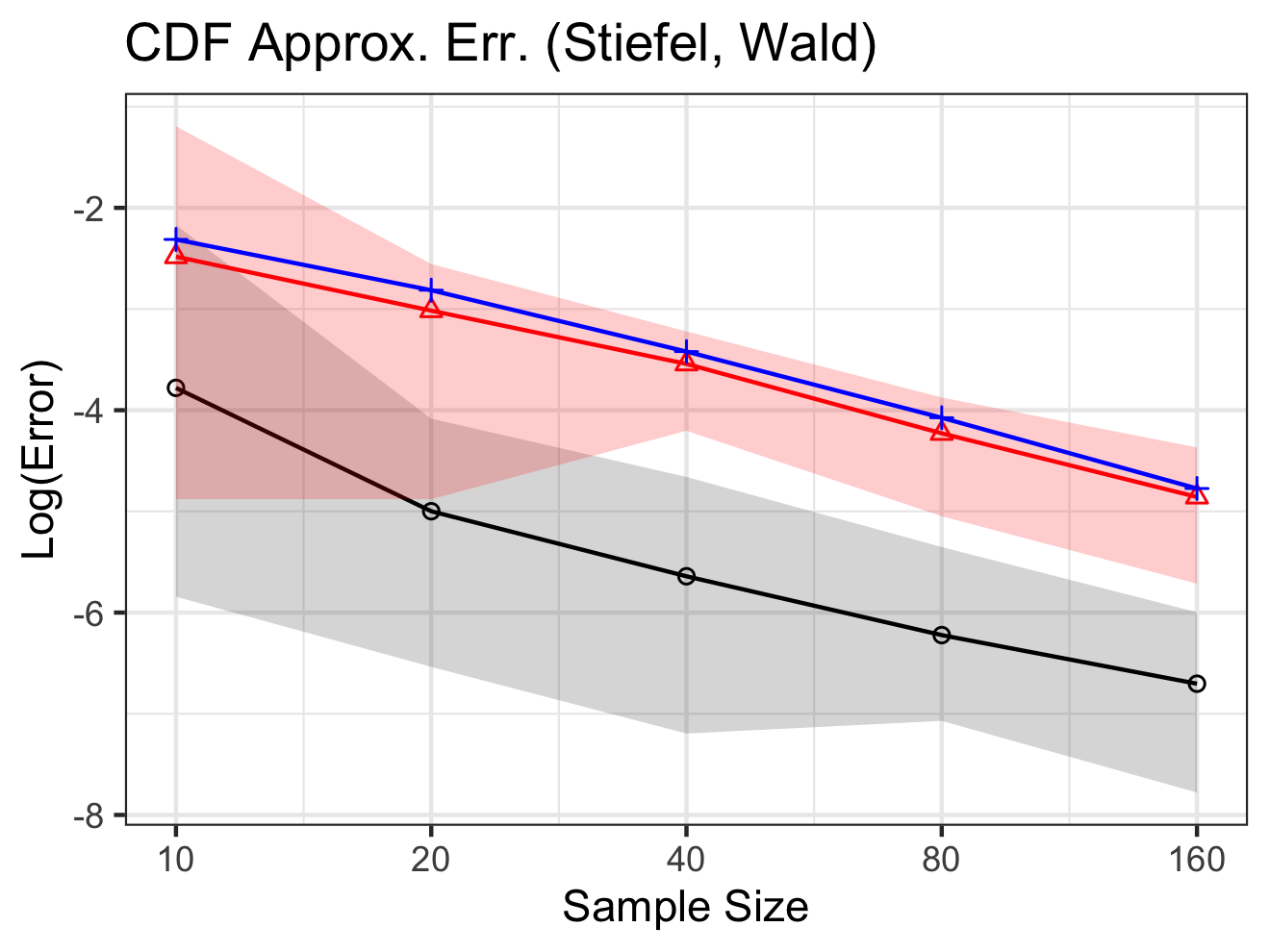}
}
	\subfigure{
		\includegraphics[width = 0.31\textwidth]{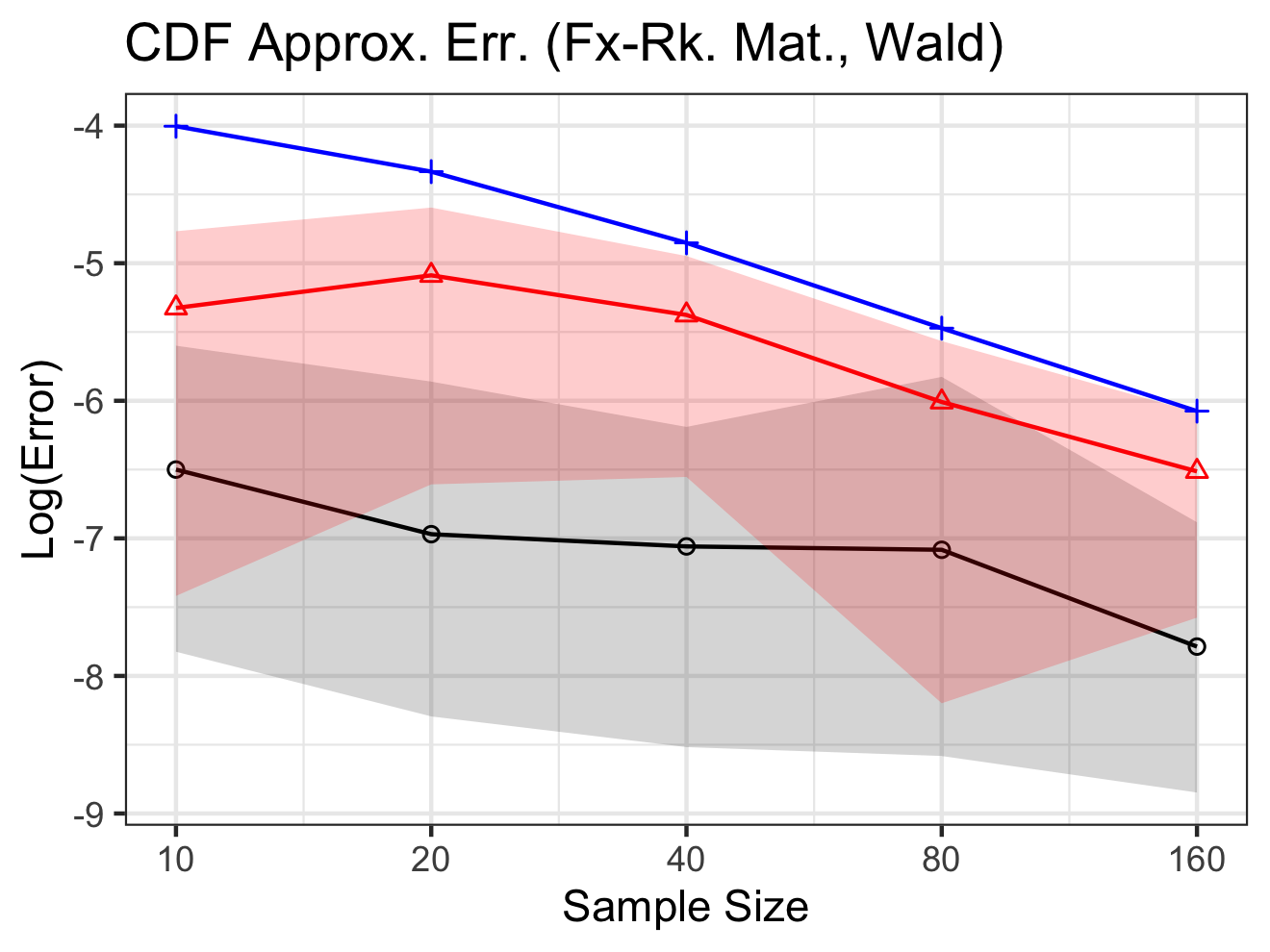}
}
	\subfigure{
		\includegraphics[width = 0.31\textwidth]{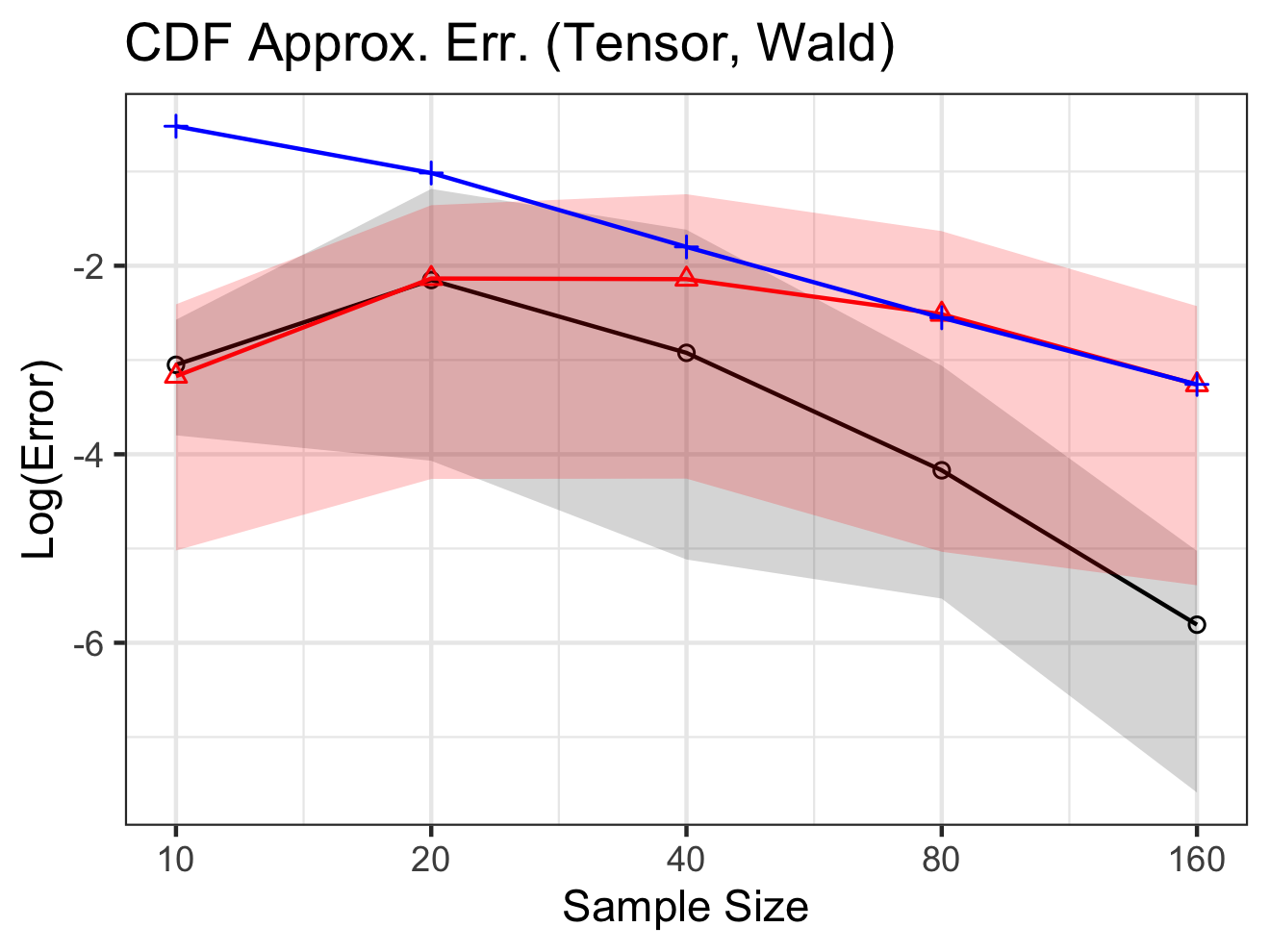}
}
	\subfigure{
		\includegraphics[width = 0.31\textwidth]{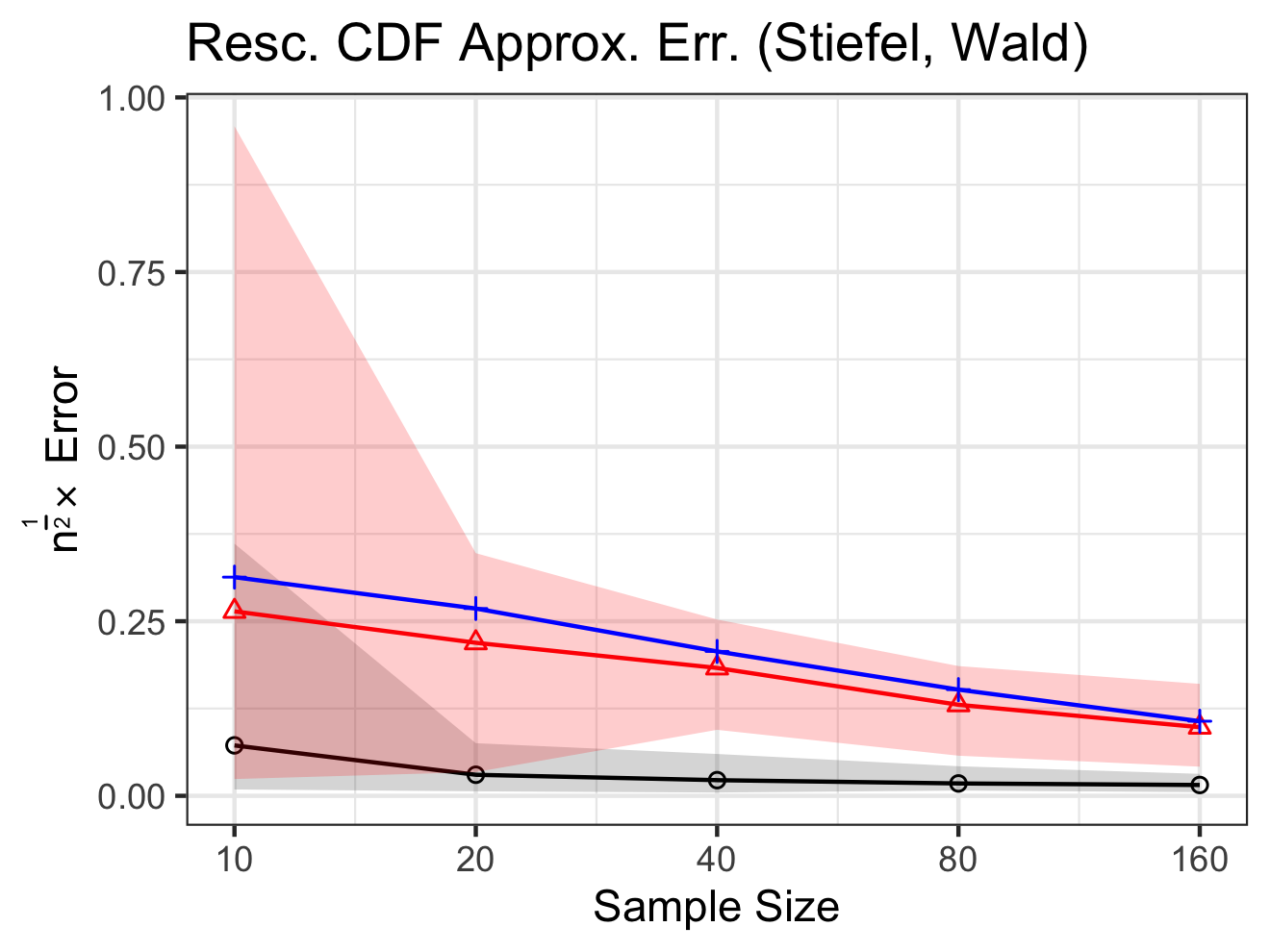}
}
	\subfigure{
		\includegraphics[width = 0.31\textwidth]{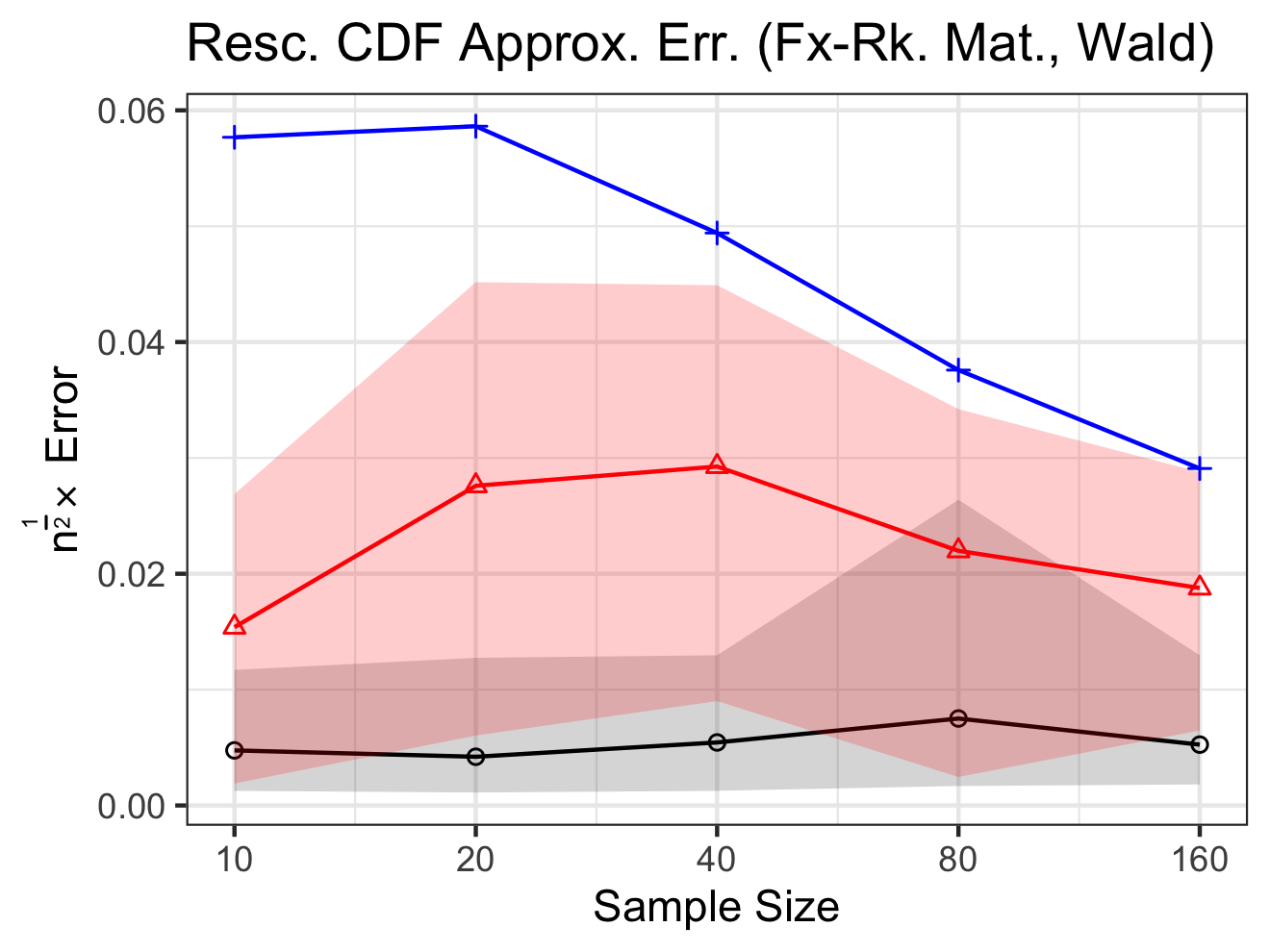}
}
	\subfigure{
		\includegraphics[width = 0.31\textwidth]{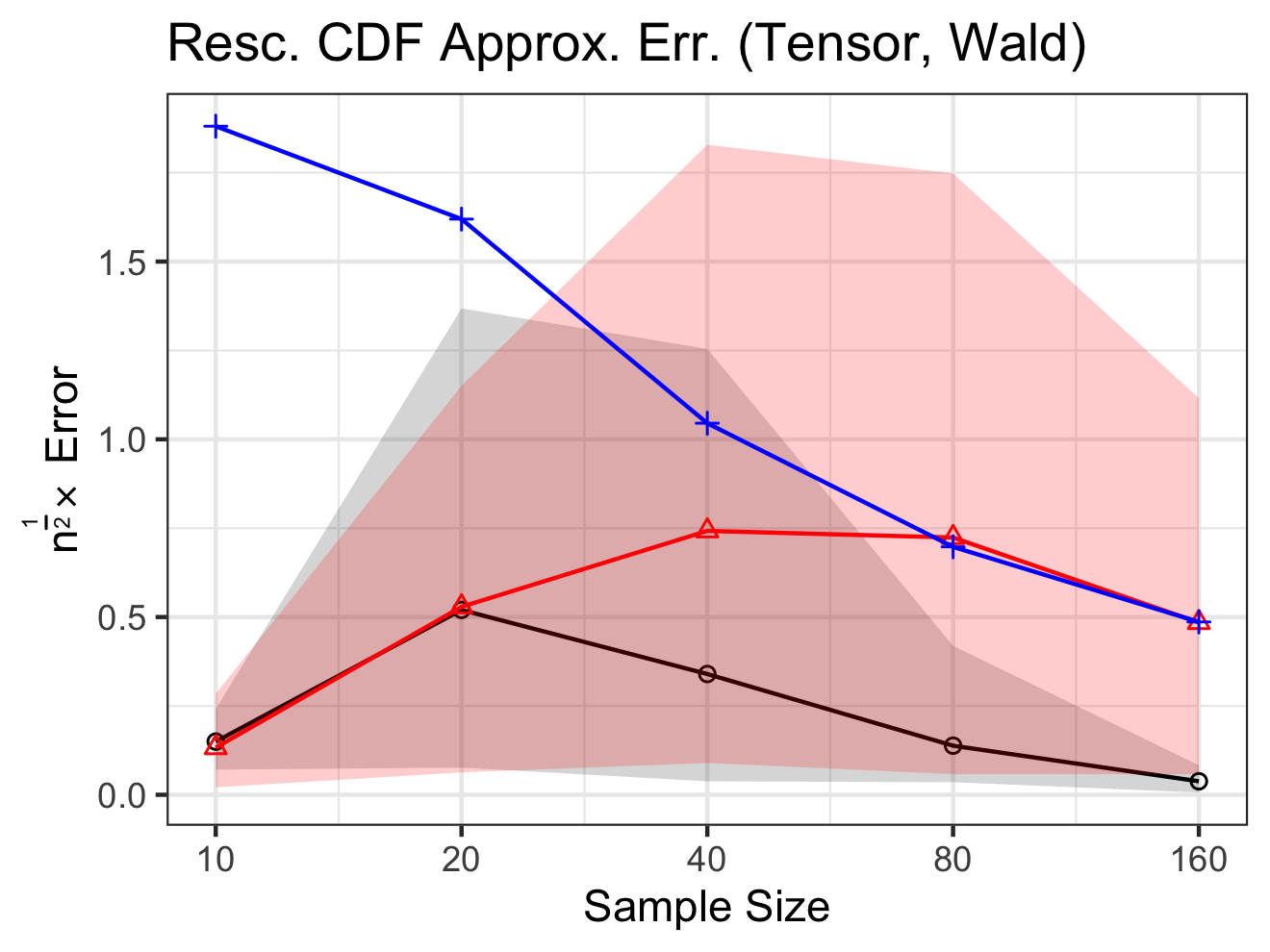}
}
	\subfigure{
		\includegraphics[width = 0.31\textwidth]{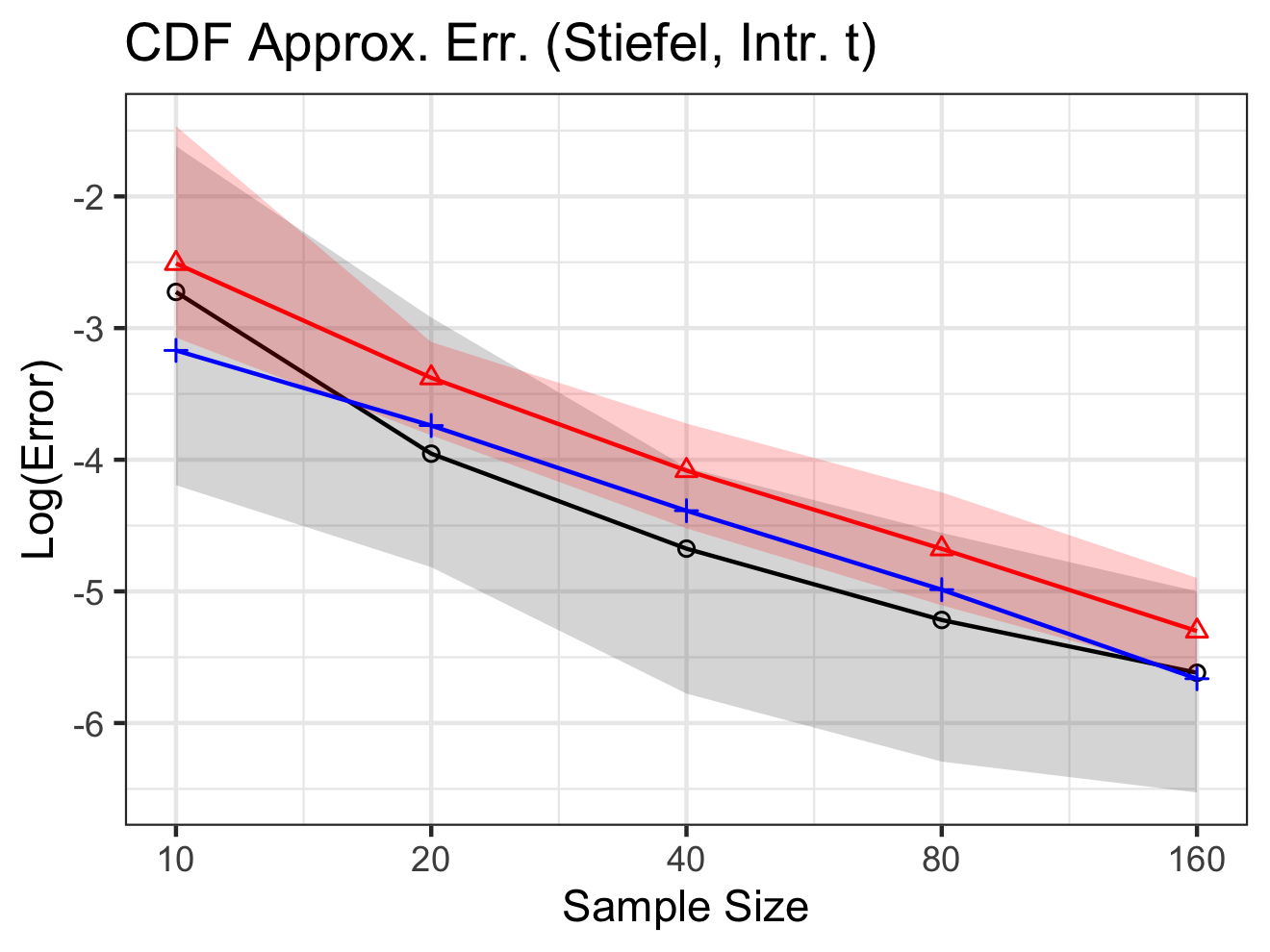}
}
	\subfigure{
		\includegraphics[width = 0.31\textwidth]{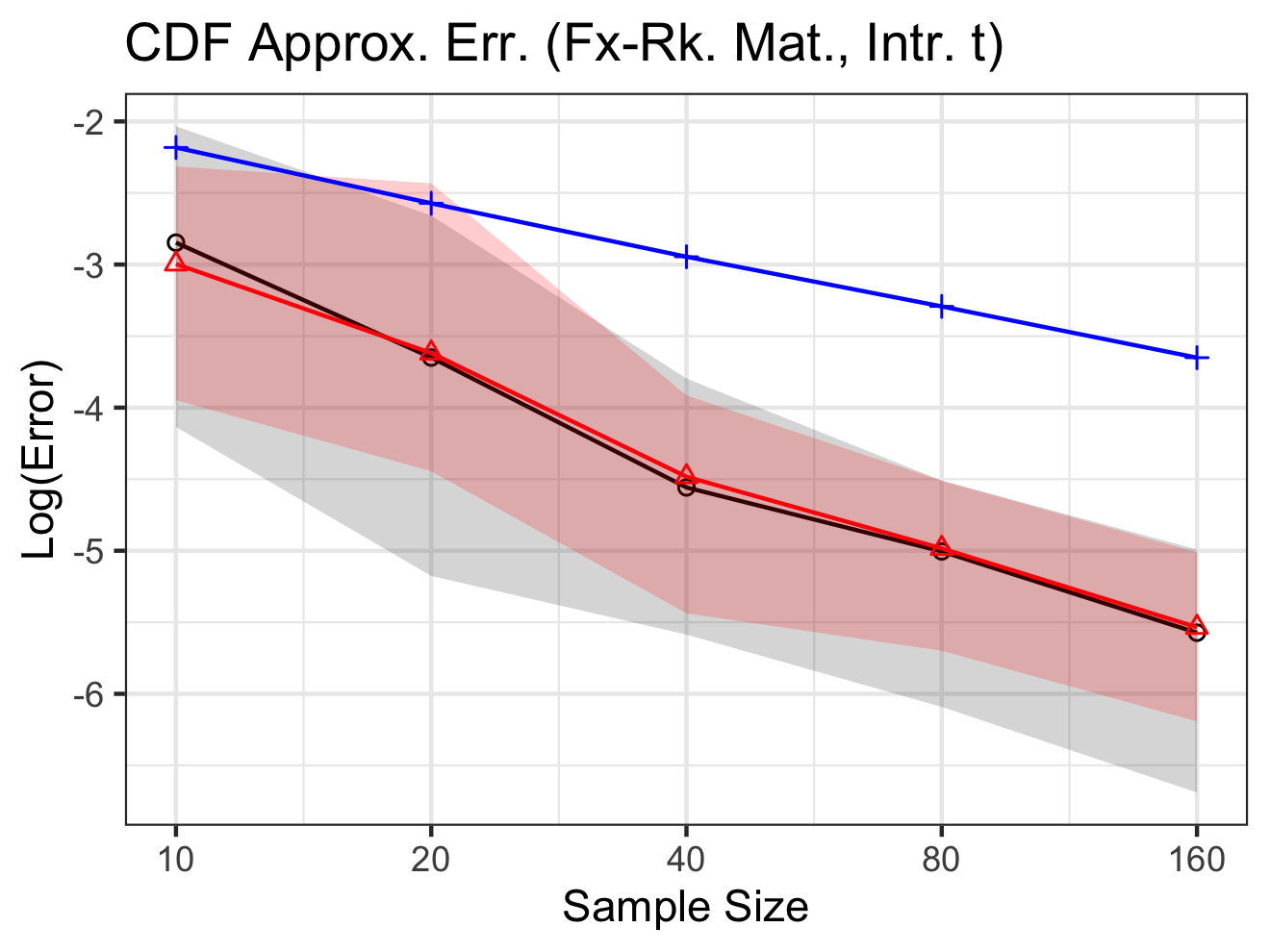}
}
	\subfigure{
		\includegraphics[width = 0.31\textwidth]{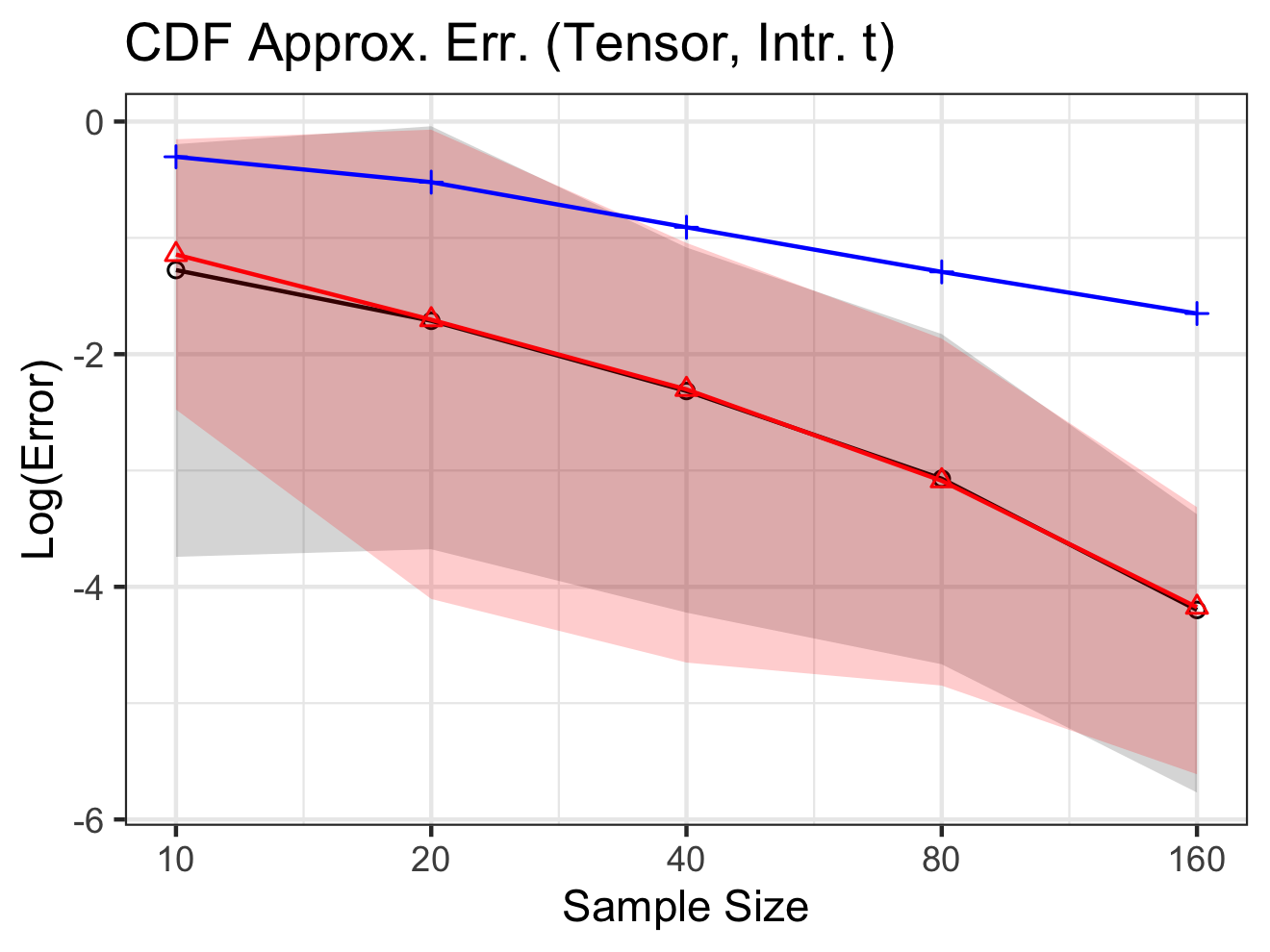}
}
	\subfigure{
		\includegraphics[width = 0.31\textwidth]{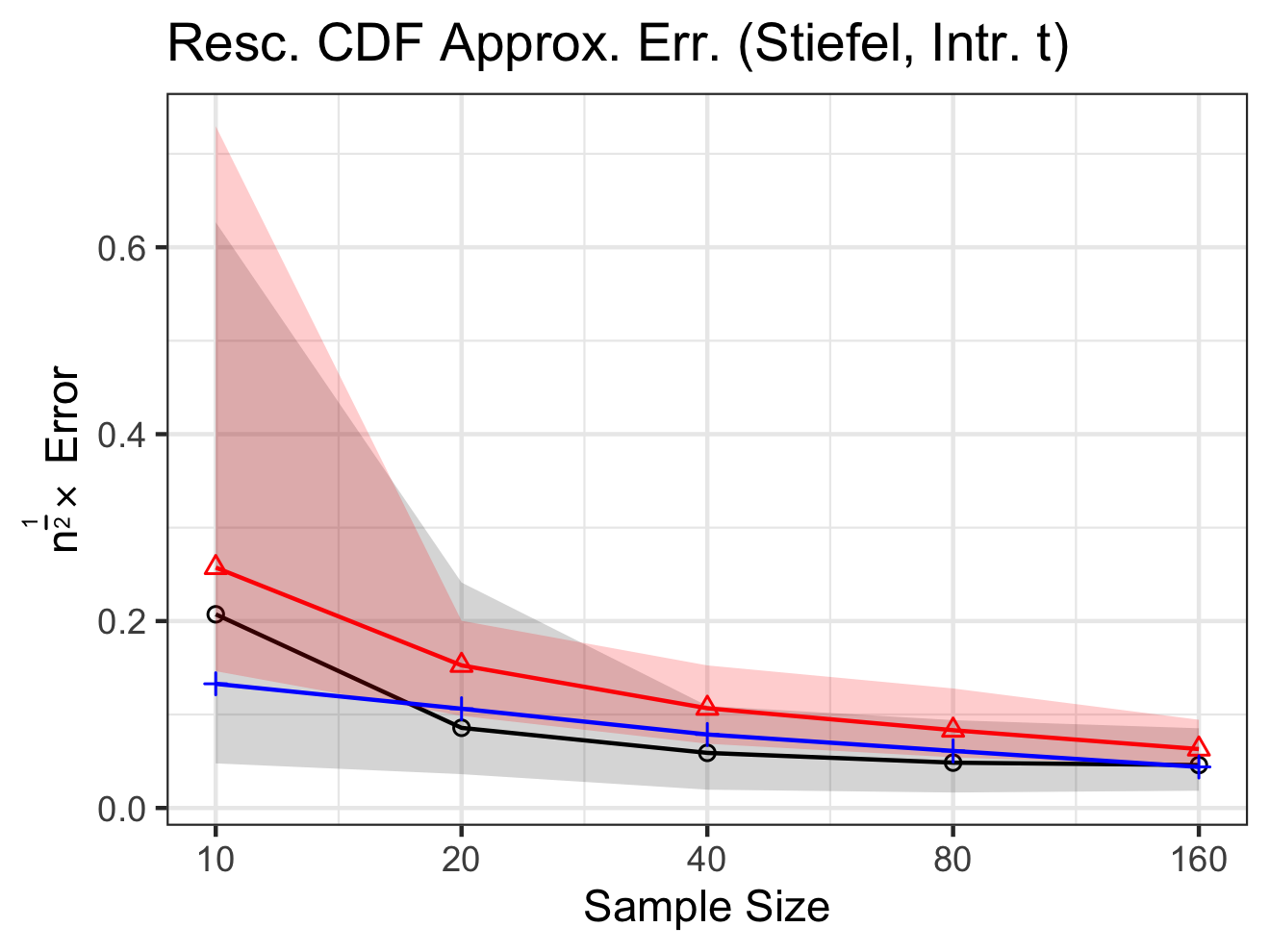}
}
	\subfigure{
		\includegraphics[width = 0.31\textwidth]{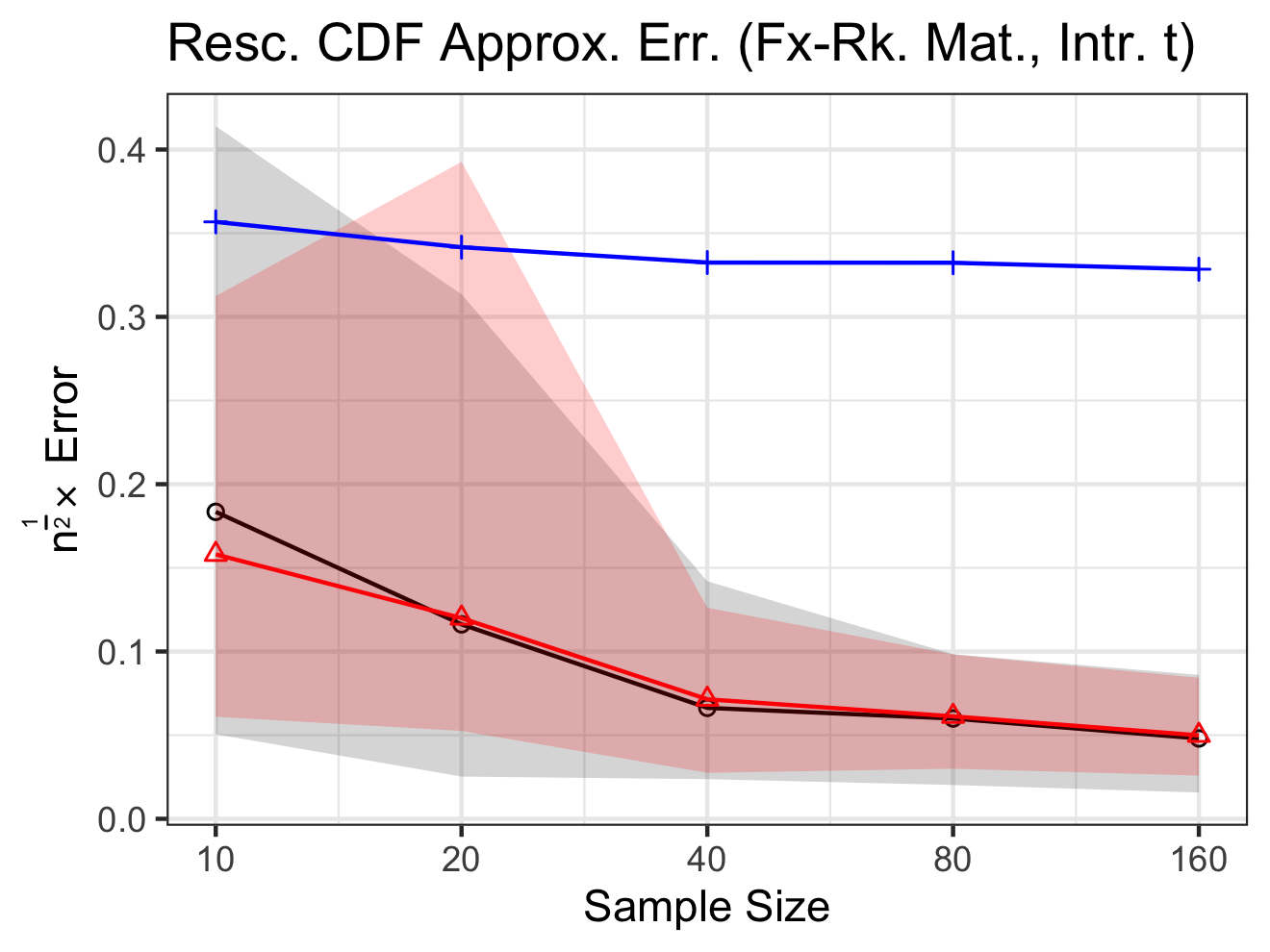}
}
	\subfigure{
		\includegraphics[width = 0.31\textwidth]{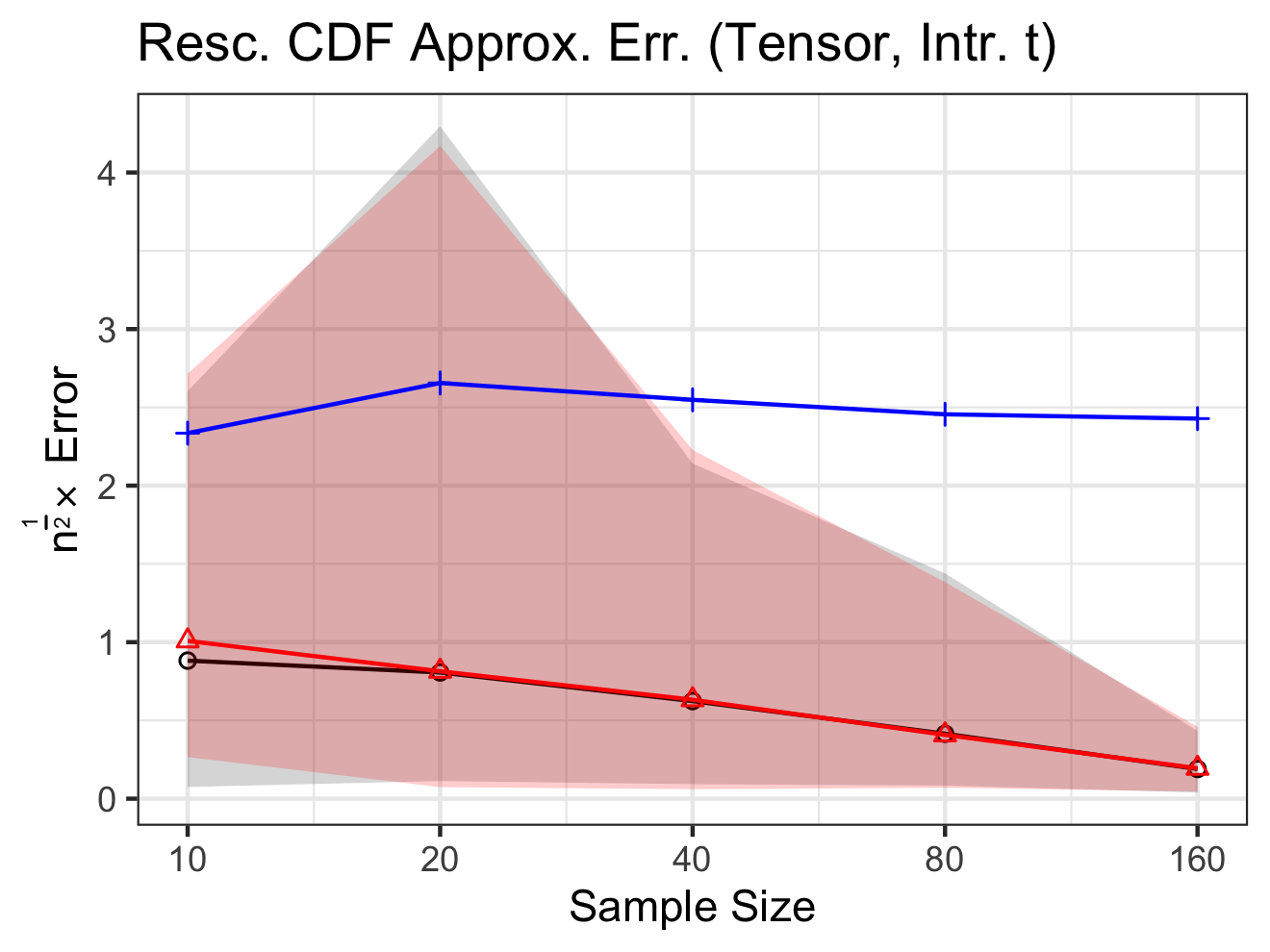}
}
\caption{Scaled Cumulative Distribution Function Difference. column $1$: \emph{Stiefel Manifold}; column $2$: \emph{fixed-rank matrix manifold}; column $3$: \emph{Rank-One Tensor Manifold}; row $1$: log-scaled Wald Statistic; row $2$: $\sqrt{n}$-scaled Wald Statistic; row $3$: log-scaled intrinsic $t$-Statistic; row $4$: $\sqrt{n}$-scaled intrinsic $t$-Statistic; 
	 {\emph{Black}}: 
	 resampled statistic (our method); {\color{red}\emph{Red}}: non-studentized resampled statistic; {\color{blue}\emph{Blue}}: chi-square distribution / standard normal distribution. Shaded regions indicate the range between the best and worst performance across replicates.}
\label{fig:CDF dif}
\end{figure}

\begin{figure}[hbtp]
	\centering
	\subfigure{\includegraphics[width = 0.45\textwidth]{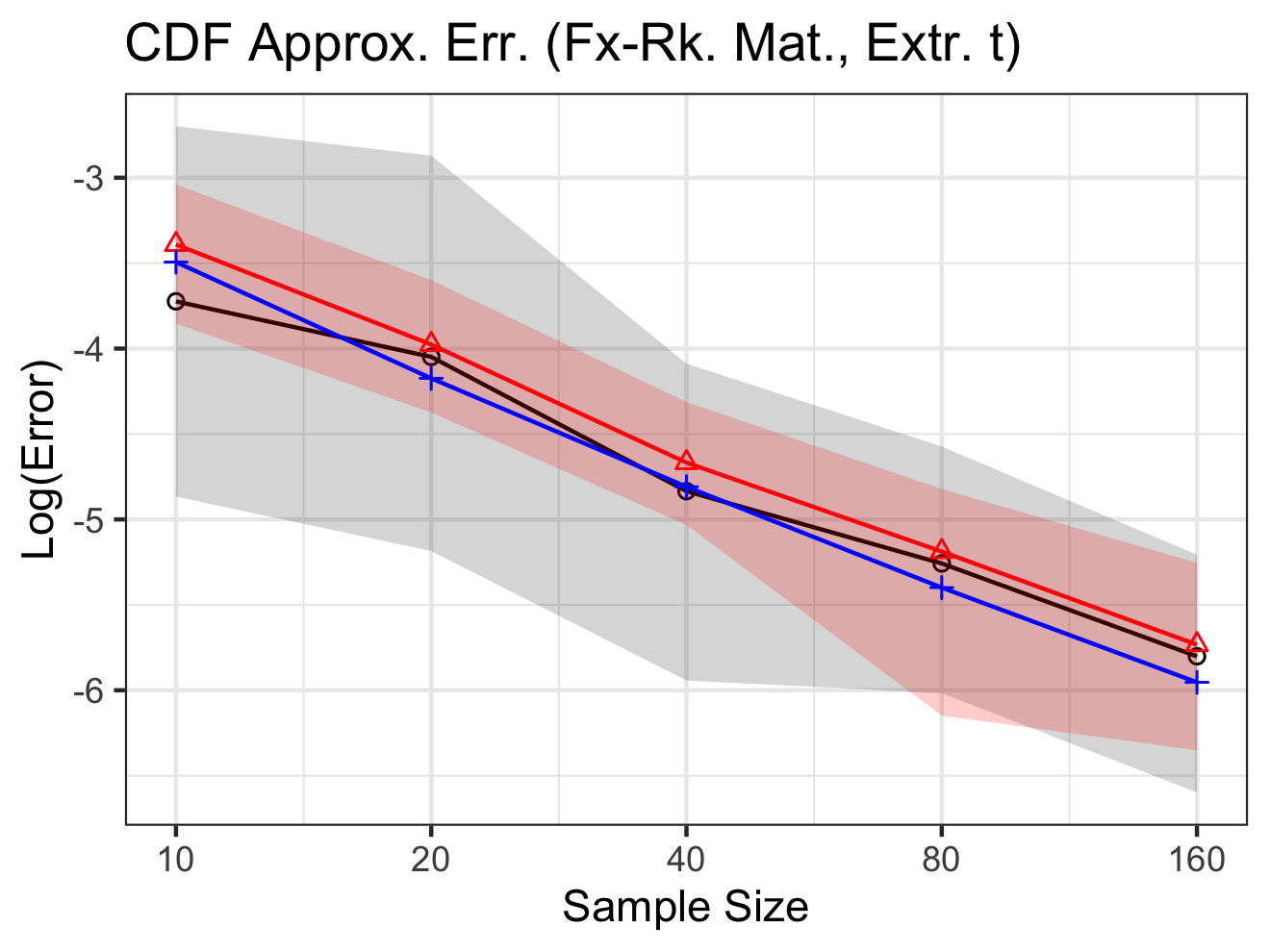}}
	\subfigure{\includegraphics[width = 0.45\textwidth]{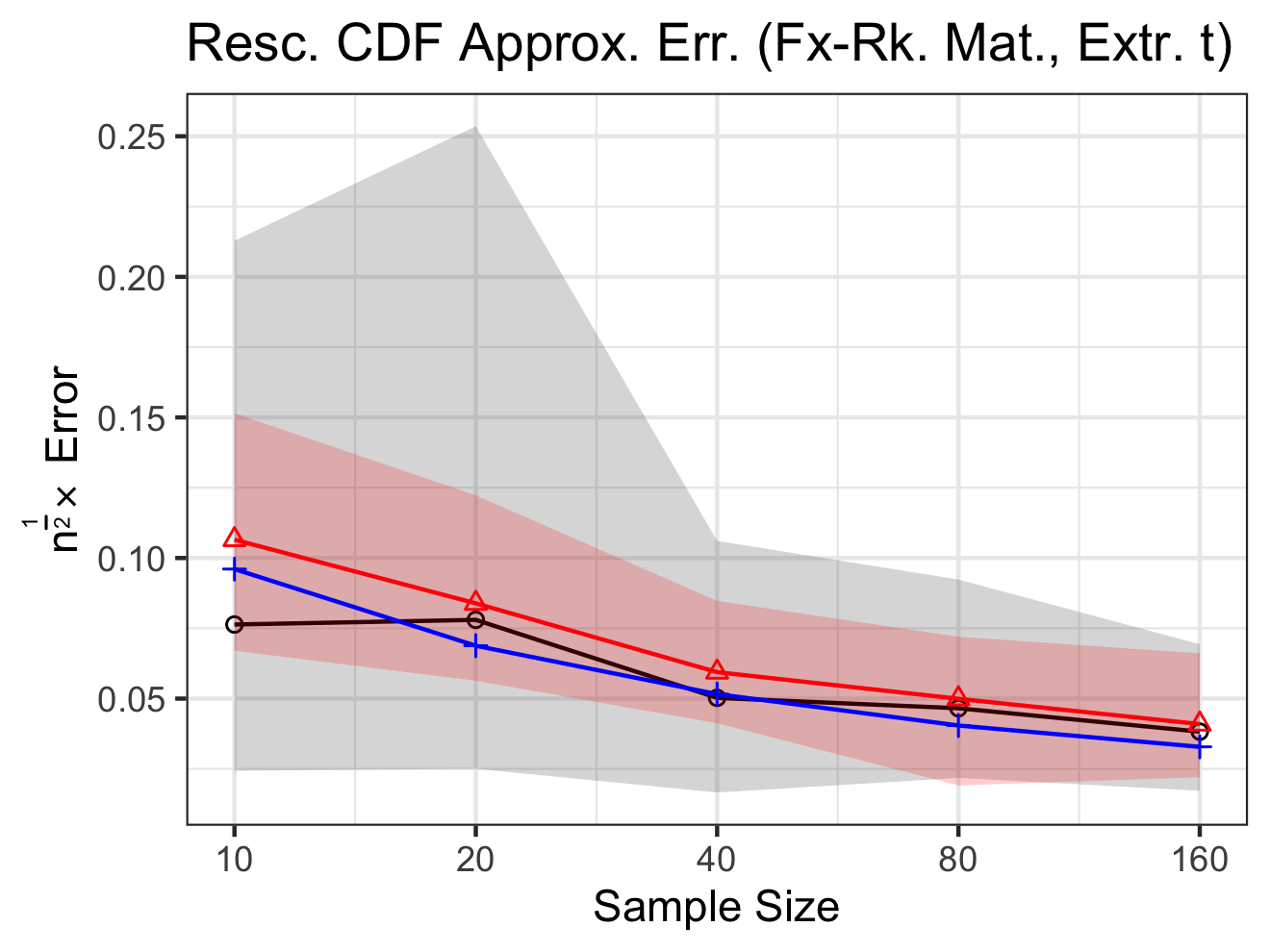}}	
	\caption{Scaled Cumulative Distribution Function Difference for Extrinsic $t$-statistic for fixed-rank matrix manifold. {\emph{Black}}: 
	 the empirical distribution of the resampled statistic; {\color{red}\emph{Red}}: the empirical distribution of the non-studentized resampled statistic; {\color{blue}\emph{Blue}}: the distribution of the standard normal distribution. Shaded regions indicate the range between the best and worst performance across replicates.}
	 \label{fig: intrinsic t}
\end{figure}

    Using the gradient and Hessian formulas in Section~\ref{subsec: stats applications}, we implemented the estimating and resampling procedures in Algorithm~\ref{algorithm: resampled newton iteration}. We validated the distributional guarantees for these procedures as detailed in Section~\ref{subsection: hypothesis testing} and Algorithms~\ref{algorithm:wald based statistic},\ref{algorithm:t based statistic},~and~\ref{algorithm:extrinsic statistic}. We repeated the bootstrapping for the resampled estimators $b$ times, where $b$ is taken to be $1000\times n$, with $n$ being the corresponding sample size. Increasing the number of bootstrap iterations with $n$ helps reduce the noise in the empirical cumulative distribution function (CDF), thereby better revealing differences from the true limiting distribution. We conduct four simulations as follows:
    \begin{enumerate}
    \item \textbf{Type-I/II Error of Tests. }To evaluate the \emph{type-I error} control as discussed in Section~2.2, we conduct an illustrative simulation under Setting 1; the samples are independently drawn from $\mathcal N(\boldsymbol \theta_0, 2\mb I_3)$ and we let $\boldsymbol \theta_1 = \boldsymbol \theta_0 = (0,1,0)\t$. For sample sizes $n=\{40, 80, 160\}$, we performed $30$ repeated experiments and summarize the observed Type-$I$ error in Table~\ref{table: test on sphere} for the $t$-statistic based on $\phi(\boldsymbol \theta)_1$ and the Wald statistic. The results in Table~\ref{table: test on sphere} demonstrate that our method provides stable approximations of the targeted Type-I error as the sample size increases. To further compare our approach with the Euclidean counterpart on the sphere, we define $\boldsymbol{\theta}_0' = (\delta, \sqrt{1 - \delta^2}, 0) \in \mathbb{S}^2$ and draw samples from $\mathcal{N}(\boldsymbol{\theta}_0', \mathbf{I}_3)$. Varying the sample size within $\{10,20,40\}$ and the deviation level $\delta\in \{0.2,0.4,0.6\}$, we compute the empirical power of both methods. Here, the Euclidean method refers to the test procedure described in Section~2.2, but disregarding the manifold structure (i.e., viewing $\mathcal{M} = \mathbb{R}^3$). Observe that our method achieves substantially higher power than the Euclidean approach. 
    
    \item \textbf{Convergence Rate of Intrinsic Statistics.} This experiment aims to validate our theoretical results on statistics intrinsic to Riemannian manifolds—specifically, the Wald and intrinsic $t$-statistics—by comparing them with benchmark statistics. We consider the additive noise setting in Section~4.2.1, where the parameter space is the Stiefel, fixed-rank matrix, or rank-$1$ tensor manifold. Samples are generated as 
\eq{
\mathbf{X}_i = \boldsymbol{\theta}_0 + \boldsymbol{\epsilon}_i, \quad \boldsymbol{\epsilon}_i \sim \mathcal{N}(0, \mathbf{I}_d),
\label{eq: additive model}
}
where $d$ is the dimension of the corresponding Euclidean representation.
    To assess the approximation accuracy to the cumulative distribution function of a targeted quantity $\hat T_n$ denoted by $F_{\hat T_n}$, we consider the following measure:
\begin{equation}
\mathrm{Error}(\hat{F}_{\hat{T}_n}) = \sup_{x \in S} \left| \hat{F}_{\hat{T}_n}(x) - F_{\hat{T}_n}(x) \right|,
\end{equation}
where $\hat F_{\hat T_n}$ is an approximation to the CDF $F_{\hat T_n}$ of $\hat T_n$. 
For the Wald statistic, we let $\hat T_n = R_{\hat \theta_n}^{-1}(\theta_0)^{\top} \hat{\mb \Sigma}^{-1}R_{\hat \theta_n}^{-1}(\theta_0)$ and \( S = [8] \), whereas for the intrinsic $t$ statistic, we let $\hat T_n = R_{\hat \theta_n}^{-1}(\theta_0)_{1} / (\hat\Sigma_{1,1})^{\frac12}$ and $ S = \{ 0.2 \cdot s : s \in \{-10, -9, \ldots, 9, 10\} \}$. For the choice of $\hat F_{\hat T_n}$, we propose using the Monte-Carlo approximation of the CDF of $R_{\hat \theta_n^*}^{-1}(\hat \theta_n)^{\top} \check{\mb \Sigma}^{-1}R_{\hat \theta_n^*}^{-1}(\hat \theta_n)$ for the Wald statistic, and $ R_{\hat \theta_n^*}^{-1}(\hat \theta_n)_{1} / (\check\Sigma_{1,1})^{\frac12}$ for the intrinsic $t$ statistic. For comparison purposes, we also consider alternative options for $\hat F_{\hat T_n}$ including the empirical CDF of the non-studentized version $R_{\hat \theta_n^*}^{-1}(\hat \theta_n)^{\top} \hat{\mb \Sigma}^{-1}R_{\hat \theta_n^*}^{-1}(\hat \theta_n)$ and $ R_{\hat \theta_n^*}^{-1}(\hat \theta_n)_{1} / (\hat\Sigma_{1,1})^{\frac12}$, as well as the CDFs of the normal/chi-squared distributions. 
Then, the \( \mathrm{Error}(\hat{F}_{\hat{T}_n}) \) is scaled by the logarithm function and the square root of \( n \), respectively. In this simulation, we are interested in the distributional consistency of Wald statistics and the intrinsic $t$-statistics, in comparison with the other benchmark methods. 
Figure~\ref{fig:CDF dif} presents the average, maximum, and minimum scaled approximation errors using different methods over $30$ epochs. The first and third rows of Figure~\ref{fig:CDF dif} indicate the (almost) linear relation between $\Log(\text{Error})$ and the sample size, with our proposed statistics often exhibiting the steepest slope. Meanwhile, the second and fourth rows of Figure~\ref{fig:CDF dif} corroborate our theoretical expectation from Theorem~3.6 that $n^{\frac12}\cdot \text{Error}$ decreases with increasing sample size, with our methods consistently outperforming or closely matching the benchmarks.

    \item \textbf{Convergence Rate of Extrinsic $t$ Statistics.} Moreover, we conducted further simulations to examine the extrinsic $t$-statistic for the fixed-rank matrix manifold under the model \eqref{eq: additive model}. The simulations follow the same setup and measurement approach as described earlier for the fixed-rank matrix manifold. Our focus is on a single coordinate of the fixed-rank matrix, specifically on $\hat T_n$, which refers to the studentized version of $\big((\boldsymbol \theta_0)_{1,1} - (\hat{\boldsymbol\theta}_n)_{1,1}\big)$. 
    The resulting errors using Algorithm~5 are shown in Figure~\ref{fig: intrinsic t}. From the log-scaled plot and the $n^{\frac12}$-scaled plot, it shows that the distributional deviation $\mathrm{Error}(\hat F)$ of the extrinsic $t$-statistic is of the order $o(n^{-\frac12})$ employing Algorithm~\ref{algorithm:extrinsic statistic}, thereby supporting our theoretical results. 


    \item \textbf{Barycenter Problem. } Lastly, we aim to demonstrate the performance of our methods on Barycenter problems discussed in Section~\ref{subsubsection: barycenter}. We consider the sphere \( \mathbb{S}^2 \) restricted to the northern hemisphere (\( x_2 \geq 0 \)). We independently generate \( n \) samples \( \{ X_i \}_{i=1}^n = \{ \cos \theta_i \sin \phi_i, \cos \phi_i, \sin \theta_i \sin \phi_i \} \), where \( \theta_i \) is independently drawn from \( \mathsf{Unif}([0, 2\pi)) \) and \( \phi_i \) is drawn from the beta distribution \( \mathsf{Beta}(2, 2) \). We then compute the barycenter using Algorithm~2. Given that the population barycenter in this case is \( (0, 1, 0) \), we analyze the distributional consistency of the corresponding statistics over $30$ repetitions, with results summarized in Figure~\ref{fig:barycenter}. These results reveal that the Wald statistic, intrinsic 
$t$-statistic, and extrinsic $t$-statistic all achieve a convergence rate of $o(n^{-\frac12})$ in terms of the CDF, aligning with our theoretical predictions.

    \end{enumerate}

    In summary, the simulations validate the effectiveness and theoretical soundness of our proposed methods across diverse statistical challenges, demonstrating superior performance and consistency compared to benchmark approaches.
 
{
\setlength{\tabcolsep}{0.15em}
\begin{table}
 \centering
\begin{tabular}{|c|c||c|c|c|c|c|c||}
\hline
  & Hypo. &\multicolumn{6}{c|}{$H_0: \boldsymbol \theta_1 = \boldsymbol\theta_0$ v.s. $H_1: \boldsymbol\theta_1 \neq \boldsymbol \theta_0$} \\
   \hline
    & Test & \multicolumn{3}{c|}{ Test Based on $T_1$} & \multicolumn{3}{|c|}{ Test Based on $W$} 
   \\
   \hline\hline
   n&$1 - \alpha$ & 0.9 & 0.95 & 0.975 & 0.9 & 0.95 & 0.975 \\ 
  \hline 
  $40$&& 0.923 & 0.971 & 0.988 & 0.961 & 0.986 & 0.996
\\
  $80$&Type-I & 0.908 & 0.958 & 0.982 & 0.940 &0.978 & 0.992 \\ 
  $160$& & 0.906 & 0.956 & 0.981 & 0.921 & 0.967 & 0.987\\
  \hline
\end{tabular}
\hspace{10pt}
\begin{tabular}{|c|c||c|c|c|c|c|c|c||}
\hline
  &Hypo. &\multicolumn{6}{c|}{$H_0: (\boldsymbol\theta_0)_1 = 0$ v.s. $H_1: (\boldsymbol\theta_0)_1 \neq 0$} \\
   \hline
   & Test & \multicolumn{3}{c|}{Manifold Test} & \multicolumn{3}{|c|}{Euclidean Test} 
   \\
   \hline
  n & $\delta$ & 0.2 & 0.4 & 0.6 & 0.2 & 0.4 & 0.6 \\ 
  \hline \hline
  $10$& & 0.188 & 0.395 & 0.550 & 0.124 & 0.271 & 0.543
\\
 $20$& Power & 0.245 & 0.554 & 0.900 & 0.220 & 0.547 & 0.808 \\ 
  $40$& & 0.391 & 0.843 & 0.990 & 0.355 & 0.794 & 0.983\\
  \hline
\end{tabular}
  \caption{
  Type-I error and Power of Hypothesis Testing on Sphere. The left-hand table denotes the type-I errors with various significance levels $\alpha$ and sample sizes. The right-hand table collects the power ($1 - $type-II error) given varying level of $\bo \theta_0$'s deviation $\delta$ and different sample sizes. 
  }
  \label{table: test on sphere}
\end{table}
}

    \begin{figure}[hbtp]
            \centering 
            \subfigure{\includegraphics[width = 0.3\textwidth]{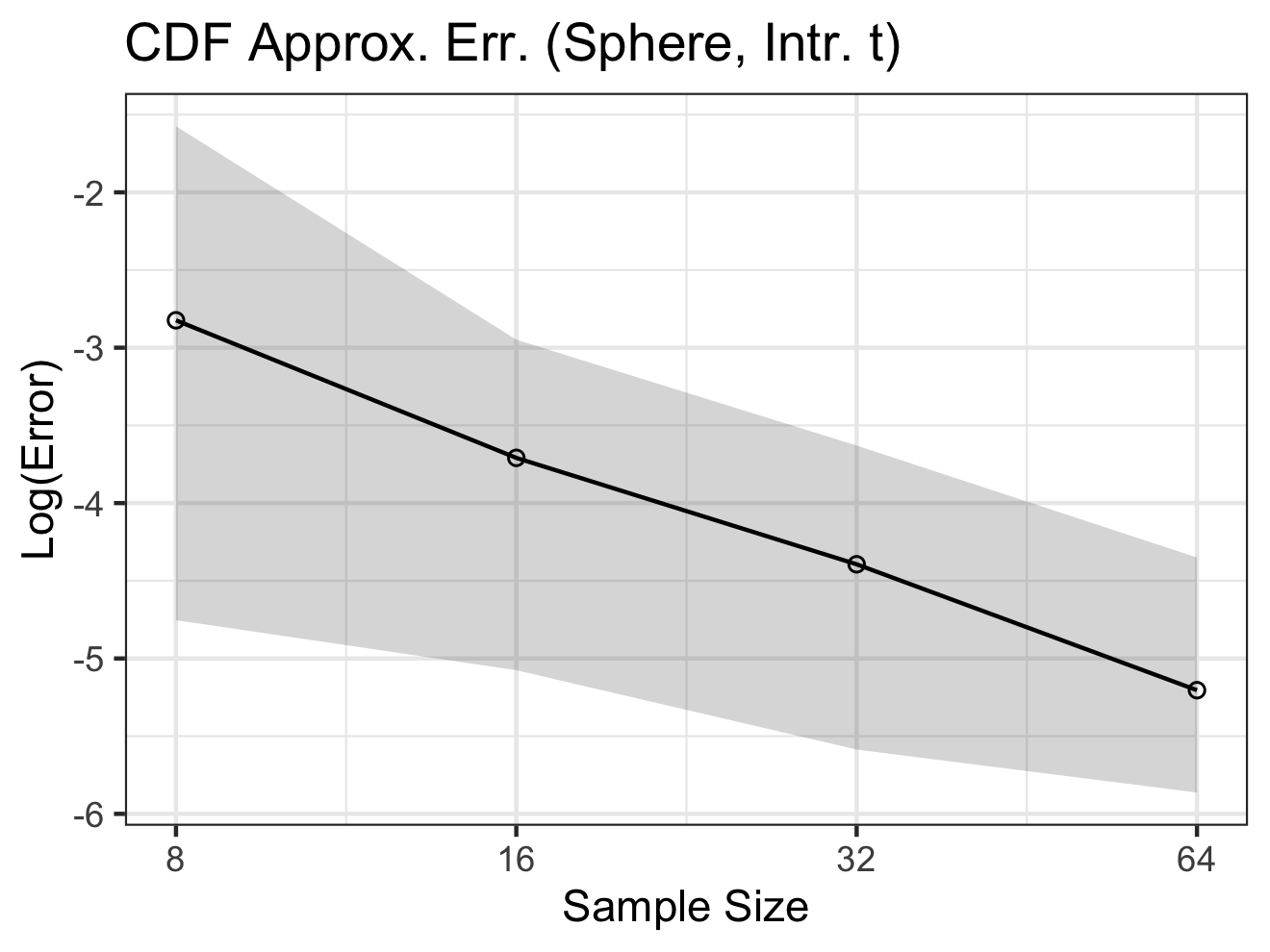}}
            \subfigure{\includegraphics[width = 0.3\textwidth]{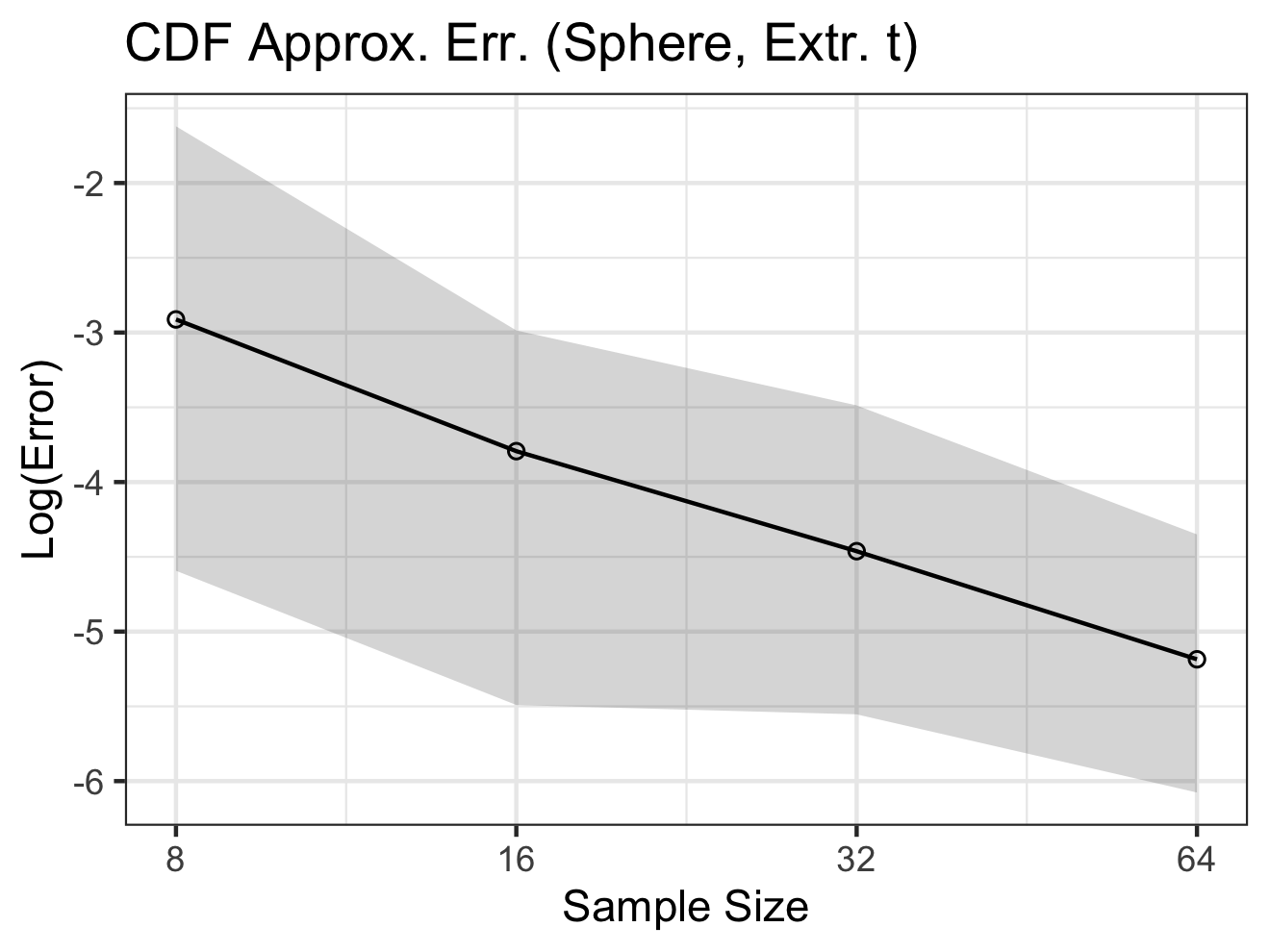}}
            \subfigure{\includegraphics[width = 0.3\textwidth]{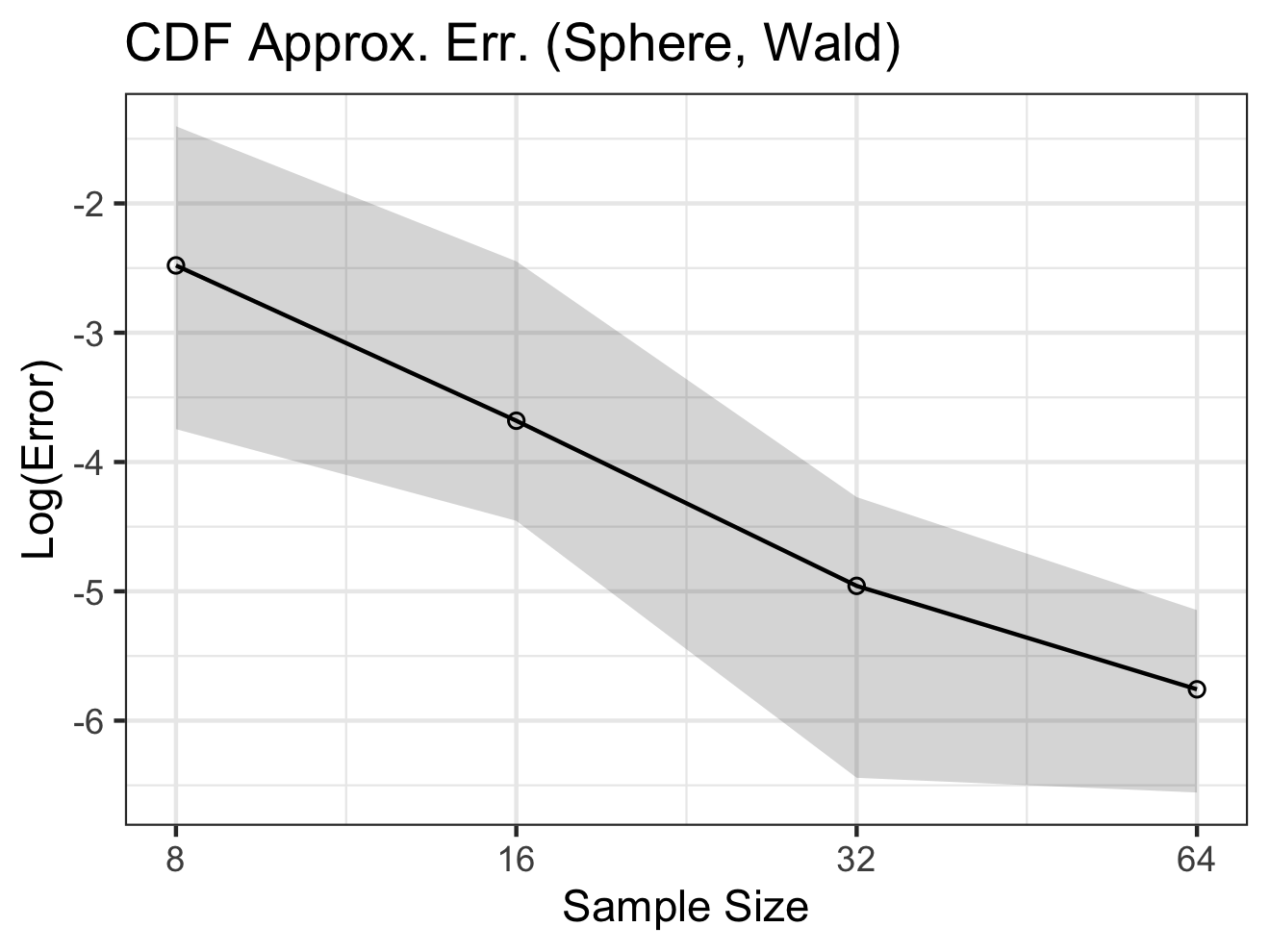}}
            \subfigure{\includegraphics[width = 0.3\textwidth]{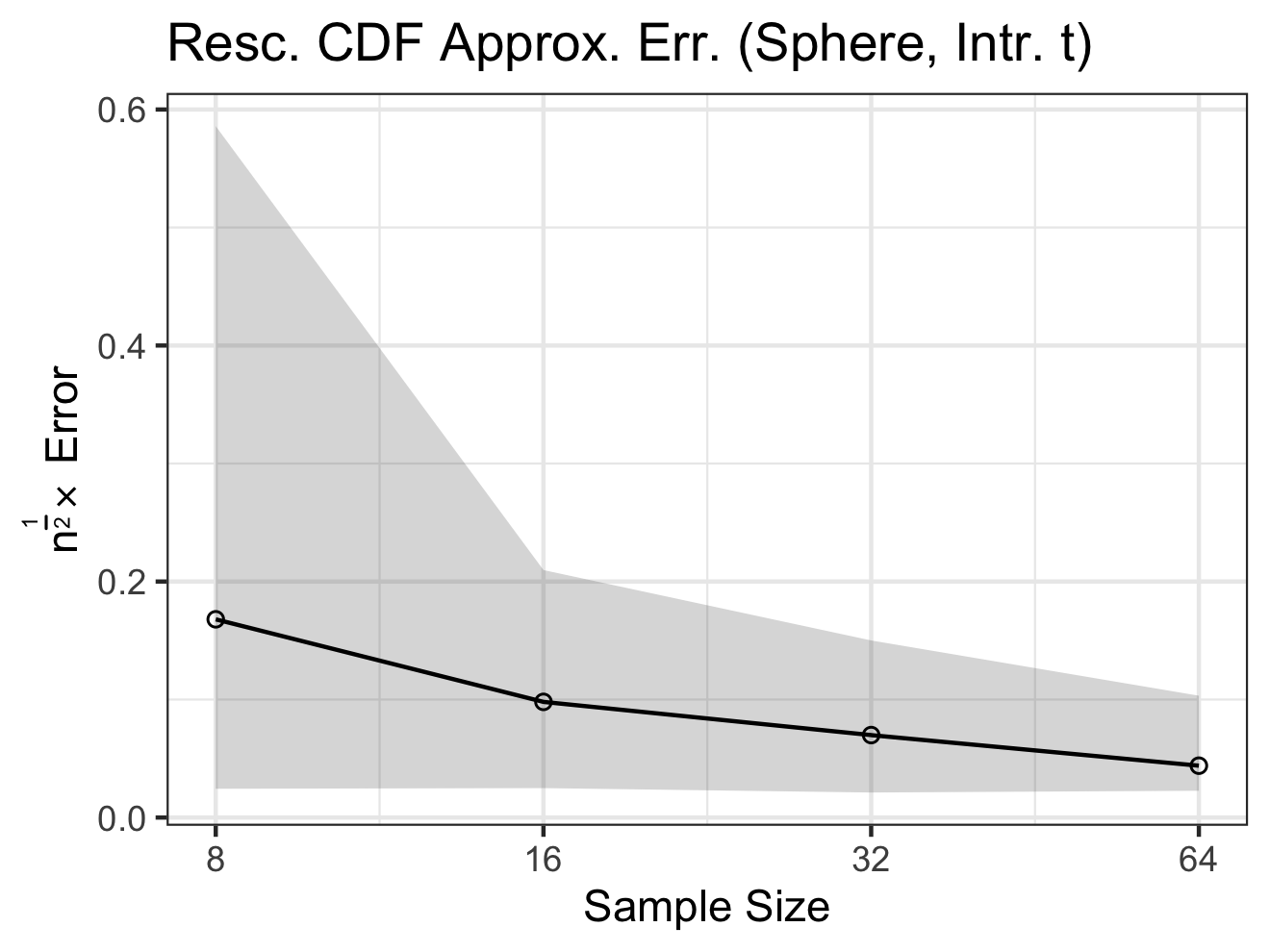}}
            \subfigure{\includegraphics[width = 0.3\textwidth]{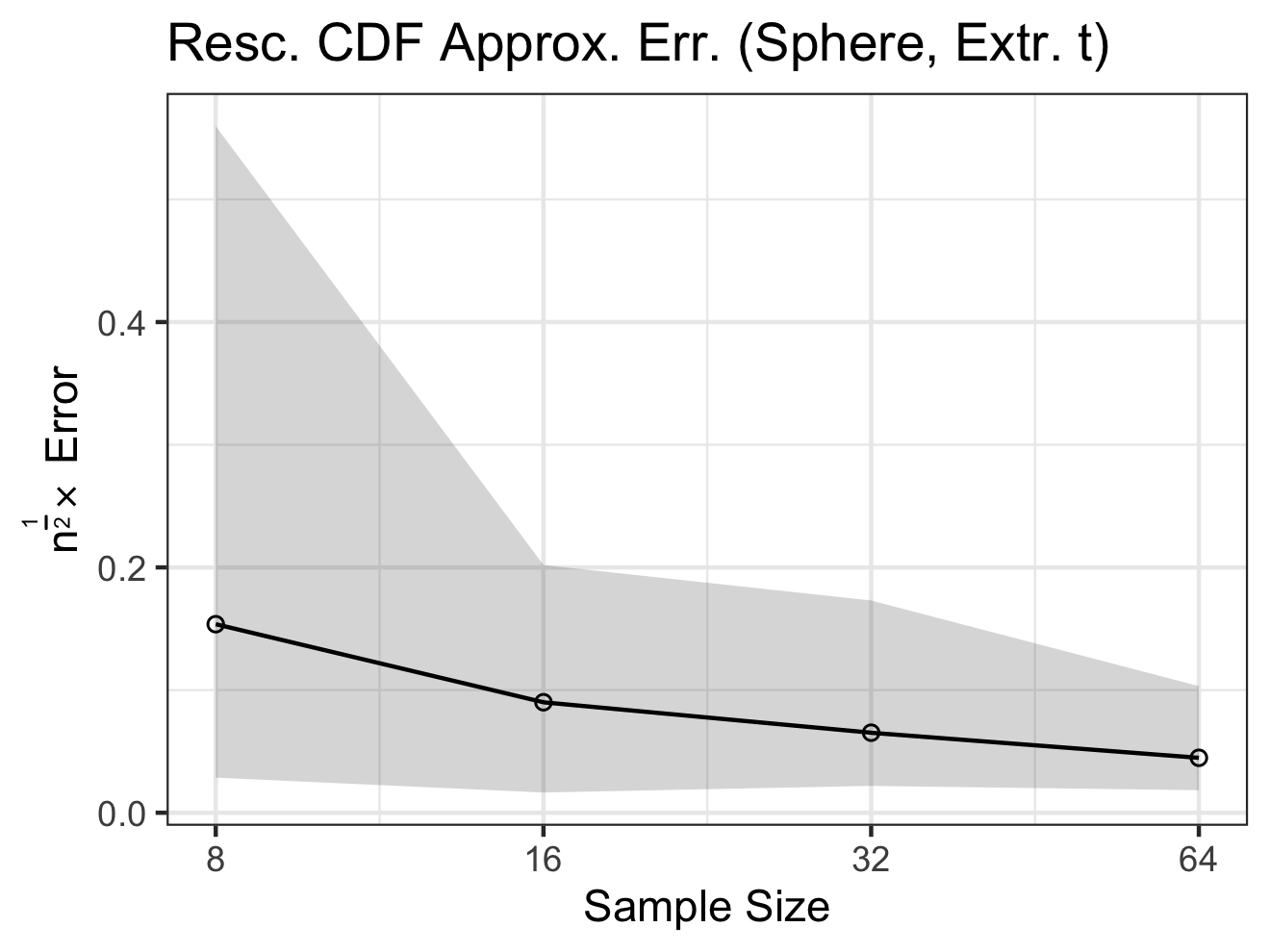}}
            \subfigure{\includegraphics[width = 0.3\textwidth]{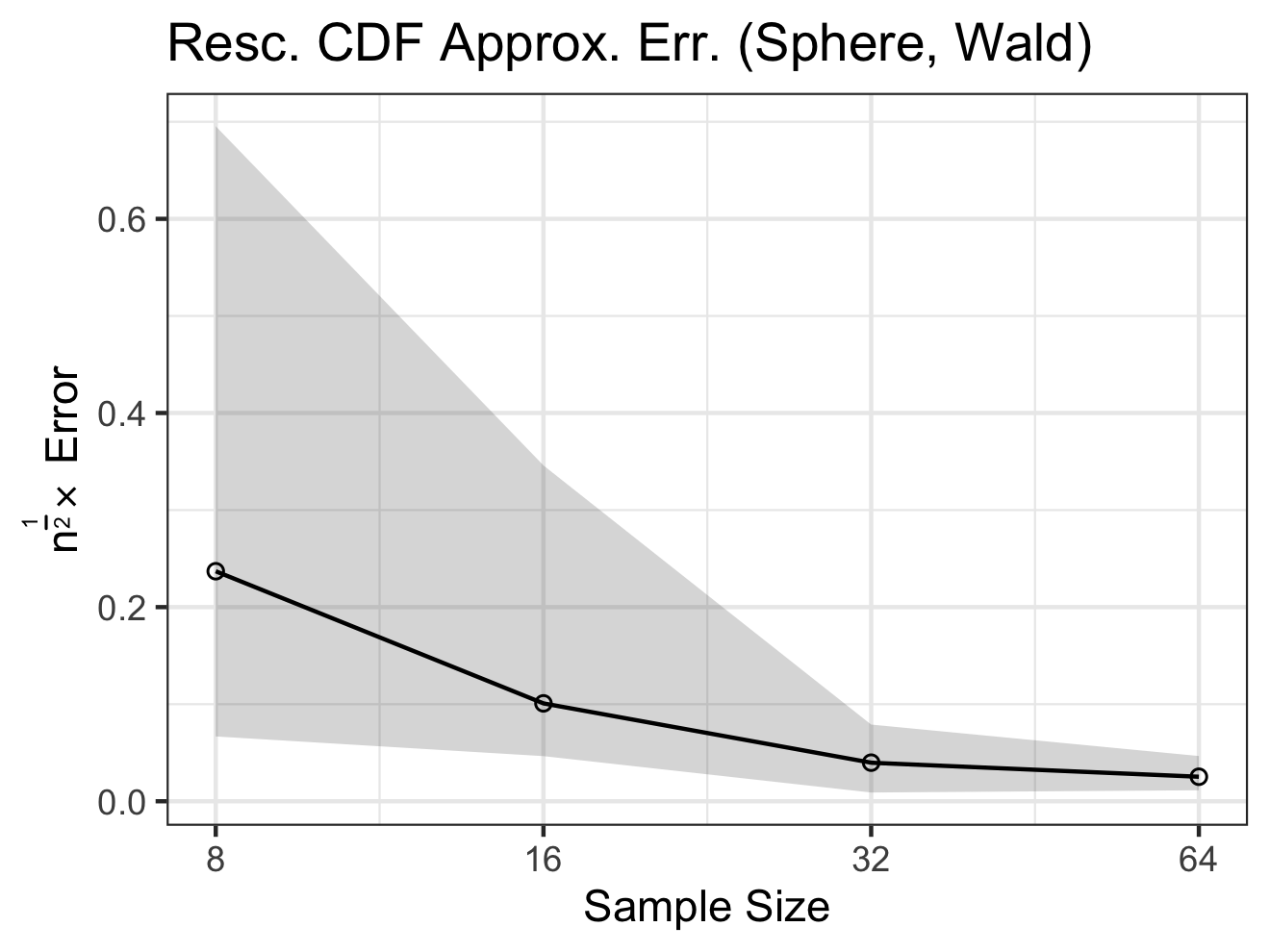}}
            \caption{Cumulative Distribution Function Values for Sphere Barycenter Problem. row $1$: Log-Scaled Error; row $2$: $\sqrt{n}$-Scaled Error; Column $1$: Intrinsic $t$; column $2$: Extrinsic $t$; column $3$: Wald. Shaded regions indicate the range between the best and worst performance across replicates.}
            \label{fig:barycenter}
            \end{figure}

\section{Proof of Main Result}

\subsection{Proof of Theorem~\ref{theorem:bootstrapmanifold}}
Instead of getting entangled in the complexities across various charts, our strategy throughout the paper is to concentrate on a fixed normal coordinate chart centered at $\theta_0$ to simplify the discussions. Our primary analytical tool is the Taylor expansion on the coordinate chart $\pi_{\{e_i\}} \circ \Log_{\theta_0}(\cdot)$, where $\{e_i\}$ is a normal basis of $\mathrm T_{\theta_0}\mc M$ to be specified later, along with the Edgeworth expansion results developed in the Euclidean case. 

Firstly, to ensure that the geodesics and exponential mappings are well-defined, we restrict the region to a geodesically complete subset of $\mc M$. By Theorem~6.17 in \cite{lee2018introduction}, there exists a positive $\rho$ such that for arbitrary $\theta \in \Exp_{ \theta_0}(B(0,\rho))$, $\Log_{\theta}$  defined on ${\Exp_{ \theta_0}(B(0,\rho))}$ is a differmorphisim onto its codomain in $\mathrm T_{\theta}\mc M$. We denote this region by $N$ in the following context. 

In Algorithm~\ref{algorithm: resampled newton iteration}, each resampling process, including Newton's step on resampled datasets, is performed independently. Consequently, we treat $\hat \theta_n^*$ as a random variable, and Algorithm~\ref{algorithm: resampled newton iteration} generates $b$ independent replications of $\hat \theta_n^*$ conditioned on the dataset $\mc X_n$. We denote the generic resampled dataset by $\mc X_n^* \coloneqq \{X_i^*\}_{i\in[n]}$ and the corresponding objective function by $L_n^*(\cdot ) \coloneqq \sum_{ i\in[n]}L(\cdot, X_i^*) / n$.

To simplify notations, we define the orthonormal basis induced by $R_{\theta_0}$ at $\theta_0$ as $\{e_i\}_{i\in[p]}\subset \mathrm T_{\theta_0}\mc M$, where $e_i = \mathrm dR_{\theta_0}\delta_i$. We use $L^{\ret, \theta}(\cdot, x)$ to denote the function $L(R_{\theta}(\cdot ), x)$. Further, given $\theta \in \mc M$ and a basis $\{e_i\}$ of $\mathrm T_{\theta}\mc M$, we denote the function $L(\Exp_{\theta}\circ \pi^{-1}_{\{e_i\}}(\cdot ), x)$ by $L^{\Exp, \theta, \{e_i\}}(\cdot, x)$ for every $\theta\in \mc M$. Similarly, we denote by $L_n^{\Exp, \theta, \{e_i\}}(\cdot)$, $  {L_n^*}^{\Exp, \theta, \{e_i\}}(\cdot)$, $L_n^{\ret, \theta}(\cdot)$, ${L_n^*}^{\ret, \theta}$ the composite functions $L_n(\Exp_{\theta}\circ \pi^{-1}_{\{e_i\}}(\cdot))$, $L_n^*(\Exp_{\theta}\circ \pi^{-1}_{\{e_i\}}(\cdot))$, $L_n(R_\theta(\cdot))$, $L_n^*(R_\theta(\cdot))$, respectively. 

We use the following coordinated forms to represent the first-order condition's solutions and the approximate solutions returned by Algorithm~\ref{algorithm: resampled newton iteration}: 
\begin{align}\label{eq: coordinated first order condition}
&\tilde {\mb \eta}_n  \text{ is the solution to }  \bar \nabla L_n^{\Exp, \theta_0, \{e_i\}}(\mb \eta) =\mb 0 ,\\
&\tilde {\mb \eta}_n  \text{ is the solution to }  \bar \nabla {L_n^*}^{\Exp, \theta_0, \{e_i\}}(\mb \eta) =\mb 0,\\ 
& \hat {\mb \eta}_n = \pi_{\{e_i\}}\circ \Log_{\theta_0}(\hat \theta_n), \quad  \hat {\mb \eta}_n^* = \pi_{\{e_i\}}\circ \Log_{\theta_0}(\hat \theta_n^*),
\end{align}
where $\hat \theta_n$ is produced by Algorithm~\ref{algorithm: resampled newton iteration} and $\hat \theta_n^*$ is a generic version of $\hat \theta_{n}^{*[i]}$ based on the dataset $\mc X_n^*$.  

Following the spirit of expressing quantities under a fixed coordinate system, we also define an analog of $\mb \Sigma$ in \eqref{eq: definition of Sigma} as follows: 
\begin{align}
&\mb \Omega \coloneqq \Big(\bb E\big[\bar \nabla^2 L^{\Exp,\theta_0, \{e_i\}}(\mb 0, X_1)\big]\Big)^{-1}\bb E\big[\bar \nabla  L^{\Exp,\theta_0, \{e_i\}}(\mb 0, X_1)\bar \nabla L^{\Exp,\theta_0, \{e_i\}}(\mb 0, X_1)\t\big]\\ 
& \qquad \cdot \Big(\bb E\big[\bar \nabla^2  L^{\Exp,\theta_0, \{e_i\}}(\mb 0, X_1)\big]\Big)^{-1}. 
\end{align}

One subtle yet crucial point is that the previously defined exact solutions $\tilde \theta_n$ and $\tilde\theta_n^*$ may not be unique, and the approximations $\hat \theta_n$ and $\hat \theta_n^*$ may be associated with a singular covariance estimate. On the computational side, we address potential singularity issues by employing the generalized matrix inverse; on the theoretical side, we will demonstrate that the solutions are well-defined and the covariance estimates are consistent with high probability, so that, in the exceptional event, the approximate solutions can be assigned arbitrary values without affecting the overall conclusion with the help of the delta method.

\subsubsection{Implicit Function of Approximate Solutions }\label{section: implicit function 1}
We first follow the route of \cite{bhattacharya1978validity} to show in the section that $\hat{\mb \eta}_n$ can be approximated by a smooth function of the sample-averaged derivatives of $L$. 

The first step is to prove the existence of a solution to \eqref{eq: first-order condition of M estimation}.
To begin with, we utilize the Taylor expansion of $L_n^{\Exp, \theta_0, \{e_i\}}$ 
 to derive that for each $i\in[p]$, 
\longeq{\label{eq:expansionofmanifoldloglikelihood}
    0 &= \frac{\partial}{\partial x_i} L^{\Exp, \theta_0, \{e_i\}}_n(\mb \eta)\\ & = e_iL_n( \theta_0) + \sum_{j}{\bar \nabla}^{\delta_i+\delta_j}L_n^{\Exp, \theta_0, \{e_i\}}(\mb 0 )\eta_j+  \sum_{\mb \nu, |\mb \nu| =2}\frac{\bar \nabla^{\mb \nu + \delta_i}L_n^{\Exp, \theta_0, \{e_i\}}(\mb 0)\mb \eta^{\mb \nu}}{\mb \nu !} + \text{Res}_i(\mb \eta),
}
where $\text{Res}_i(\mb \eta) \coloneqq    \sum_{\mb \nu: |\mb \nu| = 3}\frac{\bar \nabla^{\mb \nu + \delta_i} L_n^{\Exp, \theta_0, \{e_i\} }(t \mb \eta )\mb \eta^{\mb \nu}}{\mb \nu !}$ with some $t \in [0,1]$ represents the residual term and $\mb \nu!$ denotes the factorial product $\prod_{i=1}^{p}\nu_i!$.

Following the convention in the previous context, we denote that
\eq{\label{eq: definition of Zbar_n}
\bar{\mb Z}_n \coloneqq \sum_{i\in[n]}\mb Z_i/ n,
}
where $\mb Z_i$ is defined in Assumption~\ref{assumption:manifold}.\ref{item:E}.  

We first set out to control the fluctuation of $\bar{\mb Z}_n$.
By the finite moment conditions (Assumption~\ref{assumption:manifold}.\ref{item:B}) and the moderate deviation result (Lemma~\ref{lemma: thm1 in von1967central}), we have
\eq{
\label{eq:finitemomentconcentration}
    \bb P_{\mc X_n}\left[\big\|n^{\frac{1}{2}}\left(\bar{\mb Z}_n - \mb \mu\right)\big\|_2\geq \left(2\Lambda \log n\right)^{\frac{1}{2}}\right]\leq d_1 n^{-1}(\log n)^{-2},
}
where $\Lambda$ is the largest eigenvalue of the covariance matrix of $\mb Z_1$ and $d_1$ is some constant that depends on $\Lambda$. Similarly, applying Lemma~\ref{lemma: thm1 in von1967central} to $C(x)$ defined in Assumption~\ref{assumption:manifold}.\ref{item:B} yields that 
	\begin{align}\label{eq: C(X) bound}
		& \bb P_{\mc X_n}\Big[\frac{\sum_{i\in[n]}C(X_i)}{n} - \bb E[C( X_1)] \geq d_2n^{-\frac{1}{2}}(\log n)^{\frac{1}{2}} \Big] \leq d_2 n^{-1}(\log n)^{-2}
	\end{align}
	for some sufficiently large constant $d_2 $.

Looking into each entry of $\bar{\mb Z}_n$ separately, \eqref{eq:finitemomentconcentration} implies that
\longeq{ 
   &    \bb P_{\mc X_n}\left[\left|e_i L_n(\theta_0) \right| \geq d_3 n^{-\frac{1}{2}}\left( \log n \right)^{\frac{1}{2}}\right] \leq d_3 n^{-1}(\log n)^{-2}, 1\leq i \leq p,  \label{eq:derivativebound}\\ 
    &\bb P_{\mc X_n}\left[\big|\bar \nabla^{\mb \nu + \delta_i} L_n^{\Exp, \theta_0, \{e_i\}}(\mb 0) - \bb E[\bar \nabla^{\mb \nu} L^{\Exp, \theta_0, \{e_i\}}(\mb 0 , X_1)]\big| \geq d_3 n^{-\frac{1}{2}} \left(\log n\right)^{\frac{1}{2}}\right] \leq d_3n^{-1}(\log n)^{-2}, \\ 
    & \bb P_{\mc X_n}\Big[\big|\sum_{i\in[n]}\bar \nabla^{\mb \nu_1}L^{\Exp, \theta_0, \{e_i\}}(\tilde{ \mb \eta}_n, X_i )\bar \nabla^{\mb \nu_2}L^{\Exp, \theta_0, \{e_i\}}(\tilde{ \mb \eta}_n, X_i ) / n \\
    &- \bb E\big[\bar \nabla^{\mb \nu_1}L^{\Exp, \theta_0, \{e_i\}}(\tilde{ \mb \eta}_n, X_1 )\bar \nabla^{\mb \nu_2}L^{\Exp, \theta_0, \{e_i\}}(\tilde{ \mb \eta}_n, X_1 )\big]\big| \geq d_3 n^{-\frac{1}{2}} \left(\log n\right)^{\frac{1}{2}}\Big] \leq d_3n^{-1}(\log n)^{-2}, 
}
simultaneously hold for $1 \leq |\mb \nu| \leq 4$, $1\leq |\mb \nu_1|, |\mb \nu_2|\leq 3$, and some sufficiently large constant $d_3$. 

Moreover, invoking the definition of $C(\mb x)$ in Assumption~\ref{assumption:manifold}.\ref{item:B}, it follows from \eqref{eq: C(X) bound} that 
{\allowdisplaybreaks	
\begin{align}\label{eq: bound for residual term}
	&	\bb P_{\mc X_n}\Big[\max_{\mb \eta: \Exp_{\theta_0}\circ \pi^{-1}_{\{e_i\}}(\mb \eta) \in N} \Big|\big\vert\mathrm{Res}_i(\mb \eta)\big\vert -   \Big(\big(\sum_{\mb \nu: |\mb \nu | =3 } \big|\bb E_{\mc X_n}\big[\frac{1}{\mb \nu!}\bar \nabla^{\mb \nu+ \delta_i}L_n^{\Exp, \theta_0, \{e_i\}}(\mb 0)\big]\big| + \\ &\qquad   d_4n^{-1}(\log n)^{-2}\big)\mb \eta^{\mb \nu} + d_5\rho\big(\bb E[C(X_1)] +  d_2 n^{-\frac{1}{2}}(\log n)^{\frac{1}{2}}\big)\norm{\mb \eta}_2^3 \Big)\Big| \geq 0  \Big]\\ 
	 \leq & \bb P_{\mc X_n}\Big[
	 \Big(\sum_{\mb \nu: |\mb \nu| = 3}\frac{1}{\mb \nu !}\bar \nabla^{\mb \nu + \delta_i} L_n^{\Exp, \theta_0, \{e_i\} }(\mb 0) \\ 
     & \qquad - \sum_{\mb \nu: |\mb \nu | =3 } \big|\bb E\big[\frac{1}{\mb \nu!} \bar \nabla^{\mb \nu+ \delta_i}L_n^{\Exp, \theta_0, \{e_i\}}(\mb 0)\big]\big|\Big)\geq d_4n^{-1}(\log n)^{-2}\Big]\\
	  + &  
	 \bb P_{\mc X_n}\Big[\max_{\mb \eta: \Exp_{\theta_0}\circ \pi^{-1}_{\{e_i\}}(\mb \eta) \in N} \Big(\sum_{\mb \nu: |\mb \nu| = 3}\frac{1}{\mb \nu !} \max_{t \in [0,1]}\bar \nabla^{\mb \nu + \delta_i} L_n^{\Exp, \theta_0, \{e_i\} }(t \mb \eta)\mb \eta^{\mb \nu}\\ 
     & - \sum_{\mb \nu: |\mb \nu| = 3}\frac{1}{\mb \nu !}\bar \nabla^{\mb \nu + \delta_i} L_n^{\Exp, \theta_0, \{e_i\} }(\mb 0)\mb \eta^{\mb \nu} 
     - d_5\rho\big(\bb E[C(X_1)] +  d_2n^{-\frac{1}{2}}(\log n)^{\frac{1}{2}}\big)\norm{\mb \eta}_2^3\Big) \geq 0 \Big] \\ 
	 \leq & \bb P_{\mc X_n}\Big[
	 \Big(\sum_{\mb \nu: |\mb \nu| = 3}\frac{\bar \nabla^{\mb \nu + \delta_i} L_n^{\Exp, \theta_0, \{e_i\} }(\mb 0)}{\mb \nu !}\\ 
     & \qquad - \bb E\big[\sum_{\mb \nu: |\mb \nu | =3 } \frac{\bar \nabla^{\mb \nu+ \delta_i}L_n^{\Exp, \theta_0, \{e_i\}}(\mb 0)}{\mb \nu!}\big]\Big)\geq d_4n^{-1}(\log n)^{-2}\Big] \\ 
	 + & \bb P_{\mc X_n}\Big[\max_{\mb \eta: \Exp_{\theta_0}\circ \pi^{-1}_{\{e_i\}}(\mb \eta) \in N}\Big(\rho \sum_{\mb \nu: |\mb \nu| =3 }\frac{\mb \eta^{\mb \nu}}{\mb \nu!}\frac{\sum_{i\in[n]} C(X_i)}{n}  -  d_5\rho\big(\bb E[C(X_1)] \\ 
     & \qquad +  d_2n^{-\frac{1}{2}}(\log n)^{\frac{1}{2}}\big)\norm{\mb \eta}_2^3\Big) \geq 0  \Big]\\ 
	 \leq  & d_6 n^{-1}(\log n)^{-2}, 
\end{align}
 } holds for some constants $d_4, d_5, d_6$. 
Therefore, we can find an event $\mc F_n$ to be
\longeq{\label{eq: event F}
&\mathcal F_n \coloneqq    \Big\{\norm{\bar{\mb Z}_n}_\infty < d_7 n^{-\frac{1}{2}}\left( \log n \right)^{\frac{1}{2}} ,\big| \max_{\mb \eta: \Exp_{\theta_0}\circ \pi^{-1}_{\{e_i\}}(\mb \eta) \in N} \text{Res}_i(\mb \eta)\big|< d_7\norm{\mb \eta}_2^3, \text{ for }i,j \in[p]\Big\}
}
for some sufficiently large constant $d_7$ satisfying 
\eq{\label{eq: probability of F_n}
\bb P(\mc F_n^\complement) \leq d_7 n^{-1}(\log n)^{-2}.
}

Under the event $\mc F_n$ in \eqref{eq: event F},  multiplying by $\big(\bar \nabla^2 L_n^{\Exp, \theta_0, \{e_i\}}(\mb 0)\big)^{-1}$ both sides of \eqref{eq:expansionofmanifoldloglikelihood} and rearranging the terms yields that 
\eq{\label{eq:coordinate equation}
    \mb \eta =  \big(
    \bb E\big[\bar \nabla^2 L_n^{\Exp, \theta_0,\{e_i\}}(\mb 0)\big] + \mb \zeta_n\big)^{-1}\bigg[\mb\delta_n + \sum_{\mb \nu: |\mb \nu| =2}\frac{1}{\mb \nu!}
    \bar \nabla^{\mb \nu+\mb 1} L_n^{\Exp, \theta_0, \{e_i\}}(\mb 0)\mb \eta^{\mb \nu}  + \norm{\mb \eta}_2^3  \mb \epsilon_n\bigg]
}
where $\mb \zeta_n$ is a random matrix and $\mb \delta_n, \mb \epsilon_n$ are random vectors, each with a norm ($\norm{\cdot}_2$ or $\norm{\cdot}_F$) less than $d_7 pn^{-\frac{1}{2}}(\log n)^{\frac{1}{2}}$. 

Thus, by the Brouwer Fixed Point Theorem, at least one solution $\tilde {\mb \eta}_n$ exists solving \eqref{eq:coordinate equation} (also \eqref{eq:expansionofmanifoldloglikelihood}) under the event $\mc F_n$, within the ball $B_n \coloneqq   B\left(0,d_8 n^{-1}(\log n)^{-2}\right)\subset B(\theta_0, \frac{\rho}{2}) \subseteq N$ with a constant $d_8$ for some sufficiently large $n$.

 The next step is to leverage the implicit function theorem to obtain an approximation. We consider the following equation
\eq{\label{eq:manifoldequationsforexpansion}
    \mb 0=\mb z^{\delta_i} + \sum_{\mb \nu: |\mb \nu| = 1,2}\frac{1}{\mb \nu!} \mb z^{(\delta_i + \mb \nu)}\mb u^{\mb \nu}\eqqcolon P(\mb u,\mb z).
}
 
Obviously, the equation has a solution at $\mb u = \mb 0, \mb z= \mb \mu$ where $\mb \mu = \bb E\mb Z_1$ by Assumption~\ref{assumption:manifold}.\ref{item:C}. By the implicit function theorem and the smoothness of the function $P(\mb u,\mb z)$ with respect to $\mb u$ and $\mb z$, we conclude that there uniquely exists an injective and infinitely differentiable function $\mb H$ defined in a neighborhood of $\mb \mu$ such that \eq{\label{eq:implicitfunction}
\mb u = \mb H(\mb z) = \left(H_1(\mb z), \ldots, H_p(\mb z)\right)} with $\mb H(\mb \mu) = \mb 0$. We assume that $N$ is contained in the image of $\mb H$ by shrinking $N$ if necessary. 
This in turn implies that the first-order solution $\tilde\theta_n$ is unique on $N$ with $\dist(\tilde \theta_n, \theta_0) = O(n^{-1}\log n)$ under the event $\mc F_n$. 
Therefore it is straightforward by the existence of $\tilde{\mb \eta}_n$ that 
\longeq{
& \tilde {\mb \eta}_n = \mb H(\bar{\mb Z}_n') = \mb H(\bar{\mb Z}_n) + O(n^{-\frac32}(\log n)^{\frac32}) \label{eq: implicit function of manifold}
}
solves \eqref{eq:expansionofmanifoldloglikelihood} and also implies that $\Exp\circ\pi^{-1}_{\{e_i\}}(\tilde{\mb \eta}_n )$ is a solution of \eqref{eq: first-order condition of M estimation}, where $\bar{\mb Z}_n'$ is a modified version of $\bar{\mb Z}_n$, by replacing $e_i L_n(\theta_0)$ with $e_i L_n(\theta_0) + \mathrm{Res}_i(\tilde{\mb \eta}_n)$. 

Recall that the actual output $\hat\theta_n$ of Algorithm~\ref{algorithm: resampled newton iteration} is an approximate solution to $\eqref{eq: first-order condition of M estimation}$ according to Theorem~\ref{thm:convergence rate of newton}, and $\hat {\mb \eta}_n  = \pi_{\{e_i\}} \circ \Log_{\theta_0}(\hat \theta_n)$. 
Replacing $\tilde{\mb \eta}_n$ with  $\hat{\mb \eta}_n$ together with Theorem~\ref{thm:convergence rate of newton}, \eqref{eq: event F} yields that 
\longeq{\label{eq: approximation error of H(Z)}
& \pi_{\{T_{\theta_0\rightarrow \hat \theta_n}(e_i)\}}\circ \Log_{\hat\theta_n}(\theta_0) 
= - \left(\mb  H(\bar{\mb Z}_n)\right) + d_9n^{-\frac{3}{2}} (\log n)^{\frac{3}{2}}
}
under the event $\mathcal F_n$ for some constant $d_9$.

\subsubsection{Implicit Function of Approximate Resampled Solutions}
\label{section: implicit function 2}
We now move on to prove an empirical version of the implicit approximation. Given the extremum estimator $\hat{\theta}_n$ based on the dataset $\mc X_n = \{X_i\}$, we generate a resampled dataset $\mc X^*_n = \{X_i^*\}$ by drawing $n$ samples with replacement from $\mc X_n$, leading to the corresponding resampled extremum estimator $\hat{\theta}_n^*$ in Algorithm~\ref{algorithm: resampled newton iteration}.

Similar to \eqref{eq:expansionofmanifoldloglikelihood}, by performing a Taylor expansion of ${L_n^*}^{\Exp, \theta_0, \{e_i\}}$ at $\tilde {\mb \eta}_n$, we are interested in the solution to 
\longeq{\label{eq:empiricalexpansion}
    0 &= \frac{\partial}{\partial x_i}{L_n^*}^{\Exp, \theta_0, \{e_i\}}(\mb \eta)\\
    &
=\frac{\partial }{\partial x_i}{L_n^*}^{\Exp, \theta_0, \{e_i\}}(\tilde {\mb \eta}_n )+\sum_{j}\bar \nabla^{\delta_i,\delta_j}{L_n^*}^{\Exp, \theta_0, \{e_i\}}(\tilde{\mb \eta}_n) (\eta_j - \tilde \eta_{n,j})\\ 
&  + \frac{1}{\mb \nu!}\sum_{\mb \nu, |\mb \nu| = 2}\bar \nabla^{\mb \nu + \delta_i}{L_n^*}^{\Exp, \theta_0, \{e_i\}}(\tilde{\mb \eta}_n)(\mb \eta - \tilde{\mb \eta}_n)^{\mb \nu} 
 +\text{Res}_i^*(\mb \eta - \tilde{\mb \eta}_n )}
for $1\leq i \leq p$.

Given the samples $\mc X_n$, we naturally define a resampled version $\mb Z_n^*$ of the derivative collection $\mb Z_n$ as in Assumption~\ref{assumption:manifold}.\ref{item:E} by replacing $\mc X_n$ with $\mc X_n^*$. 
Conditional on $\mc X_n$, we shall control each entry of $\bar{\mb Z}_n^*$ by decomposing it into two parts: 
\begin{align}
	& \left| \bar \nabla^{\mb \nu}{L_n^*}^{\Exp, \theta_0, \{e_i\}}(\tilde{ \mb \eta}_n ) - \bar \nabla^{\mb \nu}L_n^{\Exp, \theta_0, \{e_i\}}(\tilde{ \mb \eta}_n)\right| \\
    \leq & \big| \underbrace{\bar \nabla^{\mb \nu}{L_n^*}^{\Exp, \theta_0, \{e_i\}}(\mb 0)}_{\text{an entry of $\bar{\mb Z}_n^*$}} -\underbrace{ \bar \nabla^{\mb \nu}L_n^{\Exp, \theta_0, \{e_i\}}(\mb 0)}_{\text{an entry of $\bar{\mb Z}_n$ }} \big| +  \left|\bar \nabla^{\mb \nu}{L_n^*}^{\Exp, \theta_0, \{e_i\}}(\tilde {\mb \eta}_n ) -  \bar \nabla^{\mb \nu}{L_n^*}^{\Exp, \theta_0, \{e_i\}}(\mb 0)\right|\\ 
    + & \left|\bar \nabla^{\mb \nu}L_n^{\Exp, \theta_0, \{e_i\}}(\tilde {\mb \eta}_n ) -  \bar \nabla^{\mb \nu}L_n^{\Exp, \theta_0, \{e_i\}}(\mb 0)\right|,
    \label{eq: nabla nu difference}  \\ 
    & \hspace{250pt} \text{for $1 \leq |\mb \nu| \leq 4$}, 
\end{align}
\begin{align}
    & \Big| \sum_{i\in[n]}\bar \nabla^{\mb \nu_1}L^{\Exp, \theta_0, \{e_i\}}(\tilde{ \mb \eta}_n, X_i^* )\bar \nabla^{\mb \nu_2}L^{\Exp, \theta_0, \{e_i\}}(\tilde{ \mb \eta}_n, X_i^* ) / n \\ 
    & - \sum_{i\in[n]}\bar \nabla^{\mb \nu_1}L^{\Exp, \theta_0, \{e_i\}}(\tilde{ \mb \eta}_n, X_i )\bar \nabla^{\mb \nu_2}L^{\Exp, \theta_0, \{e_i\}}(\tilde{ \mb \eta}_n, X_i ) / n\Big| \\ 
    \leq& \Big| \sum_{i\in[n]}\bar \nabla^{\mb \nu_1}L^{\Exp, \theta_0, \{e_i\}}(\mb 0, X_i^* )\bar \nabla^{\mb \nu_2}L^{\Exp, \theta_0, \{e_i\}}(\mb 0, X_i^* ) / n\\ 
    &  - \sum_{i\in[n]}\bar \nabla^{\mb \nu_1}L^{\Exp, \theta_0, \{e_i\}}(\mb 0, X_i^* )\bar \nabla^{\mb \nu_2}L^{\Exp, \theta_0, \{e_i\}}(\mb 0, X_i ) / n \Big| \\ 
    &   + \Big| \sum_{i\in[n]}\bar \nabla^{\mb \nu_1}L^{\Exp, \theta_0, \{e_i\}}(\tilde{ \mb \eta}_n, X_i^* )\bar \nabla^{\mb \nu_2}L^{\Exp, \theta_0, \{e_i\}}(\tilde{ \mb \eta}_n, X_i^* ) / n \\ 
    &  -  \sum_{i\in[n]}\bar \nabla^{\mb \nu_1}L^{\Exp, \theta_0, \{e_i\}}(\mb 0, X_i^* )\bar \nabla^{\mb \nu_2}L^{\Exp, \theta_0, \{e_i\}}(\mb 0, X_i^* ) / n \Big|\\
    &+  \Big| \sum_{i\in[n]}\bar \nabla^{\mb \nu_1}L^{\Exp, \theta_0, \{e_i\}}(\tilde{ \mb \eta}_n, X_i )\bar \nabla^{\mb \nu_2}L^{\Exp, \theta_0, \{e_i\}}(\tilde{ \mb \eta}_n, X_i ) / n \\ 
    &  -  \sum_{i\in[n]}\bar \nabla^{\mb \nu_1}L^{\Exp, \theta_0, \{e_i\}}(\mb 0, X_i )\bar \nabla^{\mb \nu_2}L^{\Exp, \theta_0, \{e_i\}}(\mb 0, X_i ) / n \Big|,\qquad  \text{for $1 \leq |\mb \nu_1|, |\mb \nu_2| \leq 3$}. \label{eq: quadratic nabla nu difference}
\end{align}

For an arbitrary entry $(\mb Z^*_1)_k$ of $\mb Z^*_1$, we notice that
\longeq{
& \bb E\left[\left|{(\mb Z_1^*)_k} 
- (\bar{\mb Z}_n)_k\right|^4|\mc X_n\right]\\ 
\leq & 2^4 \max_{0 \leq r \leq 4} \bb E\big[(\mb Z_1^*)_k^r |\mc X_n\big](\bar{\mb Z}_n)_k^{4-r}\\ 
\stackrel{\text{ by Jensen's inequality}}{\leq} & 2^4 \max_{0 \leq r \leq 4} \bb E\big[(\mb Z_1^*)_k^4 |\mc X_n\big]^{\frac{r}{4}}(\bar{\mb Z}_n)_k^{4-r}
\\ 
\leq & 2^4 \bb E[(\mb Z^*)_k^4] + 2^4 |(\bar{\mb Z}_n)_k|^4.
}
Hence, we have
\longeq{
\label{eq: resampled Z moment bound}
    &\bb P_{\mc X_n}\Big[\bb E\big[|{(\mb Z_1^*)_k} -(\bar{\mb Z}_n)_k|^4 \big|\mc X_n\big] \geq 2^4\left(\bb E[|(\mb Z_1)_k|^4] + |\bb E[(\bar{\mb Z})_k]|^4\right) + 2^5\Big] \\
    \leq & \bb P_{\mc X_n}\Big[ \Big|
    \sum_{i=1}^n \left(|(\mb Z_i)_k|^4 - \bb E |(\mb Z_1)_k|^4\right)
    \Big| + 
    n\left||(\bar{\mb Z}_n)_k|^4 - |\bb E [(\bar{\mb Z}_n)_k]|^4 \right|\geq 2n
    \Big]\\
    \leq & \bb P_{\mc X_n}\Big[ \Big|
    \sum_{i=1}^n \left(|(\mb Z_i)_k|^4 - \bb E[ |(\mb Z_1)_k|^4]\right)
    \Big|\geq n \Big] + \bb P_{\mc X_n}\Big[\left||(\bar{\mb Z}_n)_k|^4 - |\bb E [(\bar{\mb Z}_n)_k]|^4 \right|\geq 1
    \Big]\\
    \leq &n^{-2} \bb E_{\mc X_n}\Big[ 
    \sum_{i=1}^n \big(|(\mb Z_i)_k|^4 - \bb E_{\mc X_n} [|(\mb Z_i)_k|^4]\big)^{2}
    \Big] + \bb E_{\mc X_n}\big[\big|(\bar{\mb Z}_n)_k \big|^4\big]\\ 
     \leq & c_{\mb \nu}n^{-1}\bb E_{\mc X_n}\big[| (\mb Z_1)_k|^8\big] + c_{k} n^{-2}\bb E_{\mc X_n}\big[|(\mb Z_1)_k|^4\big]\\ 
     = &O(n^{-1}) 
}
for some constant $c_{k}$, where the finite moment condition in Assumption~\ref{assumption:manifold}.\ref{item:B} comes into play. Here the penultimate line follows from the observation that, by interpreting $S_t\coloneqq \sum_{i=1}^t (\mb Z_i)_k$ for $t = 1,\cdots, n$ as a martingale with its natural filtration, and applying Lemma~\ref{lemma: Burkholder's inequality} along with Jenson's inequality, we obtain: 
\eq{
    \bb E_{\mc X_n}\big[| (\bar{\mb Z}_n)_k|^4\big] = n^{-4} \bb E_{\mc X_n}[S_n^4] \leq c_k n^{-4} \bb E_{\mc X_n}\big[ |\sum_{i\in[n]}(\mb Z_i)_k^2|^2\big] \leq c_k n^{-2} \bb E_{\mc X_n}\big[ |(\mb Z_1)_k|^4\big].
}

\paragraph*{Analysis of Eq.~\eqref{eq: nabla nu difference}}
For the term 
$\left|\bar \nabla^{\mb \nu} {L^*_n}^{\Exp, \theta_0, \{e_i\}}(\tilde{\mb \eta}_n) -  \bar \nabla^{\mb \nu}{L^*_n}^{\Exp, \theta_0, \{e_i\}}(\mb 0)\right|$ with  $1\leq| \boldsymbol \nu|\leq 4$, 
Assumption~\ref{assumption:manifold}.\ref{item:B} implies  that
\longeq{
& \big\vert \bar \nabla^{\mb \nu} L^{\Exp, \theta_0,\{e_i\}} (\tilde{\mb \eta}_n, x) - \bar \nabla^{\mb \nu} L^{\Exp, \theta_0, \{e_i\}}(\mb 0, x) \big \vert \leq  C(x)\dist(\tilde \theta_n, \theta_0),
}
which leads to 
\begin{align}
& \label{eq: nabla nu difference 1}\big| \bar \nabla^{\mb \nu} L_n^{\Exp, \theta_0, \{e_i\}}(\tilde{\mb \eta}_n) -\bar \nabla^{\mb \nu} L_n^{\Exp, \theta_0, \{e_i\}}(\mb 0) \big| \leq \frac{\sum_{i\in[n]} C(X_i) \dist(\tilde \theta_n^*, \theta_0)}{n}, 1\leq |\mb \nu| \leq 4, \\ 
& \big|\bar \nabla^{\mb \nu}{ L_n^*}^{\Exp, \theta_0, \{e_i\}}(\tilde{\mb \eta}_n) -  \bar \nabla^{\mb \nu}{ L_n^*}^{\Exp, \theta_0, \{e_i\}}(\mb 0)\big|\leq \frac{\sum_{i\in[n]}C(X_i^*)\dist(\tilde \theta_n, \theta_0)}{n}, 1\leq |\mb \nu| \leq 4. \label{eq: nabla nu difference 2}
\end{align}

Analogous to the bound in \eqref{eq: resampled Z moment bound}, we establish that
\begin{align} 
& \bb P_{\mc X_n}\bigg[\bb E\Big[|C(X_1^*) - \sum\limits_{i\in[n]}C(X_i)/ n |^4\vert\mc X_n\Big] \geq 2^3(\bb E[|C(X)|^4] + |\bb E[C(X)]|^4) + 2^5 \bigg] = O(n^{-1}).
\end{align}
This together with Lemma~\ref{lemma: thm1 in von1967central} leads to the fact that 
\begin{align}
& \bb P_{\mc X_n}\bigg[\bb P\Big[\frac{\sum_{i\in[n]}C(X_i^*)}{ n}  \geq  f_1n^{-\frac{1}{2}}(\log n)^{\frac{1}{2}} + \bb E[C(X)] \vert \mc X_n\Big] \geq  f_1n^{-1}(\log n)^{-2} \bigg]= O(n^{-1}).\label{eq: resampled C(X) bound}
\end{align}
for some sufficiently large constant $f_1$.

By \eqref{eq: nabla nu difference}, we can bound the difference $\big|\bar \nabla^{\mb \nu}{L_n^*}^{\Exp, \theta_0, \{e_i\}}(\tilde {\mb \eta}_n ) -\bar \nabla^{\mb \nu}{L_n^{\Exp, \theta_0, \{e_i\}}}(\tilde {\mb \eta}_n )   \big|$ via 
\longeq{\label{eq: decomposition of nabla L_n}
& \big|\bar \nabla^{\mb \nu}{L_n^*}^{\Exp, \theta_0, \{e_i\}}(\tilde {\mb \eta}_n ) -\bar \nabla^{\mb \nu}{L_n^{\Exp, \theta_0, \{e_i\}}}(\tilde {\mb \eta}_n )\big| \\ 
\leq & \big|\bar \nabla^{\mb \nu}{L_n^*}^{\Exp, \theta_0, \{e_i\}}(\mb 0) -\bar \nabla^{\mb \nu}{L_n^{\Exp, \theta_0, \{e_i\}}}(\mb 0)\big|+ \frac{C\sum_{i\in[n]} C(X_i)\norm{\tilde {\mb \eta}_n}_2 }{n} + \frac{C\sum_{i\in[n]} C(X^*_i)\norm{\tilde {\mb \eta}_n}_2}{n}.
}

The decomposition \eqref{eq: decomposition of nabla L_n} together with \eqref{eq: C(X) bound}, \eqref{eq:derivativebound}, \eqref{eq: resampled Z moment bound}, and \eqref{eq: resampled C(X) bound} gives that 
\begin{equation*}
\begin{split}
& \bb P_{\mc X_n}\bigg[\bb P\Big[|\bar \nabla^{\mb \nu}{L_n^*}^{\Exp, \theta_0, \{e_i\}} (\tilde {\mb \eta}_n ) -\bar \nabla^{\mb \nu}{L_n^{\Exp, \theta_0, \{e_i\}}}(\tilde {\mb \eta}_n ) | \geq f_2n^{-\frac{1}{2}}(\log n)^{\frac{1}{2}} \mid \mc X_n\Big]\leq  f_2n^{-1}(\log n)^{-2} \bigg]\\
\geq & 1- f_2n^{-1}
\end{split}
\end{equation*}
for some sufficiently large constant $f_2$.
\paragraph*{Analysis of Eq.~\eqref{eq: quadratic nabla nu difference} }
In terms of the counterparts to the quadratic terms in $\bar{\mb Z}_n$ as decomposed in \eqref{eq: quadratic nabla nu difference}, we begin by bounding each term on the right-hand side of \eqref{eq: quadratic nabla nu difference}. 
\begin{itemize}
    \item Note that \eqref{eq: resampled Z moment bound} has already implied the uniform boundedness of the fourth moment for every entry of $\mb Z_1$, except on an event with probability $O(n^{-1})$. 
Invoking Lemma~\ref{lemma: thm1 in von1967central} again, one has 
\begin{align}
    & \bb P_{\mc X_n}\bigg[ \bb P\Big[\Big| \sum_{i\in[n]}\bar \nabla^{\mb \nu_1}L^{\Exp, \theta_0, \{e_i\}}(\mb 0, X_i^* )\bar \nabla^{\mb \nu_2}L^{\Exp, \theta_0, \{e_i\}}(\mb 0, X_i^* ) / n\\ 
    & \qquad  - \sum_{i\in[n]}\bar \nabla^{\mb \nu_1}L^{\Exp, \theta_0, \{e_i\}}(\mb 0, X_i^* )\bar \nabla^{\mb \nu_2}L^{\Exp, \theta_0, \{e_i\}}(\mb 0, X_i ) / n \Big| \geq f_3 n^{-\frac{1}{2}}(\log n)^{\frac{1}{2}} \mid \mc X_n \Big]\\ 
    & \qquad \leq  f_3 n^{-1}(\log n)^{-2} \bigg]\geq 1 - f_3n^{-1}\label{eq: quadratic nabla nu difference 1}
\end{align}
for some constant $f_3$. 
\item Further, by Assumption~\ref{assumption:manifold}.\ref{item:B}, the triangle inequality, and the Cauchy inequality, the second term in \eqref{eq: quadratic nabla nu difference} can be bounded by
\begin{align}
    & \Big| \sum_{i\in[n]}\bar \nabla^{\mb \nu_1}L^{\Exp, \theta_0, \{e_i\}}(\tilde{ \mb \eta}_n, X_i^* )\bar \nabla^{\mb \nu_2}L^{\Exp, \theta_0, \{e_i\}}(\tilde{ \mb \eta}_n, X_i^* ) / n \\ 
    &  -  \sum_{i\in[n]}\bar \nabla^{\mb \nu_1}L^{\Exp, \theta_0, \{e_i\}}(\mb 0, X_i^* )\bar \nabla^{\mb \nu_2}L^{\Exp, \theta_0, \{e_i\}}(\mb 0, X_i^* ) / n \Big|\\ 
    \leq & \sum_{i\in[n]}C(X_i^*)\dist(\tilde\theta_n, \theta_0) \big(  |\bar \nabla^{\mb \nu_1} L^{\Exp, \theta_0, \{e_i\}}(\mb 0, X_i^* )| +   |\bar \nabla^{\mb \nu_2} L^{\Exp, \theta_0, \{e_i\}}(\mb 0, X_i^* )|\big) / n\\ 
    & \quad +  \sum_{i\in[n]}C(X_i^*)^2 \dist(\tilde\theta_n, \theta_0)^2 / n \\
    \leq & \dist(\tilde\theta_n, \theta_0)\big(\sum_{i\in[n]}C(X_i^*)^2 / n\big) \cdot   \Big(2 \big(\sum_{i\in[n]} \big(\bar \nabla^{\mb \nu_1} L^{\Exp, \theta_0, \{e_i\}}(\mb 0, X_i^* ) \big)^2/ n\big)^{\frac{1}{2}} \\ 
    & \qquad +  2\big(\sum_{i\in[n]} \big(\bar \nabla^{\mb \nu_2} L^{\Exp, \theta_0, \{e_i\}}(\mb 0, X_i^* ) \big)^2/ n\big)^{\frac{1}{2}} + 1 \Big),
\end{align}
where we use $\dist(\tilde{\theta}_n, \theta_0)^2 \leq \dist(\tilde\theta_n, \theta_0)$ in a sufficiently small neighborhood $N$ to obtain the last inequality. 

 Following the route toward proving \eqref{eq: resampled Z moment bound} exploiting Lemma~\ref{lemma: thm1 in von1967central}, we derive that 
    $ \sum_{i\in[n]}C(X_i^*)^2 / n $ is upper bounded by a constant with probability at least $1- f_{C} n^{-1}$ with some constant $f_{C}$, 
 where we make use of the finite $16$-th moment condition on $C(X_1)$. This together with \eqref{eq: quadratic nabla nu difference 1} and the concentration on $\dist(\tilde\theta_n, \theta_0)$ under $\mc F_n$ yields that
\longeq{
   &\bb P_{\mc X_n}\bigg[ \bb P\Big[ \Big| \sum_{i\in[n]}\bar \nabla^{\mb \nu_1}L^{\Exp, \theta_0, \{e_i\}}(\tilde{ \mb \eta}_n, X_i^* )\bar \nabla^{\mb \nu_2}L^{\Exp, \theta_0, \{e_i\}}(\tilde{ \mb \eta}_n, X_i^* ) / n \\ 
   &  -  \sum_{i\in[n]}\bar \nabla^{\mb \nu_1}L^{\Exp, \theta_0, \{e_i\}}(\mb 0, X_i^* )\bar \nabla^{\mb \nu_2}L^{\Exp, \theta_0, \{e_i\}}(\mb 0, X_i^* ) / n \Big| \geq  f_4 n^{-\frac{1}{2}}(\log n)^{\frac{1}{2}} \mid \mc X_n\Big] \leq f_4 n^{-1}(\log n)^{-2} \bigg] \\
   \geq& 1 - f_4 n^{-1}
   \label{eq: quadratic nabla nu difference 2}
}
for some constant $f_4$.
\item 
Similar to the proof of \eqref{eq: quadratic nabla nu difference 2}, for the third term in \eqref{eq: quadratic nabla nu difference} we have
\begin{align}
    &  \Big| \sum_{i\in[n]}\bar \nabla^{\mb \nu_1}L^{\Exp, \theta_0, \{e_i\}}(\tilde{ \mb \eta}_n, X_i )\bar \nabla^{\mb \nu_2}L^{\Exp, \theta_0, \{e_i\}}(\tilde{ \mb \eta}_n, X_i ) / n \\ 
    &  -  \sum_{i\in[n]}\bar \nabla^{\mb \nu_1}L^{\Exp, \theta_0, \{e_i\}}(\mb 0, X_i )\bar \nabla^{\mb \nu_2}L^{\Exp, \theta_0, \{e_i\}}(\mb 0, X_i ) / n \Big|\\ 
    \leq & f_5 n^{-\frac{1}{2}}(\log n)^{\frac{1}{2}} \label{eq: quadratic nabla nu difference 3}
\end{align}
with probaiblity at least $1 - f_5 n^{-1}(\log n)^{-2} $. 
\end{itemize}

Combining \eqref{eq: quadratic nabla nu difference 1}, \eqref{eq: quadratic nabla nu difference 2}, and \eqref{eq: quadratic nabla nu difference 3} together with \eqref{eq: quadratic nabla nu difference} gives that
\begin{align}
    &\bb P_{\mc X_n}\bigg[ \bb P\Big[ \Big| \sum_{i\in[n]}\bar \nabla^{\mb \nu_1}L^{\Exp, \theta_0, \{e_i\}}(\tilde{\mb \eta}_n, X_i^* )\bar \nabla^{\mb \nu_2}L^{\Exp, \theta_0, \{e_i\}}(\tilde{\mb \eta}_n, X_i^* ) / n\\
    &  -  \sum_{i\in[n]}\bar \nabla^{\mb \nu_1}L^{\Exp, \theta_0, \{e_i\}}(\tilde{\mb \eta}_n, X_i )\bar \nabla^{\mb \nu_2}L^{\Exp, \theta_0, \{e_i\}}(\tilde{\mb \eta}_n, X_i ) / n \Big| \geq  f_6 n^{-1}(\log n)^{-2} \mid \mc X_n\Big] \leq f_6 n^{-1}(\log n)^{-2} \bigg] \\
    \geq& 1 - f_6 n^{-1}
\end{align}
holds for some constant $f_6$.

\paragraph*{Bound for Residual Term }
Finally, by a similar argument to \eqref{eq: bound for residual term}, \eqref{eq: resampled C(X) bound}, and \eqref{eq: decomposition of nabla L_n}, we can also prove that 
\begin{align}
     & \bb P_{\mc X_n}\bigg[\bb P\Big[\max_{\eta \in N}\left| \text{Res}_i^*(\mb \eta - \tilde{\mb \eta}_n)\right| -  f_7\norm{\mb \eta - \tilde{\mb \eta}_n}_2^s \left(1 + n^{-\frac{1}{2}} \left(\log n\right)^{\frac{1}{2}}\right) \geq 0 |\mc X_n\Big] \leq f_7n^{-1}(\log n)^{-2} \bigg] \\
     & \geq 1- f_7n^{-1} \label{eq:empiricalresonmanifold}
 \end{align} 
with sufficiently large constant $f_7$. 

\paragraph*{Putting Pieces Together }
Bringing together the above concentration inequalities, we denote by $\mc F_n^*$ the event
\longeq{ \label{eq: event F_n^*}
    \mc F_n^* \coloneqq \Big\{&\max_{1 \leq |\mb \nu|\leq 4}\big\{|\bar \nabla^{\mb \nu} {L_n^*}^{\Exp, \theta_0, \{e_i\}}(\mb \eta) - \bar \nabla^{\mb\nu} L_n^{\Exp, \theta_0, \{e_i\}}(\mb \eta) |, 1\leq |\mb \nu | \leq 4\big\}< f_8n^{-\frac{1}{2}}\left( \log n \right)^{\frac{1}{2}}, \\ 
 & \max_{1\leq |\mb \nu_1|, |\mb \nu_2|\leq 3}\Big\{\Big| \sum_{i\in[n]}\bar \nabla^{\mb \nu_1}L^{\Exp, \theta_0, \{e_i\}}(\tilde{\mb \eta}_n, X_i^* )\bar \nabla^{\mb \nu_2}L^{\Exp, \theta_0, \{e_i\}}(\tilde{\mb \eta}_n, X_i^* ) / n\\
 &  -  \sum_{i\in[n]}\bar \nabla^{\mb \nu_1}L^{\Exp, \theta_0, \{e_i\}}(\tilde{\mb \eta}_n, X_i )\bar \nabla^{\mb \nu_2}L^{\Exp, \theta_0, \{e_i\}}(\tilde{\mb \eta}_n, X_i ) / n \Big| \Big\} < f_8n^{-\frac{1}{2}}\left( \log n \right)^{\frac{1}{2}},\\
 & \frac{\sum_{i\in[n]}C(X_i^*)}{ n} <   f_8n^{-1}(\log n)^{-2} + \bb E[C(X)], \norm{\bar{\mb Z}_n^* - \bar{\mb Z}_n}_\infty <   f_8n^{-\frac{1}{2}}(\log n)^{\frac{1}{2}}, 
 \\&\left|\text{Res}^*(\mb \eta^* - \tilde {\mb \eta}_n)\right|< f_8\norm{\mb \eta- \tilde{\mb \eta}_n }_2^3 \Big\},
}
and by $\mc F_n'$ the event 
 \longeq{\label{eq: event F'}
 &\mc F_n' \coloneqq    \Big\{ \bb P\big[\mathcal F_n^* \mid \mc X_n\big] \leq f_8 n^{-1}(\log n)^{-2} \Big\} \bigcap \mc F_n
 } 
with a sufficiently large $f_8$, where its complementary event is controlled by
 \eq{\label{eq: exceptional probability of F_n'}
 \bb P_{\mc X_n}\left[{\mc F_n'}^\complement\right]\leq f_8n^{-1}
}
according to the above analysis. 

 Again, using a similar argument to the first part and the Brouwer fixed point theorem but under the event $\mathcal F_n^*\cap \mc F_n$, we can prove the existence of a solution solving 
\longeq{\label{eq: resampled brouwer fixed point}
&  \mb \eta -\tilde{\mb \eta }_n =  (
    \bar \nabla^2 {L_n^*}^{\Exp, \theta_0, \{e_i\}}(\tilde{\mb \eta}_n)+ \mb \zeta_n^*)^{-1}\\ 
    & \qquad \qquad \cdot \Big[\mb\delta_n^* + \sum_{ |\mb \nu|=2}\frac{1}{\mb \nu!}
    \bar \nabla^{\mb \nu+\mb 1}{ L_n^*}^{\Exp, \theta_0, \{e_i\}} (\tilde{\mb \eta}_n)({\mb \eta} - \tilde{\mb \eta}_n)^{\mb \nu}  + \norm{\mb \eta - \tilde{\mb \eta}_n	}_2^3 \mb \epsilon_n^*\Big],
}
where $\mb \epsilon_n^*$ is a random vector related to ${\bar \nabla^{\mb \nu} {L_n^*}^{\Exp, \theta_0, \{e_i\}}}(t(\tilde{\mb \eta}_n) + (1-t) \mb \eta),| \mb \nu| = 3, t\in[0,1]$. We denote the solution to \eqref{eq: resampled brouwer fixed point} by $\tilde{\mb \eta}_n^*$. 

Recall the smooth function defined in \eqref{eq:implicitfunction}. Applying the implicit function theorem yields that
 \eq{
    \tilde {\mb \eta}_n^* - \tilde{\mb \eta}_n  = \mb H(\bar{\mb Y}'_n )
 }
 under the event $\mathcal F_n^* \cap \mc F_n$ for a sufficiently large $n$, where $\bar{\mb Y}_n'$ corresponds to $\bar{\mb Z}_n'$ in \eqref{eq: implicit function of manifold} with $\mc X_n$ replaced by $\mc X_n^*$, while the derivatives are evaluated at $\tilde{\mb \eta}_n$ instead of $\mb 0$.

Since the implicit function $\mb H$ is deterministic and smooth, replacing $\bar{\mb Y}_n'$ with $\bar{\mb Y}_n$ implies that 
\eq{ \label{eq: approximation error of H(Y)}
    \tilde {\mb \eta}_n^* = \mb H(\bar{\mb Z}_n) + \mb H(\bar{\mb Y}_n ) + f_9n^{-\frac32}(\log n)^{\frac32},  
}
with some constant $f_9$, where we invoke \eqref{eq: implicit function of manifold} to represent $\tilde {\mb \eta}_n$, and $\bar{\mb Y}_n$ is the counterpart to $\bar{\mb Z}_n$ with $\mc X_n$ replaced by $\mc X_n^*$ with the derivatives are taken at $\tilde{\mb \eta}_n$ instead of $\mb 0$. 

Invoking Theorem~\ref{thm: convergence rate of resampled newton}, \eqref{eq: approximation error of H(Y)} turns out to be
\eq{
\hat {\mb \eta}_n^* = \mb H(\bar{\mb Z}_n) + \mb H(\bar{\mb Y}_n ) + O(n^{-\frac32}(\log n)^{\frac32}) 
}
under $\mc F_n^* \cap \mc F_n$. 

Finally, we relate the coordinate difference between $\hat \theta_n^*$ and $\hat \theta_n$ under the chart $ \pi_{\{e_i\}}\circ \Log_{\theta_0}(\cdot )$ to the coordinate of $\hat \theta_n$ under the chart $\pi_{\{T_{\theta_0\rightarrow \hat \theta_n^*}(e_i)\}}\circ \Log_{\hat\theta_n^*}(\cdot )$. Specifically,
by Lemma~\ref{lemma: normal coordinate transomation}, it holds 
\longeq{\label{eq: approximation error of H(Y) 2}
&\pi_{\{T_{\theta_0\rightarrow \hat \theta_n^*}(e_i)\}}\circ \Log_{\hat\theta_n^*}(\hat \theta_n) =  -\mb H(\bar{\mb Y}_n)  + f_{10}n^{-\frac{3}{2}} (\log n)^{\frac{3}{2}}
}
under the event $\mc F_n^* \cap \mc F_n$, provided a sufficiently large constant $f_{10}$.

\subsubsection{High-Order Asymptotics of Approximate Solutions under Second-Order Retractions}
 Now we are in a position to combine the above pieces together so as to establish the high-order asymptotics of the iteration series generated by Algorithm~\ref{algorithm:wald based statistic} and Algorithm~\ref{algorithm:t based statistic}. We first introduce two covariance matrices associated with the implicit functions in the previous steps:
 \begin{align}
     &\hat{\mb \Omega} \coloneqq \Big(\bar \nabla^2  L_n^{\Exp,\theta_0, \{e_i\}}(\mb H(\bar{\mb Z}_n)) \Big)^{-1}\Big(\sum_{i\in[n]}\bar \nabla L^{\Exp,\theta_0, \{e_i\}}(\mb H(\bar{\mb Z}_n), X_i)\bar \nabla L^{\Exp,\theta_0, \{e_i\}}(\mb H(\bar{\mb Z}_n), X_i)\t / n\Big)\\ 
&\qquad \cdot \Big(\bar \nabla^2 L^{\Exp,\theta_0, \{e_i\}}(\mb H(\bar{\mb Z}_n))\Big)^{-1},\label{eq: hat Omega}
\\
&
\hat{\mb \Omega}^* \coloneqq \Big(\bar \nabla^2  {L_n^*}^{\Exp,\theta_0, \{e_i\}}(\mb H(\bar{\mb Z}_n) + \mb H(\bar{\mb Y}_n))\Big)^{-1}\\ 
& \qquad \cdot \Big(\sum_{i\in[n]}\bar \nabla L^{\Exp,\theta_0, \{e_i\}} (\mb H(\bar{\mb Z}_n) + \mb H(\bar{\mb Y}_n), X_i^*)\bar \nabla L^{\Exp,\theta_0, \{e_i\}} (\mb H(\bar{\mb Z}_n) + \mb H(\bar{\mb Y}_n), X_i^*)\t / n\Big) \\
& \qquad \cdot \Big(\bar \nabla^2{L_n^*}^{\Exp,\theta_0, \{e_i\}}(\mb H(\bar{\mb Z}_n) + \mb H(\bar{\mb Y}_n))\Big)^{-1}.  
\label{eq: hat Omega^*}
 \end{align}
In terms of the involved matrix inversion, we will shortly elaborate on the well-definedness of these quantities under $\mc F_n$ or $\mc F_n^*$.
 
 The essential idea to demonstrate the high-order asymptotics is to exploit the delta method (Lemma~\ref{lemma:deltamethod}) to prove the equivalence between the approximate solutions and the implicit function after coordinate transformations. 
 
\paragraph*{Proof of \eqref{eq: thm1.1} (Wald Statistic) }
For the quantity $R^{-1}_{\hat \theta_n}(\theta_0)\t \hat \Sigma^{\dagger} R^{-1}_{\hat \theta_n}(\theta_0)\t$, we once again note that the generalized-inverse is designed to prevent the algorithm from failure after hundred thousands of resampling in practice; in the theoretical analysis, the generalized inverse can be replaced by the matrix inverse under the event $\mathcal F_n$ defined in \eqref{eq: event F}.

Next, we are going to express $R^{-1}_{\hat \theta_n}(\theta_0)\t \hat{\mb  \Sigma}^{-1} R^{-1}_{\hat \theta_n}(\theta_0)$ and $R^{-1}_{\hat \theta^*_n}(\hat \theta_n)\t \check {\mb \Sigma}^{-1} R^{-1}_{\hat \theta_n^*}(\hat \theta_n)$ in the coordinate system of $\pi_{\{e_i\}}\circ \Log_{\theta_0}(\cdot)$. 
 
The core of this analysis is that the difference across various charts can be eliminated in the Wald-type statistic regardless of the bases across various charts. Formally, we let $\hat{\mb Q}$ denote the orthonormal matrix whose $i$-th row corresponds to the coordinate of $\mathrm d R_{\hat \theta_n}^{-1}[T_{\theta_0 \rightarrow \hat \theta_n} e_i]$. Similarly, let $\hat{\mb Q}^*$ denote an analogous orthogonal matrix, with the $i$-th row being $\mathrm d R_{\hat \theta_n^*}^{-1}[T_{\theta_0 \rightarrow \hat \theta_n^*} e_i]$. Then, under the event $\mc F_n$ and the event $\mc F_n^* \cap \mc F_n$, respectively, applying the property of second-order retractions (Lemma~\ref{lemma: uniform control for second-order retraction}) gives that 
\longeq{
\label{eq: approximation error of retraction mapping with rotations}
& \hat{\mb Q} R_{\hat\theta_n}(\theta_0) = \pi_{\{T_{\theta_0\rightarrow \hat \theta_n}e_i\}} \circ \Log_{\hat\theta_n}(\theta_0)  + O(n^{-\frac{3}{2}}(\log n)^{\frac{3}{2}}) =  - \pi_{\{e_i\}} \circ \Log_{\theta_0}(\hat\theta_n) + O(n^{-\frac{3}{2}}(\log n)^{\frac{3}{2}}),\\
&  \hat{\mb Q}^* R_{\hat\theta_n^*}(\theta_0) = \pi_{\{T_{\theta_0\rightarrow \hat \theta_n^*}e_i\}} \circ \Log_{\hat\theta_n^*}(\theta_0)  + O(n^{-\frac{3}{2}}(\log n)^{\frac{3}{2}}) =  - \pi_{\{e_i\}} \circ \Log_{\theta_0}(\hat\theta_n^*) + O(n^{-\frac{3}{2}}(\log n)^{\frac{3}{2}}).
}

Invoking the property of the second-order retraction, it holds for every $\mb x \in U$ that
	\longeq{
	& \bar \nabla L^{\Exp, \hat \theta_n, \{T_{\theta_0 \rightarrow \hat \theta_n} e_i\}}(\mb 0, \mb x) = \hat{\mb Q}\bar \nabla L^{\ret, \hat \theta_n}(\mb 0, \mb x),\\
	&\bar \nabla^2 L^{\Exp, \hat \theta_n, \{T_{\theta_0 \rightarrow \hat \theta_n} e_i\}}(\mb 0, \mb x)  = \hat{\mb Q}\bar \nabla^2 L^{\ret, \hat \theta_n}(\mb 0, \mb x)\hat{\mb Q}\t.\label{eq: equivalence between retraction and tensor operator on gradient and Hessian}
	}

Then we make the following claim to establish the consistency of studentized terms across different charts whose proof is postponed to Section~\ref{subsubsection: proof of claim 1}. 
\begin{claim}\label{claim: claim 1}
	Under the event $\mc F_n$, it holds that 
	{\allowdisplaybreaks	
	\begin{align}
		& \Big \lVert\bar \nabla^2 L_n^{\Exp, \theta_0, \{e_i\}}(\mb H(\bar{\mb Z}_n )) - \hat{\mb Q}\bar \nabla^2 L_n^{\ret, \hat \theta_n}\hat{\mb Q}\t\Big\rVert_F = O(n^{-1}\log n),\label{eq: claim 1 eq. 1} \\ 
        &		 \Big\lVert \sum_{i\in[n]} \bar \nabla L^{\Exp, \theta_0, \{e_i\}}(\mb H(\bar{\mb Z}_n), X_i)  \bar \nabla L^{\Exp, \theta_0, \{e_i\}}(\mb H(\bar{\mb Z}_n), X_i)\t /n\\ - & \hat{\mb Q}\sum_{i=1}^n \bar \nabla L^{\ret, \hat \theta_n}(\mb 0, X_i)\bar \nabla L^{\ret, \hat \theta_n}(\mb 0, X_i)\t \hat{\mb Q}\t /n \Big\Vert_F 
        =O(n^{-1}\log n). \label{eq: claim 1 eq. 3}
        \end{align}
        Similarly, under the event $\mc F_n^*$ it has
        \begin{align}
		& \Big \lVert
		\bar \nabla^2 {L_n^*}^{\Exp, \theta_0, \{e_i\}}(\mb H(\bar{\mb Z}_n ) + \mb H(\bar{\mb Y}_n ))
 - \hat{\mb Q}^*\bar \nabla^2 {L_n^*}^{\ret, \hat \theta_n^*}{{}\hat{\mb Q}^*}\t
		\Big\rVert_F = O(n^{-1}\log n), \label{eq: claim 1 eq. 2}\\ 
		& \Big\lVert \sum_{i\in[n]} \bar \nabla L^{\Exp, \theta_0, \{e_i\}}(\mb H(\bar{\mb Z}_n ) + \mb H(\bar{\mb Y}_n), X_i^*)  \bar \nabla L^{\Exp, \theta_0, \{e_i\}}(\mb H(\bar{\mb Z}_n ) + \mb H(\bar{\mb Y}_n), X_i^*)\t /n \\- & \hat{\mb Q}^*\sum_{i=1}^n \bar \nabla L^{\ret, \hat \theta^*_n}(\mb 0, X_i^*)\bar \nabla L^{\ret, \hat \theta^*_n}(\mb 0, X_i^*)\t {{}\hat{\mb Q}^*}\t /n \Big\Vert_F = O(n^{-1}\log n). \label{eq: claim 1 eq. 4}
	\end{align}
	}
\end{claim}

As a consequence of the above claim, we have proved that, given a sufficiently large $n$, the involved matrix inverse operations in the definitions of $\hat{\mb \Sigma}$ and $\check{\mb \Sigma}$ are always feasible, under $\mc F_n$ and $\mc F_n^*\cap \mc F_n$, respectively, in light of \eqref{eq: claim 1 eq. 1} and \eqref{eq: claim 1 eq. 2}. 
Moreover, by the definitions of $\hat{\mb \Omega}$ and $\hat{\mb \Sigma}_n$, Claim~\ref{claim: claim 1} immediately implies that 
\longeq{
& \norm{\hat{\mb Q}\hat{\mb \Sigma}\hat{\mb Q}\t - \hat {\mb \Omega} }= O(n^{-1}\log n), \\ 
& \norm{\hat{\mb Q}^*\check{\mb \Sigma}{{}\hat{\mb Q}^*}\t  - \hat {\mb \Omega}^*} = O(n^{-1}\log n). \label{eq: hat Omega error}
}

As a consequence of \eqref{eq: approximation error of retraction mapping with rotations} and \eqref{eq: hat Omega error}, we have 
\longeq{
& \Big\vert R^{-1}_{\hat \theta_n}(\theta_0)\t \hat {\mb \Sigma}^{-1} R^{-1}_{\hat \theta_n}(\theta_0)\t -{\pi_{\{e_i\}} \circ \Log_{\theta_0}(\hat\theta_n)}\t \hat{\mb \Omega}^{-1} \pi_{\{e_i\}} \circ \Log_{\theta_0}(\hat\theta_n)\Big\vert =  O(n^{-2}(\log n)^{2}),  \\ 
& \hspace{12cm} \text{under $\mc F_n$,}\\ 
& \Big \vert R^{-1}_{\hat \theta^*_n}(\hat{\theta}_n )\t \check {\mb \Sigma}^{-1} R^{-1}_{\hat \theta_n^*}(\hat\theta_n)\t -\big({\pi_{\{e_i\}} \circ \Log_{\theta_0}(\hat\theta_n^*)}   - {\pi_{\{e_i\}} \circ \Log_{\theta_0}(\hat\theta_n)} \big)\t \\ 
& \hspace{4cm} {{}\hat{\mb \Omega}^*}^{-1} \big({\pi_{\{e_i\}} \circ \Log_{\theta_0}(\hat\theta_n^*)} - {\pi_{\{e_i\}} \circ \Log_{\theta_0}(\hat\theta_n)} \big)\Big \vert=O(n^{-2}(\log n)^{2}) \\ 
&  \hspace{12cm} \text{under $\mc F_n^*$. }
}

Here we are going to utilize the smooth implicit function approximation obtained in Section~\ref{section: implicit function 1} and \ref{section: implicit function 2}. Substitution of $H(\bar{\mb Z}_n)$ and $H(\bar{\mb Z}^*_n)$ into $ \pi_{\{e_i\}} \circ \Log_{\theta_0}(\hat\theta_n)$ and $\pi_{\{e_i\}} \circ \Log_{\theta_0}(\hat\theta_n^*)$ respectively together with \eqref{eq: approximation error of H(Z)} and \eqref{eq: approximation error of H(Y) 2}  gives that 
\longeq{\label{eq: quadratic approximation error}
& \big\vert R^{-1}_{\hat \theta_n}(\theta_0)\t \hat {\mb \Sigma}^{-1} R^{-1}_{\hat \theta_n}(\theta_0) -\mb H(\bar{\mb Z}_n)\t \hat{\mb \Omega}\mb  H(\bar{\mb Z}_n)\big\vert=O(n^{-2}(\log n)^{2})  \\ 
& \big \vert R^{-1}_{\hat \theta^*_n}(\hat\theta_n)\t \check {\mb \Sigma}^{-1} R^{-1}_{\hat \theta_n^*}(\hat\theta_n ) -\mb H(\bar{\mb Y}_n)\t  \hat{\mb \Omega}^* \mb H(\bar{\mb Y}_n)\big \vert = O(n^{-2}(\log n)^{2}) . 
}

Now we are left with establishing the Edgeworth expansion of $\mb H(\bar{\mb Z}_n)\t \hat{\mb \Omega}^{-1}\mb  H(\bar{\mb Z}_n)$ as well as $\mb H(\bar{\mb Y}_n)\t  {{}\hat{\mb \Omega}^*}^{-1} \mb H(\bar{\mb Y}_n)$. To see the distributional consistency of quadratic terms, we first establish the consistency of the multivariate studentized terms, namely, $\hat{\mb \Omega}^{-\frac{1}{2}}\mb H(\bar{\mb Z}_n)$ and ${{}\hat{\mb \Omega}^*}^{-\frac{1}{2}} \mb H(\bar{\mb Y}_n )$.  We first note that under the event $\mc F_n$ with every sufficiently large $n$, 
\begin{align}
    & \hat{\mb \Omega}^{-\frac12} = \mb \Omega^{-\frac12}\big(\mb I + (\hat{\mb \Omega}^{\frac12} - \mb \Omega^{\frac12})\mb \Omega^{-\frac12} \big)^{-1} = \mb \Omega^{-\frac12} \sum_{k = 0}^\infty (-1)^k\big(  (\hat{\mb \Omega}^{\frac12} - \mb \Omega^{\frac12}) \mb \Omega^{-\frac12}\big)^k \\ 
    & \qquad  =\mb  \Omega^{-\frac12} - \mb \Omega^{-\frac12}(\hat{\mb \Omega}^{\frac12} - \mb \Omega^{\frac12}) \mb \Omega^{-\frac12} + O(n^{-1}\log n), \label{eq: expansion of Omega}
\end{align}
where the validity of the expansion is ensured by the conditions in $\mc F_n$. 

Invoking Proposition~\ref{proposition: thm 3.1 in bhattacharya} first yields that
\begin{align}
    & \sup_{B \in \mc B}\Big| \bb P\big[\sqrt{n}\big(\mb  \Omega^{-\frac12} - \mb \Omega^{-\frac12}(\hat{\mb \Omega}^{\frac12} - \mb \Omega^{\frac12}) \mb \Omega^{-\frac12} \big)\mb H(\bar{\mb Z}_n)  \in B\big] - \int_{B}\psi_{4,n}^{(\mb H)} \mathrm d \mb x \Big| = o(n^{-1})
\end{align}
holds for a class $\mathcal B$ of Borel set satisfying, for some $a>0$,  
\longeq{
& \Phi_{\mb \Sigma}\big((\partial B)^{\epsilon}\big) = O(\epsilon^a),
}
where $\Phi_{\mb \Sigma}$ represents the multivariate Gaussian measure with covariance $\mb \Sigma$. Here $\psi_{4,n}^{(\mb H)}$ denotes the formal Edgeworth expansion of $ \big( \mb  \Omega^{-\frac12} - \mb \Omega^{-\frac12}(\hat{\mb \Omega}^{\frac12} - \mb \Omega^{\frac12}) \mb \Omega^{-\frac12}\big)\mb H(\bar{\mb Z}_n)$ as detailed in Section~\ref{sec: classic edgeworth expansion}. 

Then applying Lemma~\ref{lemma:deltamethod} yields that 
\begin{align}
    & \sup_{B \in \mc B}\Big| \bb P\big[\sqrt{n}\hat{\mb \Omega}^{-\frac12}\mb H(\bar{\mb Z}_n)  \in B\big] - \int_{B}\psi_{4,n}^{(\mb H)} \mathrm d \mb x \Big| = O(n^{-1 + \xi})
\end{align}
holds for every $0 < \xi \leq \frac{1}{2}$, where we make use of \eqref{eq: expansion of Omega} and \eqref{eq: probability of F_n}.

On the other hand, for its resample counterpart $\mb H(\bar{\mb Y}_n)$, a similar derivation to \eqref{eq: expansion of Omega} gives that 
\begin{align}
    & {{}\hat{\mb \Omega}^*}^{-\frac{1}{2}} = \hat{\mb \Omega}^{-\frac12} - \hat{\mb \Omega}^{-\frac12}  \big({{}\hat{\mb \Omega}^*}^{-\frac{1}{2}} - \hat{\mb \Omega}^{-\frac12} \big) \hat{\mb \Omega}^{-\frac12}  + O(n^{-1}\log n) \label{eq: hat Omega * expansion}
\end{align}
under the event $\mc F_n^*$. 

Applying Proposition~\ref{proposition: modified thm 3.3 in bhattacharya}, and Proposition~\ref{proposition: equivalence of two resampled distribution} to $\big(\hat{\mb \Omega}^{-\frac12} - \hat{\mb \Omega}^{-\frac12}  \big({{}\hat{\mb \Omega}^*}^{-\frac{1}{2}} - \hat{\mb \Omega}^{-\frac12} \big) \hat{\mb \Omega}^{-\frac12}\big)\mb  H(\bar{\mb Y}_n)$, we arrive at
\begin{equation}
     \limsup\limits_{n\rightarrow \infty} n\bb P_{\mc X_n}\bigg[\sup_{B \in \mathcal B} \Big\vert\mathbb P\big[ 
	 \big(\hat{\mb \Omega}^{-\frac12} - \hat{\mb \Omega}^{-\frac12}  \big({{}\hat{\mb \Omega}^*}^{-\frac{1}{2}} - \hat{\mb \Omega}^{-\frac12} \big) \hat{\mb \Omega}^{-\frac12}\big)\mb  H(\bar{\mb Y}_n)\in B
\vert \mathcal X_n\big]  - \int_B\psi_{4,n}^{*(\mb H)}\mathrm d \mb x  \Big \vert > n^{-1 + \xi} \bigg] < \infty\label{eq: Edgeworth expansion of hat Omega* 1}
\end{equation}
for every $0< \xi\leq \frac12$, where we invoke Lemma~\ref{lemma:deltamethod} together with \eqref{eq: exceptional probability of F_n'} and $\psi_{4,n}^{*(\mb H)}$ denotes the empirical
Edgeworth expansion of the random vector $\big(\hat{\mb \Omega}^{-\frac12} - \hat{\mb \Omega}^{-\frac12}  \big({{}\hat{\mb \Omega}^*}^{-\frac{1}{2}} - \hat{\mb \Omega}^{-\frac12} \big) \hat{\mb \Omega}^{-\frac12}\big)\mb  H(\bar{\mb Y}_n)$ conditional on $\mc X_n$. 

Combining \eqref{eq: Edgeworth expansion of hat Omega* 1} with \eqref{eq: hat Omega * expansion}, Lemma~\ref{lemma:deltamethod} yields that 
\begin{align}
    & \limsup\limits_{n\rightarrow \infty} n\bb P_{\mc X_n}\bigg[\sup_{B \in \mathcal B} \Big\vert\mathbb P\big[  {{}\hat{\mb \Omega}^*}^{-\frac12}\mb  H(\bar{\mb Y}_n)\in B
\vert \mathcal X_n\big]  - \int_B\psi_{4,n}^{*(\mb H)}\mathrm d \mb x  \Big \vert > n^{-1 + \xi} \bigg] < \infty. 
\end{align}

Then, as stated in Remark 2.2, \cite{chandra1979valid}, given a random vector $\mb X \in \bb R^d$, the validity of Edgeworth expansion for $\norm{\mb X}_2^2$ can be demonstrated by \cite[Theorem 1]{chandra1979valid} once we derive a valid Edgeworth expansion for the distribution of $\mb X$. To be specific, one has
\begin{align}
 & \sup_{t\in \bb R^+}\left|\bb P\Big[  n\mb  H(\bar{\mb Z}_n)\t \hat{\mb\Omega}^{-1}\mb H(\bar{\mb Z}_n) \leq t \Big] - \int_{x \leq t}\psi^{(\mb H^2)}_{4,n}(x)\mathrm d x\right| \leq  O(n^{-1 + \xi}),  \\ 
&\limsup\limits_{n\rightarrow \infty} n \bb P_{\mc X_n}\bigg[\sup_{t\in \bb R^+}\left|\bb P\left[n \mb H(\bar{\mb Y}_n ) \t{{}\hat{\mb\Omega}^*}^{-1}\mb H(\bar{\mb Y}_n)  \leq  t  |\mc X_n \right]  -  \int_{x\leq t} \psi^{*(\mb H^2)}_{4,n}( x)\mathrm d x\right| > n^{-1+ \xi}\bigg] < \infty.
\end{align}
Here the function $\psi^{(\mb H^2)}_{4,n}( x)$ is given by $f_{\chi^2_p}(x)\cdot\big(1 + n^{-1}\pi_1^{(\mb H^2)}(x)\big)$, where $f_{\chi^2_p}$ denotes the density function of a chi-squared distribution of freedom $p$ and $\pi_1^{(H^2)}(x)$ is a polynomial whose coefficients are determined by those of $\psi_{4,n}^{(\mb H)}(\mb x)$. Similarly, the function $\psi^{*(\mb H^2)}_{4,n}(x)$ obeys the same form,  with its coefficients determined by the ones of $\psi_{4,n}^{*(\mb H)}(\mb x)$. 

Leveraging \eqref{eq: quadratic approximation error}, we apply the delta method (Lemma~\ref{lemma:deltamethod}) to the above approximations and derive that 
\begin{align}
    \label{eq: edgeworth expansion of quadratic terms}
 & \sup_{t \in \bb R^+}\left|\bb P\Big[  R^{-1}_{\hat{\theta}_n}(\theta_0)\t {{}\hat{\mb \Sigma}}^{\dagger}R^{-1}_{\hat{\theta}_n}(\theta_0)  \leq t \Big] - \int_{x\leq t}\psi^{(\mb H^2)}_{4,n}(x)\mathrm d x\right| \leq  O(n^{-1 + \xi}),  \\ 
&\limsup\limits_{n\rightarrow \infty} n \bb P_{\mc X_n}\bigg[\sup_{t\in \bb R^+}\left|\bb P\left[n R^{-1}_{\hat{\theta}_n^*}(\hat\theta_n)\t {{}\check{\mb \Sigma}}^{\dagger}R^{-1}_{\hat{\theta}_n^*}(\hat \theta_n)  \leq  t  |\mc X_n \right]  -  \int_{x\leq t} \psi^{*(\mb H^2)}_{4,n}( x)\mathrm d x\right| > n^{-1+ \xi}\bigg] < \infty.\label{eq: empirical edgeworth expansion of quadratic terms}
\end{align}

Finally, the desired result follows by combining \eqref{eq: empirical edgeworth expansion of quadratic terms} with \eqref{eq: edgeworth expansion of quadratic terms}: 
\eq{\label{eq: wald statistic distributional consistency}
\limsup\limits_{n\rightarrow \infty} n\bb P_{\mc X_n}\Big[\sup_{x\in \bb R^+} \Big| \bb P\left[R^{-1}_{\hat{\theta}_n^*}(\hat\theta_n)\t {{}\check{\mb \Sigma}}^{\dagger}R^{-1}_{\hat{\theta}_n^*}(\hat \theta_n) \leq x\mid \mc X_n\right]  - \bb  P \left[R^{-1}_{\hat{\theta}_n}(\theta_0)\t {{}\hat{\mb \Sigma}}^{\dagger}R^{-1}_{\hat{\theta}_n}(\theta_0) \leq x \right]\Big|  > n^{-1  +\xi} \Big] < +\infty, 
}  
where we make use of the fact that the coefficients of the $n^{-1}$ term in the polynomial of $\psi^{*(\mb H^2)}_{4,n}( x)$ converge to those of $\psi^{(\mb H^2)}_{4,n}(x)$ at a rate of $O(n^{-\xi})$ for every $0< \xi < \frac12$, since:
\begin{enumerate}
    \item The coefficients are determined by the first four moments of $\mb Z_1$, and their empirical counterpart, respectively, as well as the derivatives of the functions $\big( \mb  \Omega^{-\frac12} - \mb \Omega^{-\frac12}(\hat{\mb \Omega}^{\frac12} - \mb \Omega^{\frac12}) \mb \Omega^{-\frac12}\big)\mb H(\bar{\mb Z}_n)$ and $\big(\hat{\mb \Omega}^{-\frac12} - \hat{\mb \Omega}^{-\frac12}  \big({{}\hat{\mb \Omega}^*}^{-\frac{1}{2}} - \hat{\mb \Omega}^{-\frac12} \big) \hat{\mb \Omega}^{-\frac12}\big)\mb  H(\bar{\mb Y}_n)$, respectively. 
    \item By the finite moment condition Assumption~\ref{assumption:manifold}.\ref{item:B}, Lemma~\ref{lemma: thm1 in von1967central} yields that the empirical moments converge to the population ones at a rate of $O(n^{-\xi})$ with probability $1 - O(n^{-1})$, for every $0< \xi < \frac12$. 
    \item By the convergence rate of $\bar{\mb Z}_n$ in \eqref{eq:finitemomentconcentration} and Proposition~\ref{proposition: equivalence of two resampled distribution}, we can prove that the involved derivatives also converge to their counterparts at a rate of $O(n^{-\xi})$ with probability $1 - O(n^{-1})$, for every $0< \xi < \frac12$. 
    \item The function $f_{\mc X_p^2}(x)$, when multiplied by any polynomial is always integrable. 
\end{enumerate}

\paragraph*{Proof of \eqref{eq: thm1.2} ($t$-Statistic) }
In contrast to the Wald-type statistic, where orthogonal matrices corresponding to basis transformations are canceled out, analyzing the intrinsic $t$-type statistic is more complicated because of the existence of basis transformation. Without loss of generality, we let $\mb a$ be $\delta_i$ for an arbitrary $i\in[p]$. 
In what follows, we claim the existence of a precise enough approximation, provided the differentiability of the retraction:  
\begin{claim}\label{claim: claim 2}
 Under the event $\mc F_n^*$, for each $i\in[p]$, there exists a function $\tilde H_i$ such that 
	\begin{align}
	\label{eq: expansion in claim 2.1}
	& \Big \vert \frac{\tilde H_i(\bar {\mb Z}_n, \bb E[\mb Z_1])}{(\mb \Sigma)_{i,i}^{\frac{1}{2}}} - \frac{R_{\hat \theta_n}(\theta_0)_i }{(\hat {\mb \Sigma})_{i,i}^{\frac{1}{2}}}\Big\vert = O(n^{-\frac{3}{2}}(\log n)^{\frac{3}{2}}),\\ 
	& \Big \vert \frac{\tilde H_i(\bar{\mb Y}_n, \bar{\mb Z}_n)}{(\hat {\mb \Sigma})_{i,i}^{\frac{1}{2}}} - \frac{R_{\hat \theta_n^*}(\hat\theta_n)_i }{(\check {\mb \Sigma})_{i,i}^{\frac{1}{2}}}\Big\vert = O(n^{-\frac{3}{2}}(\log n)^{\frac{3}{2}}). \label{eq: expansion in claim 2.2}
	\end{align}
    Moreover, the function $\tilde H_i$ satisfies that 
\begin{enumerate}
\item $\tilde H_i(\cdot, \mb x')$ is three-times continuously differentiable on $O$ for every $\mb  x' \in O$ where $O$ is an open neighborhood of $\bb E[\mb Z_1]$. 
\item The first three-times derivatives of $\tilde H_i(\mb x_1, \mb x_2)$ with respect to the first variable $\mb x_1 \in O$ are continuously differntiable with respect to $\mb x_2\in O$. 
\end{enumerate}
\end{claim}

The lengthy proof of Claim~\ref{claim: claim 2} is postponed to Section~\ref{subsubsection: studentization analysis}. 
For the function $\tilde H_i$, it follows from Proposition~\ref{proposition: thm 3.1 in bhattacharya} that 
\longeq{
\label{eq: edgeworth expansion of g_i}
&  \sup_{B \in \mathcal{B}}\left|\bb P\left[n^{\frac{1}{2}}\tilde H_i(\bar{\mb Z}_n, \bb E[\mb Z_1]) \in B\right] - \int_B\Psi^{(\tilde H_i)}_{4,n}(v)\mathrm dv \right| = o(n^{-1}),
 }
holds for all class $\mc B \in \mc B^p$ (the collection of Borel sets in $\bb R^p$) such that, for some $a > 0$, 
    \eq{
    \sup_{B \in \mc B}  \Phi_{\mb \Sigma}((\partial B)^\epsilon) = O(\epsilon^a), \quad \text{ as $\epsilon \downarrow 0$. }
    }
Here the explicit form of $\Psi_{4,n}^{(\tilde H_i)}(x)$ is introduced in Remark~\ref{remark: coefficients in Edgeworth expansions}. 
 
Combining the fact that $\bb P(\mc F_n^\complement) \leq C n^{-1}(\log n)^{-2}$, \eqref{eq: expansion in claim 2.1}, and Lemma~\ref{lemma:deltamethod} with \eqref{eq: edgeworth expansion of g_i} yields that
\eq{\label{eq: population edgeworth expansion of intrinsic t}
    \sup_{B \in \mathcal{B}}\bigg|\bb P_{\mc X_n}\Big[n^{\frac{1}{2}}\frac{R^{-1}_{\hat \theta_n}(\theta_0)_i}{(\hat{\mb \Sigma})_{i,i}^{\frac{1}{2}}}\in B\Big] - \int_B\Psi^{(\tilde H_i)}_{4,n}(v)\mathrm dv \bigg| = O(n^{-1 + \xi})
}
for every $0< \xi \leq \frac12$.

Similarly, Proposition~\ref{proposition: modified thm 3.3 in bhattacharya} and Proposition~\ref{proposition: equivalence of two resampled distribution} implies that 
\longeq{\label{eq: empirical edgeworth expansion of g_i}
& \limsup\limits_{n\rightarrow \infty} n \bb P_{\mc X_n}\Big[\sup_{B \in \mathcal{B}}\left|\bb P\big[n^{\frac{1}{2}}\tilde H_i(\bar{\mb Y}_n, \bar{\mb Z}_n) \in B\vert \mathcal X\right] - \int_B\Psi^{*(\tilde H_i)}_{4,n}(v)\mathrm dv \big| > n^{-1}\Big] < \infty 
}
for every class $\mathcal B$ of Borel set satisfying, for some $a>0$,  
\longeq{
& \Phi_{\mb \Sigma}\big((\partial B)^{\epsilon}\big) = O(\epsilon^a),
}
where $\Psi^{*(\tilde H_i)}_{4,n}$ is an empirical analog to $\Psi^{(\tilde H_i)}_{4,n}$.

Again, invoking \eqref{eq: expansion in claim 2.1} and the delta method (Lemma~\ref{lemma:deltamethod}), one has 
\eq{\label{eq: empirical edgeworth expansion of intrinsic t}
    \limsup\limits_{n\rightarrow \infty} n\bb P_{\mc X_n} \Big[\sup_{B \in \mc B'} \bigg|\bb P\Big[n^{\frac{1}{2}}\frac{R^{-1}_{\hat \theta_n^*}(\hat \theta_n)_i }{(\check{\mb \Sigma})_{i,i}^{\frac{1}{2}}} \in B|\mc X_n\Big] - \int_B \Psi^{*(\tilde H_i)}_{4,n}(v)\mathrm d v  \bigg|  > n^{-1 +\xi } \Big] < \infty
}
for every $0 < \xi \leq \frac12$.
Here, the coefficients of $n^{-\frac{1}{2}}$ in $\Psi^{(\tilde H_i)}_{4,n}$ and $\Psi^{*(\tilde H_i)}_{4,n}$ coincides and the coefficient of $n^{-1}$ in $\Psi^{*(\tilde H_i)}_{4,n}$ converges to that in $\Psi^{(\tilde H_i)}_{4,n}$ at a rate of $O(n^{-\xi})$ with probability $1 - O(n^{-1})$ for $0< \xi< \frac12$, following similar arguments to those below \eqref{eq: wald statistic distributional consistency}.

Consequently, combining \eqref{eq: population edgeworth expansion of intrinsic t} with \eqref{eq: empirical edgeworth expansion of g_i} gives that 
\longeq{
	& \limsup\limits_{n\rightarrow \infty} n \bb P_{\mc X_n}\Big[ \sup_{x\in \mathbb R}\Big\vert \bb P_{\mc X_n}\Big[n^{\frac{1}{2}}\frac{R^{-1}_{\hat \theta_n}(\theta_0)_i}{(\hat{\mb \Sigma})_{i,i}^{\frac{1}{2}}} \leq x\Big] -   \bb P\Big[n^{\frac{1}{2}}\frac{R^{-1}_{\hat \theta_n^*}(\hat \theta_n)_i }{(\check{\mb \Sigma})_{i,i}^{\frac{1}{2}}} \leq x|\mc X_n\Big] \Big\vert > n^{-1 + \xi}\Big] <\infty
}
for every $0 < \xi \leq \frac12$. \qed
\subsubsection{Proof of Corollary~\ref{corollary: coverage}}
\label{subsubsection: proof of corollary: coverage}
The coverage probability is a consequence of the uniform controls in the proof of Theorem~\ref{theorem:bootstrapmanifold}. For clarity, we provide the proof for the Wald statistic for instance.

We define $w_{\alpha}^*$ as the exact solution $t$ to $\int_{x\leq t}\psi^{*(\mb H^2)}_{4,n}(x)\mathrm d x = 1- \alpha$, whose existence is ensured by the continuity of the function. We then recall the definition of the approximate quantile in Algorithm~\ref{algorithm:wald based statistic} with $b = \infty$, expressed as: 
\eq{
w_{\alpha,\infty}  = \argmin_{t\in \bb R^+}\big|\bb P\big[W^* \leq t \big|\mc X_n\big] - (1 - \alpha) \big| \wedge \big|\bb P\big[W^* < t \big|\mc X_n\big] - (1 - \alpha) \big|.
}
From this, we deduce the following:
\begin{align}
    & \Big|\bb P_{\mc X_n}\big[R^{-1}_{\hat{\theta}_n}(\theta_0)\t {{}\hat{\mb \Sigma}}^{\dagger}R^{-1}_{\hat{\theta}_n}(\theta_0) \leq w_{\alpha, \infty} \big] - (1 - \alpha) \Big| \\ 
    \leq & \Big|\bb P_{\mc X_n}\big[R^{-1}_{\hat{\theta}_n}(\theta_0)\t {{}\hat{\mb \Sigma}}^{\dagger}R^{-1}_{\hat{\theta}_n}(\theta_0) \leq w_{\alpha, \infty}\big] - \bb P_{\mc X_n}\Big[\bb P\big[R^{-1}_{\hat{\theta}_n^*}(\hat\theta_n)\t {{}\check{\mb \Sigma}}^{\dagger}R^{-1}_{\hat{\theta}_n^*}(\hat \theta_n) \leq w_{\alpha, \infty} \mid \mc X_n\big]\Big]\Big|
    \\ 
    & + \bigg[ \Big|\bb P_{\mc X_n}\Big[\bb P\big[R^{-1}_{\hat{\theta}_n^*}(\hat\theta_n)\t {{}\check{\mb \Sigma}}^{\dagger}R^{-1}_{\hat{\theta}_n^*}(\hat \theta_n) \leq w_{\alpha, \infty} \mid \mc X_n\big]\Big]- (1 - \alpha) \Big| \bigg] \\ 
    \wedge &  \bigg[  \Big| \bb P_{\mc X_n}\Big[\bb P\big[R^{-1}_{\hat{\theta}_n^*}(\hat\theta_n)\t {{}\check{\mb \Sigma}}^{\dagger}R^{-1}_{\hat{\theta}_n^*}(\hat \theta_n) \leq  w_{\alpha, \infty} \mid \mc X_n\big]\Big] \\ 
    & - \bb P_{\mc X_n}\Big[\bb P\big[R^{-1}_{\hat{\theta}_n^*}(\hat\theta_n)\t {{}\check{\mb \Sigma}}^{\dagger}R^{-1}_{\hat{\theta}_n^*}(\hat \theta_n) <  w_{\alpha, \infty} \mid \mc X_n\big]\Big]\Big|\\ 
    & + \Big|\bb P_{\mc X_n}\Big[\bb P\big[R^{-1}_{\hat{\theta}_n^*}(\hat\theta_n)\t {{}\check{\mb \Sigma}}^{\dagger}R^{-1}_{\hat{\theta}_n^*}(\hat \theta_n) <  w_{\alpha, \infty} \mid \mc X_n\big]\Big] - (1- \alpha ) \Big|\bigg] 
    \\ 
    \leq &  \Big|\bb P_{\mc X_n}\big[R^{-1}_{\hat{\theta}_n}(\theta_0)\t {{}\hat{\mb \Sigma}}^{\dagger}R^{-1}_{\hat{\theta}_n}(\theta_0) \leq w_{\alpha, \infty}\big] - \bb P_{\mc X_n}\Big[\bb P\big[R^{-1}_{\hat{\theta}_n^*}(\hat\theta_n)\t {{}\check{\mb \Sigma}}^{\dagger}R^{-1}_{\hat{\theta}_n^*}(\hat \theta_n) \leq w_{\alpha, \infty} \mid \mc X_n\big]\Big]\Big| \\ 
    & + \Big|\bb P_{\mc X_n}\Big[\bb P\big[R^{-1}_{\hat{\theta}_n^*}(\hat\theta_n)\t {{}\check{\mb \Sigma}}^{\dagger}R^{-1}_{\hat{\theta}_n^*}(\hat \theta_n) \leq w_\alpha^* \mid \mc X_n\big]\Big]- (1 - \alpha) \Big| + O(n^{-1+ \xi}) \\ 
    \leq &O(n^{-1+ \xi}),
\end{align}
for every $0< \xi < \frac12 $. The penultimate inequality follows from the definition of $w_{\alpha, \infty}$ along with \eqref{eq: edgeworth expansion of quadratic terms} and \eqref{eq: empirical edgeworth expansion of quadratic terms} to derive the penultimate inequality, while the last inequality applies Theorem~\ref{theorem:bootstrapmanifold} to upper bound the first term and uses \eqref{eq: empirical edgeworth expansion of quadratic terms} to control the second term. 

To establish coverage consistency for the practical threshold $w_\alpha$
 used in the algorithm, it suffices to show that $\limsup_{b \rightarrow \infty} |w_{\alpha} - w_{\alpha, \infty} | \leq C n^{-1 } $ under an event with probability exceeding $O(n^{-1})$. To emphasize dependence on $b$, we write it as $w_{\alpha,b}$ in what follows. 
Invoking \eqref{eq: empirical edgeworth expansion of quadratic terms}, one has from the Edgeworth expansion that 
\begin{align}
& \bb P\Big[ \big|\bb P[W^* \leq w_{\alpha, \infty} - Cn^{-1} \big|\mc X_n\big] - \lim_{t\uparrow w_{\alpha,\infty}} \bb P\big[W^* \leq t \big|\mc X_n] \big| \\ 
& \qquad \vee \big|\bb P[W^* \leq w_{\alpha, \infty} + Cn^{-1} \big|\mc X_n\big] - \lim_{t\downarrow w_{\alpha,\infty}} \bb P\big[W^* \leq t \big|\mc X_n] \big| \geq n^{-1}\Big] \geq 1- O(\frac{1}{n})
\label{eq: perturbed w alpha infty}
\end{align}
holds for some positive $C$. Conditioning on $\mc X_n$ and applying the Bernstein inequality, one has with probability at least $1- O(n^{-1})$ that
\begin{align}
    & \Big|\bb P[W^* \leq w_{\alpha, \infty} - Cn^{-1} \big|\mc X_n\big] - \sum_{i\in[b]} 1\{W^{*[i]} \leq w_{\alpha, \infty} - Cn^{-1}\}/ b\Big| \leq \sqrt{\frac{\log n}{b}}, \\ 
    & \Big|\bb P[W^* < w_{\alpha, \infty} \big|\mc X_n\big] - \sum_{i\in[b]} 1\{W^{*[i]} < w_{\alpha, \infty}\}/ b\Big| \leq \sqrt{\frac{\log n}{b}},  \\ 
    & \Big|\bb P[W^* \leq w_{\alpha, \infty} \big|\mc X_n\big] - \sum_{i\in[b]} 1\{W^{*[i]} \leq w_{\alpha, \infty}\}/ b\Big| \leq \sqrt{\frac{\log n}{b}}, \\
    & \Big|\bb P[W^* \leq w_{\alpha, \infty} + Cn^{-1} \big|\mc X_n\big] - \sum_{i\in[b]} 1\{W^{*[i]} \leq w_{\alpha, \infty} + Cn^{-1}\}/ b\Big| \leq \sqrt{\frac{\log n}{b}}. 
\end{align}
 Choosing $b \geq C' (\log n) n^2$ with some sufficiently large $C'$, one has 
$
    \sqrt{\frac{\log n}{b}} \leq \frac{1}{4}n^{-1}
$,
which leads to 
\begin{align}
    &  |w_{\alpha, b} - w_{\alpha} | \leq C n^{-1}
    \label{eq: convergence of w_alpha b}
\end{align}
for all $b \geq C' (\log n) n^2$ under the event $\Big\{  \big|\bb P[W^* \leq w_{\alpha, \infty} - Cn^{-1} \big|\mc X_n\big] - \lim_{t\uparrow w_{\alpha,\infty}} \bb P\big[W^* \leq t \big|\mc X_n] \big| \vee \big|\bb P[W^* \leq w_{\alpha, \infty} + Cn^{-1} \big|\mc X_n\big] - \lim_{t\downarrow w_{\alpha,\infty}} \bb P\big[W^* \leq t \big|\mc X_n] \big| \geq n^{-1}\Big\}$. Finally, combining \eqref{eq: convergence of w_alpha b} with the empirical Edgeworth expansion
\eqref{eq: empirical edgeworth expansion of quadratic terms} completes the proof for the Wald statistic.

The coverage guarantees for the intrinsic $t$-statistic and the extrinsic $t$-statistic can be derived using analogous arguments, which we omit for brevity.

\subsubsection{Proof of Claim~\ref{claim: claim 1}} \label{subsubsection: proof of claim 1}
The basic idea is to pull the gradient and the Hessian back to the coordinate chart induced by $\pi_{\{e_i\}}\circ \Log_{\theta_0}$ by controlling the accumulated error in parallel transports. With regard to the gradient term $\bar \nabla L^{\Exp, \hat \theta_n, \{T_{\theta_0\rightarrow \hat \theta_n } e_i\}}(X_i)$, it follows by Lemma~\ref{lemma:parallel transport error} and the definition of $\mc F_n'$ that
{
\allowdisplaybreaks
	\begin{align}
		& \Big\lVert \sum_{i\in[n]} \bar \nabla L^{\Exp, \hat \theta_n, \{T_{\theta_0\rightarrow \hat \theta_n } e_i\}}(X_i)\bar \nabla L^{\Exp, \hat \theta_n, \{T_{\theta_0\rightarrow \hat \theta_n } e_i\}}(X_i)\t /n   \\
		-& \sum_{i\in[n]}  \bar \nabla L^{\Exp, \theta_0, \{e_i\}}(\pi_{\{e_i\}}\circ \Log_{\theta_0}(\hat \theta_n), X_i)\bar \nabla L^{\Exp, \theta_0, \{e_i\}}(\pi_{\{e_i\}}\circ \Log_{\theta_0}(\hat \theta_n), X_i)\t  / n\Big\Vert_F  \\ 
		\leq & 2\sqrt{p} \cdot \max_{j \in [p]}\norm{\mathrm d\Exp_{\theta_0}\circ \pi_{\{e_i\}}^{-1} [\delta'_j] - T_{\theta_0 \rightarrow\hat \theta_n}e_j} \\ \cdot &  \Big\lVert \sum_{i\in[n]}  \bar \nabla L^{\Exp, \theta_0, \{e_i\}}(\pi_{\{e_i\}}\circ \Log_{\theta_0}(\hat \theta_n), X_i)\bar \nabla L^{\Exp, \theta_0, \{e_i\}}(\pi_{\{e_i\}}\circ \Log_{\theta_0}(\hat \theta_n), X_i)\t  / n\Big\rVert_F \\ 
        & + p \max_{j \in [p]}\norm{\mathrm d\Exp_{\theta_0}\circ \pi_{\{e_i\}}^{-1} [\delta'_j] - T_{\theta_0 \rightarrow\hat \theta_n}e_j}^2 \\ 
        & \cdot \Big\lVert \sum_{i\in[n]}  \bar \nabla L^{\Exp, \theta_0, \{e_i\}}(\pi_{\{e_i\}}\circ \Log_{\theta_0}(\hat \theta_n), X_i)\bar \nabla L^{\Exp, \theta_0, \{e_i\}}(\pi_{\{e_i\}}\circ \Log_{\theta_0}(\hat \theta_n), X_i)\t  / n\Big\rVert_F\\ 
		= & O(n^{-1}\log n),\label{eq: proof of claim 1 gradient square control}
	\end{align}}
	where $\delta_j'$ denotes the $j$-th canonical tangent vector of $\bb R^p$ at $\pi^{-1}_{\{e_i\}}\circ \Log_{\theta_0}(\hat \theta_n)$ with a slight abuse of notation. Here the last line arises since the term 
    \eq{
    \Big\lVert \sum_{i\in[n]}  \bar \nabla L^{\Exp, \theta_0, \{e_i\}}(\pi_{\{e_i\}}\circ \Log_{\theta_0}(\hat \theta_n), X_i)\bar \nabla L^{\Exp, \theta_0, \{e_i\}}(\pi_{\{e_i\}}\circ \Log_{\theta_0}(\hat \theta_n), X_i)\t  / n\Big\rVert_F 
    } 
    is bounded by a constant as a consequence of \eqref{eq: C(X) bound} and $\max_{j \in [p]}\norm{\mathrm d\Exp_{\theta_0}\circ \pi_{\{e_i\}}^{-1} \delta'_j - T_{\theta_0 \rightarrow\hat \theta_n}e_j}\leq O(\dist(\theta_0, \hat \theta_n)^2)$ by Lemma~\ref{lemma:parallel transport error}	. 
	
	Following similar arguments, we can also prove that the inequality \eqref{eq: claim 1 eq. 4} holds as a result of Lemma~\ref{lemma:parallel transport error} and  \eqref{eq: resampled C(X) bound} under the event $\mc F_n'$ 
    \begin{align}
        & \Big\lVert\sum_{i\in[n]}\bar \nabla L^{\Exp, \hat \theta_n^*, \{T_{\theta_0 \rightarrow \hat \theta_n^*} \} }(X_i^*)\bar \nabla L^{\Exp, \hat \theta_n^*, \{T_{\theta_0 \rightarrow \hat \theta_n^*} \} }(X_i^*)\t / n\\ 
        - &\sum_{i\in[n]}  \bar \nabla L^{\Exp, \theta_0, \{e_i\}}(\pi_{\{e_i\}}\circ \Log_{\theta_0}(\hat \theta_n^*), X_i^*)\bar \nabla L^{\Exp, \theta_0, \{e_i\}}(\pi_{\{e_i\}}\circ \Log_{\theta_0}(\hat \theta_n^*), X_i^*)\t  / n \Big \rVert_F   \\
    \leq &  2\sqrt{p} \cdot \max_{j\in[p]}\norm{\mathrm d \Exp_{\theta_0}\cdot \pi^{-1}_{\{e_i\}} \delta_j''  - T_{\theta_0 \rightarrow \hat \theta_n^*}e_j} \\ 
    \cdot & \Big\lVert \sum_{i\in[n]}  \bar \nabla L^{\Exp, \theta_0, \{e_i\}}(\pi_{\{e_i\}}\circ \Log_{\theta_0}(\hat \theta_n^*), X_i^*)\bar \nabla L^{\Exp, \theta_0, \{e_i\}}(\pi_{\{e_i\}}\circ \Log_{\theta_0}(\hat \theta_n^*), X_i^*)\t  / n\Big\rVert_F \\ 
     & + p \max_{j\in[p]}\norm{\mathrm d \Exp_{\theta_0}\cdot \pi^{-1}_{\{e_i\}} \delta_j''  - T_{\theta_0 \rightarrow \hat \theta_n^*}e_j}^2 \\ 
    & \cdot  \Big\lVert \sum_{i\in[n]}  \bar \nabla L^{\Exp, \theta_0, \{e_i\}}(\pi_{\{e_i\}}\circ \Log_{\theta_0}(\hat \theta_n^*), X_i^*)\bar \nabla L^{\Exp, \theta_0, \{e_i\}}(\pi_{\{e_i\}}\circ \Log_{\theta_0}(\hat \theta_n^*), X_i^*)\t  / n\Big\rVert_F \\ 
    = & O(n^{-1} \log n).
    \end{align}

    Finally, invoking the relations \eqref{eq: approximation error of H(Z)} and \eqref{eq: approximation error of H(Y)} as well as the uniform bound yields the desired results \eqref{eq: claim 1 eq. 3} and \eqref{eq: claim 1 eq. 4}.

	Regarding the Hessian terms, Lemma~\ref{lemma: Hessian consistency} demonstrates that the difference is not only related to the distances among $\theta_0$, $\hat\theta_n$, and $\hat \theta_n^*$ but also to the magnitude of averaged gradients, that is, 
	\begin{align}
		& \Big \lVert\hat{\mb Q}\bar \nabla^2 L_n^{\ret, \hat \theta_n} (\mb 0)\hat{\mb Q}\t - \bar \nabla^2 L_n^{\Exp, \hat \theta_n, \{T_{\theta_0\rightarrow \hat \theta_n}e_i\}} (\pi_{\{e_i\}}\circ \Log_{\theta_0}(\hat \theta_n)) \Big\rVert  \\ 
		\leq & O(n^{-1}\log n) + O(n^{-\frac{1}{2}} (\log n)^{\frac{1}{2}} ) \norm{\nabla L_n(\hat \theta_{n})} \\ 
  = & O(n^{-1}\log n) ,\label{eq: proof of claim 1 error control on the hessian}\\ 
		&\Big \lVert
		\hat{\mb Q}^*\bar \nabla^2 {L_n^*}^{\ret, \hat \theta_n^*}(\mb 0){{}\hat{\mb Q}^*}\t
 - \bar \nabla^2 {{}L_n^*}^{\Exp, \hat \theta_n, \{T_{\theta_0\rightarrow \hat \theta^*_n}e_i\}} ( \pi_{\{e_i\}}\circ \Log_{\theta_0}(\hat \theta_n^*))
		\Big\rVert \\ 
	 \leq & O(n^{-1}\log n) + O(n^{-\frac{1}{2}} (\log n)^{\frac{1}{2}} ) \norm{\nabla L_n^*(\hat \theta_{n}^*)} \\ 
  = & O(n^{-1}\log n) \label{eq: proof of claim 1 error control on the hessian 2},
	\end{align}
 where $\norm{ \nabla L_n(\hat \theta_{n})}$ and $\norm{\nabla L_n^*(\hat \theta_{n}^*)}$ are of the order $o(n^{-\frac{1}{2}}(\log n)^{\frac{1}{2}})$ by Theorem~\ref{thm:convergence rate of newton}, Theorem~\ref{thm: convergence rate of resampled newton}, \eqref{eq: nabla nu difference 1}, and \eqref{eq: nabla nu difference 2}. 

Analogously, replacing $\pi_{\{e_i\}}\circ \Log_{\theta_0}(\hat \theta_n)$ and $\pi_{\{e_i\}}\circ \Log_{\theta_0}(\hat \theta_n^*)$ with $\mb H(\bar{\mb Z}_n)$ and $\mb H(\bar{\mb Y}_n)$, respectively, implies \eqref{eq: claim 1 eq. 1} and \eqref{eq: claim 1 eq. 2}. \qed

\subsubsection{Proof of Claim~\ref{claim: claim 2}}
\paragraph*{Step 1: Basis Transformation Analysis}\label{section: Basis Transformation Analysis}
	The goal of this part is to establish explicit approximations for the rotated coordinates as functions of $\bar{\mb Z}_n$ and $\bar{\mb Y}_n$ respectively.
    Before continuing, we remind that in Section~\ref{section: implicit function 1} and Section~\ref{section: implicit function 2}, we established that $\mb H(\bar{\mb Z}_n)$, $\mb H(\bar{\mb Y}_n)$ serve as nice approximations to $-\pi_{\{T_{\theta_0 \rightarrow \hat \theta_n}e_i\}} \circ \Log_{\hat \theta_n} (\theta_0)$ and $\pi_{\{T_{\theta_0 \rightarrow \hat \theta_n^*}e_i\}} \circ \Log_{\hat \theta_n^*} (\hat \theta_n)$, respectively. The remaining task is to address how to reconcile the bases induced by parallel transports and those arising from the adopted retraction $\{R_\theta\}$. 
    
    The subsequent analysis is a natural exploitation of the differentiability of $\mathrm{d}R_{\theta}^{-1}\big(T_{\theta_0 \rightarrow \theta}(\mathrm{d} R_{\theta_0}\delta_i )\big)$. For ease of notations, we denote the composite function $\mathrm{d}R_{\Exp_{\theta}\circ \pi^{-1}_{\{e_i\}}(\mb x)}^{-1}(T_{\theta_0 \rightarrow \Exp_{\theta}\circ \pi^{-1}_{\{e_i\}}(\mb x)}(e_j))$ by $\mb t_j(\mb x)= (t_{i,k}(\mb x))_{k\in[p]}$ for every $j \in[p]$, which is twice-differentiable as assumed in the second part of Theorem~\ref{theorem:bootstrapmanifold}. 
    
    In the context of this section and the next section, we proceed with our analysis under the event $ \mc F_n^* \cap \mc F_n$ for every $n$, which heuristically means that, $\hat \theta_n, \hat \theta_n^*, \tilde \theta_n$, and $\tilde \theta_n^*$ are always contained in a sequence of shrinking regions as illustrated in Section~\ref{section: implicit function 1} and Section~\ref{section: implicit function 2}, and the deviations of derivatives are appropriately bounded. 
    
    We begin with the Taylor expansion of $\mb t_i(\mb x) = (t_{i,k}(\mb x))_{k\in[p]}$ with respect to $\mb x$. One has 
	\longeq{\label{eq: expansion of t_i,k}
	& t_{i,k}(\mb x) =t_{i,k}(\mb 0) +  \sum_{j\in[p]}\frac{\partial t_{i,k}}{\partial x_j}(\mb 0) x_j + 
	 O(\norm{\mb x}_2^2).
	}

	Recall that the $i$-th canonical basis of $\bb R^p$ is denoted by $\delta_i$ for $i \in[p]$. 
	For the quantity $R_{\hat \theta_n}^{-1}(\theta)_i$ with $\theta \in \mc M$, it then holds that 
	\longeq{\label{eq:basis transformation 1}
	& \big(R_{\hat\theta_n}^{-1}(\theta)\big)_i = \pi_{\{\mathrm dR_{\hat \theta_n} \delta_i\}}\circ \Log_{\hat \theta_n}(\theta) + O(\dist(\theta, \hat \theta_n)^3) \\ 
	= & \sum_{j\in[p]} (\mb t_j (\hat{\mb \eta}_n ))_i\Big(\pi_{\{T_{\theta_0 \rightarrow \Exp_{\theta}\circ \pi^{-1}_{\{e_i\}}(\hat{\mb \eta}_n)}(e_i)\}} \circ \Log_{\hat \theta_n}(\theta) \Big)_j + O(\dist(\theta, \hat \theta_n)^3)\\ 
	= & \sum_{j\in[p]} (\mb t_j (\hat{\mb \eta}_n ))_i\Big(\pi_{\{T_{\theta_0 \rightarrow \Exp_{\theta}\circ \pi^{-1}_{\{e_i\}}(\hat{\mb \eta}_n )}(e_i)\}} \circ \Log_{\hat \theta_n}(\theta) \Big)_j + O(n^{-\frac{3}{2}}(\log n)^{\frac{3}{2}})\\ 
	= & \sum_{j\in[p]} \big((\mb t_j (\mb H(\bar{\mb Z}_n)))_i + O(n^{-1} \log n)\big)\Big(\pi_{\{T_{\theta_0 \rightarrow \Exp_{\theta}\circ \pi^{-1}_{\{e_i\}}(\hat{\mb \eta}_n )}(e_i)\}} \circ \Log_{\hat \theta_n}(\theta) \Big)_j + O(n^{-\frac{3}{2}}(\log n)^{\frac{3}{2}}),
	}
		where the first equality holds by Lemma~\ref{lemma: uniform control for second-order retraction}, the second one follows by the definition of $\mb t_i(\mb x)$, and the last equality arises by \eqref{eq: approximation error of H(Z)} and the fact that $\mb t_j$ is continuously differentiable. 
	
	Taking $\theta = \theta_0$, we invoke the coordinate invariance under parallel transports to derive that 
	\longeq{\label{eq: basis transformation 1.1}
	&-\big( R^{-1}_{\hat \theta_n}(\theta_0)\big)_i = \sum_{j\in[p]} -(\mb t_j(\mb H(\bar{\mb Z}_n) ))_i (\hat{\mb \eta}_n)_j +  O(n^{-\frac{3}{2}}(\log n)^{\frac{3}{2}}) \\ 
	= &  \sum_{j\in[p]}t_{j,i}(\mb H(\bar{\mb Z}_n )) H(\bar{\mb Z}_n )_j +  O(n^{-\frac{3}{2}}(\log n)^{\frac{3}{2}}) \\ 
	= & H(\bar{\mb Z}_n )_i +\sum_{j,k,l\in[p]} \frac{\partial t_{j,i}}{\partial x_k}(\mb 0) H(\bar{\mb Z}_n )_k H(\bar{\mb Z}_n )_j +  O(n^{-\frac{3}{2}}(\log n)^{\frac{3}{2}}),
	}
	where the last line holds by the fact that $t_{i,j}(0)= \mathds 1\{i = j\}$ and \eqref{eq: expansion of t_i,k}. 

	We now move on to the empirical counterpart $R_{\hat\theta_n} (\hat \theta_n^*)$. Again, by jointly applying Lemma~\ref{lemma: uniform control for second-order retraction}, the property of the double exponential mapping (Lemma~\ref{lemma: normal coordinate transomation}), and a similar expansion of $\mb t_j$ as above, we find: 
	\longeq{\label{eq:basis transformation 2}
	& - R^{-1}_{\hat \theta_n^*}(\hat \theta_n)_i = \pi_{\{dR_{\hat \theta_n^*} e_i\} }\circ \Log_{\hat \theta_n^*}(\hat \theta_n) + O(\dist(\hat\theta_n, \hat \theta_n^*)^3) \\ 
	= & \bigg(\pi_{\{dR_{\hat \theta_n^*} e_i\} }\circ T_{\theta_0 \rightarrow \hat \theta_n^*}\Big( 
	\Log_{\theta_0}(\hat \theta_n) - \Log_{\theta_0}(\hat \theta_n^*)
	\Big) \bigg)_i+  O(n^{-\frac{3}{2}}(\log n)^{\frac{3}{2}})\\ 
	= & \sum_{j\in[p]} t_{j,i}(\mb H(\bar{\mb Y}_n ))H(\bar{\mb Y}_n)_j + \sum_{j,k,l\in[p]} \frac{\partial t_{j,i}}{\partial x_k}(\mb H(\bar{\mb Z}_n )) H(\bar{\mb Y}_n )_k H(\bar{\mb Y}_n )_l +  O(n^{-\frac{3}{2}}(\log n)^{\frac{3}{2}}). 
	}  
\paragraph*{Step 2: Studentization Analysis}\label{subsubsection: studentization analysis}
	Now we are placed to prove Claim~\ref{claim: claim 2} associated with the studentization quantities in Algorithm~\ref{algorithm:t based statistic}. We aim to show that the numerators and the denominators of the studentized statistics can be replaced by clearer forms related to $\bar{\mb Z}_n$ and $\bar{\mb Y}_n $. For clarity, we herein prove \eqref{eq: expansion in claim 2.2} for $i = 1$, as the other cases similarly follow in a similar manner. 
    
    We start with the analysis of the covariance terms. We collect the coordinate-transformation functions in a matrix $\mb T(\mb x) = \big(\mb t_1(\mb x), \ldots, \mb t_p(\mb x)\big) \in \bb R^{p\times p}$. By the orthogonality assumption in the definition of a retraction, it follows that $\mb T(\mb x)$ is an orthogonal matrix for every $\mb x\in \bb R^p$. 
		By replacing $\mb T(\pi_{\{e_i\}} \circ \Log_{\theta_0}(\hat\theta_n))$ with the rotation matrix $\hat{\mb T} \coloneqq \mb T\big(\mb H(\bar{\mb Z}_n)\big)$ at $\mb H(\bar{\mb Z}_n)$ and invoking \eqref{eq: approximation error of H(Z)}, $(\hat {\mb \Sigma})_{1,1}$ can be approximated as follows:  
	\begin{align}
	&  (\hat {\mb \Sigma})_{1,1} = \bigg((\bar \nabla^2 L_n^{\ret, \hat \theta_n})^{-1}\Big(\sum_{i\in[n]}\bar \nabla L^{\ret, \hat \theta_n}(X_i)\bar \nabla L^{\ret, \hat \theta_n}(X_i)\t /n  \Big)(\bar \nabla^2 L_n^{\ret, \hat \theta_n})^{-1}\bigg)_{1,1}  \\ 
	= & \Big(\big(\hat {\mb T} \bar \nabla^2L^{\Exp, \hat \theta_n, \{T_{\theta_0\rightarrow \hat \theta_n } e_i\}}_n  \hat {\mb T}\t \big)^{-1}\\ 
    & \cdot \big(\hat{\mb T}\big(\sum_{i\in[n]}\bar \nabla L^{\Exp, \hat \theta_n, \{T_{\theta_0\rightarrow \hat \theta_n } e_i\}}(X_i)\bar \nabla L^{\Exp, \hat \theta_n, \{T_{\theta_0\rightarrow \hat \theta_n } e_i\}}(X_i)\t  /n \big) \hat{\mb T}\t\big) \\ 
	\cdot & \big(\hat {\mb T} \bar \nabla^2 L^{\Exp, \hat \theta_n, \{T_{\theta_0\rightarrow \hat \theta_n } e_i\}}_n \hat {\mb T}\t \big)^{-1}\Big)_{1,1} +O(n^{-1}\log n)\label{eq: cov dif 1}.
	\end{align} 

    In what follows, we shall switch back to the chart centered at $\theta_0$ using intermediate results established in the proof of Claim~\ref{claim: claim 1}. Invoking the error controls \eqref{eq: proof of claim 1 gradient square control} and \eqref{eq: proof of claim 1 error control on the hessian} on the sum of squared gradient terms  and on the sample Hessian matrix, \eqref{eq: cov dif 1} turns into the form that 
    {\allowdisplaybreaks
    \begin{align}
         (\hat {\mb \Sigma})_{1,1}  
	= & \Big(\big(\hat {\mb T} \bar \nabla^2 L_n^{\Exp, \theta_0, \{e_i\}}(\mb H(\bar{\mb Z}_n)) \hat {\mb T}\t \big)^{-1} \\ 
    & \cdot \big(\sum_{i\in[n]}\hat{\mb T}\bar \nabla L^{\Exp, \theta_0, \{e_i\}} (\mb H(\bar{\mb Z}_n),X_i) \bar \nabla L^{\Exp, \theta_0, \{e_i\}}(\mb H(\bar{\mb Z}_n),X_i)\t\hat{\mb T}\t  /n  \big) \\ 
	&\cdot  \big(\hat {\mb T} \bar \nabla^2 L_n^{\Exp, \theta_0, \{e_i\}}(\mb H(\bar{\mb Z}_n))\hat {\mb T}\t \big)^{-1}\Big)_{1,1}  + O(n^{-1}\log n) \\
    = &\Big(\hat {\mb T}\big( \bar \nabla^2 L_n^{\Exp, \theta_0, \{e_i\}}(\mb H(\bar{\mb Z}_n))  \big)^{-1}\\ 
    &\cdot\big(\sum_{i\in[n]}\bar \nabla L^{\Exp, \theta_0, \{e_i\}} (\mb H(\bar{\mb Z}_n),X_i) \bar \nabla L^{\Exp, \theta_0, \{e_i\}}(\mb H(\bar{\mb Z}_n),X_i)\t /n  \big) \\ 
	\cdot & \big( \bar \nabla^2 L_n^{\Exp, \theta_0, \{e_i\}}(\mb H(\bar{\mb Z}_n)) \big)^{-1} \hat {\mb T}\t\Big)_{1,1}  + O(n^{-1}\log n),
	\label{eq: proof of claim 2 decomposition of hat Sigma 1}
	\end{align}
    }
 where the last line holds since $\hat {\mb T }$ is always an orthogonal matrix. 

We proceed to parse the difference $(\mb \Sigma)_{1,1} - (\hat {\mb \Sigma})_{1,1}$ by applying Taylor's expansion to the second term in \eqref{eq: proof of claim 2 decomposition of hat Sigma 1}. To facilitate this, we define a function $\mb \phi_{\mb A}(\mb X, \mb Z): \bb R^{p\times p} \times \bb R^{p\times p} \rightarrow \bb R^{p\times p}$ as 
\eq{
\mb \phi_{\mb A}(\mb X, \mb Z) \coloneqq \mb Z\mb X^{-1}\mb A\mb X^{-1} \mb Z\t,
}
where $\mb A$ is a fixed $p$-by-$p$ matrix. 

We assume that the smallest singular value of $\mb X_0, $ is bounded away from $0$ and the largest singular values of $\mb X, \mb Z$ are bounded by a constant $c_\sigma$. By the Taylor expansion, for every $\boldsymbol \triangle_1, \boldsymbol \triangle_2 \in \bb R^{p\times p}$ with $\norm{\boldsymbol \triangle_i}_F \leq c_\triangle\leq c_\sigma, i = 1,2$, we have 
\begin{align}
    &\mb  \phi_{\mb A}(\mb X_0 + \boldsymbol \triangle_1, \mb Z_0+ \boldsymbol \triangle_2) \\
    = &\mb  \phi_{\mb A}(\mb X_0, \mb Z_0) -  \mb Z_0\mb X_0^{-1}\boldsymbol \triangle_1\mb X_0^{-1} \mb A\mb X_0^{-1} \mb Z_0\t - \mb  Z_0\mb X_0^{-1}\mb A\mb X_0^{-1}\boldsymbol \triangle_1\mb X_0^{-1}\mb Z_0\t \\ 
    + & \boldsymbol \triangle_2\mb X_0^{-1}\mb A \mb X_0^{-1} \mb Z_0\t + \mb Z_0\mb X_0^{-1}\mb A \mb X_0^{-1} \boldsymbol \triangle_2\t+  c_\delta \max\{\norm{\boldsymbol \triangle_1}_F^2, \norm{\boldsymbol \triangle_2}_F^2\},
    \label{eq: fact of taylor expansion}
\end{align} 
where the constant $c_\delta$ is determined by $c_\sigma$, $c_\triangle$, and $\norm{\mb A}_F$.

Moreover, we introduce the following function: 
\eq{
\label{eq: G_n definition}
\mb G(\mb\eta, \mc X_n) \coloneqq\sum_{i\in[n]}\bar \nabla L^{\Exp, \theta_0, \{e_i\}} (\mb \eta,X_i) \bar \nabla L^{\Exp, \theta_0, \{e_i\}}(\mb \eta,X_i)\t /n,
} and
denote the matrix $\mb G(\mb H(\bar{\mb Z}_n), \mc X_n)$ by $\hat{\mb G}_n$. 
Applying the above fact \eqref{eq: fact of taylor expansion} to  \eqref{eq: proof of claim 2 decomposition of hat Sigma 1} yields that 
{\allowdisplaybreaks
\begin{align}
    &  \bigg(\hat {\mb T}\Big( \bar \nabla^2 L_n^{\Exp, \theta_0, \{e_i\}}(\mb H(\bar{\mb Z}_n))  \Big)^{-1}\hat{\mb G}_n\Big( \bar \nabla^2 L_n^{\Exp, \theta_0, \{e_i\}}(\mb H(\bar{\mb Z}_n)) \Big)^{-1} \hat {\mb T}\t\bigg)_{1,1} \\ 
    = & \bigg( \Big( \bar \nabla^2 L_n^{\Exp, \theta_0, \{e_i\}}(\mb 0)  \Big)^{-1}\hat{\mb G}_n\Big( \bar \nabla^2 L_n^{\Exp, \theta_0, \{e_i\}}(\mb 0) \Big)^{-1}\bigg)_{1,1} \\
    + & \bigg(\big(\hat{\mb T} - \mb I\big)\Big( \bar \nabla^2 L_n^{\Exp, \theta_0, \{e_i\}}(\mb 0)  \Big)^{-1}\hat{\mb G}_n\Big( \bar \nabla^2 L_n^{\Exp, \theta_0, \{e_i\}}(\mb 0) \Big)^{-1} \bigg)_{1,1}\\ 
    +& \bigg(\Big( \bar \nabla^2 L_n^{\Exp, \theta_0, \{e_i\}}(\mb 0)  \Big)^{-1}\hat{\mb G}_n\Big( \bar \nabla^2 L_n^{\Exp, \theta_0, \{e_i\}}(\mb 0) \Big)^{-1} \big(\hat{\mb T} - \mb I\big)\t\bigg)_{1,1} \\ 
    +&  \bigg( \Big( \bar \nabla^2 L_n^{\Exp, \theta_0, \{e_i\}}(\mb 0)  \Big)^{-1}\Big( \bar \nabla^2 L_n^{\Exp, \theta_0, \{e_i\}}(\mb 0) -  \bar \nabla^2 L_n^{\Exp, \theta_0, \{e_i\}}(\mb H(\bar{\mb Z}_n))   \Big)\\ 
    \cdot & \Big( \bar \nabla^2 L_n^{\Exp, \theta_0, \{e_i\}}(\mb 0)  \Big)^{-1} \hat{\mb G}_n\Big( \bar \nabla^2 L_n^{\Exp, \theta_0, \{e_i\}}(\mb 0) \Big)^{-1}\bigg)_{1,1}\\ 
    + &  \bigg( \Big( \bar \nabla^2 L_n^{\Exp, \theta_0, \{e_i\}}(\mb 0)  \Big)^{-1}\hat{\mb G}_n\Big( \bar \nabla^2 L_n^{\Exp, \theta_0, \{e_i\}}(\mb 0) \Big)^{-1} \\
    \cdot &\Big( \bar \nabla^2 L_n^{\Exp, \theta_0, \{e_i\}}(\mb 0) -  \bar \nabla^2 L_n^{\Exp, \theta_0, \{e_i\}}(\mb H(\bar{\mb Z}_n))  \Big) \Big( \bar \nabla^2 L_n^{\Exp, \theta_0, \{e_i\}}(\mb 0)  \Big)^{-1}\bigg)_{1,1}\\ 
    + & c'\max\Big\{\norm{\hat {\mb T} - \mb I}_F^2, \norm{ \bar \nabla^2 L_n^{\Exp, \theta_0, \{e_i\}}(\mb 0) -  \bar \nabla^2 L_n^{\Exp, \theta_0, \{e_i\}}(\mb H(\bar{\mb Z}_n))}_F^2
    \Big\}\label{eq: proof of claim2 cov decomposition expansion}
\end{align}
}
holds for some constant $c'$ since $\big\lVert\hat{\mb G}_n\big\rVert$ is upper bounded by some constant under the event $\mc F_n'$. 

Furthermore, applying a Taylor expansion again to $\hat {\mb T} - \mb I$ and $ \bar \nabla^2 L_n^{\Exp, \theta_0, \{e_i\}}(\mb 0) -  \bar \nabla^2 L_n^{\Exp, \theta_0, \{e_i\}}(\mb H(\bar{\mb Z}_n))$ provides us with a linear approximation in terms of  $\mb H(\bar{\mb Z}_n)$ as follows
\begin{align}
&\hat{\mb T} - \mb I = \sum_{j\in[p]}\frac{\partial }{\partial x_j}\mb T(\mb 0)(\mb H(\bar{\mb Z}_n) )_j+ c_T\norm{\mb H(\bar{\mb Z}_n)}_2^2,  \label{eq: hat T taylor}\\ 
& \bar \nabla^2 L_n^{\Exp, \theta_0, \{e_i\}}(\mb 0) -  \bar \nabla^2 L_n^{\Exp, \theta_0, \{e_i\}}(\mb H(\bar{\mb Z}_n)) \\ 
 = & - \sum_{j\in[p]} (\mb H(\bar{\mb Z}_n) )_j \frac{\partial }{\partial x_j}\bar \nabla^2 L_n^{\Exp, \theta_0, \{e_i\}}  + c_h \norm{\mb H(\bar{\mb Z}_n)}_2^2 \label{eq: Hessian matrix taylor}
\end{align}
for some constants $c_T$ and $c_h$ since the fourth-time derivatives of $L_n^{\Exp, \theta_0, \{e_i\}}$ is locally bounded by a constant under the event $\mc F_n'$. 

What remains to be justified is the quantity $\hat{\mb G}_n$. In view of Taylor's expansion of $\mb G_n$, $\hat{\mb G}_n$ can be written as 
\begin{align}
     \hat{\mb G}_n  
    = & \sum_{i\in[n]}\bar \nabla L^{\Exp, \theta_0, \{e_i\}} (\mb 0, X_i) \bar \nabla L^{\Exp, \theta_0, \{e_i\}}(\mb 0,X_i)\t /n \\ 
    + & \sum_{j\in[p]} \Big(\sum_{i\in[n]}\frac{\partial }{\partial x_j}\bar \nabla L^{\Exp, \theta_0, \{e_i\}} (\mb 0, X_i) \bar \nabla L^{\Exp, \theta_0, \{e_i\}}(\mb 0,X_i)\t /n\Big)(\mb H(\bar{\mb Z}_n) )_j\\ 
    + &\sum_{j\in[p]} \Big(\sum_{i\in[n]}\bar \nabla L^{\Exp, \theta_0, \{e_i\}}(\mb 0,X_i) \frac{\partial }{\partial x_j}\bar \nabla L^{\Exp, \theta_0, \{e_i\}} (\mb 0, X_i)\t  /n\Big)(\mb H(\bar{\mb Z}_n) )_j\\ 
    + & \sum_{j\in[p]} \sum_{k \in [p]}\int_{0}^1(1- t)\frac{\partial^2 }{\partial x_i \partial x_j}\mb G_n(t\mb H(\bar{\mb Z}_n)) \mathrm dt (\mb H(\bar{\mb Z}_n))_i(\mb H(\bar{\mb Z}_n) )_j  \\ 
    = & \sum_{i\in[n]}\bar \nabla L^{\Exp, \theta_0, \{e_i\}} (\mb 0, X_i) \bar \nabla L^{\Exp, \theta_0, \{e_i\}}(\mb 0,X_i)\t /n \\ 
    + & \sum_{j\in[p]} \Big(\sum_{i\in[n]}\frac{\partial }{\partial x_j}\bar \nabla^2 L^{\Exp, \theta_0, \{e_i\}} (\mb 0, X_i) \bar \nabla L^{\Exp, \theta_0, \{e_i\}}(\mb 0,X_i)\t /n\Big)(\mb H(\bar{\mb Z}_n) )_j \\ 
    + &\sum_{j\in[p]} \Big(\sum_{i\in[n]}\bar \nabla L^{\Exp, \theta_0, \{e_i\}}(\mb 0,X_i) \frac{\partial }{\partial x_j}\bar \nabla L^{\Exp, \theta_0, \{e_i\}} (\mb 0, X_i)\t  /n\Big)(\mb H(\bar{\mb Z}_n) )_j \\ 
    + & c_G \norm{\mb H(\bar{\mb Z}_n)}_2^2\label{eq: Ghat taylor}
\end{align}
for some constant $c_G$ where the last equality holds since $\frac{\partial^2 }{\partial x_i \partial x_j}\mb G_n(t\mb H(\bar{\mb Z}_n))$ for $t\in[0,1]$ is uniformly upper bounded by some constant $c_G$ under the event $\mc F_n$. 

Armed with the above approximations, we substitute  \eqref{eq: proof of claim2 cov decomposition expansion}, \eqref{eq: hat T taylor}, \eqref{eq: Hessian matrix taylor}, and \eqref{eq: Ghat taylor} into \eqref{eq: proof of claim 2 decomposition of hat Sigma 1}, and derive that 
{\allowdisplaybreaks
\begin{align}
    &(\hat {\mb \Sigma})_{1,1} \\
    = &  \bigg(\hat {\mb T}\Big( \bar \nabla^2 L_n^{\Exp, \theta_0, \{e_i\}}(\mb H(\bar{\mb Z}_n))  \Big)^{-1}\hat{\mb G}_n\Big( \bar \nabla^2 L_n^{\Exp, \theta_0, \{e_i\}}(\mb H(\bar{\mb Z}_n)) \Big)^{-1} \hat {\mb T}\t\bigg)_{1,1}+ O(n^{-1}\log n) \\ 
    =& (\mb \Sigma)_{1,1} +  \hat{\mb \xi}_R(\bar{\mb Z}_n, \mb \mu) +O(n^{-1}\log n),\label{eq: proof of claim2 cov decomposition 2}
\end{align}
}
where $\hat{\mb \xi}_R(\bar{\mb Z}_n, \mb \mu)$ is a quantity determined by $\bar{\mb Z}_n$ and can be expressed as 
\begin{align}
    & \hat{\mb \xi}_R(\bar{\mb Z}_n, \mb \mu)    \coloneqq   \Big(\mb T(\mb H(\mb \mu)) \hat{\mb \xi}(\bar{\mb Z}_n) \mb T(\mb H(\mb \mu))\t\Big)_{1,1} \\ 
    &+  \bigg(\big(\sum_{j\in[p]}\frac{\partial }{\partial x_j}\mb T(\mb H(\mb \mu))(\mb H(\bar{\mb Z}_n))_j\big)\Big( \bar \nabla^2 L_n^{\Exp, \theta_0, \{e_i\}}(\mb 0)  \Big)^{-1}\mb G(\mb 0, \mc X_n)\Big( \bar \nabla^2 L_n^{\Exp, \theta_0, \{e_i\}}(\mb 0) \Big)^{-1} \bigg)_{1,1}\\ 
    &+ \bigg(\Big( \bar \nabla^2 L_n^{\Exp, \theta_0, \{e_i\}}(\mb 0)  \Big)^{-1}\mb G(\mb 0, \mc X_n)\Big( \bar \nabla^2 L_n^{\Exp, \theta_0, \{e_i\}}(\mb 0) \Big)^{-1} \big(\sum_{j\in[p]}\frac{\partial }{\partial x_j}\mb T(\mb H(\mb \mu))(\mb H(\bar{\mb Z}_n))_j\big)\t \bigg)_{1,1}, 
\end{align}
with 
\begin{align}
    \hat{\mb \xi}(\bar{\mb Z}_n) \coloneqq & \Big( \big( \bar \nabla^2 L_n^{\Exp, \theta_0, \{e_i\}}(\mb 0)  \big)^{-1}\mb G(\mb 0, \mc X_n) \big( \bar \nabla^2 L_n^{\Exp, \theta_0, \{e_i\}}(\mb 0) \big)^{-1} - \mb \Sigma\Big)_{1,1} \\
    + & \bigg( \Big( \bar \nabla^2 L_n^{\Exp, \theta_0, \{e_i\}}(\mb 0)  \Big)^{-1}\Big(  \sum_{j\in[p]} \Big(\sum_{i\in[n]}\frac{\partial }{\partial x_j}\bar \nabla L^{\Exp, \theta_0, \{e_i\}} (\mb 0, X_i) \bar \nabla L^{\Exp, \theta_0, \{e_i\}}(\mb 0,X_i)\t /n\Big)H(\bar{\mb Z}_n)_j \Big)\\ 
    \cdot & \Big( \bar \nabla^2 L_n^{\Exp, \theta_0, \{e_i\}}(\mb 0) \Big)^{-1}\bigg)_{1,1} \\ 
    + & 
    \bigg( \Big( \bar \nabla^2 L_n^{\Exp, \theta_0, \{e_i\}}(\mb 0)  \Big)^{-1}\Big(\sum_{j\in[p]} \Big(\sum_{i\in[n]}\bar \nabla L^{\Exp, \theta_0, \{e_i\}}(\mb 0,X_i) \frac{\partial }{\partial x_j}\bar \nabla L^{\Exp, \theta_0, \{e_i\}} (\mb 0, X_i)\t  /n\Big)H(\bar{\mb Z}_n)_j  \Big) \\ 
    \cdot & \Big( \bar \nabla^2 L_n^{\Exp, \theta_0, \{e_i\}}(\mb 0) \Big)^{-1}\bigg)_{1,1}
    \\ 
    -&  \bigg( \Big( \bar \nabla^2 L_n^{\Exp, \theta_0, \{e_i\}}(\mb 0)  \Big)^{-1}\Big( \sum_{j\in[p]} H(\bar{\mb Z}_n)_j \frac{\partial }{\partial x_j}\bar \nabla^2 L_n^{\Exp, \theta_0, \{e_i\}}   \Big)\\ 
    \cdot & \Big( \bar \nabla^2 L_n^{\Exp, \theta_0, \{e_i\}}(\mb 0)  \Big)^{-1}\mb G(\mb 0, \mc X_n)\Big( \bar \nabla^2 L_n^{\Exp, \theta_0, \{e_i\}}(\mb 0) \Big)^{-1}\bigg)_{1,1}\\ 
    - &  \bigg( \Big( \bar \nabla^2 L_n^{\Exp, \theta_0, \{e_i\}}(\mb 0)  \Big)^{-1}\mb G(\mb 0, \mc X_n)\Big( \bar \nabla^2 L_n^{\Exp, \theta_0, \{e_i\}}(\mb 0) \Big)^{-1} \\
    \cdot &\Big( \sum_{j\in[p]} H(\bar{\mb Z}_n)_j \frac{\partial }{\partial x_j}\bar \nabla^2 L_n^{\Exp, \theta_0, \{e_i\}}  \Big) \Big( \bar \nabla^2 L_n^{\Exp, \theta_0, \{e_i\}}(\mb 0)  \Big)^{-1}\bigg)_{1,1}. 
\end{align}
Here we make use of the fact that $\norm{\mb H(\bar{\mb Z}_n)}_2$ is of the order $O(n^{-\frac{1}{2}}(\log n)^{\frac{1}{2}})$ under the event $\mc F_n$. We hereby make note of the quantity $\hat{\mb \xi}_R(\bar{\mb Z}_n, \mb \mu)$ that all involved derivatives of $L^{\Exp, \theta_0, \{e_i\}}$ or $L_n^{\Exp, \theta_0, \{e_i\}}$ can be viewed as the components of $\bar{\mb Z}_n$. Thus $\hat{\mb \xi}_R: \bb R^{\mathrm{dim}(\bar{\mb Z}_n)} \times \bb R^{\mathrm{dim}(\bar{\mb Z}_n)} \rightarrow \bb R^{p \times p}$ is a function independently of $L^{\Exp, \theta_0, \{e_i\}}$ or $L_n^{\Exp, \theta_0, \{e_i\}}$ and will be utilized again in the analysis of the resampled quantities.

An immediate consequence of \eqref{eq: proof of claim2 cov decomposition 2} is that 
\eq{
\big|(\mb \Sigma)_{1,1} - (\hat {\mb \Sigma})_{1,1}\big| = O(n^{-\frac{1}{2}}(\log n)^{\frac{1}{2}})
}
under the event $\mc F_n$.

    Now we are ready to decompose the studentized statistic $-\frac{R^{-1}_{\hat\theta_n}(\theta_0)_i}{(\hat{\mb \Sigma})_{1,1}^{\frac12}}$. Combining \eqref{eq: basis transformation 1.1} and the fact that $\norm{\mb H(\bar{\mb Z}_n)}_2$ is of the order $O(n^{-\frac{1}{2}}(\log n)^{\frac{1}{2}})$ under the event $\mc F_n$ 
    gives that 
    {\allowdisplaybreaks
    \begin{align}
	 &-\frac{R^{-1}_{\hat \theta_n}(\theta_0)_1}{(\hat{\mb \Sigma})_{1,1}^{\frac{1}{2}}}\\
	 = & \frac{H(\bar{\mb Z}_n )_1 +\sum_{j,k\in[p]} \frac{\partial t_{j,1}}{\partial x_k}(\mb 0) H(\bar{\mb Z}_n )_k H(\bar{\mb Z}_n )_j}{(\hat{\mb \Sigma})_{1,1}^{\frac{1}{2}}} +  O(n^{-\frac{3}{2}}(\log n)^{\frac{3}{2}})
	 \\ = & \frac{H(\bar{\mb Z}_n)_1 +\sum_{j,k\in[p]} \frac{\partial t_{j,1}}{\partial x_k}(\mb 0) H(\bar{\mb Z}_n)_k H(\bar{\mb Z}_n)_j }{(\mb \Sigma)_{1,1}^{\frac{1}{2}}}\\ 
     &  + \frac{H(\bar{\mb Z}_n)_1 +\sum_{j,k\in[p]} \frac{\partial t_{j,1}}{\partial x_k}(\mb 0) H(\bar{\mb Z}_n)_k H(\bar{\mb Z}_n)_j }{(\mb \Sigma)_{1,1}}\left((\mb \Sigma)_{1,1}^{\frac{1}{2}} - (\hat {\mb \Sigma})_{1,1}^{\frac{1}{2}}\right)\\ 
	 &+ \frac{H(\bar{\mb Z}_n)_1 + \sum_{j,k\in[p]} \frac{\partial t_{j,1}}{\partial x_k}(\mb 0) H(\bar{\mb Z}_n)_k H(\bar{\mb Z}_n)_j }{(\mb \Sigma)_{1,1}(\hat {\mb \Sigma})_{1,1}^{\frac{1}{2}}}\left((\mb \Sigma)_{1,1}^{\frac{1}{2}} - (\hat {\mb \Sigma})_{1,1}^{\frac{1}{2}}\right)^2 + O(n^{-\frac{3}{2}}(\log n)^{\frac{3}{2}})\\ 
	= & \frac{H(\bar{\mb Z}_n)_1 +\sum_{j,k\in[p]} \frac{\partial t_{j,1}}{\partial x_k}(\mb 0) H(\bar{\mb Z}_n)_k H(\bar{\mb Z}_n)_j }{(\mb \Sigma)_{1,1}^{\frac{1}{2}}}\\ 
    &  + \frac{H(\bar{\mb Z}_n)_1 +\sum_{j,k\in[p]} \frac{\partial t_{j,1}}{\partial x_k}(\mb 0) H(\bar{\mb Z}_n)_k H(\bar{\mb Z}_n)_j }{(\mb \Sigma)_{1,1}}\left((\mb \Sigma)_{1,1}^{\frac{1}{2}} - (\hat {\mb \Sigma})_{1,1}^{\frac{1}{2}}\right)+ O(n^{-\frac{3}{2}}(\log n)^{\frac{3}{2}}) \\ 
	= & \frac{H(\bar{\mb Z}_n)_1 +\sum_{j,k\in[p]} \frac{\partial t_{j,1}}{\partial x_k}(\mb 0) H(\bar{\mb Z}_n)_k H(\bar{\mb Z}_n)_j }{(\mb \Sigma)_{1,1}^{\frac{1}{2}}}\\ 
    &  + \frac{H(\bar{\mb Z}_n)_1 +\sum_{j,k\in[p]} \frac{\partial t_{j,1}}{\partial x_k}(\mb 0) H(\bar{\mb Z}_n)_k H(\bar{\mb Z}_n)_j }{2(\mb \Sigma)_{1,1}^{\frac{3}{2}} }\left((\mb \Sigma)_{1,1} - (\hat {\mb \Sigma})_{1,1}\right)+ +O(n^{-\frac{3}{2}}(\log n)^{\frac{3}{2}}).\label{eq: studentized term expansion}
 	\end{align}
	}

Substituting \eqref{eq: proof of claim2 cov decomposition 2} into \eqref{eq: studentized term expansion} establishes the first part of the claim:
{\allowdisplaybreaks
	\begin{align}\label{eq: f_1}
	& -\frac{R_{\hat \theta_n}(\theta_0)_1}{(\hat{\mb \Sigma})_{1,1}^{\frac{1}{2}}} 
 \\= & 
\frac{H(\bar{\mb Z}_n)_1 +\sum_{j,k\in[p]} \frac{\partial t_{j,1}}{\partial x_k}(\mb 0) H(\bar{\mb Z}_n)_k H(\bar{\mb Z}_n)_j }{(\mb \Sigma)_{1,1}^{\frac{1}{2}}}
\\ 
& -\frac{H(\bar{\mb Z}_n)_1 +\sum_{j,k\in[p]} \frac{\partial t_{j,1}}{\partial x_k}(\mb 0) H(\bar{\mb Z}_n)_k H(\bar{\mb Z}_n)_j }{2(\mb \Sigma)_{1,1}^{\frac{3}{2}} }\hat{\mb \xi}_R(\bar{\mb Z}_n, \mb \mu)
+ O(n^{-\frac{3}{2}}(\log n)^{\frac{3}{2}}) \\
 = & \frac{1}{(\mb \Sigma)_{1,1}^{\frac{1}{2}}}\Big(H(\bar{\mb Z}_n)_1 +\sum_{j,k\in[p]} \frac{\partial t_{j,1}}{\partial x_k}(\mb 0) H(\bar{\mb Z}_n)_k H(\bar{\mb Z}_n)_j 
\\ 
& -\frac{H(\bar{\mb Z}_n)_1 +\sum_{j,k\in[p]} \frac{\partial t_{j,1}}{\partial x_k}(\mb 0) H(\bar{\mb Z}_n)_k H(\bar{\mb Z}_n)_j }{2(\mb \Sigma)_{1,1} }\hat{\mb \xi}_R(\bar{\mb Z}_n, \mb \mu)\Big)+ O(n^{-\frac{3}{2}}(\log n)^{\frac{3}{2}})\\
= & \frac{\tilde H_1(\bar{\mb Z}_n, \mb \mu)}{(\mb \Sigma)_{1,1}^{\frac{1}{2}}} + O(n^{-\frac{3}{2}}(\log n)^{\frac{3}{2}}),
\end{align}
}
	where we define the differentiable function $\tilde H_i(\mb x, \mb y)$ appearing in \eqref{eq: expansion in claim 2.1} as 
 \begin{align}
     & \tilde H_1(\mb x, \mb y) \coloneqq \sum_{j \in[p]}t_{j,1}(\mb H(\mb y)) H(\mb x)_j +  \sum_{j,k\in[p]} \frac{\partial t_{j,1}}{\partial x_k}(\mb H(\mb y)) H(\mb x)_k H(\mb x)_j \\
     & - \frac{\sum_{j \in[p]}t_{j,1}(\mb H(\mb y)) H(\mb x)_j +\sum_{j,k\in[p]} \frac{\partial t_{j,1}}{\partial x_k}(\mb H(\mb y) ) H(\mb x)_k H(\mb x)_j }{2(\mb \Sigma)_{1,1} }\hat{\mb \xi}_R(\mb x, \mb y).
     \end{align}

The proof of \eqref{eq: expansion in claim 2.2} closely follows the arguments presented earlier, but proceeds under the event $\mc F_n \cap \mc F_n^*$. Below, we provide a brief outline of the main steps.

    Recalling the function $\mb G$ defined in \eqref{eq: G_n definition}, we introduce the quantity:
    \eq{
    \mb G(\mb \eta, \mc X_n^*)  \coloneqq\sum_{i\in[n]}\bar \nabla L^{\Exp, \theta_0, \{e_i\}} (\mb \eta,X_i^*) \bar \nabla L^{\Exp, \theta_0, \{e_i\}}(\mb \eta,X_i^*)\t /n.
    }
    We also denote that $\hat {\mb T}^* \coloneqq \mb T(\mb H(\bar{\mb Z}_n) +  
    \mb H(\bar{\mb Y}_n))$ and $\hat{\mb G}_n^*\coloneqq \mb G(\mb H(\bar{\mb Z}_n) + \mb H(\bar{\mb Y}_n), \mc X_n^*)$. 

Notice that a similar fact to \eqref{eq: hat T taylor} and \eqref{eq: Hessian matrix taylor} holds under the event $\mc F_n \cap \mc F_n^*$ that 
\begin{align}
    & \hat {\mb T}^* - \hat{\mb T}  = \sum_{j\in[p]} \frac{\partial }{\partial x_j}T(\bar{\mb Z}_n)H(\bar{\mb Y}_n)_j + c_T'\norm{\mb H(\bar{\mb Y}_n)}_2^2, \\ 
    &    \bar \nabla^2 {L_n^*}^{\Exp, \theta_0, \{e_i\}}(\mb H(\bar{\mb Z}_n))- \bar \nabla^2 {L_n^*}^{\Exp, \theta_0, \{e_i\}}(\mb H(\bar{\mb Z}_n) + \mb H(\bar {\mb Y}_n))  \\
    = &- \sum_{j\in[p]} H(\bar{\mb Y}_n)_j \frac{\partial }{\partial x_j}\bar \nabla^2 {L_n^*}^{\Exp, \theta_0, \{e_i\}}(\mb H(\bar{\mb Z}_n) )  + c_h' \norm{\mb H(\bar{\mb Y}_n)}_2^2.
\end{align}

With the idea of Taylor's expansion in mind together with Claim~\ref{claim: claim 1}, we can deduce that under the event $\mc F_n \cap \mc F_n^*$
{\allowdisplaybreaks
\begin{align} 
    &(\check{\mb \Sigma})_{1,1} \\ 
    = & \bigg(\hat {\mb T}^*\Big( \bar \nabla^2 {L_n^*}^{\Exp, \theta_0, \{e_i\}}(\mb H(\bar{\mb Z}_n) + \mb H(\bar{\mb Y}_n)) \Big)^{-1}\\ 
    \cdot & \hat{\mb G}_n^* \Big( \bar \nabla^2 {L_n^*}^{\Exp, \theta_0, \{e_i\}}(\mb H(\bar{\mb Z}_n) + \mb H(\bar{\mb Y}_n)) \Big)^{-1} {{}\hat {\mb T}^*}\t\bigg)_{1,1} + O(n^{-1}\log n)\\ 
    = & (\hat {\mb \Sigma})_{1,1} \\ 
      + &\bigg( \hat {\mb T}\Big( \bar \nabla^2 {L^*_n}^{\Exp, \theta_0, \{e_i\}}(\mb H(\bar{\mb Z}_n))  \Big)^{-1}\hat{\mb G}_n^*\Big( \bar \nabla^2 {L^*_n}^{\Exp, \theta_0, \{e_i\}}(\mb H(\bar{\mb Z}_n)) \Big)^{-1} \hat {\mb T}\t - (\hat {\mb \Sigma})_{1,1}\bigg)_{1,1} \\
     + &  \bigg( \Big( \bar \nabla^2 {L_n^*}^{\Exp, \theta_0, \{e_i\}}(\mb H(\bar{\mb Z}_n) )  \Big)^{-1} \\
    & \hspace{-0.5cm}\cdot  \Big(  \sum_{j\in[p]} \Big(\sum_{i\in[n]}\frac{\partial }{\partial x_j}\bar \nabla L^{\Exp, \theta_0, \{e_i\}} (\mb H(\bar{\mb Z}_n) , X_i^*) \bar \nabla L^{\Exp, \theta_0, \{e_i\}}(\mb H(\bar{\mb Z}_n) , X_i^*)\t /n\Big)H(\bar{\mb Y}_n)_j \Big)
    \\ 
    &\cdot  \Big( \bar \nabla^2 {L_n^*}^{\Exp, \theta_0, \{e_i\}}(\mb H(\bar{\mb Z}_n)) \Big)^{-1}\bigg)_{1,1} \\ 
    +&   
    \bigg( \Big( \bar \nabla^2 {L_n^*}^{\Exp, \theta_0, \{e_i\}}(\mb H(\bar{\mb Z}_n) )  \Big)^{-1}\\
     & \hspace{-0.5cm} \cdot \Big(\sum_{j\in[p]} \Big(\sum_{i\in[n]}\bar \nabla L^{\Exp, \theta_0, \{e_i\}}(\mb H(\bar{\mb Z}_n) ,X_i^*) \frac{\partial }{\partial x_j}\bar \nabla L^{\Exp, \theta_0, \{e_i\}} (\mb H(\bar{\mb Z}_n) , X_i^*)\t  /n\Big)H(\bar{\mb Y}_n)_j  \Big) \\ 
    \cdot & \Big( \bar \nabla^2 {L_n^*}^{\Exp, \theta_0, \{e_i\}}(\mb H(\bar{\mb Z}_n) ) \Big)^{-1}\bigg)_{1,1}
    \\ 
    + & \bigg(\big(\sum_{i\in[p]}\frac{\partial }{\partial x_i}\mb T(\mb H(\bar{\mb Z}_n) )H_i(\bar{\mb Y}_n)\big) \\ 
    \cdot & \Big( \bar \nabla^2 {L_n^*}^{\Exp, \theta_0, \{e_i\}}(\mb H(\bar{\mb Z}_n) )  \Big)^{-1}\mb G_n^*(\mb H(\bar{\mb Z}_n) )\Big( \bar \nabla^2 {L_n^*}^{\Exp, \theta_0, \{e_i\}}(\mb H(\bar{\mb Z}_n) ) \Big)^{-1} \bigg)_{1,1}\\ 
    +& \bigg(\Big( \bar \nabla^2 {L_n^*}^{\Exp, \theta_0, \{e_i\}}(\mb H(\bar{\mb Z}_n) )  \Big)^{-1}\mb G_n^*(\mb H(\bar{\mb Z}_n) )\\ 
    \cdot & \Big( \bar \nabla^2 {L_n^*}^{\Exp, \theta_0, \{e_i\}}(\mb H(\bar{\mb Z}_n) ) \Big)^{-1} \big(\sum_{i\in[p]}\frac{\partial }{\partial x_i}\mb T(\mb H(\bar{\mb Z}_n) )H_i(\bar{\mb Y}_n)\big)\t \bigg)_{1,1} \\ 
    -&  \bigg( \Big( \bar \nabla^2 {L_n^*}^{\Exp, \theta_0, \{e_i\}}(\mb H(\bar{\mb Z}_n) )  \Big)^{-1}\Big( \sum_{j\in[p]} H(\bar{\mb Y}_n)_j \frac{\partial }{\partial x_j}\bar \nabla^2 {L_n^*}^{\Exp, \theta_0, \{e_i\}} (\mb H(\bar{\mb Z}_n) )   \Big)\\ 
    \cdot & \Big( \bar \nabla^2 {L_n^*}^{\Exp, \theta_0, \{e_i\}}(\mb H(\bar{\mb Z}_n) )  \Big)^{-1}\mb G_n^*(\mb H(\bar{\mb Z}_n) )\Big( \bar \nabla^2 {L_n^*}^{\Exp, \theta_0, \{e_i\}}(\mb H(\bar{\mb Z}_n) ) \Big)^{-1}\bigg)_{1,1}\\ 
    - &  \bigg( \Big( \bar \nabla^2 {L_n^*}^{\Exp, \theta_0, \{e_i\}}(\mb H(\bar{\mb Z}_n ))   \Big)^{-1}\mb G_n^*(\mb H(\bar{\mb Z}_n ))\Big( \bar \nabla^2{ L_n^*}^{\Exp, \theta_0, \{e_i\}}(\mb H(\bar{\mb Z}_n ) ) \Big)^{-1} \\
    \cdot &\Big( \sum_{j\in[p]} H(\bar{\mb Y}_n)_j \frac{\partial }{\partial x_j}\bar \nabla^2 {L_n^*}^{\Exp, \theta_0, \{e_i\}} (\mb H(\bar{\mb Z}_n) )  \Big) \Big( \bar \nabla^2 {L_n^*}^{\Exp, \theta_0, \{e_i\}}(\mb H(\bar{\mb Z}_n) )  \Big)^{-1}\bigg)_{1,1}\\
     + &  O(n^{-1}\log n) \\ 
    = & (\hat {\mb \Sigma})_{1,1} +  \hat{\mb \xi}_R(\bar{\mb Y}_n, \bar{\mb Z}_n) +O(n^{-1}\log n).
    \label{eq: proof of claim 2 resampled cov decompostion}
\end{align}
}

Again, for the resampled studentized term $-\frac{R^{-1}_{\hat\theta_n^*}(\hat \theta_n)_1 }{(\check{\mb \Sigma})_{1,1}^{\frac{1}{2}}}$ under the event $\mc F_n \cap \mc F_n^*$,  we can decompose it as follows 
{\allowdisplaybreaks
	\begin{align}
	& - \frac{R_{\hat \theta_n^*}(\hat \theta_n)_1}{(\check{\mb \Sigma})_{1,1}^{\frac{1}{2}}} \\ 
	= & \frac{\sum_{j\in[p]} t_{j,1}(\mb H(\bar{\mb Z}_n))H(\bar{\mb Y}_n)_j + \sum_{j,k\in[p]} \frac{\partial t_{j,1}}{\partial x_k}(\mb H(\bar{\mb Z}_n)) H(\bar{\mb Y}_n)_k H(\bar{\mb Y}_n)_j}{(\check{\mb \Sigma})_{1,1}^{\frac{1}{2}}} \\ + &  \frac{\sum_{j\in[p]} t_{j,1}(\mb H(\bar{\mb Z}_n))H(\bar{\mb Y}_n)_j + \sum_{j,k\in[p]} \frac{\partial t_{j,1}}{\partial x_k}(\mb H(\bar{\mb Z}_n)) H(\bar{\mb Y}_n)_k H(\bar{\mb Y}_n)_j}{(\hat {\mb \Sigma})_{1,1}}\\ 
    & \cdot  \left((\hat {\mb \Sigma})_{1,1}^{\frac{1}{2}} - (\check{\mb \Sigma})_{1,1}^{\frac{1}{2}}\right)\\ 
	+ & \frac{\sum_{j\in[p]} t_{j,1}(\mb H(\bar{\mb Z}_n))H(\bar{\mb Y}_n)_j + \sum_{j,k\in[p]} \frac{\partial t_{j,i}}{\partial x_k}(\mb H(\bar{\mb Z}_n)) H(\bar{\mb Y}_n)_k H(\bar{\mb Y}_n)_j}{(\hat {\mb \Sigma})_{1,1}(\check{\mb \Sigma})_{1,1}^{\frac{1}{2}}} \\ 
    & \cdot \left((\hat {\mb \Sigma})_{1,1}^{\frac{1}{2}} - (\check{\mb \Sigma})_{1,1}^{\frac{1}{2}}\right)^2 + O(n^{-\frac{3}{2}}(\log n)^{\frac{3}{2}}) \\ 
	= & \frac{\sum_{j\in[p]} t_{j,1}(\mb H(\bar{\mb Z}_n))H(\bar{\mb Y}_n)_j + \sum_{j,k\in[p]} \frac{\partial t_{j,1}}{\partial x_k}(\mb H(\bar{\mb Z}_n)) H(\bar{\mb Y}_n)_k H(\bar{\mb Y}_n)_j}{(\hat {\mb \Sigma})_{1,1}} \\ + &  \frac{\sum_{j\in[p]} t_{j,1}(\mb H(\bar{\mb Z}_n))H(\bar{\mb Y}_n)_j + \sum_{j,k\in[p]} \frac{\partial t_{j,1}}{\partial x_k}(\mb H(\bar{\mb Z}_n)) H(\bar{\mb Y}_n)_k H(\bar{\mb Y}_n)_j}{2(\hat {\mb \Sigma})_{1,1}^{\frac{3}{2}}}\\ 
    & \cdot \left((\hat {\mb \Sigma})_{1,1} - (\check{\mb \Sigma})_{1,1}\right) +  O(n^{-\frac{3}{2}}(\log n)^{\frac{3}{2}}),
    \label{eq: proof of claim 2 resampled cov decomposition 2}
	\end{align}
    }
    where we leverage the concentration conditions under $\mc F_n\cap \mc F_n^*$ again.

    Plugging \eqref{eq: proof of claim 2 resampled cov decompostion} into \eqref{eq: proof of claim 2 resampled cov decomposition 2} deduces that under the event $\mc F_n \cap \mc F_n^*$
    {\allowdisplaybreaks
    \longeq{
        &  - \frac{R_{\hat \theta_n^*}^{-1}(\hat \theta_n)_1}{(\check{\mb \Sigma})_{1,1}^{\frac{1}{2}}} \\ 
        = & \frac{\sum_{j\in[p]} t_{j,1}(\mb H(\bar{\mb Z}_n))H(\bar{\mb Y}_n)_j + \sum_{j,k\in[p]} \frac{\partial t_{j,1}}{\partial x_k}(\mb H(\bar{\mb Z}_n)) H(\bar{\mb Y}_n)_k H(\bar{\mb Y}_n)_j}{(\hat {\mb \Sigma})_{1,1}^{\frac{1}{2}}} \\ - &  \frac{\sum_{j\in[p]} t_{j,1}(\mb H(\bar{\mb Y}_n))H(\bar{\mb Y}_n)_j + \sum_{j,k\in[p]} \frac{\partial t_{j,1}}{\partial x_k}(\mb H(\bar{\mb Z}_n)) H(\bar{\mb Y}_n)_k H(\bar{\mb Y}_n)_l}{2(\hat {\mb \Sigma})_{1,1}^{\frac{3}{2}}}\hat{\mb \xi}_R(\bar{\mb Y}_n, \bar{\mb Z}_n)  \\ 
        & + O(n^{-\frac{3}{2}}(\log n)^{\frac{3}{2}}) \\ 
        = & \frac{1}{(\hat {\mb \Sigma})_{1,1}^{\frac{1}{2}}} \Big( \sum_{j\in[p]} t_{j,1}(\mb H(\bar{\mb Z}_n))H(\bar{\mb Y}_n)_j + \sum_{j,k\in[p]} \frac{\partial t_{j,1}}{\partial x_k}(\mb H(\bar{\mb Z}_n)) H(\bar{\mb Y}_n)_k H(\bar{\mb Y}_n)_j \\ 
        -&  \frac{\sum_{j\in[p]} t_{j,1}(\mb H(\bar{\mb Y}_n))H(\bar{\mb Y}_n)_j + \sum_{j,k\in[p]} \frac{\partial t_{j,1}}{\partial x_k}(\mb H(\bar{\mb Z}_n)) H(\bar{\mb Y}_n)_k H(\bar{\mb Y}_n)_j}{2\Sigma_{j,j}}\hat{\mb \xi}_R(\bar{\mb Y}_n, \bar{\mb Z}_n) \Big) \\ 
        & + O(n^{-\frac{3}{2}}(\log n)^{\frac{3}{2}}) \\
        \eqqcolon & \frac{\tilde H_1(\bar{\mb Y}_n, \bar{\mb Z}_n)}{(\hat {\mb \Sigma})_{1,1}^{\frac{1}{2}}} + O(n^{-\frac{3}{2}}(\log n)^{\frac{3}{2}}),
    }
    }
    which completes the proof. \qed

\subsection{Proof of Proposition~\ref{proposition:edgeworth for extrinsic t}}

    The proof of this result closely parallels that of Theorem~\ref{theorem:bootstrapmanifold}, and we therefore provide only a brief sketch. Throughout, we adopt the notation introduced in the proof of Theorem~\ref{theorem:bootstrapmanifold}. Compared with Theorem~\ref{theorem:bootstrapmanifold}, the extrinsic $t$-statistic returns to the concise representation under a fixed coordinate system. The condition that $ f\circ \Exp\circ \pi^{-1}_{\{e_i\}}$ is smooth in a neighborhood of $\mb 0$, together with  \eqref{eq: implicit function of manifold}, \eqref{eq: approximation error of H(Z)}, Theorem~\ref{thm:convergence rate of newton}, and Theorem~\ref{thm: convergence rate of resampled newton}, reveals that
\eq{\label{eq: approximation error of extrinsic t statistic}
    f(\hat{\theta}_n^*) - f(\hat{ \theta}_n) = f(\Exp_{ \theta_0}\circ \pi^{-1}_{\{e_i\}}(\mb H(\bar{\mb Z}_n)) )- f(\Exp_{ \theta_0}\circ \pi^{-1}_{\{e_i\}}(\mb H(\bar{\mb Z}_n^*))) + O(n^{-\frac{3}{2}}(\log n)^{\frac{3}{2}}).
} 

Note that $f(\Exp_{\mb \theta_0}\circ \pi^{-1}_{\{e_i\}}(\mb H (\mb x))$ is a sufficiently smooth function. By combining the Edgeworth expansions with respect to $f\big(\Exp\circ \pi^{-1}_{\{e_i\}}(\mb H(\bar {\mb Z}_n )\big)$ and $f\big(\Exp\circ \pi^{-1}_{\{e_i\}}(\mb H(\bar {\mb Z}_n^* )\big)$, together with an argument analogous to the proof of Claim~\ref{claim: claim 2} and the approximation error bound in \eqref{eq: approximation error of extrinsic t statistic}, we obtain
\longeq{
    & \limsup\limits_{n\rightarrow \infty} n\bb P_{\mc X_n}\bigg[\sup_{-\infty< x<\infty}\Big| \bb P\Big[\frac{f(\hat{\theta}_n^*) - f(\hat{\theta}_n)}{\big(\bar \nabla f(R_{\hat \theta_n^*}(\cdot))\mid_{\mb 0} \t 
    \check {\mb \Sigma}\bar \nabla f(R_{\hat \theta_n^*}(\cdot))\mid_{\mb 0}  \big)^{\frac{1}{2}}}\leq x \mid \mc X_n\Big] \\ 
    & \qquad - \bb P_{\mc X_n}\big[\frac{f(\hat{\theta}_n) - f(\theta_0)}{\big(\bar \nabla f(R_{\hat \theta_n}(\cdot))\mid_{\mb 0}\t  
    \hat {\mb \Sigma}\bar \nabla f(R_{\hat \theta_n}(\cdot))\mid_{\mb 0} \big)^{\frac{1}{2}}} \leq x\Big]\Big| \geq n^{-1+\xi} \bigg] \leq \infty, 
}
 again, by Proposition~\ref{proposition: thm 3.1 in bhattacharya}, Proposition~\ref{proposition: modified thm 3.3 in bhattacharya}, and Lemma~\ref{lemma:deltamethod}. \qed

\subsection{Proof of Proposition~\ref{proposition: hypothesis testing}}
\label{sec: proof of hypothesis testing}
The assumptions we make in Proposition~\ref{proposition: hypothesis testing} is an i.i.d. version of the assumptions in \cite[Theorme 2.(a) and Theorme 2.(b)]{andrews2002higher}. More precisely, in the context of the GMM estimator discussion in \cite{andrews2002higher} we let the function $g(X_i, \mb \eta)$ be $\bar \nabla L(\phi^{-1}(\mb \eta), X_i)$ and $\Omega$ be an identity matrix. To apply Theorem 2.(a) and Theorem 2.(b) in \cite{andrews2002higher} to the $t$-statistic and the Wald statistic, respectively, we let $d_1$ and $d_2$ be $5$ and $4$, respectively. Then those results immediately imply that 
\begin{align}
    & \bb P\big[T_{\phi,j} >  z^*_{T_{\phi,j},2, \alpha}\big]  = \alpha + o(n^{-1 + \xi}),\\ 
            & \bb P\big[T_{\phi,j} <  z^*_{T_{\phi,j}, 1-\alpha}\big] = \alpha + o(n^{-1 + \xi}), \\ 
			& \bb P \big[W_\phi > z^*_{W_\phi,\alpha}\big] = \alpha + o(n^{-\frac{3}{2} + \xi}),
\end{align}
where $z^*_{T_{\phi,j},2, \alpha}$, $z^*_{T_{\phi,j}, 1-\alpha}$, and $z^*_{W_\phi,\alpha}$ corresponds to the quantile analogs of the ``exact-version'' statistics  by replacing $\hat \theta_n$ and $\hat \theta_n^*$ with $\tilde\theta_n$ and $\tilde\theta_n^*$, respectively. It remains to prove the distributional proximity between the ``exact'' versions and the ``two-step-Newton'' version, which was originally provided by Theorem 1 in \cite{andrews2002higher} for the Euclidean Newton updates. Nonetheless, their arguments solely rely on the convergence rate of the approximation instead of in which way it is generated, which can also apply to $\hat \theta_n$ generated in Algorithm~\ref{algorithm: resampled newton iteration}. This is because the assumptions we make in Proposition~\ref{proposition: hypothesis testing} can imply the assumptions in Theorem~\ref{thm:convergence rate of newton} and Theorem~\ref{thm: convergence rate of resampled newton}, and the convergence rate provided in these theorems is sufficient for proving the distributional proximity. Consequently, by Lemma~\ref{lemma:deltamethod}, we conclude the claimed results by a similar argument to the one in Theorem 2 of \cite{andrews2002higher}. \qed

    \section{Normal Coordinate Transformation}
    \subsection{Double Exponential Mapping}
    In this part, we are going to present one of the essential tools for coordinate analysis, which was first studied in Riemannian geometry and general relativity theory and further developed in \cite{gavrilov2006algebraic,gavrilov2007double} in the sense of high-order expansion. Here we introduce the formula derived from the result in \cite{gavrilov2007double}, which fulfills our third-order requirements.
    \begin{lemma}    \label{lemma: normal coordinate transomation}
	Let $(\mc M,\nabla)$ be a smooth manifold with a smooth, symmetric connection. Given $x \in \mc M$, we let $U\subseteq \mc M$ be a geodesically complete neighborhood of $x$. Then for $v, w\in T_x\mc M$ with $\norm{v},\norm{w}\leq \epsilon$, it has 
	\longeq{
	& \Log_x\big(\Exp_{\Exp_x(v)}(T_{x\rightarrow \Exp_x{v}}(w))\big)\\ 
	=&  v + w + \frac{1}{6} R(w, v)v + \frac{1}{3}R(w,v)w + \frac{1}{12}\nabla_v R(w,v)v + \frac{1}{24}\nabla_w R(w, v) v \\+ &  \frac{5}{24}\nabla_v R(w, v)w + \frac{1}{12} \nabla_w R(w,v)w + C\epsilon^5
	}
	for a sufficiently small $\epsilon$ and a sufficiently large constant $C$. 
    \end{lemma}
    \begin{proof}
    	By Equation (3) and the main theorem in \cite{gavrilov2007double}, it can be inferred that the terms on the right hand are exactly the derivatives up to the fourth order. Noticing that $v, w$ are parameters in the system of the second-order ODE, the Picard-Lindel\''of theorem yields that $h(v,w)\coloneqq    \Exp_{\Exp_x(v)}(T_{x\rightarrow \Exp_x{v}}(w))$ is a smooth function with respect to $v,w$. By the integral form of the remainder and the boundness of the fifth-order derivative of $h(v,w)$ around $x$, the term $r(v,w)$ can be replaced by $C\epsilon^5$ with some constant $C$. 

    \end{proof}
    \subsection{Consistency of Derivative Operators}
    This part is aimed at proving the desired consistency between $ \nabla^k f(\hat \theta_n)$ and $\bar \nabla^k f^{\Exp, \theta, \{e_i\}}\big(\pi_{\{e_i\}} \circ \Log(\hat \theta_n)\big)$ with $k=2, 3$. Invoking that $\nabla^k f$ is a $(0,2)$ tensor, we consider their matrix/tensor forms. 
    \begin{lemma}\label{lemma: Hessian consistency}
    	Let $\mathcal M$ be a dimension-$p$ Riemannian manifold and $\{e_i\}$ is a basis of $\mathrm T_x\mc M$. Given $\theta\in \mc M$ and a third-order differentiable function $f: \mc M \rightarrow \bb R^p$, it holds in a neighborhood of $\theta$ that 
    	\longeq{
    	& \norm{ \bar \nabla^2 f^{\Exp, \theta', \{T_{\theta \rightarrow \theta'}(e_i)\}}(\mb 0) - \bar \nabla^2 f^{\Exp, \theta, \{e_i\}}\big(\pi_{\{e_i\}}\circ \Log_\theta(\theta')\big)}_F\\\leq
        &  C \Big(\norm{ \bar \nabla^2 f^{\Exp, \theta, \{e_i\}}\big(\pi_{\{e_i\}}\circ \Log_\theta(\theta')\big)}_F\dist(\theta, \theta')^2 + \norm{\bar \nabla f(\theta')}_2 \dist(\theta, \theta')\Big)
    	} 
    	for some positive constants $C$.
    \end{lemma}
    \begin{proof}
    	Given $\{e_i\} \in\mathrm T_{\theta}\mc M$ and $\theta'$ in a sufficient small neighborhood of $\theta$, we denote the vector bundle induced by the function $\Exp\circ \pi^{-1}_{\{e_i\}}$ at $\theta'$ by $\{e_i(\theta')\}$. To begin with, noting that the exponential mapping is a second-order retraction, Lemma~\ref{lemma:parallel transport error} allows us to relate $T_{\theta\rightarrow \theta'}e_j$ and $e_j(\theta')$ with a second-order error:
    	\longeq{\label{eq: consistency of hessian 1}
    	& \norm{e_j(\theta') - T_{\theta\rightarrow \theta'}e_j} =\norm{\frac{\mathrm d}{\mathrm dx_j}\Exp_\theta\circ \pi^{-1}_{\{e_i\}}(\theta') - T_{\theta\rightarrow \theta' } e_j} \leq C_1\dist(\theta, \theta')^2
    	}
    	with some constant $C_1$ for $j = 1, \dots, p$. Thus, for the Hessian matrices given two different bases, the tensor property of the Hessian operator yields that
        \begin{align}\label{eq: consistency of hessian 3}
            &  \norm{ \bar \nabla^2 f^{\Exp, \theta', \{T_{\theta \rightarrow \theta'}(e_i)\}}(\mb 0) - \bar \nabla^2 f^{\Exp, \theta', \{\frac{\mathrm d}{\mathrm dx_j}\Exp_x\circ \pi^{-1}_{\{e_i\}}(\theta')\}}(\mb 0)}_F\\ 
             \leq &  2 \sqrt{p}C_1 \dist(x, \theta')^2  \norm{\bar \nabla^2 f^{\Exp, \theta', \{\frac{\mathrm d}{\mathrm dx_j}\Exp_x\circ \pi^{-1}_{\{e_i\}}(\theta')\}}  (\mb 0)}_F \\ 
             + & p C_1^2 \dist(x, \theta')^4 \norm{\bar \nabla^2 f^{\Exp, \theta', \{\frac{\mathrm d}{\mathrm dx_j}\Exp_x\circ \pi^{-1}_{\{e_i\}}(\theta')\}}  (\mb 0)}_F \\ 
             \leq&C_2 \dist(x, \theta')^2  \norm{\bar \nabla^2 f^{\Exp, \theta', \{\frac{\mathrm d}{\mathrm dx_j}\Exp_x\circ \pi^{-1}_{\{e_i\}}(\theta')\}}  (\mb 0)}_F
        \end{align}
        for some constant $C_2$.

    	On the other hand, the compatibility of a Riemannian connection yields that 
    	\longeq{\label{eq: consistency of hessian 2}
    	& \frac{\mathrm d}{\mathrm dx_{j_1}}\frac{\mathrm d}{\mathrm dx_{j_2}}f(\Exp_x\circ \pi^{-1}_{\{e_i\}}(y)) =\frac{\mathrm d}{\mathrm dx_{j_1}}\langle \nabla f, e_{j_2}(y)\rangle \\ 
        =&  \langle \nabla_{e_{j_1}(y)}\nabla f, e_{j_2}(y) \rangle + \langle\nabla f(\Exp_x\circ \pi^{-1}_{\{e_{i}\}}(y)), \nabla_{e_{j_1}(y)} e_{j_2}(y) \rangle \\ 
        = & \nabla^2 f(e_{j_1}(y), e_{j_2}(y)) +  \langle\nabla f(\Exp_x\circ \pi^{-1}_{\{e_i\}}(y)), \nabla_{e_{j_1}(y)} e_{j_2}(y) \rangle. 
    	}
        By the property of geodesics, we have $\nabla_{e_i(\theta)} e_j(\theta) = 0$ and thus $\norm{\nabla_{e_i(\theta')}e_j(\theta')} = O(\dist(\theta, \theta'))$ by the smoothness of the connection. Rearranging \eqref{eq: consistency of hessian 2} gives 
    	\eq{
    	\Big\vert\frac{\mathrm d}{\mathrm dx_{j_1}}\frac{\mathrm d}{\mathrm dx_{j_2}}f(\theta') - \bar \nabla^2 f(e_i(\pi_{\{e_i\}} \circ \Log_x (\theta')), e_j(\pi_{\{e_i\}} \circ \Log_x (\theta')))\Big\vert \leq C_3\norm{ \nabla f(\theta')}_2 \dist(\theta, \theta').
    	}
        for some constant $C_3$. 
        Therefore, controlling the entrywise difference of the Hessian matrices yields a bound on the Frobenius norm of the difference:
        \begin{equation}
            \norm{ \bar \nabla^2 f^{\Exp, x, \{e_i\}}\big(\pi_{\{e_i\}}\circ \Log_x(\theta')\big) - \bar \nabla^2 f^{\Exp, \theta', \{\frac{\mathrm d}{\mathrm dx_{j_1}}\Exp_x\circ \pi^{-1}_{\{e_i\}}(\theta')\}}(\mb 0)}_F  \leq C_4\norm{ \nabla f(\theta')}_2 \dist(\theta, \theta') \label{eq: consistency of hessian 4}
        \end{equation}
        holds for some constant $C_4$.

    	Then combining \eqref{eq: consistency of hessian 1} and \eqref{eq: consistency of hessian 2} yields that
        \begin{align}
            & \norm{ \bar \nabla^2 f^{\Exp, \theta', \{T_{x\rightarrow \theta'}(e_i)\}}(\mb 0) - \bar \nabla^2 f^{\Exp, x, \{e_i\}}\big(\pi_{\{e_i\}}\circ \Log_x(\theta')\big)}_F \\ 
            \leq & \norm{\bar \nabla^2 f^{\Exp, \theta', \{\frac{\mathrm d}{\mathrm dx_i}\Exp_x\circ \pi^{-1}_{\{e_i\}}(\theta')\}}(\mb 0)- \bar \nabla^2 f^{\Exp, x, \{e_i\}}\big(\pi_{\{e_i\}}\circ \Log_x(\theta')\big)}_F \\ 
            & + C_2 \dist(x, \theta')^2  \norm{\bar \nabla^2 f^{\Exp, \theta', \{\frac{\mathrm d}{\mathrm dx_j}\Exp_x\circ \pi^{-1}_{\{e_i\}}(\theta')\}}  (\mb 0)}_F\\ 
            \leq & C_2 \dist(x, \theta')^2  \norm{\bar \nabla^2 f^{\Exp, \theta', \{\frac{\mathrm d}{\mathrm dx_j}\Exp_x\circ \pi^{-1}_{\{e_i\}}(\theta')\}}  (\mb 0)}_F + C_4\norm{\nabla f(\theta')}_2 \dist(\theta, \theta')\\ 
            \leq &C \big( \norm{\bar \nabla^2 f^{\Exp, x, \{e_i\}}\big(\pi_{\{e_i\}}\circ \Log_x(\theta')\big)}_F \dist(x, \theta')^2 + \norm{\nabla f(\theta')}_2 \dist(\theta, \theta') \big)
        \end{align}
        holds for some constant $C$ where we make use of \eqref{eq: consistency of hessian 4} again. 
    	
    \end{proof}

    \begin{lemma}\label{lemma: third order tensor consistency}
        Let $\mathcal M$ be a Riemannian manifold. Given $\theta\in \mc M$, a basis $\{e_i\}$ of $\mathrm T_\theta\mc M$, and a third-order differentiable function $f: \mc M \rightarrow \bb R^p$, it holds in a neighborhood of $\theta$ that 
    	\longeq{
    	& \norm{ \bar \nabla^3 f^{\Exp, \theta', \{T_{x\rightarrow \theta'}(e_i)\}}(\mb 0) - \bar \nabla^3 f^{\Exp, x, \{e_i\}}\big(\pi_{\{e_i\}}\circ \Log_x(\theta')\big)}_F \\ 
        \leq & C  \Big( \norm{ \bar \nabla^3 f^{\Exp, x, \{e_i\}}\big(\pi_{\{e_i\}}\circ \Log_x(\theta')\big)}_F \dist(x, \theta')^2 \\ 
        & + \norm{\bar \nabla^2 f^{\Exp, x, \{e_i\}}\big(\pi_{\{e_i\}}\circ \Log_x(\theta')\big)}_F \dist(x, \theta') + \norm{ \nabla f(\theta')} \Big)
        }
        for some constant $C$. 
    \end{lemma}
    \begin{proof}
        Similar to  \eqref{eq: consistency of hessian 1} and \eqref{eq: consistency of hessian 3}, we have 
        \longeq{
              &\norm{ \bar \nabla^3 f^{\Exp, \theta', \{T_{\theta\rightarrow \theta'}(e_i)\}}(\mb 0) - \bar \nabla^3 f^{\Exp, \theta', \{\frac{\mathrm d}{\mathrm dx_i}\Exp_\theta\circ \pi^{-1}_{\{e_i\}}(\theta')\}}(\mb 0)}_F\\
               \leq &C_1 \norm{\bar \nabla^3 f^{\Exp, \theta', \{\frac{\mathrm d}{\mathrm dx_i}\Exp_\theta\circ \pi^{-1}_{\{e_i\}}(\theta')\}}(\mb 0)}_F\dist(\theta, \theta')^2
        }
        for some constant $C_1> 0 $. 

        Recall the definition of $e_k(\theta')$ in the proof of Lemma~\ref{lemma: Hessian consistency}. By the definition of $\nabla^3 f$, we have 
        \longeq{
            & \nabla^3 f[e_i(\theta'), e_j(\theta'), e_k(\theta')] \\ 
            =&  e_i(\theta')\ip{\nabla^2 f(e_j(\theta')), e_k(\theta')} -  \ip{\nabla^2 f(e_j(\theta')), \nabla_{e_i(\theta')} e_k(\theta')} - \ip{\nabla^2 f(\nabla_{e_i(\theta')} e_j(\theta') ), e_k(\theta') }\\ 
            = & e_i(\theta')\big(\frac{\mathrm d}{\mathrm d x_j}\frac{\mathrm d}{\mathrm d x_k}f(\Exp_\theta\circ \pi^{-1}_{\{e_i\}}(y))- \big\langle\nabla f, \nabla_{e_j(\theta')}e_k(\theta')\big\rangle \big)-\ip{\nabla^2 f(e_j(\theta')), \nabla_{e_i(\theta')} e_k(\theta')}\\ 
            &  - \ip{\nabla^2 f(e_k(\theta') ), \nabla_{e_i(\theta')} e_j(\theta') }\\ 
            = &\frac{\mathrm d}{\mathrm d x_i} \frac{\mathrm d}{\mathrm d x_j}\frac{\mathrm d}{\mathrm d x_k}f(\Exp_\theta\circ \pi^{-1}_{\{e_i\}}(y)) - \ip{\nabla f^2(e_i(\theta')), \nabla_{e_j(\theta')}e_k(\theta')} -\ip{\nabla^2 f(e_j(\theta')), \nabla_{e_i(\theta')} e_k(\theta')} \\ & - \ip{\nabla^2 f(e_k(\theta') ), \nabla_{e_i(\theta')} e_j(\theta') }- \ip{\nabla f, \nabla_{e_i(\theta')} \nabla_{e_j(\theta')} e_k(\theta')}.
        }
        Again, for $i,j\in[p]$,  we have $\norm{\nabla_{e_i(\theta')}e_j(\theta')} = O(\dist(\theta, \theta'))$ and $\norm{\nabla_{e_i(\theta')}\nabla_{e_j(\theta')}e_k(\theta')}$ is bounded away from $\infty$ by smoothness. 

        Therefore, invoking Lemma~\ref{lemma: Hessian consistency}, we arrive at the conclusion that
        \begin{align}
            & \norm{ \bar \nabla^3 f^{\Exp, \theta', \{T_{\theta\rightarrow \theta'}(e_i)\}}(\mb 0) - \bar \nabla^3 f^{\Exp, \theta, \{e_i\}}\big(\pi_{\{e_i\}}\circ \Log_\theta(\theta')\big)}_F \\ 
            \leq & C_1 \norm{\bar \nabla^3 f^{\Exp, \theta', \{\frac{\mathrm d}{\mathrm dx_i}\Exp_\theta\circ \pi^{-1}_{\{e_i\}}(\theta')\}}(\mb 0)}_F\dist(\theta, \theta')^2 \\ 
            & + C_2\norm{\bar \nabla^2 f^{\Exp, \theta', \{\frac{\mathrm d}{\mathrm dx_i}\Exp_\theta\circ \pi^{-1}_{\{e_i\}}(\theta')\}}(\mb 0)}_F\dist(\theta, \theta') + C_3\norm{ \nabla f(\theta')}\\ 
            \leq & C \Big( \norm{ \bar \nabla^3 f^{\Exp, \theta, \{e_i\}}\big(\pi_{\{e_i\}}\circ \Log_\theta(\theta')\big)}_F\dist(\theta, \theta')^2 \\ 
            &  + \norm{\bar \nabla^2 f^{\Exp, \theta, \{e_i\}}\big(\pi_{\{e_i\}}\circ \Log_\theta(\theta')\big)}_F\dist(\theta, \theta') + \norm{\nabla f(\theta')} \Big)
        \end{align}
        for some constants $C_2, C_3, C> 0$.
    \end{proof}

        \section{Convergence Rates of (Resampled) Newton Iterations}
        In this section, we will work on the convergence rates of Newton's iteration (Algorithm~\ref{algorithm: update}) with the original data $\mc X_n$ and the resampled data $\mc X_n^{[i]}, i = 1,\cdots, b$, respectively. 
     We also remind that the desired quadratic convergence indeed doesn't rely on the second-order property of retractions; this condition is necessary only for constructing the confidence regions.
        
        \paragraph*{Proof of Theorem~\ref{thm:convergence rate of newton}} 
        We begin by focusing on a neighborhood $N \coloneqq \Exp_{\theta_0}(B(0,\rho))$, ensuring that $N$ is geodescially convex, i.e., for every $p,q$ in $N$ there is a unique minimizing geodesic from $p$ to $q$ lying in $N$, by Theorem~6.17 in \cite{lee2018introduction}. We also assume that $\theta^{\text{initial}} \in \Exp_{\theta_0}(B(0,c_\rho \rho))$ for some $0 <c_\rho <1$. 
        We carry out  the following analysis under the event $\mc F_n$ defined in \eqref{eq: event F} with $\bb P[\mc F_n^{\complement}] =  O(n^{-1}(\log n)^{-2})$. Given the conditions under $\mc F_n$ and the relation \eqref{eq: implicit function of manifold}, it follows that $\tilde \theta_n$ shrinks to $\theta_0$ and thus is contained in $\Exp_{\theta_0}(B(0,c_\rho \rho))$ for sufficiently large $n$. 
        
        Recap that $\hat \theta^{(0)}$ represents the initial estimate $\theta^{\text{initial}}$ and $\theta^{(k)}, k= 1,\ldots, t$ be the updated point after $k$ steps of Newton's update in Algorithm~\ref{algorithm: update}. We aim to demonstrate the targeted convergence rate by inductively establishing that, for $k =0, \cdots, b-1$, 
        \begin{align}
            & \dist(\theta^{(k+1)}, \tilde{\theta}_n) \leq c_l \dist(\theta^{(k+1)}, \tilde{\theta}_n)^2.\label{eq: quadratic convergence rate of newton}
        \end{align}
        holds, given $\theta^{(k+1)} \in \Exp_{\theta_0}(B(0,c_\rho \rho))$.

        Recall that $\tilde \theta_n$ denotes the exact solution of $\nabla	 L_n(\theta) = 0$ which lies closest to $\theta_0$. 
           Applying Lemma~\ref{lemma: taylor expansion of tensor} to $\nabla L_n$ at $\tilde \theta_n$ gives that 
        \longeq{\label{eq: taylor expansion of gradient}
        &0 =  T_{ \tilde \theta_n\rightarrow \theta^{(k)}} \nabla L_n(\tilde \theta_n) = \nabla L_n(\theta^{(k)}) + \nabla^2 L_n(\theta^{(k)})[\Log_{\theta^{(k)} }( \tilde \theta_n)] \\ 
        & + \int_{0}^{1}(1 - t) T_{\Exp_{\theta^{(k)}}(t \Log_{\theta^{(k)} }( \tilde \theta_n)) \rightarrow \theta^{(k)}}\Big(\nabla^3 L_n\big(\Exp_{\theta^{(k)}}(t \Log_{\theta^{(k)} }( \tilde \theta_n))\big)\big[\Log_{\theta^{(k)} }( \tilde \theta_n), \Log_{\theta^{(k)} }( \tilde \theta_n)\big] \Big)\mathrm dt.
        }
        Each term on the right-hand side of equation \eqref{eq: taylor expansion of gradient} requires a delicate analysis.

        For the Hessian term, a basic fact is that \eq{\nabla^2 L_n(\theta^{(k+1)})[T_{\theta_0\rightarrow \theta}e_i, T_{\theta_0\rightarrow \theta}e_j] =\big(\bar \nabla^2 L_n^{\Exp,\hat \theta^{ (k)}, \{T_{\theta_0\rightarrow \theta^{(k+1)}}(e_i )\}}(\mb 0)\big)_{i,j}.
    \label{eq: Hessian equivalence}    
    }
        Further, by Lemma~\ref{lemma: Hessian consistency}, we replace the right hand side of \eqref{eq: Hessian equivalence}  with the Hessian matrix of the function $L_n^{\Exp, \theta_0 , \{e_i\}}$ at $\pi_{\{\pi_i\}} \circ \Log_{\theta_0}(\theta^{(k+1)})$ to obtain that 
        \longeq{
            & \Big\|\bar \nabla^2 L_n^{\Exp,\hat \theta^{ (k)}, \{T_{\theta_0\rightarrow \theta^{(k+1)}}(e_i )\}}(\mb 0) - \bar \nabla^2 L_n^{\Exp, \theta_0 , \{e_i\}}(\pi_{\{\pi_i\}} \circ \Log_{\theta_0}(\theta^{(k+1)}))\Big\|_F \\ 
            \leq&  C\rho^2 \norm{ \bar \nabla^2 L_n^{\Exp, \theta_0 , \{e_i\}}(\pi_{\{\pi_i\}} \circ \Log_{\theta_0}(\theta^{(k+1)}))}_F  + C\rho \norm{\nabla L_n(\theta^{(k+1)})}. \label{eq: application of lemma: Hessian consistency}
        }
        for some constant $C$. 
    
        In order to control $ \norm{ \bar \nabla^2 L_n^{\Exp, \theta_0 , \{e_i\}}(\pi_{\{\pi_i\}} \circ \Log_{\theta_0}(\theta^{(k+1)}))}_F$, we combine \eqref{eq: C(X) bound} and \eqref{eq:derivativebound} to obtain that 
        \longeq{\label{eq: bound of hessian in convergence analysis}
            &  \norm{ \bar \nabla^2 L_n^{\Exp, \theta_0 , \{e_i\}}(\pi_{\{\pi_i\}} \circ \Log_{\theta_0}(\theta^{(k+1)}))}_F \\
            \leq &  \Big(\norm{ \bar \nabla^2 L_n^{\Exp, \theta_0 , \{e_i\}}(\pi_{\{\pi_i\}} \circ \Log_{\theta_0}(\theta^{(k+1)})) - \nabla^2 L_n^{\Exp, \theta_0 , \{e_i\}}(\mb 0)}_F + \norm{\bar \nabla^2 L_n^{\Exp, \theta_0 , \{e_i\}}(\mb 0 )}_F  \Big) \\ 
             \leq&  C_p\rho \big( \bb E[C(X_1)]+ \big|\sum_{i\in[n]} C(X_i) / n - \bb E[C(X_1)]\big| \big) \\
              & + \Big( \norm{ \bar \nabla^2 L_n^{\Exp, \theta_0 , \{e_i\}}(\mb 0) - \bb E[\bar \nabla^2 L^{\Exp, \theta_0 , \{e_i\}}(\mb 0, X_1 )] }_F + \norm{\bb E[\bar \nabla^2 L^{\Exp, \theta_0 , \{e_i\}}(\mb 0, X_1)]}_F \Big)\\ 
              \leq &  C'
        }
        for some constant $C_p,C'$ under the event $\mc F_n$. 
    
        Analogously, for the gradient term $\nabla L_n(\theta^{(k+1)})$ we have 
        \longeq{\label{eq: bound of gradient in convergence analysis}
            &\norm{\nabla L_n(\theta^{(k+1)})}\\ 
             =&  \Big\|\bar \nabla L_n^{\Exp, \theta^{(k+1)},\{T_{\theta_0 \rightarrow \theta^{(k+1)}}e_i\}}(\mb 0)\Big\|_2  \\
            \leq & 2 \Big\|\bar \nabla L_n^{\Exp, \theta_0, \{e_i\}}(\pi_{\{e_i\}} \circ \Log_{\theta_0}(\theta^{(k+1)}))\Big\|_2 \\ 
            \leq& c_p\rho \big( \bb E[C(X_1)]+ \big|\sum_{i\in[n]} C(X_i) / n - \bb E[C(X_1)]\big| \big) + c n^{-\frac{1}{2}}(\log n)^{\frac{1}{2}} \leq c',
        }
        for some constants $c, c', c_p$ under the event $\mc F_n$, where the first inequality follows by Lemma~\ref{lemma:parallel transport error}.
    
        With an appropriate choice of $\rho$, this ensures that the minimum singular value of $\bar \nabla^2 L_n^{\Exp,\hat \theta^{ (k)}, \{T_{\theta_0\rightarrow \theta^{(k+1)}}(e_i )\}}(\mb 0)$ exceeds half of the smallest eigenvalue of $\bb E[\nabla^2 L^{\Exp, \theta_0, \{e_i\}}(\mb 0, X_1)])$, making it invertible with probability at least $1 - O(n^{-1}(\log n)^{-2})$.
    
        For the last term in \eqref{eq: taylor expansion of gradient}, we decompose it using Lemma~\ref{lemma: third order tensor consistency} as follows: 
        \longeq{
            & \sup_{ t\in[0,1]}\sup_{ u_1,  u_2,  u_2 \in \mathrm T_{\Exp_{\theta^{ (k)}}(t \Log_{\theta^{ (k)} }( \tilde \theta_n))}\mc M} \big|\nabla^3 L_n\big(\Exp_{\theta^{ (k)}}(t \Log_{\theta^{ (k)} }( \tilde \theta_n))\big)[u_1, u_2, u_3]\big| \\ 
            \leq & \sup_{t\in[0,1]}\Big\| \bar \nabla^3 L_n^{\Exp, \Exp_{\theta^{ (k)}}(t \Log_{\theta^{ (k)} }( \tilde \theta_n)), \{T_{\theta_0 \rightarrow \Exp_{\theta^{ (k)}}(t \Log_{\theta^{ (k)} }( \tilde \theta_n))}(e_i)\}}(\mb 0) \\ &  - \bar \nabla^3 f^{\Exp, \theta_0 , \{e_i\}}\big(\pi_{\{e_i\}}\circ \Log_{\theta_0} (\Exp_{\theta^{ (k)}}(t \Log_{\theta^{ (k)} }( \tilde \theta_n)))\big)\Big\|_F \\ 
            & +\sup_{  t\in[0,1]} \norm{  \bar \nabla^3 f^{\Exp, \theta_0 , \{e_i\}}\big(\pi_{\{e_i\}}\circ \Log_{\theta_0} (\Exp_{\theta^{ (k)}}(t \Log_{\theta^{ (k)} }( \tilde \theta_n)))\big)}_F \\ 
            \stackrel{\text{(i)}}{\leq} & C_1 \big(\sup_{  t\in[0,1]} \norm{  \bar \nabla^3 f^{\Exp, \theta_0 , \{e_i\}}\big(\pi_{\{e_i\}}\circ \Log_{\theta_0} (\Exp_{\theta^{ (k)}}(t \Log_{\theta^{ (k)} }( \tilde \theta_n)))\big)}_F  + 1 \big) \\ 
            \stackrel{\text{(ii)}}{\leq} & C_2 \label{eq: bound of third time derivatives in convergence analysis}
        }
        for some constants $C_1$ and $C_2$ under the event $\mc F_n$, 
        where (i) holds by Lemma~\ref{lemma: third order tensor consistency}, \eqref{eq: bound of hessian in convergence analysis}, and \eqref{eq: bound of gradient in convergence analysis}, and (ii) holds by arguments similar to \eqref{eq: bound of hessian in convergence analysis}. 

        To the end, we are required to ensure the well-definedness of the geodesics between $\theta^{(k+1)}$ and the other elements in $N$. Notice that under the event $\mc F_n$,  
        \longeq{
            & \norm{\nabla^2 L_n(\theta^{(k)})^{-1} \nabla L_n(\theta^{(k)})}\\ 
            \stackrel{\text{by \eqref{eq: taylor expansion of gradient}}}{ \leq} & \norm{\Log_{\theta^{(k+1)}}(\tilde \theta_n) }\\ 
            & +  \norm{\Log_{\theta^{(k+1)}}(\tilde \theta_n)}^2\sup_{ t\in[0,1]}\Big\| \bar \nabla^3 L_n^{\Exp, \Exp_{\theta^{ (k)}}(t \Log_{\theta^{ (k)} }( \tilde \theta_n)), \{T_{\theta_0 \rightarrow \Exp_{\theta^{ (k)}}(t \Log_{\theta^{ (k)} }( \tilde \theta_n))}(e_i)\}}(\mb 0)\Big\|_F\\ 
            & \cdot  \Big\| \big(\bar \nabla^2 L_n^{\Exp, \Exp_{\theta^{ (k)}}(t \Log_{\theta^{ (k)} }( \tilde \theta_n)), \{T_{\theta_0 \rightarrow \Exp_{\theta^{ (k)}}(t \Log_{\theta^{ (k)} }( \tilde \theta_n))}(e_i)\}}(\mb 0)\big)^{-1}\Big\|  \\ 
            \leq &  C_3\norm{\Log_{\theta^{(k)}}(\tilde \theta_n)}\Big(1 + C_3 \norm{\Log_{\theta^{(k)}}(\tilde \theta_n)}\Big) \leq C_4 \dist(\theta^{(k)}, \tilde\theta_n),\\ 
            & \dist(\hat\theta^{(k)}, \theta^{(k+1)})=\norm{ \Log_{\hat\theta^{(k)}}(\hat\theta^{(k+1)})}  \stackrel{\text{by Definition \ref{definition:second order retraction}}}{\leq }c_R'\norm{\nabla^2 L_n(\theta^{(k)})^{-1} \nabla L_n(\theta^{(k)})}
        }
        hold for some constants $C_3$, $C_4$, and $c_R'$. Therefore, a proper choice of $c_\rho$ and $\rho$ ensures that the distance $\dist(\theta_0, \theta^{(k+1)}) \leq \dist(\theta_0, \tilde\theta_n) + \dist(\tilde\theta_n, \theta^{(k+1)}) + \dist(\theta^{(k+1)}, \theta^{(k+1)}) < \rho$. By the definition of a geodesically convex neighborhood, a geodesic that lies in $N$ and connects $\theta^{(k+1)}$ with a point in $N$ is always unique. 
        
        Now we are prepared to prove the quadratic convergence rate of Newton's iteration. 
        By Lemma~\ref{lemma: normal coordinate transomation} and the well-definedness of logarithmic mapping related to $\theta^{(k+1)}$, there exist some constants $c_t,c_R$ such that 
        \begin{align}
            & \norm{\Log_{\tilde \theta_n}(\theta^{(k+1)})} = \norm{\Log_{\theta^{(k)}}(\tilde \theta_n) - \Log_{\theta^{(k)}}(\theta^{(k+1)})}\\ 
            &  + c_t\big(\norm{\Log_{\theta^{(k+1)}}(\tilde \theta_n)}^2 \norm{ \Log_{\theta^{(k+1)}}(\theta^{(k+1)})} +\norm{\Log_{\theta^{(k+1)}}(\tilde \theta_n)} \norm{ \Log_{\theta^{(k+1)}}(\theta^{(k+1)})}^2  \big) \\ 
            & \stackrel{(i)}{\leq } \norm{\Log_{\theta^{(k)}}(\tilde \theta_n) + \nabla^2 L_n(\theta^{(k+1)})^{-1} \nabla L_n(\theta^{(k+1)})}\\ 
            &  + c_t\Big(\norm{\Log_{\theta^{(k+1)}}(\tilde \theta_n)}^2 \norm{ \Log_{\theta^{(k+1)}}(\theta^{(k+1)})} +\norm{\Log_{\theta^{(k+1)}}(\tilde \theta_n)} \norm{ \Log_{\theta^{(k+1)}}(\theta^{(k+1)})}^2  \Big) \\ 
            & + c_R \norm{\nabla^2 L_n(\theta^{(k+1)})^{-1} \nabla L_n(\theta^{(k+1)})}^3, 
        \end{align}
        where $(i)$ holds by Lemma~\ref{lemma: uniform control for second-order retraction} and \eqref{eq: taylor expansion of gradient}. 
    
        Since the operator $\nabla^2 L_n(\theta^{(k)})$ is invertible with probability at least $1 - O(n^{-1}(\log n)^{-2})$, rearranging \eqref{eq: taylor expansion of gradient} and plugging the above inequalities gives the desired quadratic convergence rate: 
        \longeq{
           & \dist(\tilde\theta_n,\theta^{(k+1)} ) = \norm{\Log_{\tilde \theta_n}(\theta^{(k+1)})} \\ 
            = &  \norm{\int_{0}^{1}(1 - t) T_{\Exp_{\theta^{ (k)}}(t \Log_{\theta^{ (k)} }( \tilde \theta_n)) \rightarrow \theta^{ (k)}}\Big(\nabla^3 L_n\big(\Exp_{\theta^{ (k)}}(t \Log_{\theta^{ (k)} }( \tilde \theta_n))\big)\big(\Log_{\theta^{ (k)} }( \tilde \theta_n), \Log_{\theta^{ (k)} }( \tilde \theta_n)\big)\Big) \mathrm dt} \\ 
            &  + c_t\Big(\norm{\Log_{\theta^{(k)}}(\tilde \theta_n)}^2 \norm{ \Log_{\hat\theta^{(k)}}(\hat\theta^{(k+1)})} +\norm{\Log_{\theta^{(k)}}(\tilde \theta_n)} \norm{ \Log_{\hat\theta^{(k)}}(\hat\theta^{(k+1)})}^2  \Big) \\ 
            & + c_R \norm{\nabla^2 L_n(\theta^{(k)})^{-1} \nabla L_n(\theta^{(k)})}^3 \\ 
            \leq& c_l \dist(\tilde \theta_n, \theta^{(k)})^2\label{eq: convergence rate of newton} \leq c_\rho \rho 
        }
        holds for some constant $c_l$ and sufficiently large $n$, where we make use of the fact that $\norm{\Log_{\theta^{(k)}}(\tilde \theta_n)} = \dist(\theta^{(k)}, \tilde\theta_n)$.

        Finally, the conclusion is a direct consequence of the quadratic convergence rate \eqref{eq: convergence rate of newton} that we have proved. \qed

            \paragraph*{Proof of Theorem~\ref{thm: convergence rate of resampled newton} }
    Since this proof closely follows the structure of Theorem~\ref{thm:convergence rate of newton}, we provide only a sketch of the argument below. Let $\mc X_n^*$ be a generic data collection obtained from resampling $\mc X_n$ with replacement. Throughout the proof, 
    we focus on the event $\mc F_n \cap \mc F_n^*$ defined in \eqref{eq: event F} and \eqref{eq: event F_n^*}. Our goal is to investigate the distance  $\dist(\theta^{(k+1)}, \tilde \theta_n^*)$, provided $\theta^{(k)}\in\Exp_{\theta_0}( B(\theta_0, c_\rho\rho ))$ for $k =0, 1$ and the constants $\rho$ and $c_\rho$ in the proof of Theorem~\ref{thm:convergence rate of newton}.  

        A similar application of Lemma~\ref{lemma: taylor expansion of tensor} yields that 
        \longeq{\label{eq: taylor expansion of gradient 2}
        &0 =  T_{ \tilde \theta_{n}^*\rightarrow \theta^{(k)}} \nabla L_n(\tilde \theta_{n}^*) = \nabla L_n(\theta^{ (k)}) + \nabla^2 L_n(\theta^{ (k)})\big(\Log_{\theta^{ (k)} }( \tilde \theta_{n}^*)\big) \\ 
        & + \int_{0}^{1}(1 - t) T_{\Exp_{\theta^{ (k)}}(t \Log_{\theta^{ (k)} }(\tilde \theta_{n}^*)) \rightarrow \theta^{ (k)}}\Big(\nabla^3 L_n\big(\Exp_{\theta^{ (k)}}(t \Log_{\theta^{ (k)} }(\tilde \theta_{n}^*))\big)\big(\Log_{\theta^{ (k)} }(\tilde \theta_{n}^*), \Log_{\theta^{ (k)} }(\tilde \theta_{n}^*)\big)\Big) \mathrm dt. 
        } 
    
        In terms of the control of $\nabla^2 L_n(\theta^{(k)})$ and $\nabla^3 L_n(\theta^{(k)})$, similar treatments to \eqref{eq: application of lemma: Hessian consistency}, \eqref{eq: bound of hessian in convergence analysis}, \eqref{eq: bound of gradient in convergence analysis}, \eqref{eq: bound of third time derivatives in convergence analysis} but replacing the concentrations based on \eqref{eq: C(X) bound} and \eqref{eq:derivativebound} with applying \eqref{eq: resampled C(X) bound} and the conditions in \eqref{eq: event F'} gives that
        \begin{enumerate}
        \item the minimum singular value of $\bar \nabla^2 L_n^{\Exp,\hat \theta^{ (k)}, \{T_{\theta_0\rightarrow \theta^{(k+1)}}(e_i )\}}(\mb 0)$ is greater than half of the smallest singular value of $\bb E[\nabla^2 L^{\Exp, \theta_0, \{e_i\}}(\mb 0, X_1)])$,
            \item $\sup_{ t\in[0,1]}\sup_{ u_1,  u_2,  u_2 \in \mathrm T_{\Exp_{\theta^{ (k)}}(t \Log_{\theta^{ (k)} }( \tilde \theta_n^*))}\mc M} \big|\nabla^3 L_n\big(\Exp_{\theta^{ (k)}}(t \Log_{\theta^{ (k)} }( \tilde \theta_n^*))\big)[u_1, u_2, u_3]\big|$ is upper bounded by some constant. 
            \item $\hat\theta^{(k+1)}\in \Exp_{\theta_0}(B(0,c_\rho \rho))$. 
        \end{enumerate}

        Thus $\norm{\Log_{\tilde\theta_{n}^*}(\theta^{(k+1)})}$ is upper bounded by 
        \begin{align}
            & \dist(\tilde \theta_{n}^*, \theta^{(k+1)}) =  \norm{\Log_{\tilde \theta_{n}^*}(\theta^{(k+1)})} \\ 
            \leq & \norm{\Log_{\theta^{(k+1)}}(\tilde \theta_{n}^*) + \nabla^2 L_n(\theta^{(k+1)})^{-1} \nabla L_n(\theta^{(k+1)})}\\ 
            &  + c_t\Big(\norm{\Log_{\theta^{(k+1)}}(\tilde \theta_{n}^*)}^2 \norm{ \Log_{\hat\theta^{(k)}}(\hat\theta^{(k+1)})} +\norm{\Log_{\theta^{(k+1)}}(\tilde \theta_{n}^*)} \norm{ \Log_{\hat\theta^{(k)}}(\hat\theta^{(k+1)})}^2  \Big) \\ 
            & + c_R \norm{\nabla^2 L_n(\theta^{(k+1)})^{-1} \nabla L_n(\theta^{(k+1)})}^3 \\
             \leq & c_l^*\dist(\tilde\theta_{n}^*, \theta^{(k+1)})^2
        \end{align}
        for some constant $c_l^*$, where we use the facts: 
        \begin{align}
            & \norm{\nabla^2 L_n(\theta^{(k+1)})^{-1} \nabla L_n(\theta^{(k+1)})} \leq \norm{\Log_{\theta^{(k+1)}}(\tilde \theta_{n}^*)}\big(1 + C_2 \norm{\Log_{\theta^{(k+1)}}(\tilde \theta_{n}^*)}\big) \leq C_4^* \dist(\theta^{(k)}, \tilde\theta_{n}^*),\\ 
            & \norm{ \Log_{\hat\theta^{(k)}}(\hat\theta^{(k+1)})}  \leq c_R^*\norm{\nabla^2 L_n(\theta^{(k+1)})^{-1} \nabla L_n(\theta^{(k+1)})} 
        \end{align}
        for some constants $C_4^*$ and $c_R^*$.

       Finally, the conclusion that $\dist(\theta^{(2)}, \tilde\theta_{n}^*)= O(n^{-2}(\log n)^2)$ follows from the quadratic convergence under the event $\mc F_n \cap \mc F_n$ since $\dist(\tilde \theta_n^*, \theta^{\text{initial}}) = O(n^{-\frac12}(\log n)^{\frac12})$.  \qed

\section{Classic Edgeworth Expansion Results}
\label{sec: classic edgeworth expansion}
In this section, we restate the settings considered in \cite{bhattacharya1978validity,bhattacharya1987some,denker1990asymptotic} and provide a version with some slight but necessary modifications. Consider $p$-dimentional i.i.d. observations $\mb X_j$, $j \in[n]$ with common distribution $G$.  Write 
\longeq{
& \mb \mu \coloneqq \mathbb E\mb X_j,\\ 
& \mb V = \big(  v_{rr'}\big)_{r,r' \in [p]} \coloneqq \mathrm{Cov}(\mb X_1). 
}
Let $\mb H(\mb x) = \big( H^{(1)}(\mb x), H^{(2)}(\mb x), \ldots, H^{(k)}(\mb x) \big) $ be a measurable function defined on $\bb R^k$. We consider the asymptotic distribution of $H(\bar{\mb X})$, where 
\longeq{
& \bar{\mb X}_n \coloneqq  \big(\mb X_1 + \cdots + \mb X_n\big) / n . 
}

The assumptions are as follows: 
\begin{itemize}
    \item[($A_1$)]$\bb E\norm{\mb X_1}_2^s < \infty$ for an integer $s\geq 3$. 
    \item [($A_2$)] Derivatives of $H^{(i)}(\mb x)$, $i \in [p]$  of order up to and including $s-1$ exist and are continuous in a neighborhood of $\mb \mu$. 
    \item[($A_3$)] $\bar \nabla   H^{(i)}(\mb \mu) \coloneqq    \big( \partial / \partial x^{(1)} H^{(i)} (\mb x), \partial/ \partial x^{(2)} H^{(i)} (\mb x), \ldots, \partial / \partial x^{(k)} H^{(i)}(\mb x) \big)\big\vert_{\mb x = \mb \mu}$, $ 1\leq i\leq p$, are linearly independent elements of $\bb R^p$. 
    \item[($A_4$)] The distribution of $\mb X_1$ satisfies Cram\`{e}r's condition, namely, $\limsup\limits_{\norm{t}_2 \rightarrow \infty } \big|\bb E\big[e^{i\mb t\t\mb X_1} \big]\big| <1$. 
\end{itemize}

In view of the population, the following result from \cite{denker1990asymptotic} provides us with the consistency guarantee of the population Edgeworth expansion.

\begin{proposition}[Theorem 2.1 in \cite{denker1990asymptotic}]
\label{proposition: thm 3.1 in bhattacharya}
    Under the assumptions ($A_1$) - ($A_4$), one has 
    \longeq{
    & \sup_{B \in \mc B}\Big| \bb P\big[ \sqrt{n}\big(\mb H(\bar{\mb X}_n) - \mb H(\mb \mu) \big) \in B] - \int_B \psi_{s,n}(\mb x)\mathrm d\mb x \Big| = o(n^{-\frac{(s-2)}{2}}),
    }
    for all class $\mc B \in \mc B^p$ (the collection of Borel sets in $\bb R^p$) such that, for some $a > 0$, 
    \eq{
    \sup_{B \in \mc B}  \Phi_{\mb \Sigma}((\partial B)^\epsilon) = O(\epsilon^a), \quad \text{ as $\epsilon \downarrow 0$. }
    }
\end{proposition}
\begin{remark}\label{remark: coefficients in Edgeworth expansions}
    Note that $\psi_{s,n}(\mb x)$ takes the form of $e^{-\frac{\mb x\t \mb \Sigma^{-1} \mb x}{ 2}}$ multiplied by a polynomial of $\mb x$ where the coefficients in the polynomial are determined by the moments of $\mb  X$ and the derivatives of $H$ as discussed in \cite{bhattacharya1978validity}. For the sake of completeness,
    we restate the details as follows: 
    
    We first reduce the target function $\mb H$ to its $(s - 1)$-order Taylor expansion $\mb H_{s-1}$ at $\mb \mu$ by the fact that $\sqrt{n}\big(\mb {H}(\bar{\mb X}_n)  - \mb H_{s-1}(\bar{\mb X}_n)\big) = o_{\bb P}(n^{-\frac{s-2}{2}})$ and the delta method (Lemma~\ref{lemma:deltamethod}). We denote the $\mb\nu $-th cumulant of $\sqrt{n}\big(\mb H_{s-1}(\bar{\mb X}_n) -\mb H_{s-1}(\mb \mu) \big)$ by $ \kappa_{\mb \nu,n}$\footnote{
    Since $H_k$ is a polynomial function, it is straightforward to see that $\kappa_{\mb \nu,n} $ is of the form $\kappa_{\mb \nu, n} = \sum_{k=0}^{(s-2)|\mb \nu|} a_{\mb \nu, k}n^{-\frac{k}{2}}$.
    }, and denote the approximate $\mb \nu$-th cumulant with respect to $\kappa_{\mb \nu,n}$ by
    \begin{align}
        & \tilde\kappa_{\mb \nu,n} = \left\{ 
        \begin{matrix}
         \Sigma_{i,j} + \sum_{l=1}^{s - 2  }a_{\mb \nu, l}n^{-\frac{l}{2}} & \text{ if $\mb \nu = \{i,j\}$; }\\ 
         \sum_{l=1}^{s - 2 }a_{\mb \nu, l}n^{-\frac{l}{2}} & \text{ otherwise}
        \end{matrix}
        \right.,
    \end{align}
after removing the terms related to $n^{-{\frac{k}{2}}}$ for $l \geq s-  1$
    where the coefficients $\{a_{\mb \nu, l}\}$ depend on the cumulants of $\mb X_1$ and the coefficients of $\mb H_k$; in fact, by \cite{james1958moments} we have $\tilde \kappa_{\mb \nu, n} = O(n^{-\frac{|\mb \nu| - 2}{2}})$ for $|\mb \nu| \geq 3$. 
    
    Expanding the charateristic function of $\sqrt{n}\big(\mb H(\bar{\mb X}_n) - \mb H(\mb \mu)\big)$ gives that for each $\mb t$ one has 
    \begin{align}
        & \bb E\big[\exp\big(\sqrt{n}i\mb t\t \big(\mb H(\bar{\mb X}_n) - \mb H(\mb \mu)\big) \big)\big]\\
        -& \exp\big(-\frac{\mb t\t \mb \Sigma \mb t}{2}\big) \exp\Big(i\sum_{i\in[k]}t_i\tilde \kappa_{\{i\},n} + \frac{-1}{2}\sum_{i,j\in[k]}(\tilde \kappa_{\{i,j\},n} - \Sigma_{i,j}) t_it_j + \sum_{|\mb \nu|\leq s}\frac{i^{|\mb \nu|}}{\mb \nu!}\tilde \kappa_{\mb \nu, n}\mb t^{\mb \nu}\Big) \\ 
        = & o(n^{-\frac{s- 2}{2}}), 
    \end{align}
    where $\mb \Sigma = (\Sigma_{i,j})_{i,j\in[p]}$. 
    Then expanding the exponential term $\exp\Big(i\sum_{i\in[k]}t_i\tilde \kappa_{\{i\},n} + \frac{-1}{2}\sum_{i,j\in[k]}(\tilde \kappa_{\{i,j\},n} - \Sigma_{i,j}) t_it_j + \sum_{|\mb \nu|\leq s}\frac{1}{\mb \nu!}\tilde \kappa_{\mb \nu, n}\mb t^{\mb \nu}\Big) $ implies the existence of polynomials $\{\pi_i\}_{i=1}^{s-2}$ whose coefficients do no depend on $n$ but depend on the coefficients $\{a_{\mb \nu, l}\}$ in $\{\tilde \kappa_{\mb \nu, n}\}$ such that 
    \begin{align}
        \bb E\big[\exp\big(\sqrt{n}i\mb t\t \big(\mb H(\bar{\mb X}_n) - \mb H(\mb \mu)\big) \big)\big] = \exp\big(-\frac{\mb t\t \mb \Sigma \mb t}{2}\big)\Big(1+ \sum_{i = 1}^{s-2}n^{-\frac{i}{2}}\pi_i(i\mb t) \Big)  +  o(n^{-\frac{s- 2}{2}}).
    \end{align}
        In the end, we define the \emph{formal edgeworth expansion} $\psi_{s,n}$ as 
    \begin{align}
        & \psi_{s,n}(x) = \Big(1+ \sum_{i = 1}^{s-2}n^{-\frac{i}{2}}\pi_i(-\frac{\mathrm d}{\mathrm d\mb x}) \Big)\exp\big(-\frac{\mb x\t \mb\Sigma^{-1}\mb x }{2}\big). 
    \end{align}
\end{remark}

Regarding the bootstrap procedure, we shall provide a modified version of Theorem 5.3 in \cite{bhattacharya1987some} to accommodate the impact of coordinate chart adjustment. We first consider a continuous function $\tilde {\mb H}: \bb R^k \times \bb R^k \rightarrow \bb R^p$ satisfying that
\begin{enumerate}
\item $\tilde {\mb H}(\mb x, \mb \mu) = \mb H(\mb x)$, and $\tilde {\mb H}(\cdot, \mb x')$ is $(s-1)$-times continuously differentiable on $O$ for every $\mb  x' \in O$ where $O$ is an open neighborhood of $\mb \mu$. 
\item The first $(s-1)$-times derivatives of $\tilde{\mb  H}(\mb x_1, \mb x_2)$ with respect to the first variable $\mb x_1 \in O$ are continuous with respect to $\mb x_2\in O$. 
\end{enumerate}
\begin{proposition}[A Modified Version of Theorem~5.3 in \cite{bhattacharya1987some}]
\label{proposition: modified thm 3.3 in bhattacharya}
    If Assumption~($A_1$) is strengthened to  $\bb E\norm{\mb X_1}_2^{2s} < \infty$ for some integer $s\geq 3$, and the assumptions ($A_2$) - ($A_4$) are satisfied, then one has 
    \eq{
     \limsup\limits_{n\rightarrow \infty} n\bb P_{\mc X_n}\Big[\sup_{B \in \mc B}\Big| \bb P\big[\sqrt{n} \big(\tilde{\mb H}(\bar{\mb X}_n^*, \bar{\mb X}_n) - \tilde{\mb H}(\bar{\mb X}_n, \bar{\mb X}_n) \big) \in B\vert \mc X_n\big] - \int_B \tilde \psi_{s,n}(\mb x) \mathrm d \mb x\Big| \geq cn^{-\frac{(s-1)}{2}}\Big] < \infty    \label{eq: modified resampled edgeworth expansion}
    }
    for every constant $c>0$ as $n\rightarrow \infty$, for every class $\mc B$ of Borel sets satisfying, for some $a>0$, 
    \longeq{
    & \sup_{B\in \mc B}\Psi_{\mb \Sigma}\big( (\partial B)^\epsilon \big) = O( \epsilon^a), \text{\quad as $\epsilon \downarrow 0$. }\label{eq: set family}
    }
    Here the coefficients of the polynomial in $\tilde{\psi}_{s,n}$ depends on the empirical moments of $\{\mb X_i\}_{i \in [n]}$ as well as the derivatives up to $s-1$ times of $\tilde{\mb H}(\cdot, \bar{\mb X}_n)$ evaluated at $\bar{\mb X}_n$. Moreover, the coefficients of $\tilde \psi_{s,n}$ converge to their population counterparts at a rate of $O(n^{-\frac12}(\log n)^{\frac12})$ with probability at least $1- O(n^{-\frac{s-2}{2}}(\log n)^{-\frac{s}{2}})$. 
\end{proposition} 

\begin{proof}
    Following the proof of Theorem~5.2 in \cite{bhattacharya1987some} together with moment control similar to \eqref{eq: resampled Z moment bound}, we have 
    \eq{
         \limsup\limits_{n\rightarrow \infty} n\bb P_{\mc X_n} \Big[\sup_{B \in \mc B}\Big| \bb P\big[\sqrt{n}(\bar{\mb X}_n^* - \bar{\mb X}_n) \in B\mid \mc X_n \big] - \int_B \big[ 1 + \sum_{r=1}^{s-2}n^{-\frac{r}{2}}\hat\psi_r(\mb x)\big] \phi_{\bar {\mb V}}(\mb x)\mathrm d \mb x \Big| \geq cn^{-\frac{s-1}{2}}\Big] <\infty    \label{eq: 3.18 in bhattacharya1987some}
    }
    for every constant $c>0$ as $n\rightarrow \infty$,
    for every class $\mc B$ satisfying, for some $a > 0$, $\sup_{B \in \mc B} \Psi_{\mb V}((\partial B)^\epsilon) = O( \epsilon^a)$, as $\epsilon \downarrow 0$,  where $\{\hat \phi_r\}$ are polynomials that do not depend on $n$ but depend on the moments of $\mb X_1^*$, and $\bar{\mb V}$ denotes the sample covariance matrix of $\{\mb X_i\}_{i=1}^n$. Then the reduction step from \eqref{eq: 3.18 in bhattacharya1987some} to \eqref{eq: modified resampled edgeworth expansion} follows by the same arguments as the one in \cite{bhattacharya1978validity}, with an observation from Remark~\ref{remark: coefficients in Edgeworth expansions} that the coefficients in the formal Edgeworth expansion only depends on the moments of $\mb X_1^*$ conditional on $\mc X_n$ and the derivatives up to $s-1$ times of $\tilde{\mb H}(\cdot, \bar{\mb X}_n)|_{\bar{\mb X}_n}$. The second part follows from the continuity of $\tilde{\mb H}$ and the moderate deviation (Lemma~\ref{lemma: thm1 in von1967central}). 
\end{proof}

Now we turn back to the preceding analysis in Section~\ref{section: implicit function 1} and Section~\ref{section: implicit function 2}. One subtle aspect of this analysis involves approximating $\tilde{\mb \eta}_n^* - \tilde{\mb \eta}$ by a smooth function of the derivatives of $L_n^{*\Exp, \theta_0, \{e_i\}}$ evaluated at $\tilde{\mb \eta}_n$. To leverage the assumptions we made on the derivatives of $L$ at $\theta_0$, we shall utilize the next proposition to facilitate the replacement of derivatives. 
\begin{proposition}\label{proposition: equivalence of two resampled distribution}
    Instate the assumptions in Theorem~\ref{theorem:bootstrapmanifold} and consider $\bar{\mb Z}_n$ in \eqref{eq: definition of Zbar_n}, its resampled counterpart $\bar{\mb Z}_n^*$, and $\bar{\mb Y}_n$ defined in \eqref{eq: approximation error of H(Y)}. Given a function $\tilde{\mb H}$ as introduced in Proposition~\ref{proposition: modified thm 3.3 in bhattacharya},  it holds almost surely that 
    \longeq{
    & \limsup\limits_{n\rightarrow \infty} n\bb P_{\mc X_n}\Big[\bb P\big[ \big\|\big(\tilde{\mb H}(\bar{\mb Z}_n^*, \bar{\mb Z}_n) - \tilde{\mb H}(\bar{\mb Z}_n, \bar{\mb Z}_n) \big) - \big( \tilde{\mb H}(\bar{\mb  Y}_n, \bar{\mb Z}_n) -\tilde{\mb H}(\bb E[\bar{\mb Y}_n\mid \mc X_n],  \bar{\mb Z}_n)\big) \big\|_2 
    \\ 
    & \hspace{7cm} \leq  cn^{-1}\log n \mid \mc X_n\big]\leq 1- cn^{-1}(\log n)^{-2}\Big] < \infty
    }
    for some constant $c$. 
\end{proposition}
\begin{proof}
    We shall discuss a univariate $\tilde H$ since the proof for a multivariate $\tilde{\mb H}$ simply follows a similar argument. We denote the (random) function $H(\cdot, \bar{\mb Z}_n)$ by $g(\cdot)$ herein. 
    Expanding the function $g$ at $\bar{\mb Z}_n$ and at $\bar{\mb Y}_n$ respectively yields that 
    \longeq{\label{eq: g(Z) and g(Y)}
        &  \big(g(\bar{\mb Z}_n^*) - g(\bar{\mb Z}_n) \big) - \big( g(\bar{\mb  Y}_n) -g(\bb E[\bar{\mb Y}_n\mid \mc X_n] \big)\\
        =& (\bar \nabla g(\bar{\mb Z}_n))\t \big(\bar{\mb Z}_n^* - \bar{\mb Z}_n\big) -
        (\bar \nabla g(\bb E[\bar{\mb Y}_n\mid \mc X_n]))\t \big(\bar{\mb Y}_n - \bb E[\bar{\mb Y}_n\mid \mc X_n]\big) + \Xi_1 - \Xi_2 \\ 
        = & (\bar \nabla g(\bar{\mb Z}_n))\t \Big(\big(\bar{\mb Z}_n^* - \bar{\mb Z}_n\big) -
        \big(\bar{\mb Y}_n - \bb E[\bar{\mb Y}_n\mid \mc X_n]\big) \Big) + \Xi_1 - \Xi_2 \\ 
        & + \big(\bar \nabla g(\bar{\mb Z}_n) - \bar \nabla g(\bb E[\bar{\mb Y}_n\mid \mc X_n]) \big)\big(\bar{\mb Y}_n - \bb E[\bar{\mb Y}_n\mid \mc X_n]\big)
    }
    where $\Xi_1, \Xi_2$ are the residual terms of the Taylor expansions. By \eqref{eq: exceptional probability of F_n'} and similar reasoning in its proof, it holds that 
    \begin{align*}
    & \limsup\limits_{n\rightarrow \infty} n\bb P_{\mc X_n}\Big[\bb P\big[\norm{\bar{\mb Y}_n - \bb E[\bar{\mb Y}_n\mid \mc X_n] }_2 \leq c_1(\log n)^{\frac{1}{2}} n^{-\frac{1}{2}}, \quad \norm{\bar{ \mb Z}_n^* - \bar{\mb Z}_n }_2 \leq c_1(\log n)^{\frac{1}{2}} n^{-\frac{1}{2}}, \\ 
    & \qquad  |\Xi_1| \leq c_1(\log n) n^{-1}, \quad |\Xi_2| \leq c_1 (\log n) n^{-1}, \\ 
    & \qquad \big(\nabla g(\bar{\mb Z}_n) - \nabla g(\bb E[\bar{\mb Y}_n|\mc X_n]) \big)\big(\bar{\mb Y}_n - \bb E[\bar{\mb Y}_n\mid \mc X_n ]\big) \leq c_1 (\log n) n^{-1} |\mc X_n\big] \leq   1- c_1 n^{-1}(\log n)^{-2}\Big] <+ \infty 
    \end{align*}
    for some constant $c_1$. 

    Then we look into the term $\big(\bar{\mb Z}_n^* - \bar{\mb Z}_n\big) -
        \big(\bar{\mb Y}_n - \bb E[\bar{\mb Y}_n\mid \mc X_n]\big)$. By definition and Taylor's expansion, we have 
        \longeq{
        & \big(\bar{\mb Z}_n^* - \bar{\mb Z}_n\big) -
        \big(\bar{\mb Y}_n - \bb E[\bar{\mb Y}_n\mid \mc X_n]\big) \\ 
        = & \big(\bar{\mb Z}_n^* - \bar{\mb Y}_n\big) -
        \big(\bar{\mb Z}_n - \bb E[\bar{\mb Y}_n\mid \mc X_n]\big) \\ 
        = & -\big(\mb \Pi_1 - \mb \Pi_2 \big)\t \tilde{\mb \eta}_n  + \mb \Xi_3 + \mb \Xi_4,  
        }
        where $\mb \Pi_i, i =1,2$ are the vectors that collect the derivatives of $L_n^{\Exp, \theta_0, \{e_i\}}$ and ${L_n^*}^{\Exp, \theta_0, \{e_i\}}$ up to the fifth time evalated at $\mb 0$, respectively,  
        and $\Xi_3$ and $\Xi_4$ are the residual terms. By a similar moderate deviation argument invoking Lemma~\ref{lemma: thm1 in von1967central}, we arrive at 
        \longeq{\label{eq: bound for Pi1 - Pi2}
        & \limsup\limits_{n\rightarrow \infty} n \bb P_{\mc X_n}\Big[\bb P\big[\norm{\mb \Pi_1 - \mb \Pi_2}_2 \leq  c_2 (\log n)^{\frac12}n^{-\frac12}, \\ 
        & \qquad \qquad  \norm{\mb \Xi_1}_2 \leq c_2 (\log n)^{\frac12}n^{-\frac12}, \norm{\mb \Xi_2}_2 \leq c_2 (\log n)^{\frac12}n^{-\frac12}|\mc X_n \big] \leq 1 - c_2n^{-1}(\log n)^{-2}  \Big] <\infty }
        for some constant $c_2$. Together with the fact that $\norm{\tilde {\mb \eta}_n}_2  = O(n^{-\frac{1}{2}}(\log n)^{\frac{1}{2}})$ under $\mc F_n$, we derive that 
        \longeq{
        & \limsup\limits_{n\rightarrow \infty} n\bb P_{\mc X_n}\Big[\bb P\big[\norm{\big(\bar{\mb Z}_n^* - \bar{\mb Z}_n\big) -
        \big(\bar{\mb Y}_n - \bb E[\bar{\mb Y}_n\mid \mc X_n]\big)}_2 \leq c_9(\log n)n^{-1}|\mc X_n \big] \leq 1-c_3n^{-1}(\log n)^{-2}\Big] < \infty.
        }
        
        Plugging these bounds into \eqref{eq: g(Z) and g(Y)} concludes that 
        \longeq{
             &\limsup\limits_{n\rightarrow \infty} n \bb P_{\mc X_n}\Big[\bb P\big[\big\vert \big(g(\bar{\mb Z_n}^*) - g(\bar{\mb Z}_n) \big) - \big( g(\bar{\mb  Y}_n) -g(\bb E[\bar{\mb Y}_n\mid \mc X_n] \big)  \big\vert \leq  c(\log n)n^{-1} \mid \mc X_n\big]  \\ 
             & \qquad \qquad \qquad \qquad \leq 1 - cn^{-1}(\log n)^{-2} \Big] < \infty    
        }
        holds for some constant $c$.  
\end{proof}

\section{Auxiliary Lemmas}

An elementary tool in this paper is called the delta method in the distribution theory, which allows us to neglect small terms in expansions. For the sake of completeness, we restate Lemma 5 in \cite{andrews2002higher} in the following. 
\begin{lemma}[Lemma 5 in \cite{andrews2002higher}]\label{lemma:deltamethod}
\begin{enumerate}
    \item     Let $\{A_n\}\in \bb R^{L_A}$ be a sequence of random vectors with an Edgeworth expansion with coefficients of order $O(1)$ satisfying that 
    \eq{
        \lim_{n\rightarrow \infty} n^a\sup_{\mb x \in\bb R^{L_A}}  \Big|\bb P\big[(A_n)_i \leq x_i, i \in [L_A]\big] - \int_{z_i \leq x_i}\big[ 1 + \sum_{i=1}^{2a}n^{-\frac{i}{2}} \pi_i(\partial/\partial z)\big]\phi_{\mb \Sigma_n}(\mb z)\mathrm d \mb z\Big| = 0
    }
    for some $a$ with $2a$ being a positive integer, 
where $\phi_{\mb \Sigma_n}$ denotes the density function of $\mc N(\mb 0, \mb \Sigma_n)$ and the eigenvalues of $\mb \Sigma_n$'s are bounded away from $0$ and $\infty$. 
    Let $\{\xi_N\}\in \bb R^{L_A}$ be a sequence of random vectors with $\bb P[\norm{\xi_n}_2 > \varrho_n ]  = o(n^{-b})$ for $\varrho_n = o(n^{-b})$ and some $0< b \leq a$. Then 
    \eq{
    \lim_{n\rightarrow \infty} \sup_{x\in \bb R}n^b\left|\bb P[A_n + \xi_n\leq z] - \bb P[A_n \leq z] \right| = 0.
    }
    \item Let $\{A_n^*, n \geq 1\}\in \bb R^{L_A}$ be a sequence of bootstrap random vectors that possesses
    \begin{equation}
        \lim_{n\rightarrow}n^{a} \bb P_{\mc X_n}\Big[n^a  \big|\bb P[A_n^*\leq z\mid \mc X_n] - \big[ 1 + \sum_{i=1}^{2a}n^{-\frac{i}{2}} \pi_i^*(\partial/\partial z)\big]\Phi_{\Sigma^*_n}(z)\big|\Big] =0 
    \end{equation}
    for some $a$ with $2a$ being a positive integer,
    where the coefficients of polynomials $\pi_i^*(\partial/\partial z)$ are of order $O(1)$ and the smallest eigenvalue of $\Sigma_n^*$ is bounded away from $0$ and $\infty$, except in a sequence of sets with probability $o(N^{-a})$. Let $\{\xi_n^*,n\geq1\}\in \bb R^{L_a}$ be a sequence of random vectors with $\lim_{n\rightarrow} n^b \bb P\big[ \bb P^*[\norm{\xi_n}_2 > \varrho_n] > n^{-b}\big] = 0$ for some sequence $\varrho_n = o(n^{-b})$ with $b \leq a$. Then 
    \eq{
    \lim_{n\rightarrow \infty} n^b\bb P_{\mc X_n}\Big[\sup_{z\in R^{L_A}} \big|\bb P[A_n^* + \xi_n^* \leq z|\mc X_n] - \bb P[A_n^* \leq z |\mc X_n] \big| > n^{-b}\Big] = 0. 
    }
\end{enumerate}
\end{lemma}

    \begin{lemma}[Uniform Consistency of a Second-Order Retraction]
        Let $\mc M$ be a dimension-$p$ Riemannian manifold and $R\coloneqq \{R_\theta\}$ be a second-order retraction defined in Definition~\ref{definition:second order retraction}. Then, given an arbitrary $\theta_0 \in \mc M$, there exists a neighborhood $U$ of $\mb 0 \in \bb R^m$ such that the following uniformly holds for every $\theta,\theta' \in U$ and $\mb y \in \bb R^p$ with $\norm{\mb y} \leq c$ that 
        \begin{align}
            & \norm{\Log_{\theta}(R_\theta(\mb y)) - \mathrm d R_\theta (\mb y)} \vee \dist(R_{\theta}(\mb y), \Exp_{\theta}(\mathrm d R_{\theta}(\mb y))) \leq C\norm{\mb y}_2^3, \\ 
            & \norm{R^{-1}_\theta(\theta') - \mathrm d R_\theta^{-1} \Log_{\theta}(\theta')}_2  \leq C \dist(\theta, \theta')^3 
        \end{align}
        for some constant $C$. 
        \label{lemma: uniform control for second-order retraction}
    \end{lemma}
    \begin{proof}

        By Lemma~6.14 in \cite{lee2018introduction}, there exists a neighborhood $U_0$ of $\theta_0$ such that it is contained in a geodesic ball of radius $\epsilon$ around each of its points. Moreover, from the proof of the lemma, there exists a diffeomorphism $\phi$ from $\pi^{-1}(U_0) \subset \mathrm T \mc M$ to $\mc U_0 \times \bb R^p$ ($\pi$ denotes the canonical projection), and 
        the mapping $E(\theta, v) = (\theta, \Exp_\theta(v))$ is diffeomorphic on $\phi^{-1}(U_0 \times \{\mb x: \norm{\mb x}_2 \leq \epsilon'\})\subset \pi^{-1}(U_0)$ for some $\epsilon'$. Let $U_1 \subset U_0$ be a precompact neighborhood containing $\theta_0$. 

        Then, one has uniformly for $\theta \in U_0$ and $\norm{v} \leq c_1 \epsilon'$ with some constant $c$ that 
        \begin{align}
           &\dist(R_{\theta}(\mb y), \Exp_{\theta}(\mathrm d R_{\theta}(\mb y))) \leq \int_{0}^1 \norm{\frac{\mathrm d }{\mathrm d t} \Exp_{\theta}\big((1-t)\mathrm d R_{\theta}(\mb y) + t \Log_{\theta}(R_\theta(\mb y)) \big)} \mathrm d t \\ 
           \leq&  C_1 \norm{\Log_{\theta}(R_\theta(\mb y)) - \mathrm d R_{\theta}(\mb y)}\leq C_2 \norm{\mb y}_2^3
        \end{align}
        for some constant $C_1$ and $C_2$. Here, the first inequality follows from the fact that the geodesic provides the shortest connecting curve. The second uses the precompactness of $U_1$ together with the smoothness of $E$. And the last inequality follows from the second-order retraction property, the smoothness of $\Log_\theta\circ R_\theta \circ \mathrm dR^{-1}_\theta$, and the Taylor expansion that
        \begin{align}
            & \Log_{\theta}(R_\theta(\mb y)) - \mathrm d R_{\theta}(\mb y)\\ 
            = &\frac{\mathrm d}{\mathrm d t}\big[\Log_{\theta}(R_\theta(t\mb y)) - \mathrm d R_{\theta}(t\mb y) \big]\Big\vert_{t=0} + \frac{1}{2}\frac{\mathrm d^2}{\mathrm d t^2}\big[\Log_{\theta}(R_\theta(t\mb y)) - \mathrm d R_{\theta}(t\mb y) \big]\Big\vert_{t=0} \\ 
            & + \frac{1}{6}\int_0^1 \frac{\mathrm d^3}{\mathrm d t^3}\big[\Log_{\theta}(R_\theta(t\mb y)) - \mathrm d R_{\theta}(t\mb y) \big]\Big\vert_{t=s}\mathrm d s \\ 
            = & \frac{1}{6}\int_0^1 \frac{\mathrm d^3}{\mathrm d t^3}\Log_{\theta}(R_\theta(t\mb y)) \Big\vert_{t=s}\mathrm d s. 
        \end{align}

        On the other hand, applying Taylor's expansion to $\mathrm d R_\theta R^{-1}_\theta(\theta') -  \Log_{\theta}(\theta')$ that 
        \begin{align}
            & \mathrm d R_\theta R^{-1}_\theta(\theta') -  \Log_{\theta}(\theta') = \frac{\mathrm d }{\mathrm dt } \big[ \mathrm d R_\theta R^{-1}_\theta(\Exp_\theta(t\Log_\theta(\theta'))) - t \Log_{\theta}(\theta') \big] \Big\vert_{t = 0} \\ 
            & + \frac{1}{2}\frac{\mathrm d^2}{\mathrm dt^2}\big[ \mathrm d R_\theta R^{-1}_\theta(\Exp_\theta(t\Log_\theta(\theta'))) - t \Log_{\theta}(\theta') \big]\Big\vert_{t = 0} \\ 
            & + \frac{1}{6}\int_0^1\frac{\mathrm d^3}{\mathrm dt^3} \big[ \mathrm d R_\theta R^{-1}_\theta(\Exp_\theta(t\Log_\theta(\theta'))) - t \Log_{\theta}(\theta') \big] \Big\vert_{t = s } \mathrm d s \\ 
            = & \frac{1}{6}\int_0^1\frac{\mathrm d^3}{\mathrm dt^3} \mathrm d R_\theta R^{-1}_\theta(\Exp_\theta(t\Log_\theta(\theta'))) \big\vert_{t = s } \mathrm d s. 
            \label{eq: dRR-1 expansion}
        \end{align}
         Here, since the second-order retraction property ensures that $\nabla \mathrm d R \vert_0 = 0$, we have $\nabla \mathrm d R^{-1}|_\theta = 0$. Thus, $\frac{\mathrm d^2}{\mathrm dt^2}\mathrm d R_\theta R^{-1}_\theta(\Exp_\theta(t\Log_\theta(\theta')))\big\vert_{t=0}=0$. 
         
         Using \eqref{eq: dRR-1 expansion}, along with the smoothness of $R_\theta \circ \mathrm d R^{-1}_\theta$ and the precompactness of $U_1$, we obtain that $\norm{R^{-1}_\theta(\theta') - \mathrm d R_\theta^{-1} \Log_{\theta}(\theta')}_2  \leq C_3 \dist(\theta, \theta')^3 $ for some constant $C_3$. 
         
    \end{proof}

\begin{lemma}[Taylor Expansion of a Real-valued Function]
\label{lemma: taylor expansion}
	Let $\mathcal M$ be a dimension-$p$ Riemannian manifold and $f: \mc M \rightarrow \bb R$ is a $k$-times differentiable function on $\mc M$. Then for $\theta \in \mc M$, there exists a positive $\rho$ such that for arbitrary $\theta \in \Exp_{\mb \theta}(B(0,\rho))$, $\Log_{\theta}$  defined on ${\Exp_{\mb \theta_0}(B(0,\rho))}$ is a differmorphisim onto its codomain in $\mathrm T_{\theta}\mc M$. Moreover, given a fixed normal basis of $\mathrm T_\theta \mc M$ and $\theta' \in B(\theta, \rho)$ it holds that 
    \longeq{
    & f(\theta') - f(\theta) = \sum_{i=1}^{k-1}\sum_{\mb \nu: |\mb \nu| =i}\frac{1}{\mb \nu!}\nabla^{\mb \nu} f\vert_{x = \theta}[h^{\mb \nu}] + \frac{1}{k!}\frac{\mathrm d}{\mathrm d t}f\circ \Exp(t\Log_{\theta}(\theta'))\big\vert_{t = t_0}\norm{\Log_\theta(\theta')}^k,
    }
    where $t_0 \in [0, \norm{\Log_\theta(\theta')}]$. 
\end{lemma}
\begin{proof}
    The first part is a well-known conclusion in Riemannian geometry. And the Taylor expansion is an immediate consequence of the Taylor expansion of $f\circ \Exp\circ \pi^{-1}_{\{e_i\}}$ on $\bb R^p$. 
\end{proof}

\begin{lemma}[Taylor Expansion of a $0$-$1$ Tensor]\label{lemma: taylor expansion of tensor}
    Let $V$ be a $0$-$1$ tensor on a Riemannian manifold $\mc M$ and $c:[a,b]\rightarrow\mc M$ be a geodesic curve connecting $c(a)$ and $c(b)$. Then it holds that 
    \begin{equation}
        T_{c(b) \rightarrow  c(a)} V(c(b)) = V(c(a)) + \nabla V(c(a))(\dot c(t)) + \int_a^b (1 - t)T_{c(b) \rightarrow  c(a)}\big(\nabla^2 T(c(t))[ \dot c(t), \dot c(t)]\big)  \mathrm d t.
    \end{equation}
\end{lemma}

\begin{proof}
    The proof is a direct consequence of the Taylor expansion of $V\circ c$ at $a$ and the property of parallel transport. 
\end{proof}

    \begin{lemma}
    \label{lemma:parallel transport error}
    Let $\mc M$ be a Riemannian manifold of dimension $p$, $R$ be a second-order retraction which is defined on $U \subset \mathbb R^p$ centered at $x$, and $x'$ be a point satisfying $x' \in R_{x}(U)$ with $\dist(x,x') = \epsilon $. Then given the canonical basis $\{\delta_i\}$ of $\bb R^p$, it holds that
        \eq{
        \norm{\mathrm dR_{x}( R_x^{-1}(x'))\delta_i - T_{x\rightarrow x'}(\mathrm d R_x (0)\delta_i)} \leq O(\epsilon^2)
        ,\quad  \forall i \in [p]. 
        }
    \end{lemma}
    \begin{proof}
        Given a sufficiently small $\epsilon>0$, we ssume that $\gamma(t), t \in [0,1]$ is a unique shortest constant-speed geodesic between $x$ and $x'$ and $\{\tilde e_i(t)\} \subset  \mathrm T_{\gamma(t)}\mc M$ is the basis derived from the parallel transport of $\{\mathrm dR(0)\delta_i \}$ along $\gamma$. We decompose $\mathrm dR_a(\gamma(t))\delta_i$ as $\mathrm dR_x(\gamma(t))[\delta_i] = \sum_{j=1}^pa_{i,j}(t)\tilde e_j$. Then we have 
        \eq{
        \norm{\mathrm dR_{x}(R_x^{-1}(x'))\delta_i - T_{x\rightarrow x'}(\mathrm d R_x(0)\delta_i)}^2  \leq \big( |a_{i,i}(1)-1|^2 + \sum_{j\neq i}^p |a_{i,j}(1)|^2\big).
        }

        On the other hand, for $v \coloneqq    \Log_x(x') \in  \mathrm T_x\mc M$, the property of second-order retractions yields that 
        \longeq{
        &\nabla_v (\mathrm d R_x[\delta_i])(0) = 0 = \sum_{j=1}^p (v  a_j)\tilde e_j(0) 
        }
        which leads to the fact that $va_i = 0$ for $i = 1,\ldots, p$. 

        Therefore, it follows by the Taylor expansion and the smoothness of $R$ that 
        \longeq{
        & |a_1-1| \leq \frac{c\epsilon^2}{2}\\ 
        & |a_i| \leq \frac{c\epsilon^2}{2}, \qquad \forall i = 2, \ldots, p
        }
        with a constant $c$ in a local neighborhood of $x$ which leads to the conclusion. 
    \end{proof}

    For completeness, we present the classical results on moderate deviation and moment control in the following. 
    \begin{lemma}[Theorem~1 in \cite{von1967central}]
        \label{lemma: thm1 in von1967central}
            Suppose that we have $n$ i.i.d. random vector $\{\mb X_i\}_{i\in[n]}$ of dimension $k$. If $\bb E\big[ \norm{\mb X_1}_2^r \big]$ holds for an integer $r \geq 3$, and $m$ is the largest singular value of the covariance matrix of $\mb X_1$, then 
            \begin{equation}
                \Big| \bb P\big[\|\sum_{i\in[n]}\mb X_i / \sqrt{n}\|_2 \leq a\big] - \int_{\norm{\mb x}_2 \leq a}\mathrm d \Phi(\mb x)\Big| \leq C\cdot a^{-r}\cdot n^{-\frac{r-2}{2}}
        \end{equation}
        holds for some constant $C$ and every $a\geq (\frac54 m(r-2)\log n)^{\frac12}$, 
        where $\Phi$ denotes the standard $k$-dimensional normal distribution. 
        \end{lemma}
	\begin{lemma}[Theorem 2.10 in \cite{hall2014martingale}]
	\label{lemma: Burkholder's inequality}
		If $\{S_i, \mathcal F_i,1\leq i \leq n \}$ is a martingale and $1< p < \infty$, $X_1 = S_1$, and $X_i = S_i - S_{i-1}$, $2\leq i\leq n$, then there exist constants $C_1$ and $C_2$ depending only on $p$ such that 
		\longeq{
		& C_1\bb E\big| \sum_{i\in[n]} X_i^2 \big|^{\frac{p}{2}} \leq \bb E|S_n |^p \leq C_2 \bb E\big| \sum_{i\in[n]} X_i^2\big|^{\frac{p}{2}}.
		}
	\end{lemma}
    
        Lastly, we provide a proof of Proposition~\ref{proposition:inverse retraction}. 
        \paragraph*{Proof of Proposition~\ref{proposition:inverse retraction}}
        This proposition is a consequence of Lemma~20 and Theorem~22 in \cite{absil2012projection}. To be specific, we consider a function $R: \mathrm T\mc M \rightarrow \mathrm M$ that maps $(x_0, u) \in \mathrm T\mc M$ to $x_0 + u + v$, where $v$ is a vector in $\mathrm N_{x_0}\mc M$ of the shortest length ensuring $x_0 + u + v\in \mc M$. By Lemma~20 in \cite{absil2012projection}, the existence and uniqueness of $v$ in a sufficiently small neighborthood of $(x_0, 0_{x_0})\in \mathrm T\mc M$ is guaranteed, and $\mathrm d R(x_0, \cdot)|_{0_{x_0}} = \mathrm{id}$. By choosing a sufficiently small neighborhood $U \in \mathrm T_{x_0}\mc M$, $R(x_0, \cdot)$ is a differmorphisim from $U$ to $R(x_0, U)$. By the construction of $L_{x_0}$, it immediately follows that $L_{x_0}(\cdot) = (R(x_0, \cdot))^{-1}$ on $R(x_0, U)$. 
        Notice that $R$ is proved to satisfy $\mc P_{x_0} \frac{\mathrm d^2}{\mathrm dt^2} R(x_0, t v )\big\vert_{t =0 } = 0$, or equivalently, $ \frac{\mathrm d^2}{\mathrm dt^2}\Log_{x_0}(R(x_0, t v ))\big\vert_{t =0 } = 0$, for every $v \in \mathrm T_{x_0}\mc M$, by Theorem~22 in \cite{absil2012projection}. We thus conclude that $L_{x_0}$ is a second-order inverse retraction as it satisfies $\frac{\mathrm d^2}{\mathrm dt^2}L_{x_0}(\Exp_{x_0}(tv)) =0$ for every $v \in \mathrm T_{x_0}\mc M$. \qed

\section{Further Details of Numerical Studies}
We provide an overview of the statistics used in the experiments presented in Section~\ref{sec: numerical simulations} to aid the reader’s understanding. Table~\ref{table: explicit expressions} summarizes the employed Euclidean representations as well as the specific statistics used to generate the figures in the main text.

\begin{table}[htbp]
    \centering
    \makebox[\textwidth][c]{%
    \begin{tabular}{c|c|c|c}
    \toprule
    Experiment & Manifold & Euclidean Representation & Used Statistic \\ 
    \midrule 
    Type-I/II Errors 
        & Sphere $\bb S^{2}$ 
        & $\phi(\boldsymbol \theta) = \frac{1}{\theta_2}(\theta_1, \theta_3)^{\top} $ 
        & \makecell{
        $    \phi (\hat\theta_n)\t    \mb \Sigma_{\phi}^{-1}   \phi (\hat\theta_n)$ \\ 
        \& $  \big(  \phi (\hat\theta_n^{*[i]}) -   \phi (\hat\theta_n)\big)\t  {{}  \mb \Sigma_{\phi}^{*[i]}}^{\dagger} \big(  \phi (\hat\theta_n^{*[i]}) -   \phi (\hat\theta_n)\big)$; \\ 
        $ \frac{  \phi(\hat \theta_n )_1}{(\mb \Sigma_\phi)_{1,1}^{1/2}}$~\&~$\frac{  \phi(\hat \theta_n^{*[i]})_1 -   \phi(\hat \theta_n)_1}{(\mb   \Sigma_\phi^{*[i]})_{1,1}^{1/2}}$
        } \\ 
    \midrule 
    \multirow{3}{*}{\makecell{\\ 
    Converg. of Intr. Stats. \\  (Wald~\&~$t$-stat)}} 
        & Stiefel $\mathrm{St}(4,2)$
        &  $R_{\mb X}^{-1}(\mb Y) \approx \pi_{\mb H}\circ \mathsf{Log}_{\mb X}(\mb Y
        )$ 
        &  \rule{0pt}{4ex} \rule[-3ex]{0pt}{0pt} \multirow{3}{*}{\makecell{
        $R_{\hat \theta_n}^{-1}(\theta_0)\t \hat{\boldsymbol \Sigma}^{\dagger} R^{-1}_{\hat \theta_n}(\theta_0)$\\ 
        \& $R_{\hat \theta^{*[i]}_n}^{-1}(\hat \theta_n)\t {{}\check{\boldsymbol \Sigma}^{[i]}}^{\dagger} R^{-1}_{\hat \theta_n^{*[i]}}(\hat \theta_n)$ \text{ for Wald}; \\
        $ \frac{R^{-1}_{\hat \theta_n}(\theta_0)_1}{\hat{\boldsymbol \Sigma}_{1,1}^{\frac12}}$ 
        \& $\frac{R^{-1}_{\hat \theta_n^{*[i]}}(\hat \theta_n)_1}{{{}\check{\boldsymbol \Sigma}^{[i]}_{1,1}}^{\frac12}}$ \text{ for $t$-stat.} 
        }} \\ 
    \cline{2-3}
    \rule{0pt}{4ex} \rule[-3ex]{0pt}{0pt}
        & Fx-Rk Matrix $\mc R_{2,4,4}$
        &  $R_{\mb X}^{-1}(\mb Y)$ in Example~4.5
        & \\ 
    \cline{2-3} \rule{0pt}{4ex} \rule[-3ex]{0pt}{0pt}
        & Rank-$1$ Tensor $\mc M^{(1)}_{3,3,3}$ 
        &  $R_{\mb X}^{-1}(\mb Y)$ in Example~4.6
        & 
          \\ 
    \midrule Converg. of Extr. Stats. & Fx-Rk Matrix $\mc R_{2,4,4}$ & $\mathrm{id}: \mc R_{2,4,4} \rightarrow \bb R^{4\times 4} \cong \bb R^{16}$ & $\frac{(\hat \theta_n - \theta_0)_{1,1}}{(\mb c_{\hat \theta_n}\t \hat{\boldsymbol \Sigma} \mb c_{\hat \theta_n})^{\frac12}}$~\&~$\frac{(\hat \theta_n^{*[i]} - \hat \theta_n)_{1,1}}{(\mb c_{\hat \theta_n^{*[i]}}\t \check{\boldsymbol \Sigma}^{[i]} \mb c_{\hat \theta_n^{*[i]}})^{\frac12}}$ \\ 
    \midrule \rule{0pt}{6ex} \rule[-6ex]{0pt}{0pt} Barycenter & Sphere $\bb S^2$ & \makecell{$R_{\mb X}^{-1}(\mb Y)$ in Example~4.1; \\
    $\mathrm{id}: \bb S^2 \rightarrow \bb R^3$
    } & \makecell{
        $R_{\hat \theta_n}^{-1}(\theta_0)\t \hat{\boldsymbol \Sigma}^{\dagger} R^{-1}_{\hat \theta_n}(\theta_0)$\\ 
        \& $R_{\hat \theta^{*[i]}_n}^{-1}(\hat \theta_n)\t {{}\check{\boldsymbol \Sigma}^{[i]}}^{\dagger} R^{-1}_{\hat \theta_n^{*[i]}}(\hat \theta_n)$ \text{ for Wald}; \\
        $ \frac{R^{-1}_{\hat \theta_n}(\theta_0)_1}{\hat{\boldsymbol \Sigma}_{1,1}^{\frac12}}$ 
        \& $\frac{R^{-1}_{\hat \theta_n^{*[i]}}(\hat \theta_n)_1}{{{}\check{\boldsymbol \Sigma}^{[i]}_{1,1}}^{\frac12}}$ \text{ for intr. $t$-stat.} ; 
        \\ 
        $\frac{(\hat \theta_n - \theta_0)_{1}}{(\mb c_{\hat \theta_n}\t \hat{\boldsymbol \Sigma} \mb c_{\hat \theta_n})^{\frac12}}$~\&~$\frac{(\hat \theta_n^{*[i]} - \hat \theta_n)_{1}}{(\mb c_{\hat \theta_n^{*[i]}}\t \check{\boldsymbol \Sigma}^{[i]} \mb c_{\hat \theta_n^{*[i]}})^{\frac12}}$ for extr. t-stat.
        }
    \\
    \bottomrule
    \end{tabular}}
    \caption{Explicit Expressions Used in Numerical Simulations.}
    \label{table: explicit expressions}
\end{table}

We further comment on the expressions listed in Table~\ref{table: explicit expressions}: 
\begin{enumerate}
    \item The computation of $\boldsymbol \Sigma_\phi$, $\boldsymbol \Sigma_{\phi}^{*[i]}$, $\hat{\boldsymbol \Sigma}$, and $\check{\boldsymbol \Sigma}^{[i]}$ is based on the gradient and Hessian of the respective underlying objective functions, as derived in Sections~\ref{subsec: riemannian manifold examples} and~\ref{subsec: stats applications}. 
    \item In the second experiment for the Stiefel manifold, following the approach of Example~4.2, the logarithmic map is approximated numerically (see, e.g., \cite{mataigne2025efficient}). The linear mapping $\pi_{\mb H}$ from $\mathrm{T}_{\mb X}\mathrm{St}(4,2)$ to $\bb R^5$ is induced by the matrix $\mb H$ introduced therein.
    \item In the third experiment, the vector $\mb c_{\hat \theta_n}$ coincides with $\bar \nabla f \circ R_{\hat \theta_n}\vert_{\mb 0}$, where $f(\theta) = (\theta)_{1,1}$ and $R_{\hat \theta_n}$ is the SVD-based retraction described in Example~4.5. Writing the singular value decomposition of $\hat \theta_n$ as $\hat \theta_n = \hat{\mb U} \left[ 
    \begin{matrix}
        \hat{\boldsymbol \Sigma} & \mb 0 \\ 
         \mb 0 & \mb 0
    \end{matrix}
    \right] \hat{\mb V}\t$, we obtain $\mb c_{\hat \theta_n} = \big(\mathsf{vec}(\hat{\mb U}_{1,\cdot}\t \hat{\mb V}_{1,1:2})\t, \mathsf{vec}(\hat{\mb U}_{1,1:2}\t \hat{\mb V}_{1,3:4})\t\big)\t$. An analogous definition applies to $\mb c_{\hat \theta_n^{*[i]}}$.
    \item For the final experiment, we employ the smooth mapping $f^{\mathsf{ortho}, \theta}$ introduced at the beginning of Section~\ref{sec:applications} to construct the retractions $R_\theta$ and their inverses, as detailed in Example~4.1. Accordingly, the vectors $\mb c_{\hat \theta_n}$ and $\mb c_{\hat \theta_n^{*[i]}}$ are defined as $(f^{\mathsf{ortho}, \hat \theta_n})_{1,\cdot}\t $ and $(f^{\mathsf{ortho}, \hat \theta_n^{*[i]}})_{1,\cdot}\t $, respectively.
\end{enumerate}

\end{appendix}

\end{sloppypar}

\end{document}